\def\DateTime{June 20th, 2022}

\def\ifams{\if00}
\def\ifsmf{\if01}

\ifsmf
\documentclass[reqno,english,draft]{smfbook}
\usepackage{smfthm,CJK}
\NumberTheoremsIn{section}
\fi
\ifams
\documentclass[a4paper,reqno,11pt,english,draft]{amsbook}
\fi
\usepackage[T2A,T1]{fontenc}
\ifsmf
\usepackage[english,francais]{babel}
\fi
\ifams
\usepackage[english]{babel}
\fi
\usepackage[utf8]{inputenc}
\usepackage{latexsym,amscd,color}
\usepackage{amsmath,amsfonts,amssymb,mathrsfs,mathtools}
\usepackage{qsymbols,euscript,enumitem}
\usepackage{multirow,bigdelim}
\usepackage{bm,comment}
\usepackage{url,xspace}
\usepackage{verbatim}
\usepackage{tikz}
\usepackage{tikz-cd}
\usetikzlibrary{positioning,arrows,shapes}
\usepackage[normalem]{ulem}
\usepackage{manfnt} % Added by Atsushi
\usepackage{etoolbox}

\usepackage{hyperref}
\usepackage{xstring}

\makeindex
\numberwithin{equation}{chapter}

\ifams
\theoremstyle{plain}
\newtheorem{theo}{Theorem}[section]%{\normalfont\bfseries}{\itshape}
\newtheorem{prop}[theo]{Proposition}%{\normalfont\bfseries}{\itshape}
\newtheorem{coro}[theo]{Corollary}%{\normalfont\bfseries}{\itshape}
\newtheorem{lemm}[theo]{Lemma}%{\normalfont\bfseries}{\itshape}
\newtheorem{clai}[theo]{Claim}%{\normalfont\bfseries}{\itshape}
\newtheorem{code}[theo]{Corollary-Definition}
\newtheorem{thde}[theo]{Theorem-Definition}

\theoremstyle{definition}
\newtheorem{defi}[theo]{Definition}%{\normalfont\bfseries}{\normalfont}
\newtheorem{rema}[theo]{Remark}%{\normalfont\bfseries}{\normalfont}
\fi

\author{Huayi CHEN}
\address{Universit\'e Paris Cit\'e, Sorbonne Universit\'e, CNRS, INRIA,  IMJ-PRG, F-75013 Paris, France}
\email{huayi.chen@imj-prg.fr}
\author{Atsushi MORIWAKI}
\address{Department of Mathematics, Faculty of Science, Kyoto University, Kyoto, 606-8502, Japan}
\email{moriwaki@math.kyoto-u.ac.jp}
\title[Hilbert-Samuel formula and positivity over adelic curves]{Hilbert-Samuel formula and positivity over adelic curves}
\date{\DateTime}

\input xy
\xyoption{all}
\newcommand{\emptyinnprod}{\langle\kern.15em,\kern-.02em\rangle}

\definecolor{indigo(web)}{rgb}{0.29, 0.0, 0.51}
\def\sbullet{{\scriptscriptstyle\bullet}}

\newcommand{\QQ}{{\mathbb{Q}}}
\newcommand{\RR}{{\mathbb{R}}}

\newcommand{\PP}{{\mathbb{P}}}
\newcommand{\NN}{{\mathbb{N}}}

\newcommand{\OO}{{\mathcal{O}}}

\newcommand{\Ker}{\operatorname{Ker}}

\newcommand{\Spec}{\operatorname{Spec}}

\newcommand{\Supp}{\operatorname{Supp}}

\newcommand{\rank}{\operatorname{rk}}

\newcommand{\adeg}{\widehat{\operatorname{deg}}}

\newcommand{\Sym}{\operatorname{Sym}}

\newcommand{\vol}{\operatorname{vol}}
\newcommand{\avol}{\widehat{\operatorname{vol}}}
\newcommand{\an}{\operatorname{an}}

\newcommand{\ndot}{\raisebox{.4ex}{.}}

\newcommand{\rest}[2]{\left.{#1}\right\vert_{{#2}}}

\newcommand{\sub}{\operatorname{sub}}
\newcommand{\quot}{\operatorname{quot}}

\newcommand{\indic}{1\hspace{-0.25em}\mathrm{l}}
\newcommand\redsout{\bgroup\markoverwith{\textcolor{mred}{\rule[0.5ex]{2pt}{0.7pt}}}\ULon}
\def\colorsout#1{\bgroup\markoverwith{\textcolor{#1}{\rule[0.5ex]{2pt}{0.7pt}}}\ULon} %[0.5ex]{2pt}{0.4pt}
\def\coloruline#1{\bgroup\markoverwith{\textcolor{#1}{\rule[-0.5ex]{2pt}{0.7pt}}}\ULon} %[0.5ex]{2pt}{0.4pt}
\DeclareSymbolFont{bbold}{U}{bbold}{m}{n}
\DeclareMathSymbol{\bbalpha}{\mathord}{bbold}{"0B}
\DeclareMathSymbol{\bbbeta}{\mathord}{bbold}{"0C}
\DeclareMathSymbol{\bbgamma}{\mathord}{bbold}{"0D}
\DeclareMathSymbol{\bbdelta}{\mathord}{bbold}{"0E}
\DeclareMathSymbol{\bbespilon}{\mathord}{bbold}{"0F}
\DeclareMathSymbol{\bbzeta}{\mathord}{bbold}{"10}
\DeclareMathSymbol{\bbeta}{\mathord}{bbold}{"11}
\DeclareMathSymbol{\bbtheta}{\mathord}{bbold}{"12}
\DeclareMathSymbol{\bbiota}{\mathord}{bbold}{"13}
\DeclareMathSymbol{\bbkappa}{\mathord}{bbold}{"14}
\DeclareMathSymbol{\bblambda}{\mathord}{bbold}{"15}
\DeclareMathSymbol{\bbmu}{\mathord}{bbold}{"16}
\DeclareMathSymbol{\bbnu}{\mathord}{bbold}{"17}
\DeclareMathSymbol{\bbxi}{\mathord}{bbold}{"18}
\DeclareMathSymbol{\bbpi}{\mathord}{bbold}{"19}
\DeclareMathSymbol{\bbrho}{\mathord}{bbold}{"1A}
\DeclareMathSymbol{\bbsigma}{\mathord}{bbold}{"1B}
\DeclareMathSymbol{\bbtau}{\mathord}{bbold}{"1C}
\DeclareMathSymbol{\bbupsilon}{\mathord}{bbold}{"1D}
\DeclareMathSymbol{\bbphi}{\mathord}{bbold}{"1E}
\DeclareMathSymbol{\bbchi}{\mathord}{bbold}{"1F}
\DeclareMathSymbol{\bbpsi}{\mathord}{bbold}{"20}
\definecolor{ruby}{rgb}{0.88, 0.07, 0.37}
\definecolor{coolblack}{rgb}{0.0, 0.18, 0.39}
\definecolor{darkspringgreen}{rgb}{0.09, 0.45, 0.27}
\definecolor{emerald}{rgb}{0.31, 0.78, 0.47}
\definecolor{lavenderindigo}{rgb}{0.58, 0.34, 0.92}
\definecolor{mred}{rgb}{0.83, 0.0, 0.0}%same as rossocorsa
\definecolor{indigo(web)}{rgb}{0.29, 0.0, 0.51}
\definecolor{darkblue}{rgb}{0.0, 0.0, 0.55}

\begin{document}
\begin{abstract}
We establish, in the setting of Arakelov geometry over adelic curves, an arithmetic Hilbert-Samuel theorem describing the asymptotic behaviour of the metrized graded linear series of an adelic line bundle in terms of its arithmetic intersection number. We then study positivity conditions of adelic line bundles.
\end{abstract}

%\begin{abstract}
%The notions of ``positivity'' and their properties are the most important subjects in algebraic geometry.
%In this article, we consider an arithmetic analogue over a general adelic curve as a generalization of the classical Arakelov geomerty
%over a number field.
%\end{abstract}

\ifams\fi
\ifsmf\fi
\thanks{The second author was supported by JSPS KAKENHI Grant-in-Aid for Scientific Research(S) Grant Number JP16H06335 and Scientific Research(C) Grant Number JP21K03203.}
\ifsmf\frontmatter\fi
\maketitle

\setcounter{tocdepth}{1}
\tableofcontents
\ifsmf\mainmatter\fi

%!TEX root = ./Hilbert_Samuel_Adelic_Curves.tex

\chapter*{Introduction}

\ifsmf
\newtheorem{introtheo}{Theorem}
\def\theintrotheo{\arabic{chapter}.\arabic{introtheo}}
\fi

In algebraic geometry, Hilbert function measures the growth of graded linear series of a line bundle on a projective variety. Let $k$ be a field, $X$ be an integral projective scheme of dimension $d\in\mathbb N\ (= \mathbb Z_{\geqslant 0})$ over $\Spec k$, and $L$ be an invertible $\mathcal O_X$-module. The Hilbert function of $L$ is defined as 
\[H_L:\mathbb N\longrightarrow\mathbb N,\quad H_L(n):=\dim_k(H^0(X,L^{\otimes n})). \]
If $L$ is ample, then the following asymptotic estimate holds:
\begin{equation}\label{Equ: Hilbert-Samuel function}H_L(n)=\frac{(L^d)}{d!}n^d+o(n^d).\end{equation}
This formula, which relates the asymptotic behaviour of the Hilbert function and the auto-intersection number of $L$, is for example a consequence of Hirzebruch-Riemann-Roch theorem and Serre's vanishing theorem. It turns out that the construction  and the asymptotic estimate of Hilbert function  have analogue in various context, such as graded algebra, local multiplicity, relative volume of two metrics, etc. 

In Arakelov geometry, an arithmetic analogue of Hilbert function has been introduced by Gillet and Soul\'{e} \cite{MR1189489} and an analogue of the asymptotic formula \eqref{Equ: Hilbert-Samuel function} has been deduced from their arithmetic Riemann-Roch theorem. This result is called an arithmetic Hilbert-Samuel theorem. Let $\mathscr X$ be a regular integral projective scheme of dimension $d+1$ over $\Spec\mathbb Z$, and $\overline{\mathscr L}=(\mathscr L,\varphi)$ be a Hermitian line bundle on $\mathscr X$, namely an invertible $\mathcal O_{\mathscr X}$-module $\mathscr L$ equipped with a smooth metric $\varphi$ on $\mathscr L(\mathbb C)$. For any integer $n\in\mathbb N$, we let $\|\ndot\|_{n\varphi}$ be the norm on the real vector space $H^0(\mathscr X,\mathscr L)\otimes_{\mathbb Z}\mathbb R$ defined as follows 
\[\forall\,s\in H^0(\mathscr X,\mathscr L)\otimes_{\mathbb Z}\mathbb R\subseteq H^0(\mathscr X_{\mathbb C},\mathscr L_{\mathbb C}^{\otimes n}),\quad \|s\|_{n\varphi}=\sup_{x\in\mathscr X(\mathbb C)}|s|_{n\varphi}(x).\]
Then the couple $(H^0(\mathscr X,\mathscr L^{\otimes n}),\|\ndot\|_{n\varphi})$ forms a lattice in a normed vector space. Recall that its arithmetic Euler-Poincar\'{e} characteristic  is   
\[\chi(H^0(\mathscr X,\mathscr L^{\otimes n}),\|\ndot\|_{n\varphi})=\ln\frac{\operatorname{vol}(\{s\in H^0(\mathscr X,\mathscr L^{\otimes n})\otimes_{\mathbb Z}\mathbb R\,:\,\|s\|_{n\varphi}\leqslant 1 \})}{\operatorname{covol}(H^0(\mathscr X,\mathscr L^{\otimes n}),\|\ndot\|_{n\varphi})}\]
where $\operatorname{vol}(\ndot)$ denotes a Haar measure on the real vector space \[H^0(\mathscr X,\mathscr L)\otimes_{\mathbb Z}\mathbb R,\] and $\operatorname{covol}(H^0(\mathscr X,\mathscr L^{\otimes n}),\|\ndot\|_{n\varphi})$ is the covolume of the lattice
$H^0(\mathscr X,\mathscr L^{\otimes n})$ with respect to the Haar measure $\operatorname{vol}(\ndot)$, namely the volume of any fundamental domain of this lattice. In this setting the arithmetic Hilbert-Samuel theorem shows that, if $\mathscr L$ is relatively ample and the metric $\varphi$ is positive, then the sequence
\[\frac{\chi(H^0(\mathscr X,\mathscr L^{\otimes n}),\|\ndot\|_{n\varphi})}{n^{d+1}/(d+1)!},\quad n\in\mathbb N,\;n\geqslant 1\]
converges to the arithmetic intersection number $(\overline{\mathscr L}^{d+1})$. In the case where $\overline{\mathscr L}$ is ample, the arithmetic Hilbert-Samuel theorem also permits to relate the asymptotic behaviour (when $n\rightarrow+\infty$) of 
\[\operatorname{card}(\{s\in H^0(\mathscr X,\mathscr L^{\otimes n})\,:\,\|s\|_{n\varphi}\leqslant 1\})\]
to the arithmetic intersection number of $\overline{\mathscr L}$.
 These results have various applications in arithmetic geometry, such as Vojta's proof of Mordell conjecture, equidistribution problem and Bogomolov conjecture, etc. The arithmetic Hilbert-Samuel theorem has then been reproved in various settings and also been generalized in works such as \cite{MR1343555,MR1971294,MR2208129}.

Recently, a new framework of Arakelov geometry has been proposed in \cite{CMArakelovAdelic}, which allows to consider arithmetic geometry over any countable field. Let $K$ be a field. A structure of proper adelic curve with underlying field $K$ is given by a family of absolute values $(|\ndot|_{\omega})_{\omega\in\Omega}$ on $K$ parametrized by a measure space $(\Omega,\mathcal A,\nu)$, which satisfies a product formula of the form
\[\forall\,a\in K^{\times},\quad \int_{\Omega}\ln|a|_\omega\,\nu(\mathrm{d}\omega)=0.\] 
We assume that, either $K$ is countable, or the $\sigma$-algebra $\mathcal A$ is discrete. This notion is a very natural generalization to any countable field of Weil's adelic approche of number theory. The fundament of height theory and Arakelov geometry for projective varieties over an adelic curve have been established in the works of Gubler \cite{MR1472498} (in a slightly different setting of $M$-fields) and Chen-Moriwaki \cite{CMArakelovAdelic}, respectively, see also the model theoretical approach of Hrushovski and Ben Yaacov \cite{Hrushovski16}.   
 More recently, the arithmetic intersection theory in the setting of adelic curves have been developed in \cite{CMIntersection}. Note that in general it is not possible to consider global integral models of an adelic curve. Several classic notions and constructions, such as integral lattice and its covolume, do not have adequate analogue over adelic curves. It turns out that a modified and generalized form of normed lattice --- adelic vector bundle --- has a natural avatar in the setting of adelic curves. An adelic vector bundle consists of a finite-dimensional vector space $V$ over $K$ equipped with a family of norms $(\|\ndot\|_{\omega})_{\omega\in\Omega}$ on vector spaces $V_\omega=V\otimes_KK_\omega$ (where $K_\omega$ denotes the completion of $K$ with respect to the absolute value $|\ndot|_\omega$), which satisfy dominancy and measurability conditions. The Arakelov degree of the adelic vector bundle 
\[\overline V=(V,(\|\ndot\|_\omega)_{\omega\in\Omega})\]
is then defined as 
\[\operatorname{\widehat{\deg}}(\overline V):=-\int_{\Omega}\ln\|s_1\wedge\cdots\wedge s_r\|_{\omega,\det}\,\nu(\mathrm{d}\omega),\]
where $(s_i)_{i=1}^r$ is an arbitrary basis of $E$ over $K$. This notion is a good candidate to replace the Euler-Poincar\'{e}
characteristic. 

Let $\pi:X\rightarrow\Spec K$ be a projective scheme over $\operatorname{Spec}K$. For any $\omega\in\Omega$, let $X_\omega=X\times_{\operatorname{Spec}K}\operatorname{Spec}K_\omega$ and let $X_\omega^{\mathrm{an}}$ be the analytic variety associated with $X_\omega$ (in the sense of Berkovich \cite{MR1070709} if $|\ndot|_\omega$ is non-Archimedean). If $E$ is a vector bundle on $X$, namely a locally free $\mathcal O_X$-module of finite rank, we denote by $E_\omega$ the pull-back of $E$ on $X_\omega$. As \emph{adelic vector bundle}\index{adelic vector bundle} on $X$, we refer to the data $\overline E=(E,(\psi_\omega)_{\omega\in\Omega})$ consisting of a vector bundle $E$ on $X$ and a family $(\psi_\omega)_{\omega\in\Omega}$ of continuous metrics on $E_\omega$ with $\omega\in\Omega$, which satisfy dominancy and measurability conditions. It turns out that, if $X$ is geometrically reduced, then the vector space of global sections $H^0(X,E)$ equipped with supremum norms $(\|\ndot\|_{\psi_\omega})_{\omega\in\Omega}$ forms an adelic vector bundle $\pi_*(\overline E)$ on the base adelic curve. 

Let $\pi:X\rightarrow\Spec K$ be an integral projective scheme of dimension $d$ over $\Spec K$ and  $\overline L=(L,\varphi)$ be an adelic line bundle on $X$, that is, an adelic vector bundle of rank $1$ on $X$. Assume that the line bundle $L$ is ample. We introduce the notion of $\chi$-volume as 
\[\widehat{\vol}_\chi(\overline L)=\limsup_{n\rightarrow+\infty}\frac{\widehat{\deg}(\pi_*(\overline L{}^{\otimes n}))}{n^{d+1}/(d+1)!}.\]
In view of the similarity between Arakelov degree and Euler-Poincar\'{e} characteristic of Euclidean lattices, the notion of $\chi$-volume is analogous to that of sectional capacity introduced in \cite{MR1677934}, or to that of volume in \cite{MR2425137}. Moreover, similarly to the number field case, we show in Theorem-Definition \ref{theo:limit:thm:vol:chi} that the above superior limit defining the $\chi$-volume is actually a limit. However, from the methodological view, we do not follow the classic approaches, which are difficultly implantable in the adelic curve setting. Our strategy consists in casting the Arakelov geometry over an adelic curve to that in the particular case where the adelic curve contains a single copy of the trivial absolute value on $K$, that is, the absolute value $|\ndot|_0$ such that $|a|_0=1$ for any $a\in K\setminus\{0\}$. More precisely, to each adelic vector bundle $\overline V=(V,(\|\ndot\|_\omega)_{\omega\in\Omega})$, we associate an ultrametric norm $\|\ndot\|_0$ on $V$ (where we consider the trivial absolute value $|\ndot|_0$) via Harder-Narasimhan theory in the form of $\mathbb R$-filtrations, such that 
\[\big|\widehat{\deg}(V,(\|\ndot\|_\omega)_{\omega\in\Omega})-\widehat{\deg}(V,\|\ndot\|_0)\big|\leqslant \frac 12\nu(\Omega_\infty)\dim_K(V)\ln(\dim_K(V)),\]
where $\Omega_\infty$ denotes the set of $\omega\in\Omega$ such that $|\ndot|_{\omega}$ is Archimedean. Then the convergence of the sequence defining $\widehat{\vol}_{\chi}(\overline L)$ follows from a limit theorem of normed graded linear series as follows (see Theorem \ref{Thm: existence of volume} and Corollary \ref{Cor: volume chi} for this result in a more general form and for more details):

\ifsmf\begin{introtheo}\fi
\ifams\begin{theo}\fi
Assume that the graded $K$-algebra $\bigoplus_{n\in\mathbb N}H^0(X,L^{\otimes n})$ is of finite type.
For any integer $n\geqslant 1$, let $\|\ndot\|_n$ be a norm on $H^0(X,L^{\otimes n})$ \textup{(}where we consider the trivial absolute value on $K$\textup{)}. Assume that
\begin{enumerate}[label=\rm(\alph*)] 
\item $\displaystyle\inf_{s\in V_n\setminus\{0\}}\ln\|s\|_n=O(n)$ when $n\rightarrow+\infty$,
\item 
for any $(n,m)\in\mathbb N_{\geqslant 1}^2$ and any $(s_n,s_m)\in V_n\times V_m$, one has
\[\|s_n\cdot s_m\|_{n+m}\leqslant \|s_n\|_n\cdot\|s_m\|_m.\]
\end{enumerate}
Then the sequence 
\[\frac{\widehat{\deg}(V_n,\|\ndot\|_n)}{n^{d+1}/(d+1)!},\quad n\in\mathbb N_{\geqslant 1}\]
converges in $\mathbb R$.
\ifsmf\end{introtheo}\fi
\ifams\end{theo}\fi

In view of the classic Hilbert-Samuel theorems in algebraic geometry and in Arakelov geometry, it is natural to compare the $\chi$-volume to the arithmetic intersection number of adelic line bundles that we have introduced in \cite{CMIntersection} (see also the work \cite{MR1472498} on heights of varieties over $M$-fields). Let $\pi:X\rightarrow\Spec K$ be a projective scheme of dimension $d\geqslant 0$ over $K$ and $\overline L=(L,\varphi)$ be an adelic line bundle on $X$ such that $L$ is ample and the metrics in the family $\varphi$ are semi-positive. Then the arithmetic self-intersection number $(\overline L{}^{d+1})$ of $\overline L$ is written in a recursive way as 
\begin{equation}\label{Equ: intersection number}\frac{1}{N}\bigg[(\overline L|_{\operatorname{div}(s)}^d)_S-\int_{\Omega}\int_{X_\omega^{\mathrm{an}}}\ln|s|_{\varphi_\omega}(x)\,c_1(L_\omega,\varphi_\omega)^d(\mathrm{d}x)\,\nu(\mathrm{d}\omega)\bigg],\end{equation}
where $N$ is a positive integer, and $s$ is a global section of $L^{\otimes N}$ which intersects properly with all irreducible components of the projective scheme $X$. The main result of the article is then the following theorem (see Theorem \ref{thm:lower:bound}).

\ifsmf\begin{introtheo}\fi
\ifams\begin{theo}\fi
\label{Thm: HS}
Assume that, either $X$ is geometrically integral, or the field $K$ is perfect. Let  $\overline L=(L,\varphi)$ be an adelic line bundle on $X$ such that $L$ is ample and that all metrics in the family $\varphi$ are semi-positive, then the following equality holds:
\[\widehat{\vol}_{\chi}(\overline L)=(\overline L{}^{d+1}).\]
\ifsmf\end{introtheo}\fi
\ifams\end{theo}\fi

Note that in the literature there exists a local version of the Hilbert-Samuel theorem which establishes an equality between the relative volume of two metrics and the relative Monge-Amp\`{e}re energy between them. We refer the readers to \cite{MR2657428} for the Archimedean case and to \cite{MR4061325,MR4192993} for the non-Archimedean case (see also \cite{BGF2020}). These results show that, for a fixed ample line bundle $L$ on $X$, the difference between $\widehat{\vol}_{\chi}(\overline L)$ and $(\overline L{}^{d+1})$ does not depend on the choice of the metric family on $L$ (see Proposition \ref{Pro: HS property particular case} and Remark \ref{Rem: difference is constant}). Moreover, by an argument of projection to a projective space (on which the arithmetic Hilbert-Samuel theorem can be proved by explicit computation, see Proposition \ref{prop:Hilbert:Samuel:Projective:Space}), one can show that the inequality $\widehat{\vol}_{\chi}(\overline L)\geqslant(\overline L{}^{d+1})$ holds (see Step 2 of the proof of Theorem \ref{thm:lower:bound}). 

In view of the recursive formula \eqref{Equ: intersection number} defining the self-intersection number, a natural idea to prove the above theorem could be an argument of induction, following the approach of \cite{MR1343555} by using an adaptation to non-Archimedean setting of some technics in complex analytic geometry developed in \cite{MR4192993,FangYanbo}. However, it seems that a refinement in the form of an asymptotic development of the function defining the local relative volume is needed to realize this strategy. Unfortunately such refinement is not yet available. Our approach consists in casting the arithmetic data of $\overline L$ to a series of metrics over a trivially valued field. This could be considered as a higher-dimensional generalization of the approach of Harder-Narasimhan $\mathbb R$-filtration mentioned above.  What is particular in the trivial valuation case is that the local geometry becomes automatically global, thanks to the trivial ``product formula''. In this case, the arithmetic Hilbert-Samuel theorem follows from  the equality between the relative volume and the relative Monge-Amp\`{e}re energy with respect to the trivial metric (see Theorem \ref{Thm: trivial valuation case}). Note that this result also shows that, in the case of a projective curve over a trivially valued field,  the arithmetic intersection number defined in \cite{CMIntersection} coincides with that constructed in a combinatoric way in \cite{Chen_Moriwaki2020} (see Remark \ref{Rem: comparison with combinatorics}).  The comparison of diverse invariants of $\overline L$ with respect to those of its casting to the trivial valuation case provides the opposite inequality $\widehat{\vol}_{\chi}(\overline L)\leqslant(\overline L{}^{d+1})$. As a sequel to the above arguments in terms of trivially valued fields, 
our way towards the arithmetic Hilbert-Samuel theorem over a adelic curve gives a new approach even for the classical case.

As an application, we prove the following higher dimensional generalization of Hodge index theorem (see Corollaries \ref{coro:lower:bound:nef} and \ref{Cor: bigness of line bundle}).

\ifsmf\begin{introtheo}\fi
\ifams\begin{theo}\fi
Assume that, either $X$ is geometrically integral, or the field $K$ is perfect. Let $\overline L=(L,\varphi)$ be an adelic line bundle on $X$. Assume that $L$ is nef and all metrics in the family $\varphi$ are semi-positive, then the inequality $\widehat{\vol}(\overline L)\geqslant(\overline L{}^{d+1})$ holds. In particular, if $(\overline L{}^{d+1})>0$, then the line bundle $L$ is big.
\ifsmf\end{introtheo}\fi
\ifams\end{theo}\fi

Theorem \ref{Thm: HS} naturally leads to the following refinement of the arithmetic Hilbert-Samuel theorem, in introducing a tensor product by an adelic vector bundle on $X$ (see Corollary~\ref{Cor:Hilbert:Samuel:with:auxiliary:locally:free:module}). As in Theorem \ref{Thm: HS}, we assume that, either $X$ is geometrically integral, or the field $K$ is perfect.

\ifsmf\begin{introtheo}\fi
\ifams\begin{theo}\fi
Let $\overline L=(L,\varphi)$ be an adelic line bundle on $X$ and $\overline E=(E,\psi)$ be an adelic vector bundle on $X$. Assume that $L$ is ample and the metrics in $\varphi$ are semi-positive. Moreover we suppose that either $\rank(E) = 1$ or $X$ is normal.
Then one has
\[\lim_{n\rightarrow+\infty}\frac{\operatorname{\widehat{\deg}}(H^0(X,L^{\otimes n}\otimes E),(\|\ndot\|_{n\varphi_\omega+\psi_\omega})_{\omega\in\Omega})}{n^{d+1}/(d+1)!}=\rank(E)(\overline L{}^{d+1}).\]
\ifsmf\end{introtheo}\fi
\ifams\end{theo}\fi

The second part of the article is devoted to the study of positivity conditions of adelic line bundles. Positivity of line bundles is one of the most fundamental and important notions in algebraic geometry. In Arakelov geometry, the analogue of ampleness and Nakai-Moishezon criterion have been studied by Zhang \cite{MR1189866,MR1254133}. The arithmetic bigness has been introduced in the works \cite{MR1779799, MR2425137, MR2496453} of Moriwaki and Yuan.
These positivity conditions and their properties have various applications in Diophantine geometry.

We assume that the underly field $K$ of the adelic curve $S$ is perfect. Let $X$ be a projective scheme over $\Spec K$. Given an adelic line bundle $\overline L$ on $X$,
we are interested in various positivity conditions of the adelic line bundle $\overline L$. We say that the adelic line bundle $\overline L$ is \emph{relatively ample}\index{ample!relatively ---} if the invertible $\mathcal O_X$-module is ample and if the  metrics of $\overline L$ are all semi-positive. The relative nefness can then been defined in a limit form of relative ampleness, similarly to the classic case in algebraic geometry. Recall that the global intersection number of relatively ample adelic line bundles (or more generally, integrable adelic line bundles) can be defined as the integral of local heights along the measure space in the adelic structure (cf. \cite{CMIntersection}). This construction is fundamental in the Arakelov height theory of projective varieties. 

We first introduce a numerical invariant ---  asymptotic minimal slope --- to describe the global positivity of an adelic line bundle $\overline L$ such that $L$ is ample. This invariant, which is denoted by $\widehat{\mu}_{\min}^{\mathrm{asy}}(\overline L)$, describes the asymptotic behaviour (when $n\rightarrow+\infty$) of the minimal slopes of the sectional spaces $H^0(X,L^{\otimes n})$ equipped with sup norms (which are adelic vector bundles on $S$). It turns out that this invariant is super-additive with respect to $\overline L$. This convexity property allows to extend the  construction of the asymptotic minimal slope to the cone of adelic line bundles with nef underlying invertible $\mathcal O_X$-module (see \S\ref{Sec: Relative ampleness} for the construction of the asymptotic minimal slope and its properties). The importance of this invariant can be shown by the following height estimate (see Theorem \ref{Thm: lower bound intersection number} for the proof and Proposition \ref{Pro: lower bound ample} for its generalization to the relatively nef case).

\ifsmf\begin{introtheo}\fi
\ifams\begin{theo}\fi
\label{Thm: main thm 1}
Assume that the field $K$ is perfect. Let $X$ be a  reduced projective scheme of dimension $d\geqslant 0$ over $\Spec K$, and $\overline L_0,\ldots,\overline L_d$ be a family of relatively ample adelic line bundles on $X$. For any $i\in\{0,\ldots,d\}$, let $\delta_i$ be the geometric intersection number \[(L_0\cdots L_{i-1}L_{i+1}\cdots L_d).\] Then the following inequality holds:
\[(\overline L_0\cdots\overline L_d)_S\geqslant\sum_{i=1}^d\delta_i\operatorname{\widehat{\mu}_{\min}^{\mathrm{asy}}}(\overline L_i),\]
where $(\overline L_0\cdots\overline L_d)_S$ denotes the arithmetic intersection number of $\overline L_0,\ldots,\overline L_d$. 
\ifsmf\end{introtheo}\fi
\ifams\end{theo}\fi

The asymptotic minimal slope always increases if one replaces the adelic line bundle by its pullback by a projective morphism (see Theorem \ref{Thm: mu min asy and pull-back}): if $g:X\rightarrow P$ is a projective morphism of reduced $K$-schemes of dimension $\geqslant 0$, then for any adelic line bundle $\overline M$ on $P$ such that $M$ is nef, one has $\widehat{\mu}_{\min}^{\mathrm{asy}}(g^*(\overline M))\geqslant\widehat{\mu}_{\min}^{\mathrm{asy}}(\overline M)$. Typical situations include a closed embedding of $X$ into a projective space, or a finite covering over a projective space, which allow to obtain lower bounds of $\widehat{\mu}_{\min}^{\mathrm{asy}}(\overline L_i)$ in the applications of the above theorem. Note that the particular case where $\overline L_0,\ldots,\overline L_d$ are all equal to the same adelic line bundle $\overline L$ gives the following inequality  
\begin{equation}\label{Equ: height inequality}\frac{(\overline L^{d+1})_S}{(d+1)(L^d)}\geqslant\widehat{\mu}_{\min}^{\mathrm{asy}}(\overline L),\end{equation}
which relates the normalized height of $X$ with respect to $\overline L$ and the asymptotic minimal slope of the latter. This inequality is similar to the first part of \cite[Theorem 5.2]{MR1254133}. However, the imitation of the devissage argument using the intersection of hypersurfaces defined by small sections would not work in the setting of adelic curves. This is mainly due to the fact that the analogue of Minkowski's first theorem fails for adelic vector bundles on a general adelic curve. Although (in the case where $X$ is an integral scheme) the inequality \eqref{Equ: height inequality} could be obtained in an alternative way by using the arithmetic Hilbert-Samuel formula of $\overline L$ together with the fact that the minimal slope of an adelic vector bundle on $S$ is always bounded from the above by its slope (see Proposition \ref{Pro: inequality height and mu min asy }), the proof of Theorem \ref{Thm: main thm 1} needs a new idea. Our approach consists in combining an analogue of the slope theory of Bost \cite{MR1423622,MR1863738} with the height of multi-resultant.

Bigness is another type of positivity condition which describes the growth of the total graded linear series of a line bundle. In Arakelov geometry of number fields, the arithmetic bigness describes the asymptotic behaviour of the number of small sections in the graded sectional algebra of adelic vector bundles. This notion can be generalized to the setting of Arakelov geometry of adelic curve in replacing the logarithm of the number of small sections by the positive degree of an adelic vector bundle (namely the supremum of the Arakelov degrees of adelic vector subbundles). In \cite[Proposition 6.4.18]{CMArakelovAdelic}, the arithmetic bigness has been related to an arithmetic sectional invariant --- asymptotic maximal slope, which is quite similar to asymptotic minimal slope: for any integral  projective $K$-scheme and any adelic line bundle $\overline L$ on $X$ such that $L$ is big, we introduce a numerical invariant $\widehat{\mu}_{\max}^{\mathrm{asy}}(\overline L)$ which describes the asymptotic behaviour (when $n\rightarrow+\infty$) of the maximal slopes of $H^0(X,L^{\otimes n})$ equipped with sup norms  (see \S\ref{Sec: asymptotic maximal slope} for its construction and properties). It turns out that this invariant is also super-additive with respect to $\overline L$, which allows to extend the function $\widehat{\mu}_{\max}^{\mathrm{asy}}(\ndot)$ to the cone of adelic line bundles $\overline L$ such that $L$ is pseudo-effective. Moreover, in the case where $L$ is nef, the inequality $\widehat{\mu}_{\min}^{\mathrm{asy}}(\overline L)\leqslant\widehat{\mu}_{\max}^{\mathrm{asy}}(\overline L)$ holds. 

Recall that Fujita's approximation theorem asserts that a big line bundle can be decomposed in a birational modification into the tensor product of two $\mathbb Q$-line bundles which are respectively ample and effective, with a good approximation of the volume function. In this article, we establish the following relative version of Fujita's approximation theorem for the asymptotic maximal slope (see Theorem \ref{Thm: relative fujita} and Remark \ref{Rem: strong relative Fujita}).

\ifsmf\begin{introtheo}\fi
\ifams\begin{theo}\fi
Assume that the field $K$ is perfect and the scheme $X$ is integral.
Let $\overline L$ be an adelic line bundle on $X$ such that $L$ is big. For any real number $t<\widehat{\mu}_{\max}^{\mathrm{asy}}(\overline L)$, there exist a positive integer $p$, a birational projective $K$-morphism $g:X'\rightarrow X$, a relatively ample adelic line bundle $\overline A$ and an effective adelic line bundle $\overline M$ on $X'$ such that $A$ is big, $g^*(\overline L^{\otimes p})$ is isomorphic to $\overline A\otimes\overline M$, and $\widehat{\mu}_{\min}^{\mathrm{asy}}(\overline A)\geqslant pt$. 
\ifsmf\end{introtheo}\fi
\ifams\end{theo}\fi

As an application, in the case where $X$ is an integral scheme, we can improve the height inequality in Theorem \ref{Thm: main thm 1} in relaxing the positivity condition of one of the adelic line bundles and in replacing the asymptotic minimal slope of this adelic line bundle by the asymptotic maximal slope (see Theorem \ref{Thm: lower bound intersection number bigness}).

\ifsmf\begin{introtheo}\fi
\ifams\begin{theo}\fi
\label{Thm: lower bound intersection number 2}
Assume that the field $K$ is perfect. Let $X$ be an integral projective scheme of dimension $d$ over $\Spec K$, and $\overline L_0,\ldots,\overline L_d$ be adelic line bundles on $X$ such that $\overline L_1,\ldots,\overline L_d$ are relatively ample and $L_0$ is big. For any $i\in\{0,\ldots,d\}$, let $\delta_i=(L_0\cdots L_{i-1}L_{i+1}\cdots L_d)$. Then the following inequality holds:
\[(\overline L_0\cdots\overline L_d)_S\geqslant \delta_0\operatorname{\widehat{\mu}}_{\max}^{\mathrm{asy}}(\overline L_0)+\sum_{i=1}^d\delta_i\operatorname{\widehat{\mu}_{\min}^{\mathrm{asy}}}(\overline L_i).\] 
\ifsmf\end{introtheo}\fi
\ifams\end{theo}\fi

In the case where $\overline L_0,\ldots,\overline L_d$ are all equal to the same  adelic line bundle $\overline L$, the above inequality leads to 
\[\frac{(\overline L^{d+1})_S}{(L^d)}\geqslant\widehat{\mu}_{\max}^{\mathrm{asy}}(\overline L)+d\operatorname{\widehat{\mu}_{\min}^{\mathrm{asy}}}(\overline L).\] 

In the case where the adelic curve $S$ comes from the canonical adelic structure of a number field, if $\overline L$ is a relatively ample adelic line bundle, then $\widehat{\mu}_{\min}^{\mathrm{asy}}(\overline L)$ is equal to the absolue minimum of the Arakelov (absolute) height function $h_{\overline L}$ on the set of closed points of $X$. This is essentially a consequence of \cite[Corollary 5.7]{MR1254133}. Similarly, the asymptotic maximal slope $\widehat{\mu}_{\max}^{\mathrm{sym}}(\overline L)$ is equal to the essential minimum of the height function $h_{\overline L}$. This is a result of Balla\"{y} \cite[Theorem 1.1]{MR4271920}. In this article, we show that these results can  be extended to the case of general adelic curves if we consider  all integral closed subschemes of $X$. More precisely, we obtain the following result (see Theorem \ref{Thm: minoration mu min asy L} and Proposition \ref{Pro: e1 is equal to mu max asy}).

\ifsmf\begin{introtheo}\fi
\ifams\begin{theo}\fi
\label{Thm: mu min asy formula}
Assume that the field $K$ is pefect. Let $X$ be a non-empty reduced projective scheme over $\Spec K$ and $\Theta_X$ be the set of integral closed subschemes of $X$. For any relatively ample adelic line bundle $\overline L$ on $X$, the following equalities hold:
\[\widehat{\mu}_{\min}^{\mathrm{asy}}(\overline L)=\inf_{Y\in\Theta_X}\frac{(\overline L|_Y^{\dim(Y)+1})_S}{(\dim(Y)+1)(L|_Y^{\dim(Y)})}=\inf_{Y\in\Theta_X}\widehat{\mu}_{\max}^{\mathrm{asy}}(\overline L|_Y).\]
Moreover, if $X$ is an integral scheme, the following equality holds:
\[\widehat{\mu}_{\max}^{\mathrm{asy}}(\overline L)=\sup_{\begin{subarray}{c}Y\in\Theta_X\\
Y\neq X
\end{subarray}}\inf_{\begin{subarray}{c}
Z\in\Theta_X\\
Z\not\subseteq Y
\end{subarray}}\widehat{\mu}_{\max}^{\mathrm{asy}}(\overline L|_Y).\]
\ifsmf\end{introtheo}\fi
\ifams\end{theo}\fi

We also show that a property similar to Minkowski's first theorem permits to recover the link between the asymptotic maximal/minimal slopes and the Arakelov height of closed points in the number field case. More precisely, we say that a relatively ample adelic line bundle $\overline L$ is \emph{strongly Minkowskian}\index{Minkowskian!strongly ---} if for any $Y\in\Theta_X$ one has 
\[\lim_{n\rightarrow+\infty}\frac{1}{n}\sup_{\begin{subarray}{c}s\in H^0(Y,L|_Y^{\otimes n})\\
s\neq 0\end{subarray}}\widehat{\deg}(s)\geqslant \frac{(\overline L|_Y^{\dim(Y)+1})_S}{(\dim(Y)+1)(L|_Y^{\dim(Y)})}.\] 
This condition is always satisfied notably when the adelic curve $S$ comes from a number field (consequence of Minkowski's first theorem) or a function field of a projective curve (consequence of Riemann-Roch theorem). We then establish the following result (see Corollary \ref{Cor: mu min equals absolute minimum Minkowskian}).

\ifsmf\begin{introtheo}\fi
\ifams\begin{theo}\fi
Assume that the field $K$ is pefect. Let $X$ be an integral projective scheme over $\Spec K$ and $\overline L$ be a relatively ample adelic line bundle on $X$ which is strongly Minkowskian. Denote by $X^{(0)}$ the set of closed points of $X$. Then the equality $\widehat{\mu}_{\min}^{\mathrm{asy}}(\overline L)=\displaystyle\inf_{x\in X^{(0)}}h_{\overline L}(x)$ holds.
\ifsmf\end{introtheo}\fi
\ifams\end{theo}\fi

Motivated by Theorem \ref{Thm: mu min asy formula}, we propose the following analogue of successive minima for relatively ample adelic line bundles. Let $f:X\rightarrow\Spec K$ be an integral projective scheme of dimension $d$ over $\Spec K$ and $\overline L$ be a relatively ample adelic line bundle on $X$. For $i\in\{1,\ldots,d+1\}$, let  
\[e_i(\overline L)=\sup_{\begin{subarray}{c}\text{$Y\subseteq X$ closed} \\
\operatorname{codim}(Y)\geqslant i\end{subarray}}\;\inf_{\begin{subarray}{c}\text{$Z\in\Theta_X$}\\
\text{$Z\not\subseteq Y$}
\end{subarray}}\operatorname{\widehat{\mu}}_{\max}^{\mathrm{asy}}(\overline L|_{Z})
.\]
With this notation, one can rewrite the assertion of Theorem \ref{Thm: mu min asy formula} as 
\[e_1(\overline L)=\widehat{\mu}_{\max}^{\mathrm{asy}}(\overline L),\quad e_{d+1}(\overline L)=\widehat{\mu}_{\min}^{\mathrm{asy}}(\overline L).\]
We show in Remark \ref{Rem: successive minima} that, in the number field case,  one has
\begin{equation}\label{Equ: eiL}\forall\,i\in\{1,\ldots,d+1\},\quad e_i(\overline L)=\sup_{\begin{subarray}{c}\text{$Y\subseteq X$ closed} \\
\operatorname{codim}(Y)\geqslant i\end{subarray}}\;\inf_{x\in (X\setminus Y)^{(0)}}h_{\overline L}(x).\end{equation}
Thus we recover the definition of successive minima in the sens of \cite[\S5]{MR1189866}.
We propose several fundamental questions about these invariants:
\begin{enumerate}[label=\rm(\arabic*)] 
\item Do the equalities \eqref{Equ: eiL} hold in the case of a general adelic curve, under the assumption that  $\overline L$ is strongly Minkowskian?
\item What is the relation between the invariants $e_2(\overline L),\ldots,e_d(\overline L)$ and the sectional algebra $\bigoplus_{n\in\mathbb N}f_*(\overline L^{\otimes n})$?
\item Does the analogue of some classic results in Diophantine geometry concerning the successive minima, such as the inequality
\[\frac{(\overline L^{d+1})_S}{(L^d)}\geqslant\sum_{i=1}^{d+1}e_i(\overline L),\]
still holds for general adelic curve?
\item In the case where $(X,L)$ is a polarized toric variety and the metrics in $\varphi$ are toric metrics, is it possible to describe in a combinatoric way the positivity conditions of $\overline L$, and express the the invariants $e_i(\overline L)$ in terms of the combinatoric data of $(X,\overline L)$, generalizing some results of \cite{MR3466351,MR3449209} for example?
\end{enumerate}

The last section of the article is devoted to the study of global positivity of adelic line bundles. Motivated by Nakai-Moishezon criterion of ampleness, we say that an adelic line bundle $\overline L$ on $X$ is \emph{ample}\index{ample} if it is relatively ample and if the normalized height with respect to $\overline L$ of integral closed subschemes of $X$ has a positive lower bound. We show that this condition is equivalent to the relative ampleness together with the positivity of the invariant $\widehat{\mu}_{\min}^{\mathrm{asy}}(\overline L)$. Therefore, we deduced from Theorem \ref{Thm: main thm 1} that, if $\overline L_0,\ldots,\overline L_d$ are ample adelic line bundles on $X$, where $d$ is the dimension of $X$, then one has (see Proposition \ref{Pro: positivity intersection number})
\[(\overline L_0\cdots\overline L_d)_S>0.\]
In the case where $\overline L$ is strongly Minkowskian, $\overline L$ is ample if and only if it is relatively ample and the height function $h_{\overline L}$ on the set of closed points of $X$ has a positive lower bound (see Proposition \ref{Pro: strongly minkowskian and ampleness}). Once the ample cone is specified, one can naturally define the nef cone as its closure. It turns out that the nefness can also be described in a numerical way: an adelic line bundle $\overline L$ is nef if and only if it is relatively nef and $\widehat{\mu}_{\min}^{\mathrm{asy}}(\overline L)\geqslant 0$ (see Proposition \ref{Pro: nef with positivity of mu min}). 

Bigness and pseudo-effectivity are also described in a numerical way by the invariant $\widehat{\mu}_{\max}^{\mathrm{asy}}(\ndot)$: an adelic line bundle $\overline L $ is big if and only if $L$ is big and $\widehat{\mu}_{\max}^{\mathrm{asy}}(\overline L)>0$ (which coincides with the bigness in \cite{CMArakelovAdelic}); it is pseudo-effective if and only if $L$ is pseudo-effective and $\widehat{\mu}_{\max}^{\mathrm{asy}}(\overline L)\geqslant 0$ (see \cite[Proposition 6.4.18]{CMArakelovAdelic} and Proposition \ref{Pro: criterion pseudo-effectivity}). We deduce from Theorem \ref{Thm: lower bound intersection number 2} that, if $\overline L_0,\ldots,\overline L_d$ are adelic line bundles on $X$ such that $\overline L_0$ is pseudo-effective and that $\overline L_1,\ldots,\overline L_d$ are nef, then the inequality $(\overline L_0,\ldots,\overline L_d)_S \geqslant 0$ holds (see Proposition \ref{Pro: nef pseudoeffective}).

The rest of the article is organized as follows. In the remaining of Introduction, we remind the the notation that we use all through the article. In the first chapter, we consider metric families on vector bundles and discuss their dominancy and measurability. In the second chapter, we study normed graded linear series over a trivially valued field and prove the limit theorem of their volumes. Then in the third chapter we deduce the limit theorem for graded algebra of adelic vector bundles over a general adelic curve, which proves in particular that the sequence defining the arithmetic volume function actually converges. We also show that the arithmetic Hilbert-Samuel theorem in the original form implies the generalized form with tensor product by an adelic vector bundle. In the fourth chapter, we prove the arithmetic Hilbert-Samuel theorem. We first prove that the difference of the $\chi$-volume and the arithmetic intersection product does not depend on the choice of the metric family. Then we prove  the arithmetic Hilbert-Samuel theorem in the particular case where the adelic curve contains a single copy of the trivial absolute value, and we use the method of casting to the trivial valuation case to prove the arithmetic Hilbert-Samuel theorem in general. The fifth and last chapter is devoted to the study of positivity conditions of adelic line bundles. We begin with a discussion on relative  ampleness and nefness and its relation with sectional arithmetic invariants. We also  and deduce the generalized Hodge index theorem from the arithmetic Hilbert-Samuel theorem. Then we focus on asymptotic maximal slope and its relation with positivity of adelic line bundles. In the end of the chapter, we discuss global positivity conditions. In order to obtain the main theorems of the article in positive characteristic case, we need to generalize some results of  \cite[Chapter~5]{Chen_Moriwaki2020} to any characteristic, which we resume in the appendices.

\section*{Notation and preliminaries}

\renewcommand{\thesubsection}{\arabic{chapter}.\arabic{subsection}}
\subsection{}\label{Sec:adelic curve} Throughout the article, we fix a proper adelic curve \[S=(K,(\Omega,\mathcal A,\nu),\phi),\] where $K$ is a commutative field, $(\Omega,\mathcal A,\nu)$ is a measure space and $\phi=(|\ndot|_{\omega})_{\omega\in\Omega}$ is a family of absolute values on $K$ parametrized by $\Omega$, such that, for any $a \in K^{\times}$, $(\omega \in \Omega) \mapsto \ln |a|_{\omega}$ is integrable on $(\Omega,\mathcal A,\nu)$, and the following ``product formula'' holds:
\[\forall\,a\in K,\quad \int_{\Omega}\ln|a|_\omega\,\nu(\mathrm{d}\omega)=0.\] For any $\omega\in\Omega$, we denote by $K_\omega$ the completion of $K$ with respect to the absolute value $|\ndot|_\omega$. Let $\Omega_\infty$ be the set of $\omega\in\Omega$ such that $|\ndot|_\omega$ is Archimedean. Let $\Omega_{\operatorname{fin}}=\Omega\setminus\Omega_\infty$ be the set of $\omega\in\Omega$ such that $|\ndot|_\omega$ is non-Archimedean.

We assume that, either 
the $\sigma$-algebra $\mathcal A$ is discrete, or the field $K$ is countable. Moreover, we denote by $\Omega_\infty$ the set of $\omega\in\Omega$ such that $|\ndot|_{\omega}$ is Archimedean. Note that $\nu(\Omega_\infty)<+\infty$.
For $\omega \in \Omega_{\infty}$, we always assume that $|a|_{\omega} = a$ for any $a \in \QQ_{\geqslant 0}$. Denote by $\Omega_{\operatorname{fin}}$ the set $\Omega\setminus\Omega_\infty$.

In the last chapter of the article, we assume in addition that $K$ is perfect, that is, either $K$ is of characteristic $0$, or $K$ is of characteristic $p>0$ and the Frobenius map $(a\in K)\mapsto a^p$ is surjective.

\subsection{}\label{Sec: distance of metric families}
Let $V$ be a finite-dimensional vector space over $K$. As \emph{norm family}\index{norm family} on $V$, we refer to a family $(\|\ndot\|_{\omega})_{\omega\in\Omega}$, where $\|\ndot\|_\omega$ is a norm on $V_\omega:=V\otimes_KK_\omega$. 

Let $\xi=(\|\ndot\|_{\omega})_{\omega\in\Omega}$ and $\xi'=(\|\ndot\|'_\omega)_{\omega\in\Omega}$ be norm families on $V$. For any $\omega\in\Omega$, we denote by $d_\omega(\xi,\xi')$ the following number
\[\sup_{s\in V\setminus\{0\}}\Big|\ln\|s\|_\omega-\ln\|s\|_{\omega}'\Big|.\]
In the case where $V=\boldsymbol{0}$, by convention $d_\omega(\xi,\xi')=0$.

\subsection{}\label{Subsec: adelic vector bundle} As \emph{adelic vector bundle}\index{adelic vector bundle} on $S$, we refer to the data $\overline V=(V,\xi)$ which consists of a finite-dimensional vector space $V$ over $K$ and a family of norms $\xi=(\|\ndot\|_\omega)_{\omega\in\Omega}$ on $V_\omega:=V\otimes_KK_\omega$, satisfying the following conditions:
\begin{enumerate}[label=\rm(\arabic*)]
\item the norm family $\xi$ is \emph{strongly dominated}\index{strongly dominated}, that is, there exist an integrable function $C:\Omega\rightarrow\mathbb R_{\geqslant 0}$ and a basis $(e_i)_{i=1}^r$ of $V$ over $K$, such that, for any $\omega\in\Omega$ and any $(\lambda_1,\ldots,\lambda_r)\in K_\omega^r\setminus\{(0,\ldots,0)\}$,
\[\Big|\ln\|\lambda_1e_1+\cdots+\lambda_re_r\|_\omega-\ln\max_{i\in\{1,\ldots,r\}}|\lambda_i|_\omega\Big|\leqslant C(\omega).\]
\item the norm family $\xi$ is \emph{measurable}\index{measurable}, that is, for any $s\in V$, the function $(\omega\in\Omega)\mapsto \|s\|_{\omega}$ is $\mathcal A$-measurable.
\end{enumerate}
In the article, we only consider adelic vector bundles which are ultrametric over non-Archimedean places, namely we assume that the norm $\|\ndot\|_\omega$ is ultrametric once the absolute value $|\ndot|_{\omega}$ is non-Archimedean. If in addition the norm $\|\ndot\|_{\omega}$ is induced by an inner product whenever $|\ndot|_\omega$ is Archimedean, we say that $\overline V$ is \emph{Hermitian}\index{Hermitian}. If $\dim_K(V)=1$, we say that $\overline V$ is an \emph{adelic line bundle}\index{adelic line bundle} (note that an adelic line bundle is necessarily Hermitian).

If $\overline V$ is an adelic vector bundle on $S$, any vector subspace (resp. quotient vector space) of $V$ together with the family of restricted norms (resp. quotient norms) forms also an adelic vector bundle on $S$, which is called an \emph{adelic vector subbundle}\index{adelic vector subbundle} (resp. \emph{quotient adelic vector bundle}\index{adelic vector bundle!quotient ---}) of $\overline V$. Note that if $\overline V$ is Hermitian, then all its adelic vector subbundles and quotient adelic vector bundles are Hermitian.

\subsection{}\label{Subsec: deg and deg+} Let $\overline V=(V,(\|\ndot\|_{\omega})_{\omega\in\Omega})$ be an adelic vector bundle on $S$, we define the \emph{Arakelov degree}\index{Arakelov degree} of $\overline V$ as
\[\widehat{\deg}(\overline V):= -\int_{\Omega}\ln\|e_1\wedge\cdots\wedge e_r\|_{\omega,\det}\,\nu(\mathrm{d}\omega),\]
where $(e_i)_{i=1}^r$ is a basis of $V$ over $K$, and $\|\ndot\|_{\omega,\det}$ denotes the determinant norm of $\|\ndot\|_\omega$, which is defined as (where $r=\dim_K(V)$)
\[\forall\,\eta\in\det(V)=\Lambda^r(V),\quad \|\eta\|_{\omega,\det}=\inf_{\eta=s_1\wedge\cdots\wedge s_r}\|s_1\|\cdots\|s_r\|.\]
Let $\widehat{\deg}_+(\overline V)$ be the \emph{positive degree}\index{positive degree} of $\overline V$, which is defined as 
\[\widehat{\deg}_+(\overline V)=\sup_{W\subseteq V}\widehat{\deg}(\overline W),\]
where $W$ runs over the set of vector subspaces of $V$, and in the adelic vector bundle structure of $\overline W$ we consider the restricted norms. 
In the case where $V$ is non-zero, we denote by $\widehat{\mu}(\overline V)$ the quotient $\widehat{\deg}(\overline V)/\dim_K(V)$, called the \emph{slope}\index{slope} of $V$. We define the \emph{minimal slope}\index{slope! minimal ---} of $\overline V$ as 
\[\widehat{\mu}_{\min}(\overline V):=\inf_{V\twoheadrightarrow W\neq\{0\}}\widehat{\mu}(\overline W),\]
where $\overline W$ runs over the set of all non-zero quotient adelic vector bundles of $\overline{V}$. 

\subsection{} Let $\overline E$ and $\overline F$ be two adelic vector bundles on $S$ and $\varphi:E\rightarrow F$ be a $K$-linear map. We define the \emph{height}\index{height} of $\varphi$ as 
\[h(\varphi):=\int_{\Omega}\ln\|\varphi\|_\omega\,\nu(\mathrm{d}\omega),\]
where $\|\varphi\|_{\omega}$ denotes the operator norm of the $K_\omega$-linear map $E_\omega\rightarrow F_\omega$ induced by $\varphi$. Moreover, if $E$ is non-zero and if $\varphi$ is injective, then the following slope inequality holds (see \cite[Proposition 4.3.31]{CMArakelovAdelic}):
\[\widehat{\mu}_{\max}(\overline E)\leqslant\widehat{\mu}_{\max}(\overline F)+h(\varphi).\]

\subsection{}\label{Sec: HN R filtration}Let $\overline V$ be a non-zero adelic vector bundle on $S$. For any $t\in\mathbb R$, we let 
\[\mathcal F^t(\overline V)=\sum_{\begin{subarray}{c}\{0\}\neq W\subseteq V\\
\widehat{\mu}_{\min}(\overline W)\geqslant t
\end{subarray}}W,\]
where $W$ runs over the set of all non-zero vector subspaces of $V$ such that the minimal slope of $W$ equipped with the family of restricted norms is $\geqslant t$. We call $(\mathcal F^t(\overline V))_{t\in\mathbb R}$ the \emph{Harder-Narasimhan $\mathbb R$-filtration}\index{Harder-Narasimhan $\mathbb R$-filtration} of $\overline V$. In the case where $\overline V$ is Hermitian, the following equality holds (see \cite[Theorem 4.3.44]{CMArakelovAdelic}):
\begin{gather*}\widehat{\deg}(\overline V)=-\int_{\mathbb R}t\,\mathrm{d}(\dim_K(\mathcal F^t(\overline V))),\\
\widehat{\deg}_+(\overline V)=-\int_0^{+\infty}t\,\mathrm{d}(\dim_K(\mathcal F^t(\overline V)))=\int_{0}^{+\infty}\dim_K(\mathcal F^t(\overline V))\,\mathrm{d}t.
\end{gather*} 
In general one has (see \cite[Propositions 4.3.50 and 4.3.51, and Corollary 4.3.52]{CMArakelovAdelic})
\begin{gather*}0\leqslant\widehat{\deg}(\overline V)+\int_{\mathbb R}t\,\mathrm{d}(\dim_K(\mathcal F^t(\overline V)))\leqslant\frac 12\nu(\Omega_\infty)\dim_K(V)\ln(\dim_K(V)),\\
0\leqslant\widehat{\deg}_+(\overline V)-\int_0^{+\infty}\dim_K(\mathcal F^t(\overline V))\,\mathrm{d}t\leqslant\frac 12\nu(\Omega_\infty)\dim_K(V)\ln(\dim_K(V)).
\end{gather*}

\subsection{} Let $\overline V=(V,(\|\ndot\|_{V,\omega})_{\omega\in\Omega})$ and $\overline W=(W,(\|\ndot\|_{W,\omega})_{\omega\in\Omega})$ be adelic vector bundles on $S$. For any $\omega\in\Omega$ such that $|\ndot|_\omega$ is non-Archimedean, let $\|\ndot\|_{\omega}$ be the $\varepsilon$-tensor product on $V_\omega\otimes_{K_\omega}W_\omega$, of the norms $\|\ndot\|_{V,\omega}$ and $\|\ndot\|_{W,\omega}$. Note that, for any $T\in V_\omega\otimes_{K_\omega}W_\omega$, the value of $\|T\|_\omega$ is  equal to
\[ \min\Big\{\max_{i\in\{1,\ldots,n\}}\|e_i\|_{V,\omega}\|f_i\|_{W,\omega}\, : \,\begin{array}{l}n\in\mathbb N,\; (e_i)_{i=1}^n\in V_\omega^n,\;(f_i)_{i=1}^n\in
W_\omega^n\\ T=e_1\otimes f_1+\cdots+e_n\otimes f_n
\end{array}
\Big\}.\]
In the case where $|\ndot|_\omega$ is Archimedean, let $\|\ndot\|_\omega$ be $\pi$-tensor product of $\|\ndot\|_{V,\omega}$ of $\|\ndot\|_{W,\omega}$. Recall that for any $T\in V_\omega\otimes_{K_\omega}W_{\omega}$, the value of $\|T\|_\omega$ is  equal to
\[\min\Big\{\sum_{i=1}^n\|e_i\|_{V,\omega}\|f_i\|_{W,\omega}\,:\,\begin{array}{l}n\in\mathbb N,\; (e_i)_{i=1}^n\in V_\omega^n,\;(f_i)_{i=1}^n\in W_\omega^n\\
T=e_1\otimes f_1+\cdots+e_n\otimes f_n
\end{array}
\Big\}.\] 
The pair \[\overline V\otimes_{\varepsilon,\pi}\overline W=(V\otimes_KW,(\|\ndot\|_\omega)_{\omega\in\Omega})\]
is called the \emph{$\varepsilon,\pi$-tensor product}\index{tensor product!$\varepsilon,\pi$- ---}  of $\overline V$ and $\overline W$. 

Assume that $\overline V$ and $\overline W$ are Hermitian. If $|\ndot|_\omega$ is non-Archimedean, let $\|\ndot\|_\omega^{H}$ be the $\varepsilon$-tensor product of $\|\ndot\|_{V,\omega}$ and $\|\ndot\|_{W,\omega}$; otherwise let $\|\ndot\|_\omega^H$ be the orthogonal tensor product of the Euclidean or Hermitian norms $\|\ndot\|_{V,\omega}$ and $\|\ndot\|_{W,\omega}$. Then the pair
\[\overline V\otimes\overline W=(V\otimes_KW,(\|\ndot\|_\omega^H)_{\omega\in\Omega})\]
is called the \emph{Hermitian tensor product}\index{tensor product!Hermitian ---} of $\overline V$ and $\overline W$.

\subsection{}\label{Subsec: metric}
Let $(k,|\ndot|)$ be a field equipped with a complete absolute value, $X$ be a projective scheme
over $\Spec k$. We denote by $X^{\mathrm{an}}$ the analytic space associated with $X$  (in the sense of Berkovich if $|\ndot|$ is non-Archimedean). Recall that a point $x$ of $X^{\mathrm{an}}$ is of the form $(j(x),|\ndot|_x)$, where $j(x)$ is a scheme point of $X$, $|\ndot|_x$ is an absolute value on the residue field of $j(x)$,  which extends the absolute value $|\ndot|$ on the base field $k$. We denote by $\widehat{\kappa}(x)$ the completion of the residue field of $j(x)$ with respect to the absolute value $|\ndot|_x$, on which $|\ndot|_x$ extends by continuity. The set $X^{\operatorname{an}}$ is equipped with the most coarse topology which makes continuous the map $j:X^{\operatorname{an}}\rightarrow X$ and all functions of the form 
\[|s|:U^{\operatorname{an}}\longrightarrow\mathbb R_{\geqslant 0},\quad x\longmapsto |s(x)|_x,\]
where $U$ is a non-empty Zariski open subset of $X$ and $s\in\mathcal O_X(U)$ is a regular function on $U$. In particular, if $U$ is a Zariski open subset of $X$, then $U^{\operatorname{an}}$ is an open subset of $X^{\mathrm{an}}$. We call such open subsets of $X^{\mathrm{an}}$ \emph{Zariski open subsets}\index{Zariski open subset}.

\subsection{} Let $\pi:X\rightarrow\Spec K$ be a  projective scheme over $\Spec K$. For any $\omega\in\Omega$, let $X_\omega$ be $X\times_{\Spec K}\Spec K_\omega$ and let $X_\omega^{\mathrm{an}}$ be the analytic space associated with $X_\omega$. If $L$ is an invertible $\mathcal O_X$-module, we call \emph{metric family}\index{metric family} on $L$ any family $\varphi=(\varphi_\omega)_{\omega\in\Omega}$, where $\varphi_\omega$  is a continuous metric on $L_\omega=L|_{X_\omega}$. In the particular case where $X$ is the spectrum of a finite extension $K'$ of $K$, the invertible $\mathcal O_X$-module $L$ is just a one-dimension vector space over $K'$ and a metric family of $L$  could be viewed as norm family if we consider the adelic curve $S\otimes_KK'$ obtained by algebraic extension of scalars (see \cite[\S3.4]{CMArakelovAdelic}). 

If $\overline E=(E,(\|\ndot\|_\omega)_{\omega\in\Omega})$ is a finite-dimensional $K$-vector space $E$ equipped with a norm family, $g:X\rightarrow\mathbb P(E)$ is a projective $K$-morphism and $L=g^*(\mathcal O_E(1))$, then, for each $\omega\in\Omega$, the norm $\|\ndot\|_\omega$ induces by passing to quotient by the universal surjective homomorphism \[(g_\omega\circ \pi_\omega)^*(E_\omega)\longrightarrow g_\omega^*(\mathcal O_{E_\omega}(1))=L_\omega\] a continuous metric $\varphi_\omega$ on $L_\omega$. The metric family $(\varphi_\omega)_{\omega\in\Omega}$ is called a \emph{quotient metric family}\index{metric family!quotient ---} induced by $(\|\ndot\|_\omega)_{\omega\in\Omega}$ (and by $g$).  

Let $L$ be an invertible $\mathcal O_X$-module and $\varphi=(\varphi_\omega)_{\omega\in\Omega}$ be a metric family of $L$. For any $\omega\in\Omega$, the metric $\varphi_\omega$ induces by passing to dual a metric on $L_\omega$, which we denote by $-\varphi_\omega$. The metric family $(-\varphi_\omega)_{\omega\in\Omega}$ on $L^\vee$ is denoted by $-\varphi$.

Let $L_1,\ldots,L_n$ be invertible $\mathcal O_X$-modules. For any $i\in\{1,\ldots,n\}$, let $\varphi_i$ be a metric family on $L_i$. Then the metric families $\varphi_1,\ldots,\varphi_n$ induces by tensor product a metric family on $L_1\otimes\cdots\otimes L_n$, which we denote by $\varphi_1+\cdots+\varphi_n$ in the additive form. In particular, if all $(L_i,\varphi_i)$ are equal to the same $(L,\varphi)$, the metric family $\varphi+\cdots+\varphi$ is denoted by $n\varphi$. 

Let $Y$ be a projective $K$-scheme and $f:Y\rightarrow X$ be a $K$-morphism. If $L$ is an invertible $\mathcal O_X$-module and $\varphi$ is a metric family of $L$, then $\varphi$ induces by pullback a metric family $f^*(\varphi)$ on $f^*(L)$: for any $\omega\in\Omega$ and any $y\in Y_\omega^{\mathrm{an}}$, the norm $|\ndot|_{f^*(\varphi)_\omega}(y)$ is induced by $|\ndot|_{\varphi_\omega}(f_\omega^{\mathrm{an}}(y))$ by extension of scalars.

\subsection{}\label{Sec: dominancy line bundle} Let $X$ be a projective $K$-scheme, $L$ be an invertible $\mathcal O_X$-module and $\varphi$ be a metric family of $L$. Assume that there exist invertible $\mathcal O_X$-modules $L_1$ and $L_2$, together with quotient metric families $\varphi_1$ and $\varphi_2$ on $L_1$ and $L_2$ respectively, which are induced by strongly dominated norm families (see \S\ref{Subsec: adelic vector bundle}), such that $L\cong L_1\otimes L_2^\vee$ and that $\varphi=\varphi_1-\varphi_2$, we say that the metric family $\varphi$ is \emph{dominated}\index{dominated}. We refer to \cite[\S6.1.2]{CMArakelovAdelic} for more details.

\subsection{}\label{Sub: exponent} Let $\Omega_0$ be the set of $\omega\in\Omega$ such that the absolute value $|\ndot|_\omega$ is trivial. Let $X$ be a projective scheme over $\Spec K$. For any triplet $x=(K_x,|\ndot|_x,P_x)$, where $(K_x,|\ndot|_x)$ is a valued extension of the trivially valued field $(K,|\ndot|_0)$ and $P_x:\Spec K_x\rightarrow X$ is a $K$-morphism, we denote by $S_x$ the adelic curve
\[(K_x,(\Omega_0,\mathcal A_0,\nu_0),(|\ndot|_x)_{\omega\in\Omega_0}),\]
where $\mathcal A_0=\mathcal A|_{\Omega_0}$ and $\nu_0$ is the restriction of $\nu$ to $(\Omega_0,\mathcal A_0)$. If $L$ is an invertible $\mathcal O_X$-module and if $\varphi$ is a metric family of $L$, we denote by $L_x$ the pullback $P_x^*(L)$ and by $x^*(\varphi)$ the norm family $(|\ndot|_{\varphi_\omega}(P_x^\omega))_{\omega\in\Omega_0}$ on $L_x$, where $P_x^\omega$ denotes the point of $X_\omega^{\mathrm{an}}$ determined by $(P_x,|\ndot|_x)$.

Assume that the transcendence degree of $K_x/K$ is $\leqslant 1$. Then $|\ndot|_x$ is a discrete absolute value on $K_x$. Let $\operatorname{ord}_x(\ndot):K_x\rightarrow\mathbb Z\cup\{+\infty\}$ be the corresponding discrete valuation, which is defined as 
\[\operatorname{ord}_x(a)=\sup\{n\in\mathbb Z\,:\,a\in\mathfrak m_x^n\},\]
where \[\mathfrak m_x=\{b\in K_x\,:\,|b|_x< 1\}.\] 
Then there is a non-negative real number $q$ such that 
\[|\ndot|_x=\exp(-q\operatorname{ord}_x(\ndot)).\]
We call it the \emph{exponent}\index{exponent} of $x$. 

\subsection{}\label{Sec: measurability line bundle}Let $X$ be a projective $K$-scheme, $L$ be an invertible $\mathcal O_X$-module, and $\varphi$ be a metric family of $L$. We say that the metric family $\varphi$ is \emph{measurable}\index{measurable} if the following conditions are satisfied (see \cite[\S6.1.4]{CMArakelovAdelic} for more details):
\begin{enumerate}[label=\rm(\arabic*)]
\item for any finite extension $K'/K$ and any $K$-morphism $P:\Spec K'\rightarrow X$, the norm family $P^*(\varphi)$ is measurable,
\item for any triplet $x=(K_x,|\ndot|_x,P_x)$, where $(K_x,|\ndot|_x)$ is a valued extension of transcendence degree $\leqslant 1$ and of rational exponent of the trivially valued field $(K,|\ndot|_0)$, and $P_x:\Spec K_x\rightarrow X$ is a $K$-morphism, the norm family $x^*(\varphi)$ is measurable.
\end{enumerate}

\subsection{} Let $X$ be a projective scheme over $\Spec K$, $L$ be an invertible $\mathcal O_X$-module and $\varphi$ be a metric family of $L$. We say that $\overline L=(L,\varphi)$ is an \emph{adelic line bundle}\index{adelic line bundle} on $X$ if the metric family $\varphi$ is dominated and measurable (see \S\ref{Sec: dominancy line bundle} and \S\ref{Sec: measurability line bundle}). 

Suppose that $X$ is geometrically reduced. Let $\overline L=(L,\varphi)$ be an adelic line bundle on $X$. We denote by $f_*(\overline L)$ the couple $(H^0(X,L),(\|\ndot\|_{\varphi_\omega})_{\omega\in\Omega})$, where for $s\in H^0(X_\omega,L_\omega)$,
\[\|s\|_{\varphi_\omega}=\sup_{x\in X_\omega^{\mathrm{an}}}|s|_{\varphi_\omega}(x).\]
It turns out that $f_*(\overline L)$ is an adelic vector bundle on $S$ (see \cite[Theorems 6.1.13 and 6.1.32]{CMArakelovAdelic}). 

\subsection{} Let $X$ be a projective scheme over $\Spec K$. 
Let $L$ be an invertible $\mathcal O_X$-module, $\varphi=(\varphi_\omega)_{\omega\in\Omega}$ and $\psi=(\psi_\omega)_{\omega\in\Omega}$ be metric families on $L$ such that $(L,\varphi)$ and $(L,\psi)$ are both adelic line bundles. Then we define the \emph{distance}\index{distance} between $\varphi$ and $\psi$ as 
\[d(\varphi,\psi):=\int_{\Omega}\sup_{x\in X_\omega}\mathopen{\bigg|}\ln\frac{|\ndot|_{\varphi_\omega}(x)}{|\ndot|_{\psi_\omega}(x)}\mathclose{\bigg|}\,\nu(\mathrm{d}\omega).\]
If $L$ is semiample and if there exist a positive integer $m$ and a sequence $(\varphi_n)_{n\in\mathbb N}$ of quotient metric families (where $\varphi_n$ is a metric family of $L^{\otimes nm}$), such that
\[\lim_{n\rightarrow+\infty}\frac{1}{mn}d(nm\varphi,\varphi_n)=0,\]
we say that the metric family $\varphi$ is \emph{semi-positive}\index{semi-positive}.

\subsection{} Let $X$ be a projective scheme over $\Spec K$ and $d$ be the dimension of $X$. We denote by $\operatorname{\widehat{Int}}(X)$ the set of adelic line bundles $\overline L$ on $X$ which can be written in the form $\overline A_1\otimes\overline A_2^\vee$, where each $\overline A_i$ is an ample invertible $\mathcal O_X$-module equipped with a semi-positive metric family. In \cite{CMIntersection}, we have constructed an \emph{arithmetic intersection product}\index{arithmetic intersection product}
\[\big((\overline L_0,\ldots,\overline L_d)\in\operatorname{\widehat{Int}}(X)^{d+1}\big)\longmapsto(\overline L_0\cdots\overline L_d)_S\in\mathbb R,\]
which is multi-linear with respect to tensor product. We have also related the arithmetic intersection number $(\overline L_0\cdots\overline L_d)_S$ to the height of the multi-resultant of $L_0,\ldots,L_d$.

%!TEX root = ./Hilbert_Samuel_Adelic_Curves.tex

\ifams
\renewcommand{\thetheo}{\arabic{chapter}.\arabic{section}.\arabic{theo}}
\renewcommand{\thesection}{\arabic{chapter}.\arabic{section}}
\fi

\chapter{Metric families on vector bundles}\label{sec:Family of metrics of locally free module}

The purpose of this chapter is to generalize dominancy and measurability conditions in \cite[Chapter~6]{CMArakelovAdelic} to metrized locally free modules of finite rank, and
to develop related topics. Let $S=(K,(\Omega,\mathcal A,\nu),\phi)$ be an adelic curve as introduced in \S\ref{Sec:adelic curve}.   

\section{Metric family}

Let $p:X\rightarrow\Spec K$ be a quasi-projective 
scheme over $\Spec K$.
Let $E$ be a vector bundle on $X$, that is, a locally free $\OO_X$-module of finite rank. For any $\omega\in\Omega$, let $E_\omega$ be the restriction of $E$ to $X_\omega=X\times_{\Spec K}\Spec K_\omega$ and $\psi_\omega$ be a metric on $E_\omega$. By definition $\psi_\omega$ is a family $(|\ndot|_{\psi_\omega}(x))_{x\in X_\omega^{\mathrm{an}}}$ parametrized by $X_{\omega}^{\mathrm{an}}$, where each $|\ndot|_{\psi_\omega}(x)$ is a norm on $E_\omega(x):=E_\omega\otimes_{\mathcal O_{X_\omega}}\widehat{\kappa}(x)$. We assume that the norm $|\ndot|_{\psi_\omega}(x)$ is ultrametric if the absolute value $|\ndot|_\omega$ is non-Archimedean. Moreover, we assume that the metric $\psi_\omega$ is continuous, namely, for any section $s$ of $E$ over a Zariski open subset $U$ of $X_\omega$, the function 
\[(x\in U^{\mathrm{an}})\longmapsto |s|_{\psi_\omega}(x)\]
is continuous.  The data $\psi=(\psi_\omega)_{\omega\in\Omega}$ is called a \emph{metric family}\index{metric family} on the vector bundle $E$.  

Assume that $X$ is projective and geometrically reduced. For any $\omega\in\Omega$, we denote by $\|\ndot\|_{\psi_\omega}$ the supremum norm on $H^0(X_\omega,E_\omega)$, which is defined as
\[\forall\,s\in H^0(X_\omega,E_\omega),\quad \|s\|_{\psi_\omega}=\sup_{x\in X_\omega^{\mathrm{an}}}|s(x)|_{\psi_\omega}(x).\]
We denote by $p_*(E,\psi)$ the couple $(H^0(X,E),(\|\ndot\|_{\psi_\omega})_{\omega\in\Omega})$.

If $\varphi$ and $\psi$ are two metric families of  $E$. For any $\omega\in\Omega$ we denote by $d_\omega(\varphi,\psi)$ the element
\[\sup_{x\in X_\omega^{\mathrm{an}}}\sup_{s\in E_\omega(x)\setminus\{0\}}\bigg|\ln|s|_{\varphi_\omega}(x)-\ln|s|_{\psi_\omega}(x)\bigg|\in [0,+\infty],\]
which is called the \emph{local distance}\index{distance!local ---} at $\omega$ between $\varphi$ and $\psi$.

We denote by $\mathcal O_E(1)$ the universal invertible sheaf on the projective bundle $\pi:\mathbb P(E)\rightarrow\Spec K$. For any $\omega\in\Omega$, the metric $\psi_\omega$ induces by passing to quotient a continuous metric on $\mathcal O_E(1)_\omega\cong\mathcal O_{E_\omega}(1)$, which we denote by $\psi_\omega^{\operatorname{FS}}$. Recall that, if $y$ is an element of $\mathbb P(E_\omega)^{\mathrm{an}}$ and $x=\pi_\omega^{\mathrm{an}}(y)$, then the norm $|\ndot|_{\psi_\omega^{\operatorname{FS}}}$ on $\mathcal O_E(1)_y$ is the quotient metric induced by the universal surjective homomorphism 
\[E_{\omega}(x)\otimes_{\widehat{\kappa}(x)}\widehat{\kappa}(y)\longrightarrow \mathcal O_E(1)_y, \]
where we consider the $\varepsilon$-extension of $|\ndot|_{\psi_\omega}(x)$ to $E_{\omega}(x)\otimes_{\widehat{\kappa}(x)}\widehat{\kappa}(y)$ if $|\ndot|_\omega$ is non-Archimedean, and $\pi$-extension of $|\ndot|_{\psi_\omega}(x)$ if $|\ndot|_{\omega}$ is Archimedean (see \cite[\S1.3 and \S2.2.3]{CMArakelovAdelic}). Note that, if $\varphi$ and $\psi$ are two metric families of $E$, then one has (see \cite[Proposition 2.2.20]{CMArakelovAdelic})
\begin{equation}\label{Equ: distance FS bounded by distance}\forall\,\omega\in\Omega,\quad d_\omega(\varphi^{\operatorname{FS}},\psi^{\operatorname{FS}})\leqslant d_\omega(\varphi,\psi).\end{equation}

\section{Dominancy and measurability} In this section, we fix a projective scheme $X$ over $\Spec K$.

\begin{defi}
Let $E$ be a vector bundle on $X$.
\begin{enumerate}[label=\rm(\arabic*)]
\item We say the metric family $\psi=(\psi_\omega)_{\omega\in\Omega}$ on the locally free $\mathcal O_X$-module $E$ is \emph{dominated}\index{dominated} (resp. \emph{measurable}\index{measurable}) if the metric family $\psi^{\operatorname{FS}}=(\psi_\omega^{\operatorname{FS}})_{\omega\in\Omega}$ on $\mathcal O_E(1)$ is dominated (resp. measurable). We refer the readers to \cite[Definitions 6.1.9 and 6.1.27]{CMArakelovAdelic} for the dominancy and measurability conditions of metrized line bundles. 

\item We say $(E, \psi)$ is an \emph{adelic locally free $\OO_X$-module}\index{adelic locally free $\OO_X$-module} or an \emph{adelic vector bundle}\index{adelic vector bundle} if $\psi$ is dominated and measurable, or equivalently, $(\mathcal O_E(1),\psi^{\mathrm{FS}})$ is an adelic line bundle on $\mathbb P(E)$.
\end{enumerate}
\end{defi}

\begin{prop}\label{Pro: direct image}
\begin{enumerate}[label=\rm(\arabic*)]
\item If $\psi$ is dominated, then the norm family $\xi_{\psi}=( \|\ndot\|_{\psi_{\omega}})_{\omega \in \Omega}$ on $H^0(X, E)$ is strongly dominated.

\item If the metric family $\psi$ on $E$ is measurable, then the norm family $\xi_\psi$ on $H^0(X, E)$ is measurable.
\end{enumerate}
\end{prop}

\begin{proof}
If we identify $H^0(X,E)$ with $H^0(\mathbb P(E),\mathcal O_E(1))$, then for any $\omega\in\Omega$ one has $\|\ndot\|_{\psi_\omega}=\|\ndot\|_{\psi_\omega^{\operatorname{FS}}}$ by \cite[Remark~2.2.14]{CMArakelovAdelic}. Therefore the assertions follow from 
\cite[Theorems~6.1.13 and 6.1.32]{CMArakelovAdelic}.
\end{proof}

\begin{prop}\label{Pro: criterion of dominancy}
Let $E$ be a vector bundle on $X$, and $\varphi$ and $\psi$ be two metric families of $E$. Suppose that $\varphi$ is dominated and that the local distance function
\[(\omega\in\Omega)\longmapsto d_\omega(\varphi,\psi)\]
is bounded from above by an integrable function. Then the metric family $\psi$ is also dominated.
\end{prop}
\begin{proof}
This is a consequence of \cite[Proposition 6.1.12]{CMArakelovAdelic} and \eqref{Equ: distance FS bounded by distance}.
\end{proof}

\begin{defi}\label{Def: pull-back}
Let $f:Y\rightarrow X$ be a projective $K$-morphism from a geometrically reduced projective $K$-scheme $Y$ to $X$. Let $E$ be a vector bundle on $X$ and $\psi=(\psi_{\omega})_{\omega\in\Omega}$ be a metric family on $E$. We denote by $f^*(\psi)$ the metric family on $E$ such that, for any $y\in Y_\omega^{\mathrm{an}}$, the norm $|\ndot|_{f^*(\psi)_\omega}(y)$ on \[f^*(E)_\omega(y)=E_\omega(x)\otimes_{\widehat{\kappa}(x)}\widehat{\kappa}(y)\]
is induced by $|\ndot|_{\psi_\omega}(f^{\mathrm{an}}(y))$  by $\varepsilon$-extension of scalars in the case where $|\ndot|_\omega$ is non-Archimedean, and by  $\pi$-extension of scalars if $|\ndot|_{\omega}$ is Archimedean. 
\end{defi}

\begin{prop}
We keep the notation and the assumptions of Definition \ref{Def: pull-back}. Suppose that the metric family $\psi$ is dominated \textup{(}resp. measurable\textup{)}, then its pull-back $f^*(\psi)$ is also dominated \textup{(}resp. measurable\textup{)}.
\end{prop}
\begin{proof}
The universal property of projective bundle induces a projective morphism $F:\mathbb P(f^*(E))\rightarrow\mathbb P(E)$ such that the following diagramme is cartesian. 
\[\xymatrix{\mathbb P(f^*(E))\ar[r]^-F\ar[d]_-{\pi_{f^*(E)}}&\mathbb P(E)\ar[d]^-{\pi_E}\\
Y\ar[r]_-{f}&X}\]
Moreover, one has $\mathcal O_{f^*(E)}(1)\cong F^*(\mathcal O_E(1))$ and $F^*(\psi^{\mathrm{FS}})=f^*(\psi)^{\mathrm{FS}}$. Hence the assertion follows from \cite[Propositions 6.1.12 and 6.1.28]{CMArakelovAdelic}.
\end{proof}

\begin{defi}
Let $E$ be a vector bundle on $X$ and $\psi=(\psi_{\omega})_{\omega\in\Omega}$ be a metric family of $E$. If $F$ is a vector subbundle of $E$, for any $\omega\in\Omega$ and any $x\in X_\omega^{\mathrm{an}}$, we denote by $|\ndot|_{\psi_{F,\omega}}(x)$ the restriction of $|\ndot|_{\psi_\omega}(x)$ to $F_\omega(x)$. Note that $\psi_F=(\psi_{F,\omega})_{\omega\in\Omega}$ forms a metric family of $F$, called the \emph{restriction of $\psi$ to $F$}\index{restriction}. Similarly,  if $G$ is a quotient vector bundle of $E$, we denote by $|\ndot|_{\psi_G,\omega}(x)$ the quotient norm of $|\ndot|_{\psi_\omega}(x)$ on $G_\omega(x)$. Then $\psi_G=(\psi_{G,\omega})_{\omega\in\Omega}$ is a metric family of $G$, called the \emph{quotient metric family}\index{metric family!quotient ---} of $\psi$ on $G$.
\end{defi}

\begin{prop}\label{Pro: quotient family dominated}
Let $E$ be a vector bundle on $X$ and $G$ be a quotient vector bundle of $E$. Let $\psi$ be a metric family on $E$. If $\psi$ is dominated \textup{(}resp. measurable\textup{)}, then the quotient metric family $\psi_G$ is also dominated \textup{(}resp. measurable\textup{)}.
\end{prop}
\begin{proof}
Let $i:\mathbb P(G)\rightarrow\mathbb P(E)$ be the closed embedding induced by the quotient homomorphism $E\rightarrow G$. Then one has $i^*(\psi^{\mathrm{FS}})=\psi_G^{\mathrm{FS}}$. Hence the assertion of the proposition follows from \cite[Propositions 6.1.12 and 6.1.28]{CMArakelovAdelic}.
\end{proof}

\begin{defi}
Let $E$ and $F$ be vector bundles on $X$, equipped with metric families $\psi_E$ and $\psi_F$, respectively. For any $\omega\in\Omega$ and any $x\in X_\omega^{\mathrm{an}}$, if $|\ndot|_{\omega}$ is non-Archimedean, we denote by $|\ndot|_{(\psi_E\otimes\psi_F)_\omega}(x)$ the $\varepsilon$-tensor product of the norms $|\ndot|_{\psi_{E,\omega}}(x)$ and $|\ndot|_{\psi_{F,\omega}}(x)$, if $|\ndot|_{\omega}$ is Archimedean, we denote by $|\ndot|_{(\psi_E\otimes\psi_F)_\omega}(x)$ the $\pi$-tensor product of the norms $|\ndot|_{\psi_{E,\omega}}(x)$ and $|\ndot|_{\psi_{F,\omega}}(x)$. Thus we obtain a metric family $\psi_E\otimes\psi_F$ on the  vector bundle $E\otimes F$, called the tensor product of metric families $\psi_E$ and $\psi_F$. In the case where one of the vector bundles $E$ and $F$ is of rank $1$, we also write the tensor product metric family of $\psi_E$ and $\psi_F$ in an additive way as $\psi_E+\psi_F$.
\end{defi}

\begin{prop}
\label{prop:tensor:line:bundle:dominated:measurable}
Let $E$ and $F$ be vector bundles on $X$, equipped with metric families $\psi_E$ and $\psi_F$ respectively. We assume that $E$ is a line bundle. If both metric families $\psi_E$ and $\psi_F$ are dominated \textup{(}resp. measurable\textup{)}, then $\psi_E+\psi_F$ is also dominated \textup{(}resp. measurable\textup{)}.
\end{prop}
\begin{proof}
Since $E$ is of rank $1$, we can identify $\mathbb P(E\otimes F)$ with $\mathbb P(F)$. Moreover, if we denote by $\pi:\mathbb P(F)\rightarrow X$ the structural morphism, one has $\mathcal O_{E\otimes F}(1)=\mathcal \pi^*(E)\otimes \mathcal O_F(1)$, and the metric family $(\psi_{E}+\psi_F)^{\operatorname{FS}}$ identifies with the tensor product of $\pi^*(\psi_E)$ and $\psi_F^{\operatorname{FS}}$. Hence the assertions follow from \cite[Propositions 6.1.12 and 6.1.28]{CMArakelovAdelic}.
\end{proof}

\begin{prop}\label{Pro: existence of dominated and measurable}
Let $E$ be a vector bundle on $X$. Then there exists a dominated and measurable 
metric family of $E$.
\end{prop}

\begin{proof}
Let $L$ be an ample line bundle on $X$ and $\varphi$ be a dominated and measurable 
metric family of $L^\vee$. Then, one can find a positive integer $m$ such that $ L^m\otimes E$ is ample and
generated by global sections. If $ L^m\otimes E$ has  a dominated and measurable 
metric family $\psi'$, then the tensor product of $m \varphi$  with $\psi'$ is a dominated and measurable 
metric family of $E$ by Proposition~\ref{prop:tensor:line:bundle:dominated:measurable}, so we may assume that $E$ is ample and generated by global sections.

Let $H^0(X, E) \otimes \OO_X \to E$ be the natural surjective homomorphism.
Fix a basis $e_1, \ldots, e_N$ of $H^0(X, E)$ and, for each $\omega \in \Omega$ and $(a_1, \ldots, a_N) \in K_{\omega}^N$, we set
\[
\| a_1 e_1 + \cdots + a_N e_N \|_{\omega} = \begin{cases}
\sqrt{|a_1|_{\omega}^2 + \cdots + |a_N|_{\omega}^2} & \text{if $\omega \in \Omega_{\infty}$}, \\[2ex]
\max \{ |a_1|_{\omega}, \ldots, |a_N|_{\omega} \} & \text{if $\omega \in \Omega \setminus \Omega_{\infty}$}, 
\end{cases}
\]
and $\xi$ be the norm family $( \|\ndot\|_{\omega})_{\omega \in \Omega}$.
Let $\psi$ be a metric family of $E$ induced by
$H^0(X, E) \otimes \OO_{X} \to E$  and $\xi$. 
Let $\pi : \PP(E) \to X$ be the projective bundle of $E$ and $\OO_E(1)$ be the tautological line bundle of $\PP(E)$.
Note that the metric family $\psi^{\mathrm{FS}}$ of $\OO_E(1)$ is induced by
$H^0(X, E) \otimes \OO_{\PP(E)} \to \OO_E(1)$ and $\xi$, so
it is dominated and measurable.
Thus the assertion follows.
\end{proof}

\section{Dual metric family}

In this section, let $X$ be a  projective scheme over $\Spec K$.

\begin{defi}
Let $E$ be a vector bundle on $X$, equipped with a metric family $\psi=(\psi_\omega)_{\omega\in\Omega}$. For any $\omega\in\Omega$ and any $x\in X_\omega^{\mathrm{an}}$, the norm $|\ndot|_{\psi_\omega}(x)$ on $E_\omega(x)$ induces a dual norm on $E_\omega(x)^\vee$, which we denote by $|\ndot|_{\psi_{\omega}^\vee}(x)$. It turns out that $\psi_\omega^\vee=(|\ndot|_{\psi_\omega^\vee}(x))_{x\in X_\omega^{\mathrm{an}}}$ forms a continuous metric on $E_\omega^\vee$. Hence $\psi^\vee=(\psi_\omega^\vee)_{\omega\in\Omega}$ is a metric family on $E^\vee$, called the \emph{dual metric family of $\psi$}\index{metric family!dual ---}.
\end{defi}

\begin{prop}\label{Pro: dual family dominated}
Let $E$ be a vector bundle on $X$ and $\psi$ be a metric family of $E$. If $\psi$ is dominated, then the dual metric family $\psi^\vee$ is also dominated.
\end{prop}
\begin{proof} Let $\pi_{E}:\mathbb P(E)\rightarrow X$ and $\pi_{E^\vee}:\mathbb P(E^\vee)\rightarrow X$ be the projective bundles associated with $E$ and $E^\vee$ respectively.
We consider the fiber product $\mathbb P(E)\times_{X}\mathbb P(E^\vee)$ of  projective bundles and denote by \[p_1:\mathbb P(E)\times_{X}\mathbb P(E^\vee)\longrightarrow\mathbb P(E)\quad\text{and}\quad p_2:\mathbb P(E)\otimes_X\mathbb P(E^\vee)\longrightarrow\mathbb P(E^\vee)\]
the morphisms of projection. Let 
\[\mathcal O_E(1)\boxtimes\mathcal O_{E^\vee}(1):=p_1^*(\mathcal O_E(1))\otimes p_2^*(\mathcal O_{E^\vee}(1))\]
and let \[s\in H^0(\mathbb P(E)\times_X\mathbb P(E^\vee),\mathcal O_E(1)\boxtimes\mathcal O_{E^\vee}(1))\] be the trace section of $\mathcal O_E(1)\boxtimes\mathcal O_{E^\vee}(1)$, which corresponds to the composition of the following universal homomorphisms
\[p_2^*(\mathcal O_{E^\vee}(-1))\longrightarrow p_2^*(\pi_{E^\vee}^*(E))\cong p_1^*(\pi_{E}^*(E))\longrightarrow p_1^*(\mathcal O_{E}(1)).\] 

\ifsmf\begin{enonce}{Claim}\fi
\ifams\begin{clai}\fi
Let $\psi_1=(\psi_{1,\omega})_{\omega\in\Omega}$ and $\psi_2=(\psi_{2,\omega})_{\omega\in\Omega}$ be metric families on $E$ and $E^\vee$ respectively. We equip $\mathcal O_E(1)\boxtimes\mathcal O_{E^\vee}(1)$ with the metric family $\varphi=(\varphi_\omega)_{\omega\in\Omega}$ which is the tensor product of the metric families  $p_1^*(\psi_1^{\mathrm{FS}})$ and $p_2^*(\psi_2^{\mathrm{FS}})$. Then, for any $\omega\in\Omega$ and $x\in X_\omega^{\mathrm{an}}$, one has
\[\sup_{f\in E_\omega^\vee(x)\setminus\{0\}}\frac{|f|_{\psi_{1,\omega}^\vee}(x)}{|f|_{\psi_{2,\omega}}(x)}\leqslant\|s\|_{\varphi_\omega}.\]
\ifsmf\end{enonce}\fi
\ifams\end{clai}\fi
\begin{proof}
Let $f$ be a non-zero element of $E_\omega^\vee(x)$. The one-dimensional $\widehat{\kappa}(x)$-vector space of $E_\omega^\vee(x)$ spanned by $f$ determines a  point $P_f$ of $ P(E_\omega)^{\mathrm{an}}$ valued in $(\widehat{\kappa}(x),|\ndot|_x)$ which lies over $x\in X_\omega^{\mathrm{an}}$. Suppose $Q$ is a point of $\mathbb P(E_\omega^\vee)^{\mathrm{an}}$ valued in $(\widehat{\kappa}(x),|\ndot|_x)$ which lies over $x$. Then $s(P_f,Q)$ corresponds to the following composition of universal homomorphisms
\begin{equation}\label{Equ: composition universal}\mathcal O_{E^\vee}(-1)(Q)\longrightarrow E_\omega(x)\longrightarrow\mathcal O_{E}(1)(P_f),\end{equation}
and $|s|_{\varphi_\omega}(P_f,Q)$ is the operator norm of this homomorphism.
We pick an arbitrary non-zero element $\ell$ of $\mathcal O_{E^\vee}(-1)(Q)$. The dual element in $\mathcal O_E(-1)(P_f)$  of the image of $\ell$ by \eqref{Equ: composition universal} is $f(\ell)^{-1}f$.  Therefore one has 
\[|s|_{\varphi_\omega}(P_f,Q)=\frac{|f(\ell)|_x}{|\ell|_{\psi_{1,\omega}}(x)\cdot|f|_{\psi_{2,\omega}}(x)}\leqslant\|s\|_{\varphi_\omega}.\] 
Taking the supremum with respect to $\ell$, we obtain the required inequality.
\end{proof}

In the above claim, if both metric families $\psi_1$ and $\psi_2$ are dominated, then the metric family $\varphi$ on $\mathcal O_E(1)\boxtimes\mathcal O_{E^\vee}(1)$ is also dominated. In particular, the function \[(\omega\in\Omega)\longmapsto \ln\|s\|_{\varphi_\omega}\]
is bounded from above by an integrable function. Then the claim shows that the function
\[(\omega\in\Omega)\longmapsto \sup_{x\in X_\omega^{\mathrm{an}}}\sup_{f\in E_\omega^\vee(x)\setminus\{0\}}\Big(\ln|f|_{\psi_{1,\omega}^\vee}(x)-\ln|f|_{\psi_{2,\omega}}(x)\Big)\]
is bounded from above by an integrable function. Therefore, the function 
\[(\omega\in\Omega)\longmapsto\sup_{Q\in\mathbb P(E^\vee)^{\mathrm{an}}}\sup_{\begin{subarray}{c}f\in\mathcal O_{E^\vee}(1)(Q)\\
f\neq 0
\end{subarray}}\Big(\ln|f|_{\psi_{1,\omega}^{\vee,\operatorname{FS}}}(Q)-\ln|f|_{\psi_{2,\omega}^{\operatorname{FS}}}(Q)\Big)\]
is bounded from above by an integrable function. For the same reason, by interchanging the roles of $E$ and $E^\vee$ we obtain that the function
\[(\omega\in\Omega)\longmapsto\sup_{P\in\mathbb P(E)^{\mathrm{an}}}\sup_{\begin{subarray}{c}t\in\mathcal O_{E}(1)(P)\\
t\neq 0
\end{subarray}}\Big(\ln|t|_{\psi_{2,\omega}^{\vee,\operatorname{FS}}}(P)-\ln|t|_{\psi_{1,\omega}^{\operatorname{FS}}}(P)\Big)\]
is also bounded from above by an integrable function. In particular, if we denote by $\widetilde{\varphi}$ the tensor product of the metric families $p_1^*(\psi_2^{\vee,\operatorname{PS}})$ and $p_2^*(\psi_1^{\vee,\operatorname{PS}})$, then the function
\[(\omega\in\Omega)\longmapsto \ln\|s\|_{\widetilde{\varphi}_\omega}\]
is still bounded from above by an integrable function. Hence the above claim implies that the function 
\[(\omega\in\Omega)\longmapsto \sup_{x\in X_\omega^{\mathrm{an}}}\sup_{f\in E_\omega^\vee(x)\setminus\{0\}}\Big(\ln|f|_{\psi_{2,\omega}}(x)-\ln|f|_{\psi_{1,\omega}^\vee}(x)\Big)\]
is bounded from above by an integrable function. Therefore we obtain that the local distance function
\[(\omega\in\Omega)\longmapsto d_\omega(\psi_1^\vee,\psi_2)\]
is bounded from above by an integrable function. By Proposition \ref{Pro: criterion of dominancy}, the metric family $\psi_1^\vee$ is dominated. By Proposition \ref{Pro: existence of dominated and measurable}, there exists at least a dominated metric family on $E^\vee$, the assertion is thus proved.  
\end{proof}

\begin{defi}
Let $E$ be a vector bundle on $X$, $\psi=(\psi_\omega)_{\omega\in\Omega}$ be a metric family of $E$. Let $K'/K$ be a finite extension and let $P:\Spec K'\rightarrow X$ be a $K$-morphism. Let
\[(K',(\Omega',\mathcal A',\nu'),\phi')=S\otimes_K K'.\]
Recall that $\Omega'$ is a disjoint union 
\[\Omega'=\coprod_{\omega\in\Omega}\Omega'_{\omega},\]
where $\Omega'_{\omega}$ denotes the set of all absolute values on $\Omega'$ extending $|\ndot|_\omega$. For any $\omega\in\Omega$ and any $x\in\Omega'_{\omega}$, we let $P_x:\Spec K'_x\rightarrow X_\omega$ be the morphism induced by \[\Spec K'_x\longrightarrow\Spec K'\stackrel{P}{\longrightarrow} X\] and the canonical morphism $\Spec K'_x\rightarrow\Spec K_\omega$.
\[\xymatrix{K'_x\ar@/_/[rdd]\ar@/^/[rrd]\ar@{.>}[rd]|-{P_x}\\
&X_\omega\ar[r]\ar[d]&X\ar[d]\\
&\Spec K_\omega\ar[r]&\Spec K }\] 
We denote by $\|\ndot\|_x$ the norm on 
\[(E\otimes_{K}K')\otimes_{K'}K'_x\cong E_\omega\otimes_{K_\omega}K'_x\]
which is induced by $|\ndot|_{\psi_\omega}(P_x)$ by $\varepsilon$-extension of scalars if $|\ndot|_\omega$ is non-Ar\-chi\-me\-de\-an, and by $\pi$-extension of scalars if $|\ndot|_\omega$ in Archimedean. Then, $(\|\ndot\|_x)_{x\in\Omega'}$ forms a norm family of $P^*(E)$, which we denote by $P^*(\psi)$.
\end{defi}

\begin{defi}
Let $\Omega_0$ be the set of $\omega\in\Omega$ such that the absolute value $|\ndot|_\omega$ is trivial. Let $x=(K_x,|\ndot|_x,P_x)$ be a triplet, where $(K_x,|\ndot|_x)$ is a  valued extension of the trivially valued field $(K,|\ndot|_0)$, and $P_x:\Spec K_x\rightarrow X$ is a $K$-morphism. Assume that $E$ is a vector bundle on $X$ and $\psi=(\psi_\omega)_{\omega\in\Omega}$ be a metric family of $E$. Denote by $E_x$ the $K_x$-vector space $P_x^*(E)$. We consider the following adelic curve
\[S_x=(K_x,(\Omega_0,\mathcal A_0,\nu_0),(|\ndot|_x)_{\omega\in\Omega_0}),\]
where $\mathcal A_0$ is the restriction of the $\sigma$-algebra $\mathcal A$ to $\Omega_0$, and $\nu_0$ is the restriction of $\nu$ to $(\Omega_0,\mathcal A_0)$. We denote by $x^*(\psi)$ the norm family $(|\ndot|_{\psi_\omega}(P_x^{\omega}))_{\omega\in\Omega_0}$ on $E_x$, where $P_x^\omega$ denotes the point of $X_\omega^{\mathrm{an}}$ determined by $(P_x,|\ndot|_x)$.

\end{defi}

\begin{prop}
Let $E$ be a vector bundle on $X$ and $\psi=(\psi_\omega)_{\omega\in\Omega}$ be a metric family of $E$. Then the metric family $\psi$ is measurable if and only if both of the following conditions are satisfied:
\begin{enumerate}[label=\rm(\arabic*)]
\item for any finite extension $K'/K$ and any $K$-morphism $P:\Spec K'\rightarrow X$, the norm family $P^*(\psi)$ is measurable,
\item for any triplet $x=(K_x,|\ndot|_x,P_x)$, where $(K_x,|\ndot|_x)$ is a valued extension of transcendence degree $\leqslant 1$ and of rational exponent \textup{(}see \textup{\S\ref{Sub: exponent}}\textup{)} of the trivially valued field $(K,|\ndot|_0)$, and $P_x:\Spec K_x\rightarrow X$ is a $K$-morphism, the norm family $x^*(\psi)$ is measurable.
\end{enumerate}
\end{prop}
\begin{proof}It suffices to treat the case where the field $K$ is countable. 
Recall that the measurability of the metric family $\psi$ signifies that the following two conditions are satisfied (see \S\ref{Sec: measurability line bundle}):
\begin{enumerate}[label=(\arabic*')]
\item for any finite extension $K'/K$ and any $K$-morphism $Q:\Spec K'\rightarrow\mathbb P(E)$, the norm family $Q^*(\psi^{\operatorname{FS}})$ is measurable,
\item for any triplet $y=(K_y,|\ndot|_y,Q_y)$, where $(K_y,|\ndot|_y)$ is a valued extension of transcendence degree $\leqslant 1$ and of rational exponent of the trivially valued field $(K,|\ndot|_0)$, and $Q_y:\Spec K_y\rightarrow\mathbb P(E)$ is a $K$-morphism, the norm family $Q_y^*(\psi^{\operatorname{FS}})$ is measurable. 
\end{enumerate}

Let $K'/K$ be a finite extension. Any $K$-morphism $Q:\Spec K'\rightarrow\mathbb P(E)$ corresponds to a $K$-morphisme $P:\Spec K'\rightarrow X$ together with a one-dimensional quotient vector space $L$ of $P^*(E)$, which identifies with $Q^*(\mathcal O_E(1))$. Moreover, the norm family $Q^*(\psi^{\operatorname{FS}})$ identifies with the quotient norm family of $P^*(\psi)$. If the norm family $P^*(\psi)$ is measurable, by \cite[Proposition 4.1.24]{CMArakelovAdelic}, we obtain that $Q^*(\psi^{\operatorname{FS}})$ is also measurable. Conversely, if for any one-dimensional quotient vector space of $P^*(E)$, the quotient norm family of $P^*(\psi)$ on it is measurable, by passing to dual we obtain from \cite[Proposition 4.1.24]{CMArakelovAdelic} that $P^*(\psi)^\vee$ is measurable and therefore $P^*(\psi)$ is also measurable. 

Let $x=(K_x,|\ndot|_x,P_x)$ be a triplet, where $(K_x,|\ndot|_x)$ is a valued extension of transcendence degree $\leqslant 1$ and rational exponent of the trivially valued field $(K,|\ndot|_0)$, and $P_x:\Spec K_x\rightarrow X$ be a $K$-morphism. Note that the field $K_x$ is countable. Similarly to the above argument, the norm family $P_x^{*}(\psi)$ is measurable if and only if all its quotient norm families on one-dimensional quotient subspaces are measurable. The proposition is thus proved.
\end{proof}

\begin{prop}\label{Pro: measurability of dual}
Let $E$ be a vector bundle on $X$ and $\psi=(\psi_\omega)_{\omega\in\Omega}$ be a metric family on $E$. If $\psi$ is measurable, then the dual metric family $\psi^\vee$ of $E^\vee$ is also measurable.
\end{prop}
\begin{proof}
Let $K'/K$ be a finite extension and $P:\Spec K'\rightarrow X$ be a $K$-morphism. If $P^*(\psi)$ is measurable, by \cite[Proposition 4.1.24]{CMArakelovAdelic} we obtain that $P^*(\psi^\vee)=P^*(\psi)^\vee$ is measurable. Similarly, for any triplet $x=(K_x,|\ndot|_x,P_x)$, where $(K_x,|\ndot|_x)$ is a valued extension of transcendence degree $\leqslant 1$ and of rational exponent of the trivially valued field $(K,|\ndot|_0)$ and $P_x:\Spec K_x\rightarrow X$ is a $K$-morphism, if the norm family $P^*_x(\psi)$ is measurable, then $P_x^*(\psi^\vee)=P_x^*(\psi)^\vee$ is also measurable. The proposition is thus proved.    
\end{proof}

\begin{coro}\label{coro:sub:metric:dominated}
Let $E$ be a vector bundle on $X$, $F$ be a vector subbundle of $E$, $\psi_E$ be a metric family of $E$, and $\psi_F$ be the restriction of $\psi_E$ to $F$. If the metric family $\psi_E$ is dominated \textup{(}resp. measurable\textup{)}, then the restricted metric family $\psi_F$ is also dominated \textup{(}resp. measurable\textup{)}.
\end{coro}
\begin{proof}
The homomorphism of inclusion $F\rightarrow E$ induces by passing to dual a surjective homomorphism $E^\vee\rightarrow F^\vee$. Thus $F^\vee$ can be considered as a quotient vector bundle of $E^\vee$. Note that $\psi_F^\vee$ identifies with the quotient metric family of $\psi_E^\vee$. Hence the assertion follows from Propositions \ref{Pro: dual family dominated}, \ref{Pro: measurability of dual} and \ref{Pro: quotient family dominated}.
\end{proof}

\section{Metric families on torsion-free sheaves} 
In this section, we assume that the $K$-scheme $X$ is geometrically integral.

\begin{defi}\label{def:globally:adelic:torsion:free}
Let $E$ be a torsion free $\OO_X$-module and $U$ be a non-empty Zariski open set of $X$ such that $E|_U$ is a vector bundle 
For any $\omega\in\Omega$, let $\psi_{\omega}$ be a continuous metric of $E_{\omega}$ over $U_{\omega}^{\an}$ such that, for any $s \in H^0(X_\omega, E_\omega)$,
\[\|s\|_{\psi_{\omega}}:=\sup \{ |s|_{\psi_{\omega}}(x) \,:\, x \in U^{\an}_{\omega} \} < +\infty.\] 
We set $\psi = ( \psi_{\omega} )_{\omega \in \Omega}$ and $\xi_{\psi} = ( \|\ndot\|_{\psi_{\omega}} )_{\omega \in \Omega}$.
We say that $(E, U,\psi)$ is a \emph{sectionally  adelic torsion free $\OO_X$-module}\index{adelic torsion free $\OO_X$-module!sectionally  ---} if $(H^0(X, E), \xi_{\psi})$ is an adelic vector bundle on $S$.
By Proposition~\ref{prop:tensor:line:bundle:dominated:measurable}, if $(E,\psi)$ is an adelic vector bundle on $X$, then, for any non-empty Zariski open set $U$ of $X$, $(E,U,\psi)$ is a sectionally adelic torsion free $\mathcal O_X$-module.
\end{defi}

\begin{defi}\label{def:birat:adelic:torsion:free}
Let $E$ be a torsion free $\mathcal O_X$-module and $U$ be a non-empty Zariski open set of $X$ such that $\rest{E}{U}$ is a vector bundle.
Let $\psi = ( \psi_{\omega} )_{\omega \in \Omega}$ be a metric family of $E|_U$.
We say $(E, U, \psi)$ is a \emph{birationally adelic torsion free $\OO_X$-module}\index{adelic torsion free $\OO_X$-module!birationally ---} if it satisfies the following properties:
\begin{enumerate}[label=\rm(\arabic*)]
\item There exist a birational morphism $\mu : X' \to X$ of geometrically integral projective schemes over $K$ such that
$\mu^{-1}(U) \to U$ is an isomorphism,
an adelic vector bundle $(E', \psi')$ on $X'$, and an injective morphism of $\mathcal O_{X'}$-modules $E \rightarrow \mu_*(E')$ whose restriction to $U$ gives an isomorphism $\rest{E}{U} \to \rest{\mu_*(E')}{U}\cong E'|_{\mu^{-1}(U)}$.

\item The isomorphism $\rest{E}{U} \to \rest{E'}{\mu^{-1}(U)}$ yields an isometry \[\rest{(E, \psi)}{U} \longrightarrow \rest{(E', \psi')}{\mu^{-1}(U)}.\]
\end{enumerate}
\end{defi}

By definition, for $s \in H^0(X, E)$ and each $\omega \in \Omega$, \[\| s \|_{\psi_{\omega}}:=\sup \{ |s|_{\psi_{\omega}}(\xi) \,:\, \xi \in U_{\omega}^{\an}\}\]
exists.
Note that $\|\ndot\|_{\psi_{\omega}}$ is the restriction of $\|\ndot\|_{\psi'_\omega}$ to $H^0(X, E)$ by using
the injective homomorphism $H^0(X, E) \to H^0(X', E')$, so that
\[(H^0(X, E), ( \|\ndot\|_{\psi_{\omega}} )_{\omega \in \Omega} )\] is an adelic vector bundle on $S$, that is,
a birationally adelic torsion free $\OO_X$-module is sectionally adelic in the sense of Definition~\ref{def:globally:adelic:torsion:free}.

\begin{lemm}\label{lem:proper:cont:max}
Let $\pi : X \to Y$ be a continuous map of locally compact Hausdorff spaces such that $\pi$ is open and proper.
Let $f : X \to \RR$ be a continuous function on $X$ and $\tilde{f} : Y \to \RR$ be a function on $Y$ given by \[\tilde{f}(y) = \max \{ f(x) \,:\, \pi(x)=y \}.\]
Then $\tilde{f}$ is continuous on $Y$.
\end{lemm}

\begin{proof}
Fix $y_0 \in Y$. Since $\pi^{-1}(y_0)$ is compact, for $\varepsilon > 0$,
there exist $x_{1}, \ldots, x_{n}  \in \pi^{-1}(y_0)$ and open subsets $U_1, \ldots, U_n$ of $X$ such that
\[\pi^{-1}(y_0) \subseteq U_1 \cup \cdots \cup U_n,\] $x_i \in U_i$ for all $i\in\{1,\ldots,n\}$ and $|f(x) - f(x_{i})| \leqslant \varepsilon$ for all $i\in\{1,\ldots,n\}$ and $x \in U_i$.
If we set $Z = X \setminus U_1 \cup \cdots \cup U_n$, then
$\pi(Z)$ is closed and $y_0 \not\in \pi(Z)$.
We choose an open set $W$ of $Y$ such that $y_0 \in W$ and \[W \subseteq \pi(U_1) \cap \cdots \cap \pi(U_n) \cap (Y \setminus \pi(Z)).\]
Note that $\pi^{-1}(W) \subseteq U_1 \cup \cdots \cup U_n$. 
Let $y \in W$ and \[\lambda_i = \sup \{ f(x) \,:\, \text{$x \in U_i$ and $y = \pi(x)$} \}.\]
Then $\tilde{f}(y) = \max \{ \lambda_1, \ldots, \lambda_n \}$ and $\lambda_i -\varepsilon \leqslant f(x_{i}) \leqslant \lambda_i + \varepsilon$ for all $i\in\{1,\ldots,n\}$,
so that \[\tilde{f}(y) - \varepsilon \leqslant \tilde{f}(y_0) \leqslant \tilde{f}(y) + \varepsilon,\] as required.
\end{proof}

Let $\pi : X \to Y$ be a generically finite morphism of geometrically integral projective schemes over $\Spec K$ and 
$(M, U,\psi)$ be a sectionally adelic torsion free $\OO_X$-module. 
Note that $\pi_*(M)$ is a torsion free $\OO_Y$-module. The pushforward $\pi_*(\psi)$ is defined as follows:
We choose a non-empty Zariski open set $V$ of $Y$ such that \[\rest{\pi}{\pi^{-1}(V)} : \pi^{-1}(V) \longrightarrow V\] is \'{e}tale and $\pi^{-1}(V) \subseteq U$.
Note that $\pi_*(M)$ is locally free over $V$.
For $y \in V_{\omega}^{\an}$ and $s \in \pi_*(M) \otimes \hat{\kappa}(y)$, $|s|_{\pi_*(\psi)_\omega}(y)$ is defined to be
\[|s|_{\pi_*({\psi_\omega})}(y) := \max \{ |s|_{\psi_{\omega}}(x) \,:\, x \in (\pi^{\an}_\omega)^{-1}(y)\}.\]
Since $\pi^{-1}(V)^{\an}_\omega \to V^{\an}_\omega$ is proper and open (c.f. \cite[Lemma~3.2.4]{MR1070709}),  
by Lemma~\ref{lem:proper:cont:max},
$\pi_*(\psi)_{\omega}$ yields a continuous metric of $\pi_*(M)_{\omega}$ over $V_{\omega}^{\an}$. We denote $( \pi_*(\psi)_\omega )_{\omega \in \Omega}$ by $\pi_*(\psi)$.
For $s \in H^0(Y, \pi_*(M)) = H^0(X, M)$, as 
\[
\sup \{ |s|_{\pi_*(\psi)_\omega}(y) \,:\, y \in V_{\omega}^{\an} \} = \sup \{ |s|_{\psi_\omega}(x) \,:\, x \in \pi^{-1}(V)_{\omega}^{\an} \},
\]
one has $\| s \|_{\pi_*(\psi)_\omega} = \| s \|_{\psi_\omega} < \infty$, so that 
$(\pi_*(M),V, \pi_*(\psi))$ forms a sectionally adelic torsion free $\OO_Y$-module and
$(H^0(Y, \pi_*(M)),  ( \|\ndot\|_{\pi_*(\psi)_\omega} )_{\omega \in \Omega})$ is isometric to 
$(H^0(X, M),  ( \|\ndot\|_{\psi_\omega} )_{\omega \in \Omega})$. We call $V$ an \emph{open subscheme of definition}\index{open subscheme of definition} of $\pi_*(\psi)$.

%!TEX root = ./Hilbert_Samuel_Adelic_Curves.tex

\chapter{Volumes of  normed graded linear series}

In this chapter, we let $k$ be a commutative field and we denote by $|\ndot|_0$ the trivial absolute value on $k$. Recall that $|a|_0=1$ for any $a\in k^{\times}$. Moreover, $S_0=(k,\{0\},|\ndot|_0)$ forms an adelic curve. 

\section{Adelic vector bundle on $S_0$}
Adelic vector bundles on $S_0$ are just finite-dimensional ultrametrically normed vector space over $k$. If $\overline V=(V,\|\ndot\|)$ is an adelic vector bundle on $S_0$, then the function $\|\ndot\|$  only takes finitely many values. Moreover, if the vector space $V$ is non-zero, then one has (see \cite[Remark 4.3.63]{CMArakelovAdelic})
\[\widehat{\mu}_{\max}(\overline V)=-\min_{s\in V\setminus\{0\}}\ln\|s\|,\quad \widehat{\mu}_{\min}(\overline V)=-\max_{s\in V}\ln\|s\|.\]
The Harder-Narasimhan $\mathbb R$-filtration of $\overline V$ is give by 
\[\forall\,t\in\mathbb R,\quad \mathcal F^t(\overline V)=\{s\in V\,:\,\|s\|\leqslant\mathrm{e}^{-t}\}.\]
Note that 
\begin{gather*}\widehat{\deg}_+(\overline V):=\sup_{W\subset V}\widehat{\deg}(\overline W)=\int_0^{+\infty}\dim_{k}(\mathcal F^t(\overline V))\,\mathrm{d}t,\\
\widehat{\deg}(\overline V)=-\int_{\mathbb R}\,t\,\mathrm{d}\dim_k(\mathcal F^t(\overline V))\,\mathrm{d}t.\end{gather*}

\section{Normed graded algebra}

\label{Def: normed graded algebra}
Let $V_\sbullet=\bigoplus_{n\in\mathbb N}V_n$ be a graded $k$-algebra. We assume that each $V_n$ is a finite-dimensional vector space over $k$. For any $n\in\mathbb N_{\geqslant 1}$, let $\|\ndot\|_n$ be an ultrametric norm on $V_n$. Then the pair $\overline V_{\!\sbullet}=(V_{\sbullet},(\|\ndot\|_n)_{n\in\mathbb N_{\geqslant 1}})$ is called a \emph{normed graded algebra}\index{normed graded algebra} over $(k,|\ndot|_0)$. Let $f:\mathbb N_{\geqslant 1}\rightarrow\mathbb R_{\geqslant 0}$ be a function. If, for all $\ell\in\mathbb N_{\geqslant 2}$, $(n_1,\ldots,n_\ell)\in\mathbb N_{\geqslant 1}^\ell$ and  $(s_1,\ldots,s_\ell)\in V_{n_1}\times\cdots\times V_{n_\ell}$, one has
\begin{equation}\label{Equ: sub multiplicative}\|s_1\cdots s_\ell\|_{n_1+\cdots+n_\ell}\leqslant \mathrm{e}^{f(n_1)+\cdots+f(n_\ell)}\|s_1\|_{n_1}\cdots\|s_{\ell}\|_{n_\ell},\end{equation}
we say that $\overline V_{\!\sbullet}$ is  \emph{$f$-sub-multiplicative}\index{sub-multiplicative!$f$- ---}. In the particular case where $f$ is the constant function taking value $0$, we just say the $\overline V_{\!\sbullet}$ is \emph{sub-mutiplicative}\index{sub-mutiplicative}.  If there exist two constant $C_1$ and $C_2$ such that, for any $n\in\mathbb N$ and any $s\in V_n\setminus\{0\}$, one has
\begin{equation}\label{Equ: bounded}\mathrm{e}^{C_1n}\leqslant\|s\|_n\leqslant \mathrm{e}^{C_2n},\end{equation}
we say that $\overline{V}_{\!\sbullet}$ is \emph{bounded}\index{bounded}.

\begin{prop}\label{Pro: graded algebra}
Let $\overline V_{\!\sbullet}$ be a normed graded algebra over $(k, |\ndot|_0)$ and $f:\mathbb N_{\geqslant 1}\rightarrow\mathbb R_{\geqslant 0}$ be a function such that 
\[\lim_{n\rightarrow+\infty}\frac{f(n)}{n}=0.\] Assume that $V_\sbullet$ is an integral domain and that $\overline V_{\!\sbullet}$ is $f$-sub-multiplicative and bounded. 
\begin{enumerate}[label=\rm(\arabic*)]
\item\label{Item: spectral norm} For any $n\in\mathbb N_{\geqslant 1}$ and any $s\in V_n$, the sequence
\[\|s^N\|_{nN}^{1/N},\quad N\in\mathbb N,\;N\geqslant 1\]
converges.
\item\label{Item: norme spectrale} For any $n\in\mathbb N_{\geqslant t 1}$, the map
\[\|\ndot\|_{\operatorname{sp},n}:V_n\longrightarrow\mathbb R_{\geqslant 0},\quad s\longmapsto\lim_{N\rightarrow+\infty}\|s^N\|_{nN}^{1/N}\]
is an ultrametric norm on $V_n$. 
\item\label{Item: sous multiplicativite} The family of norms $(\|\ndot\|_{\operatorname{sp},n})_{n\in\mathbb N}$ satisfies the following sub-mul\-ti\-pli\-ca\-tiv\-i\-ty condition: for any $(n,m)\in\mathbb N^2$ and any $(s_n,s_m)\in V_n\times V_m$,
\[\|s_ns_m\|_{\operatorname{sp},n+m}\leqslant\|s_n\|_{\operatorname{sp},n}\cdot\|s_m\|_{\operatorname{sp},m}.\]
\item\label{Item: comparison} For any $n\in\mathbb N_{\geqslant 1}$ and any $s\in V_n\setminus\{0\}$, one has 
\begin{equation}\label{Equ: upper bound}\|s\|_{\operatorname{sp},n}\leqslant \mathrm{e}^{f(n)}\|s\|_n.\end{equation}
\end{enumerate}
\end{prop}
\begin{proof}
\ref{Item: spectral norm} It suffices to treat the case where $s\neq 0$. By \eqref{Equ: sub multiplicative}, for $\ell\in\mathbb N _{\geqslant 2}$, and $(N_1,\ldots,N_\ell)\in\mathbb N_{\geqslant 1}^\ell$, one has 
\[\ln\|s^{N_1+\cdots+N_\ell}\|_{n(N_1+\cdots+N_\ell)}\leqslant \sum_{i=1}^\ell\Big(\ln\|s^{N_i}\|_{nN_1}+f(nN_i)\Big).\]
Moreover, by \eqref{Equ: bounded}, the sequence 
\[\frac{1}{N}\ln\|s^N\|_{nN},\quad N\in\mathbb N,\;N\geqslant 1\]
is bounded. Therefore this sequence converges in $\mathbb R$ (see \cite[Proposition 1.3.1]{MR2768967}), which shows that the sequence
\[\|s^N\|_{nN}^{1/N},\quad N\in\mathbb N,\;N\geqslant 1\]
converges to a positive real number.

\ref{Item: norme spectrale} It suffices to show that $\|\ndot\|_{\operatorname{sp},n}$ satisfies the strong triangle inequality. Let $s$ and $t$ be two elements of $V_n$. For any $N\in\mathbb N_{\geqslant 1}$, one has
\[(s+t)^N=\sum_{i=0}^N\binom{N}{i}s^it^{N-i}\] 
and hence
\[\|(s+t)^N\|_{nN}\leqslant\max_{i\in\{0,\ldots,N\}}\|s^it^{N-i}\|_{nN}.\]
Let 
\[M=\max_{j\in\mathbb N,\,j\geqslant 1}\frac{1}{j}\max\{\ln\|s^j\|_{nj},\ln\|t^j\|_{nj},0\}.\]
Let $(\varepsilon_j)_{j\in\mathbb N}$ be a sequence of real numbers in $[0,\frac 12]$ such that 
\[\lim_{j\rightarrow+\infty}\varepsilon_j=0,\quad \lim_{j\rightarrow+\infty}j\varepsilon_j=+\infty,\quad \lim_{j\rightarrow+\infty}(j-j\varepsilon_j)=+\infty.\]
If $i/N\leqslant\varepsilon_N$, one has
\begin{multline*}\frac{1}{N}\ln\|s^it^{N-i}\|_{nN}\leqslant \varepsilon_NM+\frac{N-i}{N}\cdot\frac{1}{N-i}\ln\|t^{N-i}\|_{n(N-i)} \\
+\frac{f(ni)}{N}+\frac{f(n(N-i))}{N}.
\end{multline*}
Similarly, if $(N-i)/N\leqslant\varepsilon_N$, one has
\[\frac{1}{N}\ln\|s^it^{N-i}\|_{nN}\leqslant \frac{i}{N}\cdot\frac{1}{i}\ln\|s^i\|_{ni}+\varepsilon_NM+\frac{f(ni)}{N}+\frac{f(n(N-i))}{N}.\]
If $N\varepsilon_N<i<N-N\varepsilon_N$, one has 
\begin{multline*}
\frac{1}{N}\ln\|s^it^{N-i}\|_{nN}\leqslant\frac{i}{N}\cdot\frac{1}{i}\ln\|s^i\|_{ni}+\frac{N-i}{N}\cdot\frac{1}{N-i}\ln\|t^{N-i}\|_{n(N-i)} \\
+\frac{f(ni)}{N}+\frac{f(n(N-i))}{N}.\end{multline*}
Taking the superior limit when $N\rightarrow+\infty$, we obtain that 
\[\limsup_{N\rightarrow+\infty}\max_{i\in\{0,\ldots, N\}}\frac{1}{N}\ln\|s^it^{N-i}\|_{nN}\leqslant\max\{\|s\|_{\operatorname{sp},n},\|t\|_{\operatorname{sp},n}\}.\]

\ref{Item: sous multiplicativite} Let $(n,m)\in\mathbb N^2$ and $(s_n,s_m)\in V_n\times V_m$. For any $N\in\mathbb N$ such that $N\geqslant 1$, one has
\[\|(s_ns_m)^{N}\|_{(n+m)N}\leqslant \mathrm{e}^{f(nN)+f(mN)}\|s_n^N\|_{nN}\cdot\|s_m^N\|_{mN}.\]
Taking the $N$-th root and letting $N\rightarrow+\infty$ we obtain
\[\|s_ns_m\|_{\operatorname{sp},n+m}\leqslant\|s_n\|_{\operatorname{sp},n}\cdot\|s_m\|_{\operatorname{sp},m}.\]

\ref{Item: comparison} For any $N\in\mathbb N_{\geqslant 1}$, the following inequality holds:
\[\|s^N\|_{nN}\leqslant \mathrm{e}^{Nf(n)}\|s\|_n^N.\]
Taking the $N$-th root and then letting $N\rightarrow+\infty$, we obtain
\[\|s\|_{\operatorname{sp},n}\leqslant \mathrm{e}^{f(n)}\|s\|_n.\]
\end{proof}

\section{Reminder on graded linear series}

In this subsection, we let $k'/k$ be a finitely generated extension of fields. As \emph{graded linear series}\index{graded linear series} of $k'/k$, we refer to a graded sub-$k$-algebra $V_\sbullet$ of 
\[k'[T]=\bigoplus_{n\in\mathbb N}k'T^n\]
such that $V_0=k$.
We denote by $\mathbb N(V_\sbullet)$ the set of $n\in\mathbb N$ such that $V_n\neq\mathbf{0}$. 
If $V_\sbullet$ is a graded linear series and $\mathbb N(V_\sbullet)\neq\{0\}$, we denote by $k(V_\sbullet)$ the sub-extension of $k'/k$ generated by 
\[\bigcup_{n\in\mathbb N(V_\sbullet) \setminus \{ 0 \}}\{f/g\,|\,(f,g)\in V_n\times(V_n\setminus\{0\})\}\] 
over $k$.
If $\mathbb N(V_\sbullet)\neq\{0\}$, then we denote by $\dim(V_\sbullet)$ the transcendence degree of the extension $k(V_\sbullet)/k$, and call it the \emph{Kodaira-Iitaka dimension}\index{Kodaira-Iitaka dimension} of $V_\sbullet$. 
In the case where $V_n=\{0\}$ for any $n\in\mathbb N_{\geqslant 1}$, by convention $\dim(V_\sbullet)$ is defined to be $-\infty$.
If $\mathbb N(V_\sbullet)\neq\{0\}$ and if the field $k(V_\sbullet)$
coincides with $k'$,
we say that the graded linear series $V_\sbullet$ is \emph{birational}\index{birational}.

  We say that $V_\sbullet$ is \emph{of sub-finite type}\index{of sub-finite type} if there exists a graded linear series $W_\sbullet$ of $k'/k$ which is a $k$-algebra of finite type and contains $V_\sbullet$ as a sub-$k$-algebra. By \cite[Theorem 3.7]{MR4104565}, there exists a graded sub-$k$-algebra of finite type $W_\sbullet$ of the polynomial ring
\[k(V_\sbullet)[T]=\bigoplus_{n\in\mathbb N}k(V_\sbullet)T^n\]
such that $k(W_\sbullet)=k(V_\sbullet)$, which contains $V_\sbullet$ as a sub-$k$-algebra. In other words, $V_\sbullet$  viewed as a graded linear series of $k(V_\sbullet)/k$ is sub-finite.

Let $V_\sbullet$ be a graded linear series of sub-finite type, and $d$ be its Kodaria-Iitaka dimension. If $\mathbb N(V_\sbullet)\neq\{0\}$, we define the \emph{volume}\index{volume} of $V_\sbullet$ as the limit (see \cite[Theorem 6.2]{MR4104565} for the convergence)
\[\operatorname{vol}(V_\sbullet):=\lim_{n\in\mathbb N(V_\sbullet),\,n\rightarrow+\infty}\frac{\dim_k(V_n)}{n^d/d!}.\]
Note that $V_\sbullet$ satisfies the Fujita approximation property, namely, one has
\[\vol(V_\sbullet)=\sup_{\begin{subarray}{c}W_\sbullet\subset V_\sbullet\\
\dim(W_\sbullet)=\dim(V_\sbullet)\end{subarray}}\vol(W_\sbullet),\]
where $W_\sbullet$ runs over the set of all graded sub-$k$-algebras of finite type of $V_\sbullet$ such that $\dim(W_\sbullet)=\dim(V_\sbullet)$.

\section{Normed graded linear series}
In this subsection, we fix a finitely generated extension $k'/k$, a graded linear series $V_\sbullet$ of $k'/k$ which is of sub-finite type, and a $f:\mathbb N_{\geqslant 1}\rightarrow\mathbb R_{\geqslant 0}$ such that 
\[\lim_{n\rightarrow+\infty}\frac{f(n)}{n}=0.\] Let $d$ be the Kodaira-Iitaka dimension of $V_\sbullet$. We assume that $d\geqslant 0$ (namely $\mathbb N(V_\sbullet)=\{n\in\mathbb N\,:\,V_n\neq\mathbf{0}\}\neq\{0\}$) and we equip the graded algebra $V_\sbullet$ with a family of norms $(\|\ndot\|_n)_{n\in\mathbb N}$ such that $\overline{V}_{\!\sbullet}=(V_\sbullet,(\|\ndot\|_n)_{n\in\mathbb N_{\geqslant 1}})$ forms a normed graded algebra which is $f$-sub-mul\-ti\-pli\-ca\-tive and bounded (see \S\ref{Def: normed graded algebra}). For any $n\in\mathbb N_{\geqslant 1}$, let $\|\ndot\|_{\operatorname{sp},n}:V_n\rightarrow\mathbb R_{\geqslant 0}$ be the map defined as 
\[\|s\|_{\operatorname{sp},n}:=\lim_{N\rightarrow+\infty}\|s^N\|_{nN}^{1/N}.\]
Then $(V_\sbullet,(\|\ndot\|_{\operatorname{sp},n})_{n\in\mathbb N_{\geqslant 1}})$ forms a normed graded algebra which is sub-mul\-ti\-pli\-ca\-tive and bounded. Moreover, we denote by $\widehat{\mu}_{\max}^{\mathrm{asy}}(\overline V_{\!\sbullet})$ the \emph{asymptotic maximal slope}\index{slope!asymptotic maximal ---} of $\overline V_{\!\sbullet}$, which is defined as
\[\begin{split}\widehat{\mu}_{\max}^{\mathrm{asy}}(\overline V_{\!\sbullet})&=-\lim_{n\in\mathbb N(V_\sbullet),\,n\rightarrow+\infty}\min_{s\in V_n\setminus\{0\}}\frac{1}{n}\ln\|s\|_n\\&=\lim_{n\in\mathbb N(V_\sbullet),\,n\rightarrow+\infty}\widehat{\mu}_{\max}(V_n,\|\ndot\|_n).
\end{split}\]
Note that the existence of the limit is ensured by the inequality \eqref{Equ: sub multiplicative}, which implies that 
\[\widehat{\mu}_{\max}(V_{n_1+\cdots+n_\ell},\|\ndot\|_{n_1+\cdots+n_\ell})\geqslant\sum_{i=1}^{\ell}\Big(\widehat{\mu}_{\max}(V_{n_i},\|\ndot\|_{n_i})-f(n_i)\Big).\]
We refer the readers to \cite[Corollary 1.3.2]{MR2768967} for a proof of the convergence.

\begin{prop}\label{Pro: pente maximale asymptotique V}
The following equality holds:
\[\widehat{\mu}_{\max}^{\operatorname{asy}}(\overline{V}_{\!\sbullet})=\lim_{n\in\mathbb N(V_\sbullet),\,n\rightarrow+\infty}\frac 1n\widehat{\mu}_{\max}(V_n,\|\ndot\|_{\operatorname{sp},n}).\]
\end{prop}
\begin{proof}
By Proposition \ref{Pro: graded algebra}, one has
\[\|\ndot\|_{\operatorname{sp},n}\leqslant \mathrm{e}^{f(n)}\|\ndot\|_n\]
and hence for $n\in\mathbb N(V_\sbullet)$ the following inequality holds
\[\widehat{\mu}_{\max}(V_n,\|\ndot\|_{\operatorname{sp},n})\geqslant\widehat{\mu}_{\max}(V_n,\|\ndot\|_{n})-f(n).\] 
This implies 
\[\lim_{n\in\mathbb N(V_\sbullet),\,n\rightarrow+\infty}\frac 1n\widehat{\mu}_{\max}(V_n,\|\ndot\|_{\operatorname{sp},n})\geqslant \lim_{n\in\mathbb N(V_\sbullet),\,n\rightarrow+\infty}\frac 1n\widehat{\mu}_{\max}(V_n,\|\ndot\|_{n}).\]
Conversely, for any fixed $n\in\mathbb N(V_\sbullet)$ and $s\in V_n\setminus\{0\}$ such that 
\[\ln\|s\|_{\operatorname{sp},n}=-\widehat{\mu}_{\max}(V_n,\|\ndot\|_{\operatorname{sp},n}),\] 
one has 
\[\begin{split}\widehat{\mu}_{\max}^{\mathrm{asy}}(\overline{V}_{\!\sbullet})&=\lim_{N\rightarrow+\infty}\frac{1}{nN}\widehat{\mu}_{\max}(V_{nN},\|\ndot\|_{nN})\\&\geqslant\lim_{N\rightarrow+\infty}\frac{-1}{nN}\ln\|s^N\|_{nN}=-\frac{1}{n}\ln\|s\|_{\operatorname{sp},n}=\frac 1n\widehat{\mu}_{\max}(V_n,\|\ndot\|_{\operatorname{sp},n}).
\end{split}\]
Taking the limit when $n\rightarrow+\infty$, we obtain 
\[\widehat{\mu}_{\max}^{\operatorname{asy}}(\overline V_{\!\sbullet})\geqslant\widehat{\mu}_{\max}^{\operatorname{asy}}(V_\sbullet,(\|\ndot\|_{\operatorname{sp},n})_{n\in\mathbb N_{\geqslant 1}}).\]
\end{proof}

\begin{defi}
We define the \emph{arithmetic volume}\index{volume!arithmetic ---} of $\overline V_{\!\sbullet}$ as (see \S\ref{Subsec: deg and deg+} for the definition of $\widehat{\deg}_+$)
\begin{equation}\label{Equ: volume arithmetique}\operatorname{\widehat{\mathrm{vol}}}(\overline V_{\!\sbullet}):=\limsup_{n\in\mathbb N(V_\sbullet),\,n\rightarrow+\infty}\frac{\widehat{\deg}_+(V_n,\|\ndot\|_n)}{n^{d+1}/(d+1)!}.\end{equation}
\end{defi}

\begin{theo}\label{Thm: existence of volume}
The superior limit in the formula \eqref{Equ: volume arithmetique} defining the arithmetic volume function is actually a limit. Moreover, the following equalities hold:
\[\operatorname{\widehat{\mathrm{vol}}}(\overline V_{\!\sbullet})=\lim_{n\in\mathbb N(V_\sbullet),\,n\rightarrow+\infty}\frac{\widehat{\deg}_+(V_n,\|\ndot\|_{\operatorname{sp},n})}{n^{d+1}/(d+1)!}=(d+1)\int_0^{+\infty}\operatorname{vol}(V_\sbullet^t)\,\mathrm{d}t,\]
where for $t\in\mathbb R$,
\[V_\sbullet^t:=k\oplus\bigoplus_{n\in\mathbb N,\,n\geqslant 1}\operatorname{Vect}_k(\{s\in V_n\,:\,\|s\|_{\operatorname{sp},n}\leqslant\mathrm{e}^{-nt}\}).\]
\end{theo}
\begin{proof}  By replacing $k'$ by $k(V_\sbullet)$, we may assume that the graded linear series $V_\sbullet$ is birational. For simplifying the notation, we let $M$ be the asymptotic maximal slope of $\overline{V}_{\!\sbullet}$. Note that $M$ is also the asymptotic maximal slope of $(V_\sbullet,(\|\ndot\|_{\operatorname{sp},n})_{n\in\mathbb N})$ (see Proposition \ref{Pro: pente maximale asymptotique V}). Moreover, since $\overline{V}_{\!\sbullet}$ is bounded, there exists a constant $A\geqslant 0$ such that $\|s\|_{n}\leqslant\mathrm{e}^{nA}$ for any $n\in\mathbb N_{\geqslant 1}$ and any $s\in V_n$.

By the same argument as  the proof of \cite[Proposition 6.6]{MR4104565}, we obtain that, for any $t<M$, one has $k(V_\sbullet^t)=k(V_\sbullet)$. Moreover, for any $t>M$ and any $n\in\mathbb N_{\geqslant 1}$, one has $V_n^t=\mathbf{0}$. Therefore, combining the construction of Newton-Okounkov bodies in  \cite[Theorem 1.1]{MR3938628} and that of the concave transform developed in \cite[\S1.3]{MR2822867}, we obtain, in a similar way as \cite[Corollary 1.13]{MR2822867} that
\begin{multline*}\widehat{\vol}(V_\sbullet,(\|\ndot\|_{\operatorname{sp},n})_{n\in\mathbb N_{\geqslant 1}})=\lim_{n\in\mathbb N(V_\sbullet),\,n\rightarrow+\infty}\frac{\widehat{\deg}_+(V_n,\|\ndot\|_{\operatorname{sp},n})}{n^{d+1}/(d+1)!}\\
=(d+1)\int_0^{+\infty}\operatorname{vol}(V_\sbullet^t)\,\mathrm{d}t.\end{multline*}
Moreover, by \eqref{Equ: upper bound} we obtain that 
\[\widehat{\deg}_+(V_n,\|\ndot\|_{\operatorname{sp},n})\geqslant \widehat{\deg}_+(V_n,\|\ndot\|_{n})-\dim_k(V_n)f(n),\]
which leads to 
\[\limsup_{n\in\mathbb N(V_\sbullet),\,n\rightarrow+\infty}\frac{\widehat{\deg}_+(V_n,\|\ndot\|_n)}{n^{d+1}/(d+1)!}\leqslant \widehat{\vol}(V_\sbullet,(\|\ndot\|_{\operatorname{sp},n})_{n\in\mathbb N_{\geqslant 1}})\]
since $\dim_{k}(V_n)=O(n^d)$ when $n\in\mathbb N(V_\sbullet)$, $n\rightarrow+\infty$.

Let $\varepsilon$ be an element of $\mathopen{]}0,M\mathclose{[}$, $t$ be an element of $\mathopen{[}\varepsilon,M\mathclose{[}$. Let $W_\sbullet^t$ be a graded sub-$k$-algebra of finite type of $V_\sbullet^t$, which is generated by a family of homogeneous elements $s_1,\ldots,s_\ell$ of homogeneous degrees $n_1,\ldots,n_\ell$ respectively. For any $i\in \{1,\ldots,\ell\}$, there exists $a_i\in\mathbb N_{\geqslant 1}$ such that the inequalities
\begin{equation}\label{Equ: majoration norm si N}\|s_i^{N}\|_{n_iN}\leqslant \mathrm{e}^{n_iN\varepsilon/2}\|s_i\|_{\operatorname{sp},n_i}^N\leqslant\mathrm{e}^{n_iN(\varepsilon/2-t)}\end{equation} 
hold for any integer $N\geqslant a_i$. Therefore, by the inequality \eqref{Equ: sub multiplicative} we obtain that, for any $(N_1,\ldots,N_\ell)\in\mathbb N_{\geqslant 1}^\ell$, one has 
\[\ln\|s_1^{N_1}\cdots s_\ell^{N_\ell}\|_{n_1N_1+\cdots+n_\ell N_\ell}\leqslant\sum_{i=1}^\ell\big(\ln\|s_i^{N_i}\|_{n_iN_i}+f(n_iN_i)\big).\]
By \eqref{Equ: majoration norm si N}, we obtain that 
\begin{multline*} 
\ln\|s_1^{N_1}\cdots s_\ell^{N_\ell}\|_{n_1N_1+\cdots+n_\ell N_\ell} \\
\leqslant\sum_{\begin{subarray}{c}i\in\{1,\ldots,\ell\}\\
N_i\geqslant a_i
\end{subarray}}n_iN_i\Big(\frac {\varepsilon}2-t\Big)+\sum_{\begin{subarray}{c}i\in\{1,\ldots,\ell\}\\
N_i<a_i\end{subarray}}n_iN_iA\\
\leqslant \Big(\frac {\varepsilon}2-t\Big)\sum_{i=1}^\ell n_i(N_i-a_i)+\sum_{i=1}^\ell n_ia_iA \\
\leqslant \Big(\frac {\varepsilon}2-t\Big)\sum_{i=1}^\ell n_iN_i+\sum_{i=1}^\ell n_ia_i(A+M).
\end{multline*}
Therefore, for $(N_1,\ldots,N_\ell)\in \mathbb N_{\geqslant 1}^\ell$ such that $n_1N_1+\cdots+n_\ell N_\ell$ is sufficiently large, one has 
\[\|s_1^{N_1}\cdots s_\ell^{N_\ell}\|_{n_1N_1+\cdots+n_\ell N_\ell}\leqslant \mathrm{e}^{(\varepsilon-t)(n_1N_1+\cdots+n_\ell N_\ell)}.\]
In particular, for $n\in\mathbb N(V_\sbullet)$ sufficiently large, one has \[W_n^t\subset\mathcal F^{(t-\varepsilon)n}(V_n,\|\ndot\|_n),\] which leads to 
\[\liminf_{n\in\mathbb N(V_\sbullet),\,n\rightarrow+\infty}\frac{\dim_k(\mathcal F^{(t-\varepsilon)n}(V_n,\|\ndot\|_n))}{n^d/d!}\geqslant\operatorname{vol}(W^t_\sbullet).\]
Taking the supremum when $W_\sbullet^t$ varies, by the Fujita approximation property of $V_\sbullet^t$ we obtain that 
\begin{equation}\label{Equ: liminf bounded by vol}\liminf_{n\in\mathbb N(V_\sbullet),\,n\rightarrow+\infty}\frac{\dim_k(\mathcal F^{(t-\varepsilon)n}(V_n,\|\ndot\|_n))}{n^d/d!}\geqslant\operatorname{vol}(V^t_\sbullet).\end{equation}
Note that 
\begin{align*}
\widehat{\deg}_+(V_{n},\|\ndot\|_n)&=\int_0^{+\infty}\dim_k(\mathcal F^t(V_{n},\|\ndot\|_n))\,\mathrm{d}t \\
& =n\int_0^{+\infty}\dim_k(\mathcal F^{nt}(V_{n},\|\ndot\|_n))\,\mathrm{d}t\\
&\geqslant n\int_{\varepsilon}^M\dim_k(\mathcal F^{n(t-\varepsilon)}(V_{n},\|\ndot\|_n))\,\mathrm{d}t.
\end{align*}
Taking the integral with respect to $t$, by Fatou's lemma we deduce from \eqref{Equ: liminf bounded by vol} that 
\begin{multline*}\liminf_{n\in\mathbb N(V_\sbullet),n\rightarrow+\infty}\frac{\widehat{\deg}_+(V_{n},\|\ndot\|_n)}{n^{d+1}/(d+1)!} \\
\geqslant \liminf_{n\in\mathbb N(V_\sbullet),\,n\rightarrow+\infty}\frac{(d+1)!}{n^d}\int_{\varepsilon}^M\dim_k(\mathcal F^{n(t-\varepsilon)}(V_n,\|\ndot\|_n))\\
\geqslant (d+1)\int_{\varepsilon}^{M}\mathrm{vol}(V_\sbullet^t)\,\mathrm{d}t=(d+1)\int_{\varepsilon}^{+\infty}\operatorname{vol}(V_\sbullet^t).
\end{multline*}
Finally, taking the supremum with respect to $\varepsilon$, we obtain the inequality 
\[\liminf_{n\in\mathbb N(V_\sbullet),n\rightarrow+\infty}\frac{\widehat{\deg}_+(V_{n},\|\ndot\|_n)}{n^{d+1}/(d+1)!}\geqslant\widehat{\vol}(V_\sbullet,(\|\ndot\|_{\operatorname{sp},n})_{n\in\mathbb N_{\geqslant 1}}).\]
The theorem is thus proved. 
\end{proof}

\begin{coro}\label{Cor: volume chi}
The sequences 
\[\frac{\widehat{\deg}(V_n,\|\ndot\|_n)}{n^{d+1}/(d+1)!},\quad n\in\mathbb N(V_\sbullet)\]
and 
\[\frac{\widehat{\deg}(V_n,\|\ndot\|_{\operatorname{sp},n})}{n^{d+1}/(d+1)!},\quad n\in\mathbb N(V_\sbullet)\]
converge to the same real number, which is equal to 
\[-\int_{\mathbb R}\,t\,\mathrm{d}\operatorname{vol}(V_\sbullet^t).\]
\end{coro}
\begin{proof}
Let $A$ be a positive constant such that $\|s\|_n\leqslant \mathrm{e}^{nA}$ for any $n\in\mathbb N_{\geqslant 1}$ and any $s\in V_n$. For any $n\in\mathbb N_{\geqslant 1}$, let $\|\ndot\|_{n}'=\mathrm{e}^{-nA}\|\ndot\|_n$. Then, $(V_\sbullet,(\mathrm{e}^{-nA}\|\ndot\|_n)_{n\in\mathbb N_{\geqslant 1}})$ forms a normed graded algebra over $(k,|\ndot|_0)$, which is $f$-sub-multiplicative and bounded. Moreover, for any $n\in\mathbb N_{\geqslant 1}$, one has 
\[\widehat{\deg}(V_n,\|\ndot\|_n')=\widehat{\deg}_+(V_n,\|\ndot\|_n')=nA\dim_k(V_n)+\widehat{\deg}(V_n,\|\ndot\|_{n}),\]
where the first equality comes from the fact that the image of $\|\ndot\|_n'$ is contained in $[0,1]$. For any $n\in\mathbb N$ one has 
\[\|\ndot\|'_{\operatorname{sp},n}=\mathrm{e}^{-nA}\|\ndot\|_{\operatorname{sp},n}.\]
By \eqref{Equ: upper bound}, for any $n\in\mathbb N_{\geqslant 1}$ and any $s\in V_n$, one has 
\[\forall\,N\in\mathbb N_{\geqslant 1},\quad \|s\|_{\operatorname{sp},n}=\|s^N\|^{1/N}_{\operatorname{sp},nN}\leqslant \mathrm{e}^{f(nN)/N}\|s^N\|_{nN}^{1/N}\leqslant \mathrm{e}^{f(nN)/N+nA}.\]
Taking the limit when $N\rightarrow+\infty$, we obtain $\|s\|_{\operatorname{sp},n}\leqslant \mathrm{e}^{nA}$ and hence $\|\ndot\|'_{\operatorname{sp},n}$ also takes value in $[0,1]$. Therefore, for any $n\in\mathbb N_{\geqslant 1}$, one has 
\[\widehat{\deg}(V_n,\|\ndot\|_{\operatorname{sp},n}')=\widehat{\deg}_+(V_n,\|\ndot\|_{\operatorname{sp},n}')=nA\dim_k(V_n)+\widehat{\deg}(V_n,\|\ndot\|_{\operatorname{sp},n}),\]
Hence Theorem \ref{Thm: existence of volume} leads to the convergence of the sequences  
\[\frac{\widehat{\deg}(V_n,\|\ndot\|_n)+nA\dim_k(V_n)}{n^{d+1}/(d+1)!},\quad n\in\mathbb N(V_\sbullet) \]
and 
\[\frac{\widehat{\deg}(V_n,\|\ndot\|_{\operatorname{sp},n})+nA\dim_k(V_n)}{n^{d+1}/(d+1)!},\quad n\in\mathbb N(V_\sbullet)\]
to the same limit, which is equal to 
\[\begin{split}&\quad\;(d+1)\int_0^{+\infty}\operatorname{vol}(V_\sbullet^{t-A})\,\mathrm{d}t=(d+1)\int_{-A}^{+\infty}\operatorname{vol}(V_\sbullet^{t})\,\mathrm{d}t\\
&=A(d+1)\operatorname{vol}(V_\sbullet)-\int_{\mathbb R}t\,\mathrm{d}\operatorname{vol}(V_\sbullet^t),
\end{split}\]
where the last equality comes from the fact that $V_\sbullet^t=V_\sbullet$ when $t\leqslant -A$. By the formula
\[\lim_{n\in\mathbb N(V_\sbullet),\,n\rightarrow+\infty}\frac{\dim_k(V_n)}{n^d/d!}=\vol(V_\sbullet),\]
we obtain the assertion. 
\end{proof}

\begin{defi}
We define the \emph{$\chi$-volume}\index{volume!$\chi$- ---} of the normed graded linear series $\overline{V}_{\!\sbullet}$ as
\[\widehat{\vol}_\chi(\overline{V}_{\!\sbullet})=\lim_{n\in\mathbb N(V_\sbullet),\,n\rightarrow+\infty}\frac{\widehat{\deg}(V,\|\ndot\|_n)}{n^{d+1}/(d+1)!}.\]
By Corollary \ref{Cor: volume chi}, we obtain that $\widehat{\vol}_{\chi}(\overline V_{\!\sbullet})=\widehat{\vol}_{\chi}(V_\sbullet,(\|\ndot\|_{\operatorname{sp},n})_{n\in\mathbb N_{\geqslant 1}})$. 
\end{defi}

%!TEX root = ./Hilbert_Samuel_Adelic_Curves.tex

\chapter{Arithmetic volumes over a general adelic curve}

In this chapter, we use the results of  the previous chapter to study the volume functions of a normed graded algebra over a general adelic curve. Let $S=(K,(\Omega,\mathcal A,\nu),\phi)$ be the adelic curve defined in \S\ref{Sec:adelic curve}.  We let $|\ndot|_0$ be the trivial absolute value on $K$, and denote by $S_0=(K,\{0\},|\ndot|_0)$ the adelic curve consisting of a single copy of the trivial absolute value $|\ndot|_0$ on $K$. 

\section{Graded algebra of adelic vector bundles}

In this section, we consider basic facts on graded algebras of adelic vector bundles.

\begin{defi}\label{Def:graded algebra in adelic vector bundle}
Let $E_\sbullet=\bigoplus_{n\in\mathbb N}E_n$ be a graded $K$-algebra. We assume that each vector space $E_n$ is finite-dimensional over $K$. For any $n\in\mathbb N$, let $\xi_n=(\|\ndot\|_{n,\omega})_{\omega\in\Omega}$ be a norm family on $E_n$ such that $\overline E_n=(E_n,\xi_n)$ forms an adelic vector bundle on $S$. We call $\overline E_{\sbullet}=(\overline E_n)_{n\in\mathbb N}$ a \emph{graded algebra of adelic vector bundles on $S$.}\index{graded algebra of adelic vector bundles}
For any $n\in\mathbb N$ such that $n\geqslant 1$,  let $(\mathcal F^t(\overline E_{n}))_{t\in\mathbb R}$ be the Harder-Narasimhan $\mathbb R$-filtration of $\overline{V}_{\!n}$ (see \S\ref{Sec: HN R filtration}). We denote by $\|\ndot\|_n^{\operatorname{HN}}$ the norm on $E_n$ (viewed as a vector space over $(K,|\ndot|_0)$) defined as 
\[\forall\,s\in E_n,\quad \|s\|_n^{\operatorname{HN}}=\exp\big(-\sup\{t\in\mathbb R\,:\,s\in\mathcal F^t(\overline E_n)\}\big).\]
Then, the couple $(E_\sbullet,(\|\ndot\|_n^{\operatorname{HN}})_{n\in\mathbb N_{\geqslant 1}})$ forms a normed graded algebra over $(K,|\ndot|_0)$ (see \S\ref{Def: normed graded algebra}). Moreover, if we view $(E_n,\|\ndot\|_{n}^{\operatorname{HN}})$ as  an adelic vector bundle on $S_0$, then its Harder-Narasimhan filtration coincides with that of $(E_n,\|\ndot\|_n)$. In particular, by the results recalled in \S\ref{Sec: HN R filtration}, the following estimates hold:
\begin{gather}0\leqslant\widehat{\deg}(E_n,\|\ndot\|_n)-\widehat{\deg}(E_n,\|\ndot\|_n^{\operatorname{HN}})\leqslant\frac 12\nu(\Omega_\infty)\dim_K(E_n)\ln(\dim_K(E_n)),\label{Equ: comparision deg}\\
0\leqslant \widehat{\deg}_+(E_n,\|\ndot\|_n)-\widehat{\deg}_+(E_n,\|\ndot\|_n^{\operatorname{HN}})\leqslant\frac 12\nu(\Omega_\infty)\dim_K(E_n)\ln(\dim_K(E_n)).\label{Equ: comparision deg+}
\end{gather}

Let $b=(b_n)_{n\in\mathbb N_{\geqslant 1}}$ be a sequence of non-negative integrable functions on $(\Omega,\mathcal A,\nu)$. 
We say that a graded 
algebra of adelic vector bundles $\overline E_{\sbullet}$ 
is \emph{$b$-sub-multiplicative}\index{sub-multiplicative!$b$- ---} if for all $\omega\in\Omega$, $\ell\in\mathbb N_{\geqslant 2}$, $(n_1,\ldots,n_\ell)\in\mathbb N_{\geqslant 1}^\ell$ and $(s_1,\ldots,s_\ell)\in E_{n_1,\omega}\times\cdots\times E_{n_\ell,\omega}$, the following inequality holds
\begin{equation}\label{Equ: submultiplicativity}\|s_1\cdots s_\ell\|_{n_1+\cdots+n_\ell,\omega}\leqslant\mathrm{e}^{b_{n_1}(\omega)+\cdots+b_{n_\ell}(\omega)}\|s_1\|_{n_1,\omega}\cdots\|s_\ell\|_{n_\ell,\omega}.\end{equation}
If for any $n$, $b_n$ is the constant function taking $0$ as its value, we simply say that $\overline E_{\sbullet}$ is \emph{sub-multiplicative}\index{sub-multiplicative}.
\end{defi}

\begin{prop}\label{Pro: f sub multiplicative}
Assume that $K$ is perfect. Let $b=(b_n)_{n\in\mathbb N_{\geqslant 1}}$ be a sequence of non-negative integrable functions on $(\Omega,\mathcal A,\nu)$, and  $\overline E_\sbullet$ be a graded algebra of adelic vector bundles on $S$, which is $b$-sub-multiplicative. Let $f:\mathbb N_{\geqslant 1}\rightarrow\mathbb R_{\geqslant 0}$ be the function defined as
\[f(n)=\frac{3}{2}\nu(\Omega_\infty)\ln(\dim_K(E_n))+\int_{\Omega}b_n(\omega)\,\nu(\mathrm{d}\omega).\]
Then the normed graded algebra $(E_\sbullet,(\|\ndot\|_n^{\operatorname{HN}})_{n\in\mathbb N_{\geqslant 1}})$ is $f$-sub-multiplicative.
\end{prop}
\begin{proof}
Let $\ell\in\mathbb N_{\geqslant 1}$ and $(n_1,\ldots,n_\ell)\in\mathbb N_{\geqslant 1}^\ell$. For any $i\in\{1,\ldots,\ell\}$, let $F_{n_i}$ be a $K$-vector subspace of $E_{n_i}$. For any $\omega\in\Omega$, we consider the $K_\omega$-linear map 
\[F_{n_1,\omega}\otimes\cdots\otimes F_{n_\ell,\omega}\longrightarrow E_{n_1+\cdots+n_\ell,\omega}\]
induced by the $K$-algebra structure of $E_\sbullet$. If we equip with $F_{n_1,\omega}\otimes\cdots\otimes F_{n_\ell,\omega}$ with the $\varepsilon$-tensor product of the norms $\|\ndot\|_{n_1,\omega},\ldots,\|\ndot\|_{n_\ell,\omega}$ when $|\ndot|_\omega$ is non-Archimedean, and with the $\pi$-tensor product when $|\ndot|_\omega$ is Archimedean, then the operator norm of the above map is bounded from above by $\exp(b_{n_1}(\omega)+\cdots+b_{n_\ell}(\omega))$. Moreover, by \cite[Corollary 5.6.2]{CMArakelovAdelic} (Although this result has been stated under the assumption that $\operatorname{char}(K)=0$, this assumption is only used in the application of  \cite[Theorem 5.4.3]{CMArakelovAdelic}, which actually applies to any perfect field. Moreover, the lifting of invariants from the symmetric power to the tensor power, that we have used to prove \cite[Proposition  5.3.1]{CMArakelovAdelic}, is valid in any characteristic. For details, see Remark \ref{Rem: positive characteristic minimal slope}.), 
one has
\[\widehat{\mu}_{\min}(\overline F_{\!n_1}\otimes_{\varepsilon,\pi}\cdots\otimes_{\varepsilon,\pi}\overline F_{\!n_\ell})\geqslant\sum_{i=1}^\ell\Big(\widehat{\mu}_{\min}(\overline F_{\!n_i})-\frac 32\nu(\Omega_\infty)\ln(E_{n_i})\Big).\] Let $F_{n_1+\cdots+n_\ell}$ be the image of the map 
\[F_{n_1}\otimes\cdots\otimes F_{n_\ell}\longrightarrow E_{n_1+\cdots+n_\ell}.\]
By \cite[Proposition 4.3.31]{CMArakelovAdelic}, we obtain that 
\begin{equation}\label{Equ: mu min F n1}\widehat{\mu}_{\min}(\overline F_{\! n_1+\cdots+n_\ell})\geqslant\sum_{i=1}^\ell\Big(\widehat{\mu}_{\min}(\overline F_{\!n_i})-\frac 32\nu(\Omega_\infty)\ln(E_{n_i})-\int_{\Omega}b_{n_i}(\omega)\,\nu(\mathrm{d}\omega)\Big).\end{equation}
Therefore, we obtain that, for any $(t_1,\ldots,t_\ell)\in\mathbb R^\ell$, one has  
\[\mathcal F^{t_1}(\overline E_{n_1})\cdots\mathcal F^{t_\ell}(E_{n_\ell})\subset\mathcal F^{t_1+\cdots+t_\ell-f(n_1)-\cdots-f(n_\ell)}(E_{n_1+\cdots+n_\ell}),\]
which shows that 
the normed graded algebra $(E_\sbullet,(\|\ndot\|_n^{\operatorname{HN}})_{n\in\mathbb N_{\geqslant 1}})$ is 
$f$-sub-mul\-ti\-pli\-ca\-tive.
\end{proof}

\ifsmf\begin{enonce}{Corollary-Definition}\fi
\ifams\begin{code}\fi
\label{Cor_Def_volchi} Assume that the field $K$ is perfect. 
Let $b=(b_n)_{n\in\mathbb N_{\geqslant 1}}$ be a sequence of non-negative integrable functions on $(\Omega,\mathcal A,\nu)$ such that 
\[\lim_{n\rightarrow+\infty}\frac{1}{n}\int_{\Omega}b_n(\omega)\,\nu(\mathrm{d}\omega)=0.\]
Let $\overline E_{\sbullet}$ be a graded algebra of adelic vector bundles on $S$, which is $b$-sub-multiplicative. Denote by $\mathbb N(E_\sbullet)$ the set of $n\in\mathbb N$ such that $E_n\neq\boldsymbol{0}$. Assume that \begin{enumerate}[label=\rm(\arabic*)]\item $E_\sbullet$ is isomorphic to a graded linear series of sub-finite type of a finitely generated extension of $K$, which is of Kodaira-Iitaka dimension $d\geqslant 0$,
\item there exists $C>0$ such that, for any $n\in\mathbb N(E_\sbullet)$,
\[-Cn\leqslant\widehat{\mu}_{\min}(\overline E_n)\leqslant\widehat{\mu}_{\max}(\overline E_n)\leqslant Cn. \]
\end{enumerate}  
Then the sequences
\[
\frac{\widehat{\deg}(\overline E_n)}{n^{d+1}/(d+1)!},\quad n\in\mathbb N(E_\sbullet)\]
and\[
\frac{\widehat{\deg}_+(\overline E_n)}{n^{d+1}/(d+1)!},\quad n\in\mathbb N(E_\sbullet)
\]
converge to two real numbers $\widehat{\vol}_{\chi}(\overline E_\sbullet)$ and $\widehat{\vol}(\overline E_\sbullet)$, which we call \emph{$\chi$-volume}\index{volume!$\chi$- ---} and \emph{volume}\index{volume} of $\overline E_\sbullet$, respectively.
\ifsmf\end{enonce}\fi
\ifams\end{code}\fi
\begin{proof}
These results follow from Proposition \ref{Pro: f sub multiplicative}, Theorem \ref{Thm: existence of volume}, Corollary \ref{Cor: volume chi} and the comparisons \eqref{Equ: comparision deg}, \eqref{Equ: comparision deg+} and the convergence of the sequence
\[\frac{\dim_K(E_n)}{n^d/d!},\quad n\in\mathbb N(E_\sbullet).\] 
\end{proof}

\begin{rema}\label{Rem: bounded mu min}
Assume that the field $K$ is perfect.  Let $b=(b_n)_{n\in\mathbb N_{\geqslant 1}}$ be a sequence of non-negative integrable functions on $(\Omega,\mathcal A,\nu)$ such that 
\[\lim_{n\rightarrow+\infty}\frac{1}{n}\int_{\Omega}b_n(\omega)\,\nu(\mathrm{d}\omega)=0.\]
Let $\overline E_{\sbullet}$ be a graded algebra of adelic vector bundles on $S$, which is $b$-sub-multiplicative. We assume that $n_1,\ldots,n_\ell$ are elements of $\mathbb N(E_\sbullet)\setminus\{0\}$ such that 
\[K\oplus\bigoplus_{n\in\mathbb N,\,n\geqslant 1}E_n\] is generated as $K$-algebra by $E_{n_1}\cup\cdots\cup E_{n_\ell}$. By \eqref{Equ: mu min F n1} we obtain that, for any $(a_1,\ldots,a_\ell)\in\mathbb N^\ell\setminus\{(0,\ldots,0)\}$, the canonical image of \[E_{n_1}^{\otimes a_1}\otimes\cdots\otimes E_{n_\ell}^{\otimes a_\ell}\]
in $E_{a_1n_1+\cdots+a_{\ell}n_\ell}$
has a minimal slope
\[\geqslant\sum_{i=1}^\ell a_i\bigg(\widehat{\mu}_{\min}(\overline E_i)-\frac 32\nu(\Omega_\infty)\ln(E_{n_i})-\int_{\Omega}b_{n_i}(\omega)\,\nu(\mathrm{d}\omega)\bigg).\]
Therefore we deduce that, for any $n\in\mathbb N(E_\sbullet)\setminus\{0\}$, the minimal slope of $\overline E_n$ is bounded from below by 
\[\min_{\begin{subarray}{c}
(a_1,\ldots,a_\ell)\in\mathbb N^\ell\\
n=a_1n_1+\cdots+a_\ell n_\ell
\end{subarray}}\sum_{i=1}^\ell a_i\bigg(\widehat{\mu}_{\min}(\overline E_i)-\frac 32\nu(\Omega_\infty)\ln(E_{n_i})-\int_{\Omega}b_{n_i}(\omega)\,\nu(\mathrm{d}\omega)\bigg).\]   Hence there exists $C>0$ such that $\widehat{\mu}_{\min}(\overline E_n)\geqslant -Cn$ holds for any $n\in\mathbb N(E_\sbullet)$.
\end{rema}

\section{Arithmetic $\chi$-volumes of adelic line bundles}

In this subsection, we introduce the arithmetic $\chi$-volume of an adelic line bundle.
\ifsmf\begin{enonce}{Theorem-Definition}\fi
\ifams\begin{thde}\fi
\label{theo:limit:thm:vol:chi}
Let $p:X\rightarrow\Spec K$ be an integral projective scheme over $\Spec K$, $d$ be the dimension of $X$, and $\overline L=(L,\varphi)$ be an  adelic line bundle on $X$. We suppose that, either $K$ is perfect, or $X$ is geometrically integral. Assume that $L$ is big and the graded $K$-algebra
\[\bigoplus_{n\in\mathbb N}H^0(X,L^{\otimes n})\]
is of finite type. We denote the adelic vector bundle \[\big(H^0(X, L^{\otimes n}), ( \|\ndot\|_{n\varphi_{\omega}} )_{\omega \in \Omega}\big)\] over $S$ by $p_*(\overline L{}^{\otimes n})$.
Then the sequence 
\begin{equation}\label{Equ: convergence vol chi}\frac{\widehat{\deg}(p_*(\overline L{}^{\otimes n}))}{n^{d+1}/(d+1)!},\quad n\in\mathbb N,\;n\geqslant 1\end{equation}
converges to a real number, which we denote by $\widehat{\vol}_{\chi}(\overline L)$ and which we call the \emph{$\chi$-volume}\index{volume!$\chi$- ---} of $\overline L$.
\ifsmf\end{enonce}\fi
\ifams\end{thde}\fi
\begin{proof}
Let $K^{\operatorname{pc}} $ be the perfect closure of $K$. Recall that, if $K^{\operatorname{ac}}$ denotes the algebraic closure of $K$, then $K^{\operatorname{pc}}$ is the intersection of all subfields of $K^{\operatorname{ac}}$ containing $K$ which are perfect fields. Note that $K^{\operatorname{pc}}/K$ is a purely inseparable algebraic extension of $K$. Therefore, for any $\omega\in\Omega$, the absolute value $|\ndot|_{\omega}$ extends in a unique way to $K^{\operatorname{pc}}/K$. In other words, the measure space in the adelic curve structure of $S\otimes_KK^{\operatorname{pc}}$ coincides with $(\Omega,\mathcal A,\nu)$.

For any $n\in\mathbb N$, let 
\[E_n=H^0(X,L^{\otimes n})\otimes_KK^{\operatorname{pc}}=H^0(X_{K^{\operatorname{pc}}},L_{K^{\operatorname{pc}}}^{\otimes n}).\]
The norm family of $p_*(\overline L^{\otimes n})$ induces by extension of scalars a norm family on  $E_n$, which we denote by $\xi_n$. By \cite[Proposition 4.3.14]{CMArakelovAdelic}, the equality 
\[\widehat{\deg}(E_n,\xi_n)=\widehat{\deg}(p_*(\overline L{}^{\otimes n}))\]
holds. Moreover, 
\[E_\sbullet=\bigoplus_{n\in\mathbb N} E_n\]
is a graded $K^{\operatorname{pc}}$-algebra of finite type, which is isomorphic to a graded linear series of the function field of $X_{K^{\operatorname{pc}}}$ over $K^{\operatorname{pc}}$. As a graded $K^{\operatorname{pc}}$-algebra of adelic vector bundles on $S\otimes_KK^{\operatorname{pc}}$, $\overline E_\sbullet=(\overline E_n)_{n\in\mathbb N}$ is sub-multiplicative. By \cite[Proposition 6.2.7]{CMArakelovAdelic}, we obtain, following the proof of \cite[Proposition 6.4.4]{CMArakelovAdelic}, that the sequence 
\[\frac{\widehat{\mu}_{\max}(\overline E_{n})}{n},\quad n\in\mathbb N,\;n\geqslant 1\] is bounded from above. Therefore the assertion follows from Corollary-Definition \ref{Cor_Def_volchi} (see also Remark \ref{Rem: bounded mu min}).
\end{proof}

\begin{rema}\label{Rem: limit theorem in average}
Under the notation and the assumption of the above theorem-definition,  the following relation holds
\[\lim_{n\rightarrow+\infty}\frac{\widehat{\deg}(p_*(\overline L{}^{\otimes n}))}{n\dim_K(H^0(X,L^{\otimes n}))}=\frac{\widehat{\operatorname{vol}}_{\chi}(\overline L)}{(d+1)\operatorname{vol}(L)}.\]
\end{rema}

\section{Normed graded module}
Let $\overline{R}_\sbullet=(\overline R_{n})_{n\in\mathbb N}$ be a graded algebra of adelic vector bundles on $S$, where $\overline R_n=(R_n,(\|\ndot\|_{n,\omega})_{\omega\in\Omega})$. Let $M_\sbullet=\bigoplus_{n\in\mathbb N}M_n$ be a graded module over $R_\sbullet=\bigoplus_{n\in\mathbb N}R_n$. If each $M_n$ is a finite-dimensional vector space over $K$ and is equipped with a norm family $(\|\ndot\|^M_{n,\omega})_{\omega\in\Omega}$ such that $\overline M_{\!n}=(M_n,(\|\ndot\|^M_{n,\omega})_{\omega\in\Omega})$ is an adelic vector bundle on $S$, we say that $\overline M_{\!\sbullet}=(\overline M_{\!n})_{n\in\mathbb N}$ is \emph{a graded $\overline R_{\sbullet}$-module of adelic vector bundles on $S$}\index{graded module of adelic vector bundles}.  

Assume that $\overline R_{\sbullet}$ is sub-multiplicative (see Definition \ref{Def:graded algebra in adelic vector bundle}). If, for all $(n,m)\in\mathbb N^2$, $\omega\in\Omega$ and $(a,s)\in R_{n,\omega}\times M_{m,\omega}$, one has
\[\|as\|_{n+m,\omega}^M\leqslant \|a\|_{n,\omega}\cdot\|s\|_{m,\omega}^M,\]
we say that $\overline M_{\!\sbullet}$ is \emph{sub-multiplicative}\index{sub-multiplicative}.

\begin{lemm}\label{lem:sub:quot:R:module}
Let $\overline M_{\!\sbullet}=((M_n,\xi_{M_n}))_{n\in\mathbb N}$ be a graded $\overline R_{\sbullet}$-module of adelic vector bundle on $S$. Let $Q = \bigoplus_{n=0}^{\infty} Q_n$ be a graded quotient $R$-module of $M$, that is, $Q_n$ is a quotient vector space of $M_n$ over $K$ for all $n$ and
$a_\ell\cdot: M_n \to M_{n+\ell}$ induces by passing to quotient $a_{\ell}\cdot : Q_n \to Q_{n+\ell}$ for $a_\ell \in R_{\ell}$. Let $\xi_{Q_n}$ be the quotient norm family of $Q_n$
induced by $M_n \to Q_n$ and $\xi_{M_n}$. Then $\overline Q_{\sbullet}=((Q_n,\xi_{Q_n}))_{n\in\mathbb N}$ is a graded $\overline R_{\sbullet}$-algebra.
\end{lemm}

\begin{proof}
Assume that $\xi_{M_n}$ and $\xi_{Q_n}$ are of the form $(\|\ndot\|_{n,\omega}^M)_{\omega\in\Omega}$ and $(\|\ndot\|_{n,\omega}^Q)_{\omega\in\Omega}$, respectively.
Let $(n,n')\in\mathbb N^2$, $\omega\in\Omega$, $a\in R_{n,\omega}$ and $q\in Q_{n',\omega}$. For any $s\in M_{n',\omega}$ which represents the class $q\in Q_{n',\omega}$, one has 
\[\|aq\|^Q_{n+n',\omega}\leqslant\|as\|^M_{n+n',\omega}\leqslant\|a\|_{n,\omega}\cdot\|s\|^M_{n',\omega}.\]
Taking the infimum with respect to $s$, we obtain
\[\|aq\|_{n+n',\omega}^Q\leqslant\|a\|_{n,\omega}\cdot\|s\|_{n',\omega}^Q,\]
as required.
\end{proof}

\begin{prop}\label{prop:positivity:torsion:sheaf}
Suppose that $R_\sbullet$ is a $K$-algebra of finite type. Let $\overline M_{\sbullet}=((M_n,\xi_{M_n}))_{n\in\mathbb N}$ be a graded $\overline R_{\sbullet}$-module of adelic vector bundle on $S$, such that $M_\sbullet$ is an $R_\sbullet$-module of finite type.
Suppose that \[\liminf\limits_{n\to\infty} \frac{\dim_K(M_n)}{n^d} = 0\] for some non-negative integer $d$, then 
\[\liminf_{n\to\infty} \frac{\adeg(M_n, \xi_{M_n})}{n^{d+1}} \geqslant 0.\]
\end{prop}

\begin{proof} Let $x_1,\ldots,x_r$ be homogeneous elements of $R$ which generate $R$ as $K$-algebra.
We choose non-zero homogeneous elements $m_1, \ldots, m_{\ell}$ of $M$
such that $M$ is generated by $m_1, \ldots, m_{\ell}$ over $R$.
We  set $e_i = \deg(x_i)$ and $f_i = \deg (m_i)$ for $i \in \{ 1, \ldots, r \}$.
For $\alpha = (a_1,\ldots, a_r) \in \mathbb N^r$, we denote $x_1^{a_1} \cdots x_r^{a_r}$ by $x^{\alpha}$.
If we set $d_n = \dim_K (M_n)$, then, for $n \geqslant \max \{ f_1, \ldots, f_r \}$,
we can find $\alpha_1, \ldots, \alpha_{d_n} \in \mathbb N^r$ and
$m_{i_1}, \ldots, m_{i_{d_n}} \in \{ m_1, \ldots, m_\ell\}$ such that
$x^{\alpha_1}m_{i_1}, \ldots, x^{\alpha_{d_n}} m_{i_{d_n}}$ form a basis of $M_n$. Note that
\begin{align*}
\| (x^{\alpha_1}m_{i_1}) \wedge\cdots \wedge (x^{\alpha_{d_n}} m_{i_{d_n}}) \|^M_{n, \omega, \det} & \leqslant
\| x^{\alpha_1}m_{i_1} \|_{n,\omega}^{M} \cdots \| x^{\alpha_{d_n}} m_{i_{d_n}} \|^M_{n, \omega} \\
& \kern-14em \leqslant \| x^{\alpha_1} \|_{n - f_{i_1}, \omega} \cdots \| x^{\alpha_{d_n}} \|_{n - f_{i_{d_n}},\omega} \cdot\| m_{i_1} \|^M_{f_{i_1},\omega} \cdots
 \| m_{i_{d_n}} \|^M_{f_{i_{d_n}},\omega} \\
& \kern-14em \leqslant \max \{ 1, \|x_1\|_{e_1,\omega}, \ldots , \|x_r\|_{e_{r}, \omega} \}^{nd_n} \max \{ 1, \|m_1\|^M_{f_1,\omega}, \ldots , \|m_\ell\|^M_{f_{\ell}, \omega} \}^{d_n},
\end{align*}
so that
\begin{multline*}
\adeg(M_n, \xi_{M_n}) \geqslant
nd_n \int_{\Omega} \min \{ 0, - \ln \|x_1\|_{e_1,\omega}, \ldots , - \ln \|x_r\|_{e_{r}, \omega} \}\, \nu(\mathrm{d}\omega)\\
+
d_n  \int_{\Omega}  \min \{ 0, - \ln\|m_1\|^M_{f_1,\omega}, \ldots , - \ln\|m_\ell\|^M_{f_{\ell}, \omega} \}\, \nu(\mathrm{d}\omega).
\end{multline*}
Thus the assertion follows.
\end{proof}

\section{Bounds of $\chi$-volume with auxiliary torsion free module}

Let us begin with the following lemma.

\begin{lemm}\label{lemma:finitely:generated:norm:space}
Let $X$ be an integral projective scheme over a field $k$, $L$
be an invertible $\OO_X$-module and $F$ be a coherent $\OO_{X}$-module.
We assume that there exist a surjective morphism $f : X \to Y$ of integral projective schemes over $k$ and an
ample invertible $\OO_Y$-module $A$ such that $f^*(A) = L$.
Then $R = \bigoplus_{n=0}^{\infty} H^0(X, L^{\otimes n})$ is a finitely generated algebra over $k$ and
$M = \bigoplus_{n=0}^\infty H^0(X, F \otimes L^{\otimes n})$ is a finitely generated $R$-module.
\end{lemm}

\begin{proof}
By \cite[\S1.8]{MArakelov}, there exist positive integers $d$ and $n_0$ such that
\[
H^0(Y, A^{\otimes d}) \otimes H^0(Y, A^{\otimes n} \otimes f_*(F)) \longrightarrow H^0(Y, A^{\otimes (d + n)} \otimes f_*(F)) 
\]
is surjective for all $n \geqslant n_0$, and hence
\[
H^0(X, L^{\otimes d}) \otimes H^0(X, L^{\otimes n} \otimes F) \to H^0(X, L^{\otimes (d + n)} \otimes F) 
\]
is surjective for all $n \geqslant n_0$ because $f_*(L^{\otimes n}) = A^{\otimes n} \otimes f_*(\OO_X)$, $f_*(L^{\otimes n} \otimes F) = A^{\otimes n} \otimes f_*(F)$,
$\OO_Y \subseteq f_*(\OO_X)$. Thus, by the arguments  in \cite[\S1.8]{MArakelov}, one can see the assertion.
\end{proof}

In the rest of the section, let $p:X\rightarrow\Spec K$ be a $d$-dimensional geometrically integral projective variety over $K$.
Let $\overline{L} = (L, \varphi)$ be an adelic invertible $\OO_X$-module.
Let $E$ be a torsion free $\OO_X$-module and $U$ be a non-empty Zariski open set of $X$ such that $E|_U$ is a vector bundle.
Let $\psi=(\psi_\omega)_{\omega\in\Omega}$ be a metric family of $E|_U$. We assume that $( L^{\otimes n}\otimes E, U, n \varphi|_U+\psi)$ is a sectionally adelic torsion free $\OO_X$-module (see Definition \ref{def:globally:adelic:torsion:free}) for all $n\in\mathbb N$. Note that, if the sectional algebra $\bigoplus_{n\in\mathbb N}H^0(X,L^{\otimes n})$ is of finite type over $K$ (this condition is true notably when $L$ satisfies the hypothesis of  Lemma  \ref{lemma:finitely:generated:norm:space}), by Theorem-Definition \ref{theo:limit:thm:vol:chi}, the sequence 
\[\frac{\widehat{\deg}(p_*(\overline L{}^{\otimes n}))}{n^{d+1}/(d+1)!},\quad n\in\mathbb N,\;n\geqslant 1\]
converges to a real number denoted by $\widehat{\vol}_{\chi}(\overline L)$.

\begin{theo}\label{thm:lower:bound:chi:volume} If there are a birational morphism $f : X \to Z$ of geometrically integral projective schemes over $\Spec K$ and an ample invertible $\OO_Z$-module
$A$ such that $L = f^{*}(A)$, then the following inequality holds:
\[ \rank (E) \operatorname{\widehat{\vol}_{\chi}}(\overline L) \leqslant \liminf_{n\to\infty} \frac{\widehat{\deg}(p_*(\overline L{}^{\otimes n}\otimes \overline E))}{n^{d+1}/(d+1)!}.\]
\end{theo}
\begin{proof}
Let $r$ be the rank of $E$.
Note that $p_*(\overline L{}^{\otimes n}\otimes \overline E)$ forms an adelic vector bundle over $S$ for any $n \in\mathbb N$.
For a sufficiently large positive integer $n_0$,  shrinking $U$ if necessarily,
we can find  
$e_1, \ldots, e_r \in H^0(X, L^{\otimes n_0} \otimes E)$ such that
$e_1, \ldots, e_r$ yield a basis of $L^{\otimes n_0} \otimes E$ over $U$.
Indeed, there is a positive integer $n_0$ such that \[H^0(Z, A^{\otimes n_0} \otimes f_*(E)) \otimes \OO_Z \longrightarrow A^{\otimes n_0} \otimes f_*(E)\]
is surjective, and hence \[H^0(X, L^{\otimes n_0} \otimes E) \otimes \OO_X \longrightarrow L^{\otimes n_0} \otimes E\] is surjective on some non-empty
Zariski open subset of $X$. Thus the assertion follows.
Let $\OO_X^{\oplus r} \to L^{\otimes n_0} \otimes E$ be the homomorphism given by \[(a_1, \ldots, a_r) \longmapsto a_1 e_1 + \cdots + a_r e_r.\]
Let $Q$ be the cokernel of $\OO_X^{\oplus r} \to L^{\otimes n_0} \otimes E$. The sequence
\[
0 \longrightarrow \OO_X^{\oplus r} \longrightarrow L^{\otimes n_0} \otimes E \longrightarrow Q \longrightarrow 0
\]
is exact, and so is
\[
0 \longrightarrow (L^{\otimes n})^{\oplus r} \longrightarrow L^{\otimes n+ n_0} \otimes E \longrightarrow L^{\otimes n} \otimes Q \longrightarrow 0.
\]
Thus 
\[
0 \longrightarrow H^0(X, L^{\otimes n})^{\oplus r} \longrightarrow H^0(X, L^{\otimes n+ n_0} \otimes E)  \longrightarrow H^0(X, L^{\otimes n} \otimes Q)
\]
is also exact for all $n\geqslant 0$. Let $Q_n$ be the image of \[H^0(X, L^{\otimes n+ n_0} \otimes E)  \longrightarrow H^0(X, L^{\otimes n} \otimes Q).\]
We equip $H^0(X, L^{\otimes n+ n_0} \otimes E)$ with the norm family \[\xi_{(n+n_0)\varphi+\psi}=( \|\ndot\|_{(n+n_0)\varphi_\omega + \psi_{\omega}})_{\omega\in \Omega}.\] 
Let $\xi_n^L=( \|\ndot\|^L_{n, \omega} )_{\omega \in \Omega}$ be its restricted norm family on  $H^0(X, L^{\otimes n})^{\oplus r}$ induced by
the injection \[H^0(X, L^{\otimes n})^{\oplus r} \longrightarrow H^0(X, L^{\otimes n+ n_0} \otimes E).\]
Let  $\xi^{Q}_n=(\|\ndot\|^Q_{n, \omega} )_{\omega \in \Omega}$
be its quotient family on $Q_n$ induced by the surjection \[H^0(X, L^{\otimes n+ n_0} \otimes E)  \longrightarrow Q_n.\]
Then, by \cite[Proposition~4.3.13, (4.26)]{CMArakelovAdelic}, 
\begin{multline*}
\adeg(H^0(X, L^{\otimes n})^{\oplus r}, \xi_n^L) + \adeg(Q_n, \xi_n^Q) 
\leqslant \adeg(H^0(X, L^{\otimes n+ n_0} \otimes E),\xi_{(n+n_0)\varphi+\psi}). 
\end{multline*}
Since $\dim \Supp(Q) < \dim X$,  by Proposition~\ref{prop:positivity:torsion:sheaf},
\[
\liminf_{n\to\infty} \frac{ \adeg(Q_n, ( \|\ndot\|^Q_{n, \omega} )_{\omega \in \Omega})}{n^{d+1}} \geqslant 0.
\]
Therefore, by the super-additivity of  inferior limit, we obtain 
\begin{equation}\label{eqn:thm:lower:bound:chi:volume:01}
\liminf_{n\to\infty} \frac{\adeg(H^0(X, L^{\otimes n})^{\oplus r}, \xi_n^L)}{n^{d+1}/(d+1)!} \leqslant \liminf_{n\to\infty} \frac{\adeg(p_*(\overline L{}^{\otimes n}\otimes E))}{n^{d+1}/(d+1)!}.
\end{equation}
Let us consider the homomorphism of identity
\[ (H^0(X, L^{\otimes n})^{\oplus r}, (\|\ndot\|_{n\varphi_\omega}^{\oplus r})_{\omega\in\Omega}) \longrightarrow (H^0(X, L^{\otimes n})^{\oplus r}, (\|\ndot\|^L_{n,\omega})_{\omega\in\Omega}),\]
where
\[
\|(s_1, \ldots, s_r) \|_{n\varphi_{\omega}}^{\oplus r} = \begin{cases}
\max\limits_{i \in \{ 1, \ldots, r \} }  \|s_i\|_{n\varphi_{\omega}}  & \text{if $\omega \in \Omega \setminus \Omega_{\infty}$}, \\[2ex]
 (\| s_1 \|^2_{n\varphi_{\omega}}+\cdots+\|s_r\|^2_{n\varphi_{\omega}})^{1/2} & \text{if $\omega \in \Omega_{\infty}$}.
\end{cases}
\]
If $\omega \in \Omega \setminus \Omega_{\infty}$, then
\begin{align*}
\| (s_1, \ldots, s_r) \|^L_{n, \omega} & \leqslant \| s_1 e_1 + \cdots + s_r e_r \|_{(n+n_0)\varphi_\omega + \psi_{\omega}} \\
& \leqslant \max_{i \in \{ 1, \ldots, r \}}  \|s_i\|_{n\varphi_\omega}\|e_i \|_{n_0\varphi_\omega + \psi_{\omega}}  \\
& \leqslant \|(s_1, \ldots, s_r) \|_{n\varphi_{\omega}}^{\oplus r}\Big( \max_{i \in \{ 1, \ldots, r \}} \|e_i \|_{n_0\varphi_\omega + \psi_{\omega}}\Big) .
\end{align*}
Moreover, if $\omega \in \Omega_{\infty}$, then by Cauchy-Schwarz inequality
\begin{align*}
\| (s_1, \ldots, s_r) \|^L_{n, \omega}& \leqslant \| s_1 e_1 + \cdots + s_r e_r \|_{(n+n_0)\varphi_\omega + \psi_{\omega}} \\
& \leqslant \sum\nolimits_{i=1}^r \|s_i\|_{n\varphi_\omega}\|e_i \|_{n_0\varphi_\omega + \psi_{\omega}} \\
& \leqslant \Big(\sum\nolimits_{i=1}^r \|s_i\|_{n\varphi_\omega}\Big) \Big(\max_{i \in \{ 1, \ldots, r \}}  \|e_i \|_{n_0\varphi_\omega + \psi_{\omega}} \Big)\\
& \leqslant  \sqrt{r} \|(s_1, \ldots, s_r) \|_{n\varphi_{\omega}}^{\oplus r} \max_{i \in \{ 1, \ldots, r \}}  \|e_i \|_{n_0\varphi_\omega + \psi_{\omega}}.
\end{align*}
Therefore,
\[
h(f_n) \leqslant \int_{ \Omega } \max_{i\in\{1,\ldots,r\}}  \log \|e_i \|_{n_0\varphi_\omega + \psi_{\omega}} \nu(d\omega) + \frac{1}{2} \log(r) \vol(\Omega_{\infty}),
\]
and hence, by \cite[Proposition~4.3.18]{CMArakelovAdelic},
\begin{multline*}
r \operatorname{\adeg} (H^0(X, L^{\otimes n}), \xi_{n\varphi}) = \adeg (H^0(X, L^{\otimes n})^{\oplus r}, \xi_{n\varphi}^{\oplus r}) \\
\kern-17em \leqslant \adeg(H^0(X, L^{\otimes n})^{\oplus r}, \xi^L_{n}) + \\
r h^0(L^{\otimes n}) \left( \int_{ \Omega } \max_{i\in\{1,\ldots,r\}}  \log \|e_i \|_{n_0\varphi_\omega + \psi_{\omega}} \, \nu(d\omega) +  \frac{1}{2} \log(r) \vol(\Omega_{\infty}) \right),
\end{multline*}
where 
\[ h^0(L^{\otimes n}) = \dim_k H^0(X, L^{\otimes n}), \quad \xi_{n\varphi}=(\|\ndot\|_{n\varphi_\omega})_{\omega\in\Omega},\quad \xi_{n\varphi}^{\oplus r}=(\|\ndot\|_{n\varphi_\omega}^{\oplus r})_{\omega\in\Omega}.\]
Thus, 
\[
r \operatorname{\avol_{\chi}}(\overline L) \leqslant \liminf_{n\to\infty} \frac{\adeg(H^0(X, L^{\otimes n})^{\oplus r}, \xi^L_{n})}{n^{d+1}/(d+1)!}.
\]
Combining this inequality with \eqref{eqn:thm:lower:bound:chi:volume:01}, we obtain the assertion.
\end{proof} 

\begin{coro}\label{coro:chi:volume:generically:finite}
Let $\pi : Y \to X$ be a generically finite morphism of geometrically integral projective schemes over $K$,
$\overline{L} = (L, \varphi)$ be an adelic invertible $\OO_X$-module and $\overline{M} = (M, \psi)$ be
an adelic invertible $\OO_Y$-module. If there are a birational morphism $f : X \to Z$ of geometrically integral projective schemes over $K$ and an ample invertible $\OO_Z$-module
$A$ such that $L = f^{*}(A)$, then
\[
\deg(\pi) \operatorname{\avol_{\chi}}(\overline{L}) \leqslant \liminf_{n\to\infty} \frac{\adeg\big((p\mbox{{\tiny $\circ$}}\pi)_*(\pi^*(\overline L)^{\otimes n}\otimes\overline M)\big)}{n^{d+1}/(d+1)!}.
\]
In particular, $\deg(\pi) \operatorname{\avol_{\chi}}(\overline{L}) \leqslant \avol_\chi(\pi^{*}(\overline{L}))$.
\end{coro}

\begin{proof}
Since $\pi^*(\overline{L}{}^{\otimes n}) \otimes \overline{M}$ is an adelic invertible $\mathcal O_Y$-module, one can see that \[(L^{\otimes n} \otimes \pi_*(M), \pi_*(n \pi^{*}(\varphi) + \psi))\] is sectionally adelic for all $n \geqslant 0$ (see the last section of Chapter~\ref{sec:Family of metrics of locally free module}).
Note that \[\pi_*(n \pi^{*}(\varphi) + \psi) = n \varphi + \pi_*(\psi) \text{\quad and\quad} \rank (\pi_*M) = \deg(\pi).\]
Thus, by Theorem~\ref{thm:lower:bound:chi:volume}, 
\[\deg(\pi) \operatorname{\avol_\chi}(\overline{L}) \leqslant \liminf_{n\rightarrow+\infty}\frac{\widehat{\deg}(H^0(X,L^{\otimes n}\otimes\pi_*(M)),(\|\ndot\|_{n\varphi_\omega+\pi_*(\psi)_\omega})_{\omega\in\Omega})}{n^{d+1}/(d+1)!}.\]
Moreover, \[(H^0(X,  L^{\otimes n}\otimes\pi_*(M)), ( \|\ndot\|_{  n \varphi_{\omega}+\pi_*(\psi)_{\omega} })_{\omega \in \Omega})\] is isometric to
\[(H^0(Y,  \pi^*(L^{\otimes n})\otimes M), ( \|\ndot\|_{  n \pi_{\omega}^*(\varphi_{\omega})+\psi_{\omega}})_{\omega \in \Omega}).\]
Thus we obtain the required inequality.
\end{proof}

\begin{theo}\label{thm:upper:estimate:chi:vol:rank:one}
Let $\overline{L}  = (L, \varphi)$ be an adelic invertible $\OO_X$-module and $ \overline{E} = (E, U,\psi)$ be a birationally adelic torsion free $\OO_X$-module.
We assume that 
there are a birational morphism $f : X \to Z$ of geometrically integral projective varieties over $K$ and an ample invertible $\OO_Z$-module
$A$ with $L = f^{*}(A)$.
If either $(E, \psi)$ is an adelic invertible $\OO_X$-module or $X$ is normal, then
the sequence \[
 \frac{\adeg(p_*(\overline L{}^{\otimes n}\otimes \overline E))}{n^{d+1}/(d+1)!}, \quad n\in\mathbb N,\; n\geqslant 1
\] is convergent to $\rank(E)\operatorname{\avol_{\chi}}(\overline{L})$.
\end{theo}

\begin{proof}
In view of Theorem \ref{thm:lower:bound:chi:volume}, it suffices to establish the following inequality
\[\limsup_{n\to\infty} \frac{\adeg(p_*(\overline L{}^{\otimes n}\otimes \overline E))}{n^{d+1}/(d+1)!}\leqslant\rank(E)\operatorname{\avol_{\chi}}(\overline{L}).\]
First we assume that $(E, \psi)$ is an adelic invertible $\OO_X$-module.
Let us begin with the following claim:

\ifsmf\begin{enonce}{Claim}\fi
\ifams\begin{clai}\fi
\phantomsection\label{Claim:rank:one:norm:estimate:01}
One has the following inequality:
\[\limsup_{n\to\infty} \frac{\adeg(p_*(\overline{L}{}^{\otimes n}\otimes\overline E))}{n^{d+1}/(d+1)!}
\leqslant
\limsup_{n\to \infty} \frac{\adeg(p_*(\overline L{}^{\otimes(n+n_0)}))}{n^{d+1}/(d+1)!}
\]
for some positive integer $n_0$. 
\ifsmf\end{enonce}\fi
\ifams\end{clai}\fi

\begin{proof}
Since $L$ is nef and big, we can choose a positive integer $n_0$ and $s_0 \in H^0(X, L^{\otimes n_0} \otimes {E}^{\vee}) \setminus \{ 0 \}$.
Note that $s_0$ gives rise to an injective homomorphism
\[
H^0(X,  L^{\otimes n}\otimes E) \longrightarrow H^0(X, L^{\otimes (n+n_0)}).
\]
Let $\xi_{\sub,n}=(\|\ndot\|_{\sub,n,\omega})_{\omega\in\Omega}$ be the restricted norm family  of $H^0(X, L^{\otimes n}\otimes E)$ induced by
the above injective homomorphism and \[\xi_{(n+n_0)\varphi}=(\|\ndot\|_{(n+n_0)\varphi_\omega})_{\omega\in\Omega}.\]
In order to show Claim~\ref{Claim:rank:one:norm:estimate:01}, it is sufficient to see the following two inequalities:
\begin{multline*}
\limsup_{n\to \infty} \frac{\adeg(H^0(X, L^{\otimes n}\otimes E),  \xi_{\sub,n} )}{n^{d+1}/(d+1)!} \\
\leqslant
\limsup_{n\to \infty} \frac{\adeg(H^0(X, L^{\otimes (n+n_0)}),  \xi_{(n+n_0)\varphi} )}{n^{d+1}/(d+1)!}
\end{multline*}
and
\[
\limsup_{n\to\infty} \frac{\adeg(p_*(\overline{L}{}^{\otimes n}\otimes\overline E))}{n^{d+1}/(d+1)!}
\leqslant
\limsup_{n\to \infty} \frac{\adeg(H^0(X, L^{\otimes n}\otimes E),  \xi_{\sub,n} )}{n^{d+1}/(d+1)!}.
\]
The first inequality is a consequence of Lemma~\ref{lem:sub:quot:R:module}, Proposition~\ref{prop:positivity:torsion:sheaf}, \cite[Lemma 1.2.16]{CMIntersection} and
\cite[Proposition~4.3.13, (4.26)]{CMArakelovAdelic}.
Let us consider the homomorphism of identity  \[f:\Big(H^0(X, L^{\otimes n}\otimes E),  (\|\ndot\|_{ n \varphi_\omega+\psi_\omega})_{\omega\in\Omega}\Big) \longrightarrow \Big(H^0(X,  L^{\otimes n}\otimes E), \xi_{ \sub,n} \Big).\] For $s \in H^0(X, E \otimes L^{\otimes n}) \setminus \{ 0 \}$,
\begin{align*}
\frac{\| s \|_{\sub,n,\omega}}{\| s \|_{n\varphi_{\omega}+\psi_{\omega}}} & = \frac{\| s s_0 \|_{(n+n_0)\varphi_{\omega}}}{\| s \|_{n\varphi_{\omega}+\psi_{\omega}}} \\
& \leqslant \frac{\| s \|_{n\varphi_{\omega}+\psi_{\omega}} \|s_0 \|_{n_0\varphi_{\omega} - \psi_{\omega}}}{\| s \|_{n\varphi_{\omega}+\psi_{\omega}}} = \|s_0 \|_{n_0\varphi_{\omega} - \psi_{\omega}},
\end{align*}
so that 
$\| f \|_{\omega} \leqslant  \|s_0 \|_{n_0\varphi_{\omega} - \psi_{\omega}}$. Therefore, by \cite[Proposition~4.3.18]{CMArakelovAdelic},
\begin{multline*}
\adeg(H^0(X,  L^{\otimes n}\otimes E),  (\|\ndot\|_{ n \varphi_\omega+\psi_\omega})_{\omega\in\Omega}) \leqslant
\adeg(H^0(X, L^{\otimes n}\otimes E),  \xi_{\sub,n} )\\
+ \dim H^0(X, L^{\otimes n}\otimes E) \int_{\Omega}
\log \|s_0 \|_{n_0\varphi_{\omega} - \psi_{\omega}} \nu(d \omega).
\end{multline*}
Thus the second inequality follows.
\end{proof}

By Lemma  \ref{lemma:finitely:generated:norm:space}, Theorem-Definition \ref{theo:limit:thm:vol:chi} and the relation 
\[\lim_{n\rightarrow+\infty}\frac{(n+n_0)^{d+1}}{n^{d+1}}=1,\]
we obtain that
\[\lim_{n\rightarrow+\infty}\frac{\adeg(p_*(\overline L{}^{\otimes(n+n_0)}))}{n^{d+1}/(d+1)!}=\widehat{\vol}_{\chi}(\overline L).\]
Hence Claim~\ref{Claim:rank:one:norm:estimate:01} leads to
\[
\avol_{\chi}(\overline{L}; \overline{E}) 
\leqslant
\avol_{\chi}(\overline{L}),
\]
as required.

\bigskip
Next we assume that $X$ is normal. We prove the assertion by induction on $r := \rank(E)$.
Let $\mu : X' \to X$,  $(E', \psi')$ and $U$ be a birational morphism, an adelic invertible $\OO_{X'}$-module and
a non-empty Zariski open set of $X$, respectively, as in Definition~\ref{def:birat:adelic:torsion:free}.
First we suppose that $r = 1$.

\ifsmf\begin{enonce}{Claim}\fi
\ifams\begin{clai}\fi
\label{Claim:rank:one:norm:estimate:02}
One has the following inequality:
\[
\limsup_{n\to\infty} \frac{\adeg(p_*(\overline{L}{}^{\otimes n}\otimes\overline E))}{n^{d+1}/(d+1)!}\leqslant \limsup_{n\to\infty} \frac{\adeg((p\mbox{{\tiny $\circ$}}\mu)_*(\mu^*(\overline{L})^{\otimes n}\otimes\overline{E'}))}{n^{d+1}/(d+1)!}
\]
\ifsmf\end{enonce}\fi
\ifams\end{clai}\fi

\begin{proof}
This is a consequence of Lemma~\ref{lem:sub:quot:R:module}, Proposition~\ref{prop:positivity:torsion:sheaf}, Lemma~\ref{lemma:finitely:generated:norm:space} and
\cite[(4.26) in Proposition~4.3.13]{CMArakelovAdelic}.
\end{proof}

By Claim~\ref{Claim:rank:one:norm:estimate:02} together with the case where $(E, \psi)$ is an adelic invertible $\OO_X$-module,
one has \[\limsup_{n\to\infty} \frac{\adeg(p_*(\overline{L}{}^{\otimes n}\otimes\overline E))}{n^{d+1}/(d+1)!} \leqslant \avol_{\chi}(\mu^*(\overline{L})).\] On the other hand, since $X$ is normal,
one can see that $\avol_{\chi}(\mu^*(\overline{L})) = \avol_{\chi}(\overline{L})$, as desired.

In the case where $r \geqslant 2$, considering a birational morphism $X'' \to X'$ if necessarily, we may assume that there exists
an exact sequence $0 \to F' \to E' \to Q' \to 0$ on $X'$ such that
$F'$ and $Q'$ are locally free, $\rank(F') = 1$ and $\rank(Q') = r - 1$.
Let $\psi_{F'}$ be the restricted metric of $F'$ over $X'$ and $\psi_{Q'}$ be the quotient metric of $Q'$ over $X'$.
Let $Q$ be the image of $E \to \mu_*(E') \to \mu_*(Q')$ and $F$ be the kernel of $E \to Q$.
Shrinking $U$ if necessarily, $\psi_{Q'}$ and $\psi_{F'}$ descent to metric families $\psi_{Q}$ and $\psi_{F}$ of $\rest{Q}{U}$ and $\rest{F}{U}$.
Note that $\overline Q=(Q, \psi_Q)$ and $\overline F=(F, \psi_F)$ are birationally adelic torsion free $\OO_X$-modules by Proposition~\ref{Pro: quotient family dominated} and Corollary~\ref{coro:sub:metric:dominated}.
Therefore, by hypothesis of induction, 
\[
\begin{cases}
{\displaystyle \limsup_{n\to\infty} \frac{\adeg(p_*(\overline L{}^{\otimes n}\otimes \overline F))}{n^{d+1}/(d+1)!}  \leqslant \avol_{\chi}(L, \varphi),} \\[2ex]
{\displaystyle \limsup_{n\to\infty} \frac{\adeg(p_*(\overline L{}^{\otimes n}\otimes\overline{Q}))}{n^{d+1}/(d+1)!}  \leqslant (r-1) \operatorname{\avol_{\chi}}(L, \varphi).}
\end{cases}
\]
For any $n \in\mathbb N$, one has an exact sequence
\begin{equation}\label{Equ: short exact sequence EQF}
0 \to H^0(X, L^n\otimes ) \to H^0(X, L^n\otimes E) \to H^0(X, L^n\otimes Q) \to H^1(X,L^{\otimes n}\otimes F).
\end{equation}
Let $Q_n$ be the image of 
\[H^0(X,L^{\otimes n}\otimes E)\longrightarrow H^0(X,L^{\otimes n}\otimes Q).\]
Let $\xi_{n, \sub}=(\|\ndot\|_{n, \sub,\omega})_{\omega\in\Omega} $ be the restricted norm family of \[\xi_{n\varphi+\psi}=(\|\ndot\|_{n\varphi_\omega+\psi_\omega})_{\omega\in\Omega}\]
on $H^0(X, L^n\otimes F)$ and 
$\xi_{n, \quot}=(\|\ndot\|_{n, \quot,\omega})_{\omega\in\Omega}$ be the quotient norm family of $\xi_{ n\varphi+\psi}$ on $H^0(X,L^n\otimes Q)$. By \cite[(4.28)]{CMArakelovAdelic},  \begin{multline*}
\adeg(H^0(X,  L^n\otimes E), \xi_{n\varphi+\psi}) - \delta(H^0(X,  L^n\otimes E), \xi_{n\varphi+\psi}) \\
\leqslant
\Big( \adeg(H^0(X, L^n\otimes F), \xi_{n, \sub}) - \delta(H^0(X, L^n\otimes F), \xi_{n, \sub}) \Big) \\
+
\Big( \adeg(H^0(X, L^n\otimes Q), \xi_{n, \quot}) - \delta(H^0(X, L^n\otimes Q), \xi_{n, \quot}) \Big),
\end{multline*}
where for any adelic vector bundle $\overline V$ on $S$, $\delta(\overline V)$ denotes the sum $\widehat{\deg}(\overline V)+\widehat{\deg}(\overline V{}^\vee)$. Let $\xi_{n\varphi+\psi_Q,\operatorname{sub}}=(\|\ndot\|_{n\varphi_\omega+\psi_{Q,\omega},\sub})_{\omega\in\Omega}$ be the restriction of \[\xi_{n\varphi+\psi_{Q}}=(\|\ndot\|_{n\varphi_\omega+\psi_{Q,\omega}})_{\omega\in\Omega}\]
to $Q_n$.
It is easy to see that, for any $\omega\in\Omega$, 
\[\|\ndot\|_{n, \sub,\omega} = \|\ndot\|_{ n\varphi_\omega+\psi_{F,\omega} },
\quad
\|\ndot\|_{n, \quot,\omega} \geqslant \|\ndot\|_{ n\varphi_\omega+\psi_{Q,\omega},\sub}.\]
Thus, by \cite[Proposition~4.3.18]{CMArakelovAdelic}, 
\[\adeg(Q_n,\xi_{n,\operatorname{quot}}) \leqslant \adeg(Q_n,\xi_{n\varphi+\psi_Q,\sub}),\]
so that
\begin{multline*}
\adeg(p_*(\overline L{}^{\otimes n}\otimes\overline E)) - \delta(p_*(\overline L{}^{\otimes n}\otimes\overline E)) \\
\leqslant
\Big( \adeg(p_*(\overline L{}^{\otimes n}\otimes\overline F )) - \delta(H^0(X,  L^n\otimes F), \xi_{n, \sub}) \Big) \\
+
\Big(\widehat{\deg}( Q_n,\xi_{n\varphi+\psi_Q,\sub}) - \delta(Q_n, \xi_{n, \quot}) \Big).
\end{multline*}
Moreover, by \cite[Proposition~4.3.10]{CMArakelovAdelic},
\[
\begin{cases}
{\displaystyle \lim_{n\to\infty} \frac{ \delta(p_*(\overline L{}^{\otimes n}\otimes\overline E))}{n^{d+1}} = 0,} \\[2ex]
{\displaystyle \lim_{n\to\infty} \frac{ \delta(H^0(X,  L^n\otimes F), \xi_{n,\sub})}{n^{d+1}} = 0,} \\[2ex]
{\displaystyle \lim_{n\to\infty} \frac{ \delta(Q_n, \xi_{n, \quot})}{n^{d+1}} = 0,}
\end{cases}
\]
so that one obtains
\[
\limsup_{n\rightarrow+\infty}\frac{\adeg(p_*(\overline L{}^{\otimes n}\otimes\overline E))}{n^{d+1}/(d+1)!}  \leqslant
\avol_{\chi}(L, \varphi) + \limsup_{n\rightarrow+\infty}\frac{\adeg(Q_n,\xi_{n\varphi+\psi_Q,\sub})}{n^{d+1}/(d+1)!},
\]
and hence it is sufficient to show that 
\begin{equation}\label{Equ: comparaison limsup}\limsup_{n\rightarrow+\infty}\frac{\widehat{\deg}(Q_n,\xi_{n\varphi+\psi_Q,\sub})}{n^{d+1}/(d+1)!}\leqslant\limsup_{n\rightarrow+\infty}\frac{\adeg(p_*(\overline L{}^{\otimes n}\otimes\overline Q))}{n^{d+1}/(d+1)!}.\end{equation}

\ifsmf\begin{enonce}{Claim}\fi
\ifams\begin{clai}\fi
\label{Cl: Leray spectral sequence}
If we set $T_n=H^0(X,L^{\otimes n}\otimes Q)/Q_n$, then \[\lim_{n\rightarrow+\infty}\dim_K(T_n)/n^d=0.\]
\ifsmf\end{enonce}\fi
\ifams\end{clai}\fi
\begin{proof}
By the Leray spectral sequence
\[E_2^{p,q}=H^p(Z,A^{\otimes n}\otimes R^qf_*(F))\Longrightarrow H^{p+q}(X,L^{\otimes n}\otimes F),\]
if $n$ is sufficiently large, then one has an injective homomorphism
\[H^1(X,L^{\otimes n}\otimes F)\longrightarrow H^0(Z,A^{\otimes n}\otimes R^1f_*(F))\]
so that 
\[\lim_{n\rightarrow+\infty}\frac{\dim_K(H^1(X,L^{\otimes n}\otimes F))}{n^d}=0\]
because $\operatorname{Supp}(R^1f_*(F))$ has Krull dimension $<d$. Thus the assertion follows by \eqref{Equ: short exact sequence EQF}.  
\end{proof}

By Lemma~\ref{lemma:finitely:generated:norm:space}, $\bigoplus_{n=0}^{\infty} H^0(X, Q \otimes L^n)$ is finitely generated over
\[\bigoplus_{n=0}^{\infty} H^0(X, L^n),\]
so that $\bigoplus_{n=0}^{\infty} T_n$ is also finitely generated over it.
Let $\xi_{T_n}$ be the quotient norm family of $\xi_{n\varphi+\psi_Q}$ on $T_n$. Then by Claim \ref{Cl: Leray spectral sequence} together with Proposition \ref{prop:positivity:torsion:sheaf}, we obtain that 
\[\liminf_{n\rightarrow+\infty}\frac{\widehat{\deg}(T_n,\xi_{T_n})}{n^{d+1}}\geqslant 0,\]
that is, for any $\varepsilon>0$,
\[\frac{\widehat{\deg}(T_n,\xi_{T_n})}{n^{d+1}}\geqslant -\varepsilon\]
for sufficiently large $n$. Moreover, by \cite[Proposition~4.3.13, (4.26)]{CMArakelovAdelic}, 
\[\frac{\widehat{\deg}(Q_n,\xi_{n\varphi+\psi_Q,\sub})}{n^{d+1}}+\frac{\widehat{\deg}(T_n,\xi_{T_n})}{n^{d+1}}\leqslant\frac{\widehat{\deg}(p_*(\overline L{}^{\otimes n}\otimes\overline Q))}{n^{d+1}},\]
so that 
\[\frac{\widehat{\deg}(Q_n,\xi_{n\varphi+\psi_Q,\sub})}{n^{d+1}}-\varepsilon\leqslant\frac{\widehat{\deg}(p_*(\overline L{}^{\otimes n}\otimes Q))}{n^{d+1}}\]
for sufficiently large $n$. Thus,
\[\limsup_{n\rightarrow+\infty}\frac{\widehat{\deg}(Q_n,\xi_{n\varphi+\psi_Q,\sub})}{n^{d+1}}-\varepsilon\leqslant \limsup_{n\rightarrow+\infty}\frac{\widehat{\deg}(p_*(\overline L{}^{\otimes n}\otimes Q))}{n^{d+1}}.\]
Since $\varepsilon$ is arbitrary, we obtain the inequality \eqref{Equ: comparaison limsup}. 
\end{proof}

\begin{coro}
Let $(E,U, \psi)$ be a birational adelic torsion free $\OO_X$-module.
If $X$ is normal and $L$ is ample, then \[\lim_{n\rightarrow+\infty}\frac{\adeg(p_*(\overline{L}{}^{\otimes n}\otimes\overline E))}{n^{d+1}/(d+1)!}  =  \rank (E)\operatorname{\avol_{\chi}}(L, \varphi).\]
\end{coro}

\begin{proof}
This is a consequence of Theorem~\ref{thm:lower:bound:chi:volume} and Theorem~\ref{thm:upper:estimate:chi:vol:rank:one}.
\end{proof}

%!TEX root = ./Hilbert_Samuel_Adelic_Curves.tex

\chapter{Hilbert-Samuel property}
Let $f:X\rightarrow\Spec K$ be an integral projective scheme over $\Spec K$, $d$ be the dimension of $X$ and $L$ be an ample invertible $\mathcal O_X$-module. We assume that, either the field $K$ is perfect, or the scheme $X$ is geometrically integral. We denote by $\mathscr M(L)$ the set of metrics families $\varphi=(\varphi_\omega)_{\omega\in\Omega}$ such that all metrics $\varphi_\omega$ are semi-positive and that $(L,\varphi)$ forms an adelic line bundle on $X$.

\section{Definition and reduction}
\begin{defi}We say that $\varphi\in\mathscr M(L)$ satisfies the \emph{Hilbert-Samuel property}\index{Hilbert-Samuel property} if the equality 
\[\widehat{\operatorname{vol}}_{\chi}(L,\varphi)=((L,\varphi)^{d+1})\]
holds, namely the $\chi$-volume and the self-intersection number of $(L,\varphi)$ coincides.
\end{defi} 

\begin{rema}\label{remark:Hilbert:Samual:power:n}
Note that Theorem-Definition \ref{theo:limit:thm:vol:chi} shows that, for any positive integer $n$, one has 
\[\widehat{\operatorname{vol}}_{\chi}(L^{\otimes n},n\varphi)=n^{d+1}\widehat{\operatorname{vol}}_{\chi}(L,\varphi).\]
Therefore, if $\varphi$ satisfies the Hilbert-Samuel property, then for any positive integer $n$, the metric family $n\varphi$ also satisfies the Hilbert-Samuel property. Conversely, if there exists a positive integer $n$ such that $n\varphi$ satisfies the Hilbert-Samuel property, then so does the metric family $\varphi$.\end{rema}

In order to show the Hilbert-Samuel property for all metrics families in $\mathscr M(L)$, it suffices to check the property for one arbitrary metric family in $\mathscr M(L)$.

\begin{lemm}\label{Lem: distance of determinant}
Let $E$ be a finite-dimensional vector space over $K$. If $\xi=(\|\ndot\|_{\omega})_{\omega\in\Omega}$ and $\xi'=(\|\ndot\|_{\omega}')_{\omega\in\Omega}$ are two norm families on $E$, then one has
\begin{equation}\label{Equ: local hadamard}d_\omega(\det(\xi),\det(\xi'))\leqslant r d_\omega(\xi,\xi').\end{equation}
In particular, if $\xi$ is strongly dominated, so is $\det(\xi)$.
\end{lemm}
\begin{proof}
Let $r$ be the dimension of $E$ over $K$. If $\eta$ is a non-zero element of $\det(E_{\omega})$, then one has 
\[\begin{split}\ln\|\eta\|_{\omega,\det}-\ln\|\eta\|_{\omega,\det}'&=\sup_{\begin{subarray}{c}
(s_1,\ldots,s_r)\in E_{\omega}^r\\
\eta=s_1\wedge\cdots\wedge s_r
\end{subarray}}\ln\|s_1\wedge\cdots\wedge s_r\|_{\omega,\det}-\sum_{i=1}^r\ln\|s_i\|_{\omega}'\\
&\leqslant\sup_{\begin{subarray}{c}
(s_1,\ldots,s_r)\in E_{\omega}^r\\
\eta=s_1\wedge\cdots\wedge s_r\end{subarray}}\sum_{i=1}^r\ln\|s_i\|_{\omega}-\ln\|s_i\|_{\omega}'\leqslant rd_\omega(\xi,\xi').
\end{split}\] 
Interchanging $\xi$ and $\xi'$, the above inequality leads to 
\[\ln\|\eta\|_{\omega,\det}'-\ln\|\eta\|_{\omega,\det}\leqslant rd_\omega(\xi,\xi').\]
Therefore, the inequality \eqref{Equ: local hadamard} holds. 
\end{proof}

\begin{prop}\label{Pro: HS property particular case}
Assume that there exists a metric family $\psi\in\mathscr M(L)$ which satisfies the Hilbert-Samuel property. Then any metric family $\varphi\in\mathscr M(L)$ satisfies the Hilbert-Samuel property.
\end{prop}
\begin{proof}
For any $n\in\mathbb N$, let $E_n$ be the $K$-vector space $H^0(X,L^{\otimes n})$ and $r_n$ be the dimension of $E_n$ of $K$. For any $\omega\in\Omega$, let $E_{n,\omega}=E_n\otimes_KK_\omega$, \[d_{n,\omega}=\sup_{s\in E_{n,\omega}\setminus\{0\}}\Big|\ln\|s\|_{n\varphi_\omega}-\ln\|s\|_{n\psi_\omega}\Big|\] 
be the distance  of $\|\ndot\|_{n\varphi_\omega}$ and $\|\ndot\|_{n\psi_\omega}$, and 
\[\delta_{n,\omega}=\sup_{\eta\in\det(E_{n,\omega})\setminus\{0\}}
\ln\|\eta\|_{n\varphi_\omega,\det}-\ln\|\eta\|_{n\psi_\omega,\det}.\]
Note that the function $(\omega\in\Omega)\mapsto\delta_{n,\omega}$ is $\nu$-integrable, and one has 
\[\int_{\Omega}\delta_{n,\omega}\,\nu(\mathrm{d}\omega)=\widehat{\deg}(p_*(L^{\otimes n},n\psi))-\widehat{\deg}(p_*(L^{\otimes n},n\varphi)).\] 
By Lemma \ref{Lem: distance of determinant}, one has 
\[|\delta_{n,\omega}|\leqslant r_n d_{n,\omega}\leqslant nr_nd_\omega(\varphi,\psi).\]
Note that the function 
\[(\omega\in\Omega)\longrightarrow d_\omega(\varphi,\psi)\]
is dominated (see \cite[Proposition 6.1.12]{CMArakelovAdelic}). 
Moreover, by \cite[Theorem 1.7]{BGF2020},  one has 
\[\lim_{n\rightarrow+\infty}\frac{\delta_{n,\omega}}{n^{d+1}/(d+1)!}=\sum_{j=0}^d\int_{X_\omega^{\mathrm{an}}}f_\omega(x)\,\mu_{(L_\omega,\varphi_\omega)^j(L_\omega,\psi_\omega)^{d-j}}(\mathrm{d}x),\]
where $f_\omega$ is the continuous function on $X_\omega^{\mathrm{an}}$ such that
\[\mathrm{e}^{f_\omega(\omega)}|\ndot|_{\psi_\omega}(x)=|\ndot|_{\varphi_\omega}\]
for any $x\in X_\omega^{\mathrm{an}}$.  Hence Theorem-Definition \ref{theo:limit:thm:vol:chi} and Lebesgue's dominated convergence theorem lead to (see Remark \ref{Rem: limit theorem in average}) 
\begin{align*}
\widehat{\operatorname{vol}}_{\chi}(L,\psi)-\widehat{\operatorname{vol}}_{\chi}(L,\varphi)
&=\lim_{n\rightarrow+\infty}\frac{1}{n^{d+1}/(d+1)!}\int_{\Omega}\delta_{n,\omega}\,\nu(d\omega)\\
&=\sum_{j=0}^d\int_{\Omega}\int_{X_\omega^{\mathrm{an}}}f_\omega(x)\,\mu_{(L_\omega,\varphi_\omega)^j(L_\omega,\psi_\omega)^{d-j}}(\mathrm{d}x)\,\nu(\mathrm{d}\omega)\\
&=((L,\psi)^{d+1})-((L,\varphi)^{d+1}).
\end{align*}
The proposition is thus proved.
\end{proof}

\begin{defi}
Let $X$ be a geometrically integral projective scheme over $\operatorname{Spec}K$ and $L$ be an ample invertible $\mathcal O_X$-module. If there exists a metric family $\varphi\in\mathscr M(L)$ which satisfies the Hilbert-Samuel property, or equivalently, any metric family $\varphi\in\mathscr M(L)$ satisfies the Hilbert-Samuel property (see Proposition \ref{Pro: HS property particular case}), we say that the ample invertible $\mathcal O_X$-module $L$ \emph{satisfies the Hilbert-Samuel property}\index{Hilbert-Samuel property}.
\end{defi}

\begin{rema}\label{Rem: difference is constant}
The proof of Proposition \ref{Pro: HS property particular case} actually shows a more precise result: the function 
\[(\varphi\in\mathscr M(L))\longrightarrow \widehat{\operatorname{vol}}_{\chi}(L,\varphi)-((L,\varphi)^{d+1})\] 
is constant.
\end{rema}

\section{Case of a projective space}

In this section, we assume that $X=\mathbb P_K^d$ is the projective case and $L=\mathcal O_{\mathbb P_K^d}(1)$ is the universal line bundle. We  show that any metric family in $\mathscr M(L)$ satisfies the Hilbert-Samuel property. Without loss of generality (by Proposition \ref{Pro: HS property particular case}), we consider a particular case as follows. Let $E$ be a $(d+1)$-dimensional vector space over $K$ and $(e_i)_{i=0}^d$ be a basis of $E$. Let $\xi=(\|\ndot\|_{\omega})_{\omega\in\Omega}$ be the Hermitian norm family on $E$ such that $(e_i)_{i=0}^d$ forms an orthonormal basis of $E$ with respect to $\|\ndot\|_\omega$. We then identify $\mathbb P_K^d$ with $\mathbb P(E)$ and let $\varphi=(\varphi_\omega)_{\omega\in\Omega}$ be the quotient metric family on $L$ induced by $\xi$. Note that, for any integer $n\in\mathbb N$, the vector space $H^0(X,L^{\otimes n})$ is isomorphic to the symmetric power $S^n(E)$. We denote by $r_n$ the dimension of $S^n(E)$. One has
\[r_n=\binom{n+d}{d}.\]

\begin{defi}\label{Def: Gauss point}Let $\omega\in\Omega$ such that $|\ndot|_\omega$ is non-Archimedean. Let $x$ be the point in $\mathbb P(E_\omega)^{\mathrm{an}}$ which consists of the generic scheme point of $\mathbb P(E_\omega)$ equipped with the absolute value \[\textstyle|\ndot|_{x}:k\big(\frac{e_0}{e_r},\ldots,\frac{e_{r-1}}{e_r}\big)\longrightarrow\mathbb R_{\geqslant 0}\]
such that, for any  
\[ P=\sum_{\boldsymbol{a}=(a_0,\ldots,a_{r-1})\in\mathbb N^d}\lambda_{\boldsymbol{a}}\Big(\frac{e_0}{e_r}\Big)^{a_0}\cdots\Big(\frac{e_{r-1}}{e_r}\Big)^{a_{r-1}}\in \textstyle{k\big[\frac{e_0}{e_r},\ldots,\frac{e_{r-1}}{e_r}\big]},\]
one has 
\[|P|_{x}=\max_{\boldsymbol{a}\in\mathbb N^d}|\lambda_{\boldsymbol{a}}|_\omega.\]
Note that the point $x$ does not depend on the choice of the orthonormal basis $(e_j)_{j=0}^r$. In fact, the norm $\|\ndot\|$ induces a symmetric algebra norm on $K_\omega[E_\omega]$ (which is often called a \emph{Gauss norm}\index{Gauss norm}) and hence defines an absolute value on the fraction field of $K_\omega[E_\omega]$. The restriction of this absolute value to the field of rational functions on $\mathbb P(E_\omega)$ identifies with $|\ndot|_x$. Hence $x$ is called the \emph{Gauss point}\index{Gauss point} of $\mathbb P(E_\omega)^{\mathrm{an}}$.
\end{defi}

\begin{lemm}\label{Lem: global section norm and symmetric norm}
Let $\omega$ be an element of $\Omega$ such that $|\ndot|_\omega$ is non-Archimedean, and $n\in\mathbb N$. Let $\|\ndot\|_{n,\omega}$ be the $\varepsilon$-tensor power of $\|\ndot\|_\omega$ on the tensor power space $E_\omega^{\otimes n}$ and let $\|\ndot\|_{n,\omega}'$ be the quotient norm of $\|\ndot\|_{n,\omega}$ by the quotient homomorphism $E_\omega^{\otimes n}\rightarrow S^n(E_\omega)$. Then the norm $\|\ndot\|_{n,\omega}'$ coincides with the supremum norm $\|\ndot\|_{n\varphi_\omega}$ of the metric $n\varphi_\omega$ on $L_\omega^{\otimes n}$.
\end{lemm}
\begin{proof} For any $\omega\in\Omega$, we denote by $E_\omega$ the $K_\omega$-vector space $E\otimes_KK_\omega$.
By \cite[Propositions 1.3.16 and 1.2.36]{CMArakelovAdelic}, if we consider the Segre embedding $\mathbb P(E_\omega)\rightarrow\mathbb P(E_\omega^{\otimes n})$, then the metric $n\varphi_\omega$ identifies with the quotient metric induced by the norm $\|\ndot\|_{n,\omega}$. Moreover, if we denote by $\mathcal O_{E_\omega^{\otimes n}}(1)$ the universal invertible sheaf of $\mathbb P(E_\omega^{\otimes n})$ and by $\psi_\omega$ the quotient metric on this invertible sheaf induced by the norm $\|\ndot\|_{n,\omega}$. By \cite[Proposition 2.2.22]{CMArakelovAdelic}, the supremum norm $\|\ndot\|_{\psi_\omega}$ on
\[H^0(\mathbb P(E_\omega^{\otimes n}),\mathcal O_{E_\omega^{\otimes n}}(1))=E_\omega^{\otimes n}\] 
of the metric $\psi_\omega$ coincides with $\|\ndot\|_{n,\omega}$. Since $L^{\otimes n}$ is the restriction of $\mathcal O_{E_\omega^{\otimes n}}(1)$ to $X$ and the restriction map
\[H^0(\mathbb P(E_\omega^{\otimes n}),\mathcal O_{E_\omega^{\otimes n}}(1))\longrightarrow H^0(\mathbb P(E_\omega),L_\omega^{\otimes n})\]
identifies with the quotient homomorphism $E_\omega^{\otimes n}\rightarrow S^n(E_\omega)$. In particular, the supremum norm $\|\ndot\|_{\varphi_\omega^{\otimes n}}$ is bounded from above by the quotient norm $\|\ndot\|_{n,\omega}'$. 

Let $x$ be the Gauss point of the Berkovich analytic space $\mathbb P(E_\omega)^{\mathrm{an}}$ (see Definition \ref{Def: Gauss point}). If 
\[F=\sum_{\begin{subarray}{c}I=(a_0,\ldots,a_d)\in \mathbb N^{d+1}\\
a_0+\cdots+a_d=n\end{subarray}}\lambda_Ie_0^{a_0}\cdots e_d^{a_d}\]
is an element of $S^n(E)$, then the relation 
\[F(x)=\bigg(\sum_{\begin{subarray}{c}I=(a_0,\ldots,a_d)\in\mathbb N^d\\
a_0+\cdots+a_d=n
\end{subarray}}\lambda_I \Big(\frac{e_1}{e_0}\Big)^{a_1}\cdots\Big(\frac{e_d}{e_0}\Big)^{a_d} \bigg)e_0(x)^{\otimes n}\]
holds. In particular, one has
\[\|F\|_{n\varphi_\omega}\geqslant|F|_{n\varphi_\omega}(x)=\max_{\begin{subarray}{c}I=(a_0,\ldots,a_d)\in\mathbb N^d\\
a_0+\cdots+a_n=d
\end{subarray}}|\lambda_I|_{\omega}.\]
Since $F$ is the image of the element 
\[\widetilde F=\sum_{\begin{subarray}{c}
I=(a_0,\ldots,a_d)\in\mathbb N^d\\
a_0+\cdots+a_d=n
\end{subarray}}\lambda_Ie_0^{\otimes a_0}\otimes\cdots\otimes e_d^{\otimes a_d}\]
by the quotient map $E_\omega^{\otimes n}\rightarrow S^n(E_\omega)$, we obtain that 
\[\|F\|_{n\varphi_\omega}\geqslant\|\widetilde F\|_{n,\omega}\geqslant \|F\|_{n,\omega}'.\]
Therefore the equality $\|\ndot\|_{n\varphi_\omega}=\|F\|_{n,\omega}'$ holds.
\end{proof}

\begin{rema}\label{Rem: supremum norm ultrametric}
As a byproduct, the proof of the above lemma shows that, for any \[F=\sum_{\begin{subarray}{c}I=(a_0,\ldots,a_d)\in \mathbb N^{d+1}\\
a_0+\cdots+a_d=n\end{subarray}}\lambda_Ie_0^{a_0}\cdots e_d^{a_d}\in S^n(E_\omega),\]
one has
\[\|F\|_{n\varphi_\omega}=\max_{\begin{subarray}{c}I=(a_0,\ldots,a_d)\in\mathbb N^d\\
a_0+\cdots+a_n=d
\end{subarray}}|\lambda_I|_{\omega}.\] 
In other words, the family 
\[(e_0^{a_0}\cdots e_d^{a_d})_{\begin{subarray}{c}(a_0,\ldots,a_d)\in\mathbb N^{d+1}\\
a_0+\cdots+a_d=n
\end{subarray}}\]
forms an orthonormal basis of $(S^d(E_\omega),\|\ndot\|_{n\varphi_\omega})$.
\end{rema}

\begin{lemm}\label{Lem: computation of integral}
For any integer $d\in\mathbb N$ and any any $x>0$, let \begin{gather*}
P_{d+1,x}=\{(t_0,\ldots,t_d)\in\mathbb R_{\geqslant 0}^{d+1}\,|\,t_0+\cdots+t_d\leqslant x\},\\
\Delta_{d,x}=\{(t_0,\ldots,t_d)\in\mathbb R_{\geqslant 0}^{d+1}\,|\, t_0+\cdots+t_d = x \}.
\end{gather*}
We denote by $\operatorname{vol}_{d+1}$ the Lebesgue measure on $\mathbb R^d$. For any affine hyperplane of $\mathbb R^d$, we denote by $\nu_d$ the translate of the Haar measure on the underlying hyperplane which is normalized with respect to the canonical Euclidean norm on $\mathbb R^{d+1}$ (namely the parallelotope spanned by an orthonormal basis has volume $1$).   
\begin{enumerate}[label=\rm(\arabic*)]
\item\label{Item: volume P} The volume of $P_{d+1,x}$ with respect to $\operatorname{vol}_{d+1}$ is $x^{d+1}/(d+1)!$.
\item\label{Item: volume delta} The volume of $\Delta_{d,x}$ with respect to $\nu_d$ is $x^d\sqrt{d+1}/d!$.
\item\label{Item: integral of entropy} Let $\mu_d$ be the uniform probability distribution on $\Delta_{d,x}$. One has
\[\int_{\Delta_{d,1}}t_0\ln(t_0)+\cdots+t_d\ln(t_d)\,\mu_d(\mathrm{d}t)=-\frac{1}{d+1}\sum_{m=1}^d\sum_{\ell=1}^m\frac{1}{\ell}.\]
\end{enumerate} 
\end{lemm}
\begin{proof}

\ref{Item: volume P} We reason by induction on $d$. The case where $d=0$ is trivial. In the following we assume the induction hypothesis that the lemma holds for $\mathbb R^d$.  By Fubini's theorem, we have 
\[\operatorname{vol}_{d+1}(P_{d+1,x})=\int_0^x\operatorname{vol}_{d}(P_{d,x-t})\,\mathrm{d}t=\int_0^x\frac{(x-t)^d}{d!}=\frac{x^{d+1}}{(d+1)!}.\]

\ref{Item: volume delta} The distance from the origin to the affine hyperplane containing $\Delta_{d,x}$ is $x/\sqrt{d+1}$. Therefore, by the equality
\[\operatorname{vol}_{d+1}(P_{d+1,x})=\frac{1}{d+1}\ndot\frac{x}{\sqrt{d+1}}\nu_{d}(\Delta_{d,x}),\]
we obtain 
\[\nu_d(\Delta_{d,x})=\sqrt{d+1}\frac{x^d}{d!}.\]

\ref{Item: integral of entropy}
By Fubini's theorem, one has
{\allowdisplaybreaks
\begin{align*} \int_{P_{d+1,x}}t_0\ln(t_0)\,\operatorname{vol}_{d+1}(\mathrm{d}t_0,\ldots,\mathrm{d}t_d) &=\int_0^x t\ln(t)\operatorname{vol}_d(P_{d,x-t})\,\mathrm{d}t\\
&\kern-12em =\frac{1}{d!}\int_0^x t(x-t)^d\ln(t)\,\mathrm{d}t=\frac{1}{d!}\sum_{i=0}^d(-1)^i\binom{d}{i}x^{d-i}\int_0^xt^{i+1}\ln(t)\,\mathrm{d}t\\
&\kern-12em =\frac{1}{d!}\sum_{i=0}^d(-1)^i\binom{d}{i}x^{d-i}\frac{1}{i+2}\bigg(x^{i+2}\ln(x)-\frac{1}{i+2}x^{i+2}\bigg)\\
&\kern-12em =\frac{x^{d+2}\ln(x)}{d!}\sum_{i=0}^d(-1)^i\binom{d}{i}\frac{1}{i+2}-\frac{x^{d+2}}{d!}\sum_{i=0}^d(-1)^i\binom{d}{i}\frac{1}{(i+2)^2}.
\end{align*}}

By a change of variables, we obtain
\[\int_{P_{d+1},x}t_0\ln(t_0)\,\operatorname{vol}_{d+1}(\mathrm{d}t_0,\ldots,\mathrm{d}_{t_{d}})=\frac{1}{\sqrt{d+1}}\int_0^{x}\int_{\Delta_{d,u}}t_0\ln(t_0)\,\nu_{d}(\mathrm{d}t)\,\mathrm{d}u.\]
Taking the derivative with respect to $x$, we obtain  
\begin{multline*}\frac{(d+2)x^{d+1}\ln(x)+x^{d+1}}{d!}\sum_{i=0}^d(-1)^i\binom{d}{i}\frac{1}{i+2} \\ -\frac{(d+2)x^{d+1}}{d!}\sum_{i=0}^d(-1)^i\binom{d}{i}\frac{1}{(i+2)^2}\\
=\frac{1}{\sqrt{d+1}}\int_{\Delta_{d,x}}t_0\ln(t_0)\,\nu_d(\mathrm{d}t)=\frac{\nu_d(\Delta_{d,x})}{\sqrt{d+1}}\int_{\Delta_{d,x}}t_0\ln(t_0)\,\mu_d(\mathrm{d}t).
\end{multline*}
In particular, one has
\[\begin{split}\int_{\Delta_{d,1}}t_0\ln(t_0)\,\mu_d(\mathrm{d}t)&=\sum_{i=0}^d(-1)^i\binom{d}{i}\frac{1}{i+2}\Big(1-\frac{d+2}{i+2}\Big)\\
&=\sum_{i=0}^d(-1)^i\frac{d!}{i!(d-i)!}\cdot\frac{i-d}{(i+2)^2}\\
&=-\frac{1}{d+1}\sum_{i=0}^{d-1}(-1)^i\frac{(d+1)!}{(i+2)!(d-i-1)!}\cdot\frac{i+1}{i+2}.
\end{split}\]
Therefore
{\allowdisplaybreaks\begin{align*}(d+1)\int_{\Delta_{d,1}}t_0\ln(t_0)\,\mu_d(\mathrm{d}t)-d\int_{\Delta_{d-1,1}}t_0\ln(t_0)\,\mu_{d-1}(\mathrm{d}t)\\
&\kern-23.5em=-\sum_{i=0}^{d-1}(-1)^i\frac{(d+1)!}{(i+2)!(d-i-1)!}\cdot\frac{i+1}{i+2} \\
+\sum_{i=0}^{d-2}(-1)^i\frac{d!}{(i+2)!(d-i-2)!}\cdot\frac{i+1}{i+2}\\
&\kern-23.5em=-\sum_{i=0}^{d-1}(-1)^i\frac{d!}{(i+2)!(d-i-1)!}\cdot\frac{i+1}{i+2}(d+1-(d-i-1))\\
&\kern-23.5em=-\sum_{i=0}^{d-1}(-1)^i\frac{d!}{(i+2)!(d-i-1)!}((i+2)-1)\\
&\kern-23.5em=-\sum_{i=0}^{d-1}(-1)^i\bigg(\binom{d}{i+1}-\frac{1}{d+1}\binom{d+1}{i+2}\bigg)\\
&\kern-23.5em=\sum_{i=1}^d(-1)^i\binom{d}{i}+\frac{1}{d+1}\sum_{i=2}^{d+1}(-1)^i\binom{d+1}{i}\\
&\kern-23.5em=-1+\frac{1}{d+1}(-1+(d+1))=-\frac{1}{d+1}.
\end{align*}}%
Combining with
\[2\int_{\Delta_{1,1}}t_0\ln(t_0)\,\mu_1(\mathrm{d}t)=2\int_0^1t\ln(t)\,\mathrm{d}t=-\int_0^1t\,\mathrm{d}t=-\frac{1}{2},\]
by induction we obtain 
\[(d+1)\int_{\Delta_{d,1}}t_0\ln(t_0)\,\mu_d(\mathrm{d}t)=-\sum_{i=1}^d\frac{1}{i+1}=-\sum_{\ell=2}^{d+1}\frac{1}{\ell}.\]
By symmetry of $(t_0,\ldots,t_d)$, we get
\[(d+1)\sum_{i=0}^d\int_{\Delta_{d,1}}t_i\ln(t_i)\,\mu_d(\mathrm{d}t)=-(d+1)\sum_{\ell=2}^{d+1}\frac{1}{\ell}.\]
Since 
\[\begin{split}\sum_{m=1}^d\sum_{\ell=1}^m\frac{1}{\ell}&=\sum_{\ell=1}^d\sum_{m=\ell}^d\frac{1}{\ell}=\sum_{\ell=1}^d\frac{d+1-\ell}{\ell}=(d+1)\sum_{\ell=1}^d\frac{1}{\ell}-d\\
&=(d+1)\sum_{\ell=2}^{d+1}\frac{1}{\ell}+(d+1)-\frac{d+1}{d+1}-d=(d+1)\sum_{\ell=2}^{d+1}\frac{1}{\ell},
\end{split}\]
we obtain the desired result.
\end{proof}

\begin{prop}\label{prop:Hilbert:Samuel:Projective:Space}
The universal invertible sheaf $\mathcal O_{\mathbb P_K^d}(1)$ satisfies the Hilbert-Samuel property.
\end{prop}
\begin{proof}
By Proposition \ref{Pro: HS property particular case}, it suffices to prove that the particular quotient metric family $\varphi=(\varphi_\omega)_{\omega\in\Omega}$ defined in the beginning of the subsection satisfies the Hilbert-Samuel property. For any $n\in\mathbb N$, let
\[\eta_n=\bigwedge_{\begin{subarray}{c}(a_0,\ldots,a_d)\in\mathbb N^{d+1}\\
a_0+\cdots+a_d=n\end{subarray}}e_0^{a_0}\cdots e_d^{a_d}\in\det(S^{n}(E)).\]
By Lemma \ref{Lem: global section norm and symmetric norm} and \cite[Proposition 1.2.23]{CMArakelovAdelic}, for any $\omega\in\Omega$ such that $|\ndot|_{\omega}$ is non-Archimedean, one has 
\[\|\eta_n\|_{n\varphi_\omega,\det}=1.\]

Let $\omega$ be an element of $\Omega$ such that $|\ndot|_\omega$ is Archimedean. Similarly to Lemma \ref{Lem: global section norm and symmetric norm}, for each $n\in\mathbb N$, we let $\|\ndot\|_{n,\omega}$ be the orthogonal tensor power norm on $E_\omega^{\otimes n}$ and $\|\ndot\|_{n,\omega}'$ be its quotient norm on $S^n(E_\omega)$. Note that 
\[(e_0^{a_0}\cdots e_d^{a_d})_{\begin{subarray}{c}(a_0,\ldots,a_d)\in\mathbb N^{d+1}\\
a_0+\cdots+a_d=n
\end{subarray}}\]
forms an orthogonal basis of $(S^d(E_\omega),\|\ndot\|_{n,\omega}')$ and 
\[\|e_0^{a_0}\cdots e_d^{a_d}\|_{n,\omega}'=\Big(\frac{a_0!\cdots a_d!}{n!}\Big)^{\frac 12}.\]
By \cite[Proposition 1.2.25]{CMArakelovAdelic}, one has
\[\|\eta_n\|_{n,\omega,\det}'=\prod_{\begin{subarray}{c}
(a_0,\ldots,a_d)\in\mathbb N^{d+1}\\
a_0+\cdots+a_d=n
\end{subarray}}\Big(\frac{a_0!\cdots a_d!}{n!}\Big)^{\frac 12}.\]
In particular, using Stirling's formula one obtains
\[\begin{split}\lim_{n\rightarrow+\infty}\frac{\ln\|\eta_n\|_{n,\omega,\det}'}{n r_n}&=-\frac 12\int_{\Delta}(t_0\ln(t_0)+\cdots+t_d\ln(t_d))\,\mathrm{d}\mu\\&=\frac{1}{2(d+1)}\sum_{m=1}^d\sum_{\ell=1}^m\frac{1}{\ell},
\end{split}\]
where $\mu$ denotes the uniform probability measure on the simplex 
\[\Delta=\{(t_0,\ldots,t_d)\in\mathbb R_{\geqslant 0}^{d+1}\,|\,t_0+\cdots+t_d=1\},\] and the second equality comes from Lemma \ref{Lem: computation of integral}.

By \cite[Lemma 4.3.6]{MR1260106} and  \cite[Lemma 30]{MR1189489}  (see also \cite[VIII.2.5 lemma 2]{MR1208731}), one has 
\[\sup_{s\in S^n(E_\omega)\setminus\{0\}}\Big|\ln(r_n^{-\frac 12}\|s\|_{n,\omega}')-\ln\|s\|_{\varphi_\omega^{\otimes n}}\Big|=O(\ln(n)).\]
Moreover, \[\ln(r_n^{-\frac 12})=-\frac 12\ln r_n=O(\ln(n)).\]
Hence by Lemma \ref{Lem: distance of determinant}
we obtain 
\[\lim_{n\rightarrow+\infty}\frac{\ln\|\eta_n\|_{\varphi_\omega^{\otimes n},\det}}{nr_n}=\lim_{n\rightarrow+\infty}\frac{\ln\|\eta_n\|'_{n,\omega,\det}}{nr_n}.\]
The proposition is thus proved.
\end{proof}

\section{Trivial valuation case}

In this section, we show the Hilbert-Samuel property in the trivial valuation case. 
Let $v = (k,|\ndot|)$ be a trivially valued field.
Let us begin with the following Lemma:

\begin{lemm}\label{lemma:trivial:norm:Fubini:st:trivial:value}
Let $X$ be a integral projective scheme of dimension $d$ over $\Spec k$ and $L$ be a very ample invertible $\OO_X$-module.
Let $\|\ndot\|$ be the trivial norm on  $H^0(X, L)$, that is, $\| e \| = 1$ for $e \in H^0(X, L) \setminus \{ 0 \}$.
Let $\varphi$ be the Fubuni-Study metric of $L$ induced by the surjective homomorphism $H^0(X, L) \otimes \OO_X \to L$ and $\|\ndot\|$.
Then we have \[\avol_{\chi}(L, \varphi) = ((L, \varphi)^{d+1})_v = 0,\]
where in the construction of $\avol_{\chi}(L, \varphi)$ we consider the adelic curve consisting of one copy of the trivial absolute value on $k$ and the counting measure.
\end{lemm}

\begin{proof}
Let $X \hookrightarrow \PP^\ell_k$ be the embbedding given by $L$, where $\ell = \dim_k H^0(X, L) - 1$.
We can find a positive integer $n_0$ such that
$H^0(\PP^{\ell}_k, \OO_{\PP^{\ell}_k}(n)) \to H^0(X, L^{\otimes n})$  is surjective for all $n \geqslant n_0$.
In order to see $\avol_{\chi}(L, \varphi) = 0$, it is sufficient to show that the norm $\|\ndot\|_{n\varphi}$ is trivial for all $n \geqslant n_0$.
As $H^0(\PP^{\ell}_k, \OO_{\PP^{\ell}_k}(n)) = \Sym^n(H^0(X, L))$,  one has that $\Sym^n(H^0(X, L))  \to H^0(X, L^{\otimes n})$  is surjective
for all $n \geqslant n_0$.
Let $(T_0, \ldots, T_\ell)$ be a homogeneous coordinate of $\PP^\ell_k$. 
For $n \geqslant n_0$ and $s \in H^0(X, L^{\otimes n})$, if 
\[
s \equiv \sum_{\substack{(i_0, \ldots, i_\ell) \in \NN^{\ell+1} \\ i_0 + \cdots + i_\ell = n}} a_{i_0, \ldots, i_\ell} T_0^{i_0} \cdots T_\ell^{i_\ell} 
\]
modulo $\Ker(\Sym^n(H^0(X, L))  \to H^0(X, L^{\otimes n}))$,
then
\[
\| s \|_{n\varphi} = \sup_{x \in (X \cap U_0)^{\an}} \frac{  \bigg|\sum_{\substack{(i_0, \ldots, i_\ell) \in \NN^{\ell+1} \\ i_0 + \cdots + i_\ell = n}} a_{i_0, \ldots, i_\ell} z_1^{i_1} \cdots z_\ell^{i_\ell} \bigg|_x}{\Big(\max \{ 1,  |z_1|_x, \ldots, |z_\ell|_x \} \Big)^n},
\]
where $z_i = T_i/T_0$ and $U_0 = \{ (T_0, \ldots, T_\ell) \in \PP^{\ell}_k \,:\, T_0 \not= 0 \}$. Note that
\begin{multline*}
\bigg|\sum\nolimits_{\substack{(i_0, \ldots, i_\ell) \in \NN^{\ell+1} \\ i_0 + \cdots + i_\ell = n}} a_{i_0, \ldots, i_\ell} z_1^{i_1} \cdots z_\ell^{i_\ell} \bigg|_x  \\
\leqslant \max \{ |z_1|_x^{i_1} \cdots |z_\ell|_x^{i_\ell} \,:\, (i_0, \ldots, i_\ell) \in \NN^{\ell+1}, i_0 + \cdots + i_\ell = n \} \\
\leqslant  \Big(\max \{ 1,  |z_1|_x, \ldots, |z_\ell|_x \} \Big)^n,
\end{multline*}
and hence $\| s \|_{n\varphi}  \leqslant 1$. Let $k^{\mathrm{ac}}$ be an algebraic closure of $k$. 
We assume $s \not= 0$. We choose $\xi = (1, \xi_1, \ldots, \xi_n) \in X(k^{\mathrm{ac}})$ such that $s(\xi) \not=0$.
Then, as $\sum_{\substack{(i_0, \ldots, i_\ell) \in \NN^{\ell+1} \\ i_0 + \cdots + i_\ell = n}} a_{i_0, \ldots, i_\ell} \xi_1^{i_1} \cdots \xi_\ell^{i_\ell} \in k^{\mathrm{ac}} \setminus \{ 0 \}$ and $\xi_1,\ldots, \xi_\ell \in  k^{\mathrm{ac}}$, one has
\[
\Big| \sum\nolimits_{\substack{(i_0, \ldots, i_\ell) \in \NN^{\ell+1} \\ i_0 + \cdots + i_\ell = n}} a_{i_0, \ldots, i_\ell} \xi_1^{i_1} \cdots \xi_\ell^{i_\ell} \Big|_{v'} = 1
\quad\text{and}\quad
\max \{ 1,  |\xi_1|_{v'}, \ldots, |\xi_\ell|_{v'} \} = 1
\]
where $v'$ is the pair of $k^{\mathrm{ac}}$ and its trivial absolute value.
Therefore, $\| s \|_{n\varphi} = 1$.

\medskip
Next let us see that $ ((L, \varphi)^{d+1})_v = 0$.
Note that \[H^0(\PP^\ell_k, \OO_{\PP^\ell_k}(1)) = H^0(X, L)\ \text{and}\ 
\Sym^n(H^0(\PP^\ell_k, \OO_{\PP^\ell_k}(1))) = H^0(\PP^\ell_k, \OO_{\PP^\ell_k}(n))\]
for $n \geqslant 1$.
Let $\psi$ be the Fubuni-Study metric of $\OO_{\PP^\ell_k}(1)$ induced by the surjective homomorphism $H^0(\PP^\ell_k, \OO_{\PP^\ell_k}(1))) \otimes \OO_{\PP^\ell_k} \to \OO_{\PP^\ell_k}(1)$ and $\|\ndot\|$. Then $\rest{\psi}{X^{\an}} = \varphi$.
In the same way as before, $\| \ndot \|_{n\psi}$ on $H^0(\PP^\ell_k, \OO_{\PP^\ell_k}(n))$ is trivial for $n \geqslant 1$.
Therefore, the induced norm on $H^0(\check{\PP}^\ell_k \times \cdots \times \check{\PP}^\ell_k, \OO_{\check{\PP}^\ell_k}(\delta) \boxtimes \cdots \boxtimes \OO_{\check{\PP}^\ell_k}(\delta))$ is also  trivial, where $\delta = (L^d)$.
Thus the assertion follows.
\end{proof}

\begin{theo}\label{Thm: trivial valuation case}
Assume that, for any $\omega\in\Omega$, $|\ndot|_\omega$ is the trivial absolute value on $K$. Then any ample line bundle $L$ on $X$ satisfies the Hilbert-Samuel property.
\end{theo}
\begin{proof}

By Remark \ref{remark:Hilbert:Samual:power:n}, we may assume that $L$ is very ample.
Let $E$ be the vector space $H^0(X,L)$. For any $\omega\in\Omega$, we denote by $\|\ndot\|_{\omega}$ the trivial norm on $E=E_\omega$.
Let $\xi=(\|\ndot\|_{\omega})_{\omega\in\Omega}$ and $\varphi=(\varphi_\omega)_{\omega\in\Omega}$ be the quotient metric family on $L$ induced by $\xi$ and the canonical closed embedding $X\rightarrow\mathbb P(E)$. 
Then, Lemma~\ref{lemma:trivial:norm:Fubini:st:trivial:value} implies
\[\operatorname{vol}_\chi(L,\varphi)=((L, \varphi)^{d+1}) = 0.\]  
Therefore, by Proposition \ref{Pro: HS property particular case} we obtain that the invertible sheaf $L$ satisfies the Hilbert-Samuel property.
\end{proof}

\begin{rema}\label{Rem: comparison with combinatorics}
In \cite{Chen_Moriwaki2020}, an intersection product of metrized divisors has been introduced in the setting of curves over a trivially valued field $(k,|\ndot|)$. Let $X$ be a regular projective curve over $\operatorname{Spec} k$. Recall that the Berkovich space $X^{\operatorname{an}}$ is an infinite tree
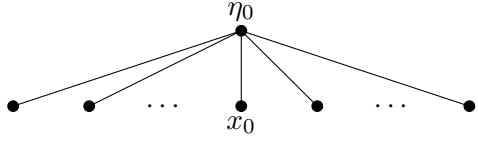
\vspace{3mm}
\begin{center}
\begin{tikzpicture}
\filldraw(0,1) circle (2pt) node[align=center, above]{$\eta_0$};
\filldraw(-3,0) circle (2pt) ;
\draw (-1,0) node{$\cdots$};
\filldraw(-2,0) circle (2pt) ;
\filldraw(-0,0) circle (2pt) node[align=center, below]{$x_0$} ;
\filldraw(1,0) circle (2pt) ;
\draw (2,0) node{$\cdots$};
\filldraw(3,0) circle (2pt) ;
\draw (0,1) -- (0,0);
\draw (0,1) -- (-3,0);
\draw (0,1) -- (1,0);
\draw (0,1) -- (-2,0);
\draw (0,1) -- (3,0);
\end{tikzpicture}
\end{center}
\vspace{3mm}
where the root point $\eta_0$ corresponds to the generic point of $X$ together with the trivial absolute value on $\kappa(\eta)$, and each leaf $x_0$ corresponds to the closed point $x$ together with the trivial absolute value on $\kappa(x)$. Moreover, each  branch $\mathopen{]}\eta_0,x_0\mathclose{[}$ is parametrized by  $\mathrm{]}0,+\infty\mathrm{[}$, where $t\in\mathrm{]}0,+\infty\mathrm{[}$ corresponds to the generic point $\eta$ together with the absolute value \[|\ndot|_{x,t}=\exp(-t\operatorname{ord}_x(\ndot)).\]
We denote by $t(\ndot):X^{\mathrm{an}}\rightarrow [0,+\infty]$ the parametrization map, where $t(\eta_0)=0$ and $t(x_0)=+\infty$.
Let $D$ be a Cartier divisor on $X$. Recall that a Green function $g$ of $D$ is of the form 
\[g=g_D+\varphi_g,\] 
where $g_D$ is the canonical Green function of $D$, which is defined as 
\[g_D(\xi)=\operatorname{ord}_x(D)t(\xi),\]
and $\varphi_g$ is a continuous real-valued function on $X^{\mathrm{an}}$ (which is hence bounded since $X^{\mathrm{an}}$ is compact). Then, the intersection number of two integrable metrized Cartier divisor $\overline D_0=(D_0,g_0)$ and $\overline D_1=(D_1,g_1)$ has been defined as
\begin{multline}\label{Equ: alternative intersection product}g_1(\eta_0)\deg(D_0)+g_0(\eta_0)\deg(D_1)
\\-\sum_{x\in X^{(1)}}[\kappa(x):k]\int_0^{+\infty}\varphi'_{g_0\circ\xi_x}(t)\varphi'_{g_1\circ\xi_x}(t)\,\mathrm{d}t,
\end{multline}
where $X^{(1)}$ is the set of closed points of $X$, $\xi_x:[0,+\infty]\rightarrow[\eta_0,x_0]$ is the map sending $t\in[0,+\infty]$ to the point in $[\eta_0,x_0]$ of parameter $t$, and the function $\varphi_{g_1\circ\xi_x}'(\ndot)$ should be considered as right-continuous version of the Radon-Nikodym density of the function $\varphi_{g_1\circ\xi_x}(\ndot)$ with respect to the Lebesgue measure.

Let $(L, \varphi_0)$ and $(L_1, \varphi_1)$ be integrable metrized invertible $\OO_X$-modules.
By \cite[Remark~7.3]{Chen_Moriwaki2020}, the above intersection number with respect to $(L, \varphi_0)$ and $(L_1, \varphi_1)$ is well-defined.
To destinguish this intersection number with the intersection number defined in \cite[Definition~3.10.1]{CMIntersection}
it is denoted by $((L_0, \varphi_0) \cdot (L_1, \varphi_2))'$.
Then one can see
\begin{equation}
 ((L_0, \varphi_0) \cdot (L_1, \varphi_1)) = ((L_0, \varphi_0) \cdot (L_1, \varphi_1))'.
\end{equation}
Indeed, by using the linearity of $(\ \cdot\ )$ and $(\ \cdot\ )'$, we may assume that
$L_0$ and $L_1$ are ample, and $\varphi_0$ and $\varphi_1$ are semipositive. Moreover, as
\[
{\small
\begin{cases}
{\displaystyle ((L_0, \varphi_0) \cdot (L_1, \varphi_1)) = \frac{ (((L_0, \varphi_0) + (L_1, \varphi_1))^2) -  ((L_0, \varphi_0)^2)  -  ((L_1, \varphi_1)^2) }{2}}, \\[2ex]
{\displaystyle ((L_0, \varphi_0) \cdot (L_1, \varphi_1))' = \frac{ (((L_0, \varphi_0) + (L_1, \varphi_1))^2)' -  ((L_0, \varphi_0)^2)'  -  ((L_1, \varphi_1)^2)' }{2}}, 
\end{cases}}
\]
we may further assume that $(L_0, \varphi_0) = (L_1, \varphi_1)$, say $(L, \varphi)$.
Then, by  \cite[Theorem~7.4]{Chen_Moriwaki2020},
\[
\lim_{n\to\infty} \frac{-\ln \| s_1 \wedge \cdots \wedge s_{r_n} \|_{n\varphi, \det}}{n^2/2} = ((L, \varphi) \cdot (L, \varphi))',
\]
where $\{ s_1, \ldots, s_{r_n} \}$ is a basis of $H^0(X, L^{\otimes n})$. On the other hand, 
\[
\lim_{n\to\infty} \frac{-\ln \| s_1 \wedge \cdots \wedge s_{r_n} \|_{n\varphi, \det}}{n^2/2} = ((L, \varphi) \cdot (L, \varphi))
\]
by Theorem~\ref{Thm: trivial valuation case} (the Hilbert-Samuel formula over a trivially valued field),  as required.
\end{rema}

\section{Casting to the trivial valuation case}

In this section, we assume that $K$ is perfect. Let $X$ be a projective $K$-scheme, $d$ be the dimension of $X$, $E$ be a finite-dimensional vector space over $K$, $f:X\rightarrow\mathbb  P(E)$ be a closed embedding, and $L$ be the restriction of the universal invertible sheaf $\mathcal O_E(1)$ to $X$. We assume that, for any positive integer $n$, the restriction map  
\[S^n(E)=H^0(\mathbb P(E),\mathcal O_E(n))\longrightarrow H^0(X,L^{\otimes n})\]
is surjective. We equip $E$ with a Hermitian norm family $\xi=(\|\ndot\|_{\omega})_{\omega\in\Omega}$ such that the couple $\overline E=(E,\xi)$ forms a strongly adelic vector bundle on the adelic curve $S$. Denote by $\varphi=(\varphi_\omega)_{\omega\in\Omega}$ the quotient metric family on $L$ induced by $\xi$ and the closed embedding $f$.

Let $\mathcal F=(\mathcal F^t(\overline E))_{t\in\mathbb R}$ be the Harder-Narasimhan $\mathbb R$-filtration of $\overline E$. Recall that 
\[\mathcal F^t(\overline E)=\sum_{\begin{subarray}{c}
0\neq F\subset E\\
\widehat{\mu}_{\min}(\overline F)\geqslant t
\end{subarray}}F\]
(cf. \cite[Corollary~4.3.4]{CMArakelovAdelic}).
Note that this $\mathbb R$-filtration actually defines an ultrametric norm $\|\ndot\|_0$ on $E$, where we consider the trivial absolute value $|\ndot|_0$ on the field $K$. More precisely, for any $s\in E$, one has
\[\|s\|_0=\exp(-\sup\{t\in\mathbb R\,:\,s\in\mathcal F^t(\overline E)\})\]
(cf. \cite[Remark~1.1.40]{CMArakelovAdelic}).
Denote by $\varphi_0$ the quotient metric on $L$ induced by $\|\ndot\|_0$. If we consider the  adelic curve $S_0$ consisting of a single copy of the trivial absolute value on $K$, then $(L,\varphi_0)$ becomes an adelic line bundle on $X$.

\begin{prop}\label{Pro: comparison of intersection numbers}
The following inequality holds:
\begin{equation}
((L,\varphi)^{d+1})\geqslant ((L,\varphi_0)^{d+1})-\nu(\Omega_\infty)\big((d+1)\delta\ln(r)+\ln(\delta!)\big),
\end{equation}
where $r$ denotes the dimension of $E$ over $K$ and $\delta$ is the degree of $X$ with respect to the line bundle $L$, that is,
$\delta=(L^d)$.
\end{prop}
\begin{proof}
For any $\omega\in\Omega$, let $\|\ndot\|_{\omega,*}$ be the dual norm on $E_\omega^\vee$ and let $\|\ndot\|_{\omega,*,\delta}$ be the $\delta$-th symmetric power of the norm $\|\ndot\|_{\omega,*}$, that is the quotient norm of the $\varepsilon$-tensor power (resp. orthogonal tensor power) of $\|\ndot\|_{\omega,*}$ by the canonical quotient map if $|\ndot|_\omega$ is non-Archimedean (resp. Archimedean). Let $\|\ndot\|_{\omega,*}'$ be the $\varepsilon$-tensor product (resp. orthogonal tensor product) of $d+1$ copies of the norm $\|\ndot\|_{\omega,*,\delta}$ if $|\ndot|_\omega$ is non-Archimedean (resp. Archimedean). By \cite[Proposition 1.2.36]{CMArakelovAdelic}, this norm also identifies with the quotient of the tensor power of $\|\ndot\|_{\omega,*}$ by the quotient map
\begin{equation}p_\omega:E_{\omega}^{\vee\otimes\delta(d+1)}\cong (E_\omega^{\vee\otimes\delta})^{\otimes(d+1)}\longrightarrow S^\delta(E_\omega^\vee)^{\otimes(d+1)}.\end{equation}
We denote by $\xi'$ the norm family $(\|\ndot\|'_{\omega,*})_{\omega\in\Omega}$. It turns out that \[(S^\delta(E^\vee)^{\otimes(d+1)},\xi')\] 
forms an
adelic vector bundle on $S$. Moreover, if we let $R\in S^\delta(E^\vee)^{\otimes(d+1)}$ be a resultant of $X$ with respect to $d+1$ copies of the closed embedding $f:X\rightarrow\mathbb P(E)$, then the following inequality holds:
\begin{equation}\label{Equ: intersection is degree}((L,\varphi)^{d+1})\geqslant -\operatorname{\widehat{\deg}}_{\xi'}(R)-\frac 12\nu(\Omega_\infty)(d+1)\ln\binom{r+\delta-1}{\delta},\end{equation}
where $r$ is the dimension of $E$ over $K$.
This is a consequence of \cite[Theorem 3.9.7]{CMIntersection} and \cite[Corollary 1.4.3, formula (1.4.10) and Lemma 4.3.6]{MR1260106}. Note that in the case where $\Omega_\infty=\varnothing$, the equality 
\begin{equation}\label{Equ: intersection is degree2}
((L,\varphi)^{d+1})=-\operatorname{\widehat{\deg}}_{\xi'}(R)\end{equation}
holds.

We now consider the trivial absolute value $|\ndot|_0$ on $K$ and we let $\xi'_0$ be the ultrametric norm on $S^{\delta}(E_\omega^\vee)^{\otimes (d+1)}$ defined as the quotient norm of the $\varepsilon$-tensor power of $\|\ndot\|_{0,*}$ by the quotient map
\[p:E^{\vee\otimes\delta(d+1)}\cong(E^{\vee\otimes\delta})^{\otimes(d+1)}\longrightarrow S^{\delta}(E^\vee)^{\otimes(d+1)}.\]
Similarly to \eqref{Equ: intersection is degree2}, the following equality holds:
\begin{equation}\label{Equ: height with respect to trivial valuation}((L,\varphi_0)^{d+1})=-\operatorname{\widehat{\deg}}_{\xi_0'}(R).\end{equation}
Note that the dual norm $\|\ndot\|_{0,*}$ corresponds to the Harder-Narasimhan $\mathbb R$-filtration of the dual  adelic vector bundle $\overline E{}^\vee=(E^\vee,\xi^\vee)$, where $\xi^\vee=(\|\ndot\|_{\omega,*})_{\omega\in\Omega}$ (see the proof of \cite[Proposition 4.3.41]{CMArakelovAdelic}). Therefore,  if we denote by $\Psi$ the one-dimensional vector sub-space of $S^{\delta}(E^\vee)^{\otimes (d+1)}$  spanned by the resultants of $X$ with respect to $d+1$ copies of $f:X\rightarrow\mathbb P(E)$,  then the dual statement of \cite[Theorem 5.6.1]{CMArakelovAdelic} (see Remark \ref{Rem: positive characteristic minimal slope}) leads to  
\[\widehat{\deg}(\Psi,\xi')\leqslant\widehat{\deg}(\Psi,\xi'_0)+\frac 12\nu(\Omega_\infty)\delta(d+1)\ln(r)+\nu(\Omega_\infty)\ln(\delta!),\]
or equivalently
\begin{equation}\label{Equ: lower bound of degree}-\widehat{\deg}_{\xi'}(R)\geqslant-\widehat{\deg}_{\xi_0'}(R)-\frac 12\nu(\Omega_\infty)\delta(d+1)\ln(r)-\nu(\Omega_\infty)\ln(\delta!).
\end{equation}
In the case where $\Omega_\infty$ is empty, we use Theorem \ref{Thm: maximal slope of symm power} to determine the Harder-Narasimhan $\mathbb R$-filtration of $S^\delta(\overline E^\vee)$ and apply the dual statement to the tensor product of $d+1$ copies of $S^\delta(\overline E^\vee)$. In the case where $\Omega_\infty$ is not empty, we uses the anti-symmetrization map (see Remark \ref{Rem: operator norm symmetri cpower}) to identify $S^\delta(E^\vee)$ with a vector subspace of $E^{\vee\otimes\delta}$ and apply the dual statement to $\delta(d+1)$ copies of $\overline E^\vee$. Note that the anti-symmetrization map $\operatorname{sym}'$ has height $\leqslant\nu(\Omega_\infty)\ln(\delta!)$ (see Propositions \ref{Pro: norm antisym p} and \ref{Pro: norm antisym 0}). By \eqref{Equ: intersection is degree},  \eqref{Equ: height with respect to trivial valuation} and \eqref{Equ: lower bound of degree}, we obtain 
\[
\begin{split}((L,\varphi)^{d+1})&\geqslant ((L,\varphi_0)^{d+1})-\frac 12\nu(\Omega_\infty)(d+1)\ln\binom{r+\delta-1}{\delta}\\
&\qquad-\frac 12\nu(\Omega_\infty)\delta(d+1)\ln(r)-\nu(\Omega_\infty)\ln(\delta!)\\
&\geqslant((L,\varphi_0)^{d+1})-\nu(\Omega_\infty)\delta(d+1)\ln(r)-\nu(\Omega_\infty)\ln(\delta!),
\end{split}\]
by using the inequality
\[\binom{r+\delta-1}{\delta}\leqslant r^{\delta}.\]
The proposition is thus proved.

\end{proof}

\section{Arithmetic Hilbert-Samuel theorem}
The purpose of this section is to prove the following theorem.

\begin{theo}\label{thm:lower:bound}
 Let $X$ be an integral projective $K$-scheme, $d$ be the dimension of $X$ and $L$ be an ample invertible $\mathcal O_X$-module. We assume that, either $K$ is perfect, or $X$ is geometrically integral. Then for any metric family $\varphi\in\mathscr M(L)$, the following equality holds
\begin{equation}\label{Equ: lower bound of chi vol}\widehat{\vol}_{\chi}(L,\varphi)= ((L,\varphi)^{d+1}).\end{equation}
\end{theo}
\begin{proof}
{\noindent\bf Step 1:} {\it We first prove the inequality $\widehat{\vol}_{\chi}(L,\varphi)\leqslant ((L,\varphi)^{d+1})$.}

Let $K'$ be the perfect closure of $K$. Note that each absolute value $|\ndot|_{\omega}$, $\omega\in\Omega$, extends in a unique way to $K'$, so that the underlying measure space of $S\otimes_KK'$ identifies with $(\Omega,\mathcal A,\nu)$.
Let $X' = X \times_{\Spec K} \Spec K'$, $L' = L \otimes_K K'$, and
$\varphi'$ be the extension of $\varphi$ to $L'$. Let $(s_1, \ldots, s_N)$ be a basis of $H^0(X, L^{\otimes n})$. 
Note that, for any $\omega\in\Omega$, then $\|\ndot\|_{n\varphi'_{\omega}}$ is an extension of $\|\ndot\|_{n\varphi_{\omega}}$ (cf. \cite[Proposition~2.1.19]{CMArakelovAdelic}). Therefore, by \cite[Proposition~1.1.66]{CMArakelovAdelic} or Appendix~\ref{subsection:determinant:norm},
one has
\[ \| s_1 \wedge \cdots \wedge s_N \|_{n\varphi_{\omega}, \det} \geqslant \| s_1 \wedge \cdots \wedge s_N \|_{n\varphi'_{\omega}, \det}, \] so that $\widehat{\vol}_{\chi}(L,\varphi) \leqslant \widehat{\vol}_{\chi}(L',\varphi')$.
Moreover, by \cite[Theorem~4.3.6]{CMIntersection}, \[((L,\varphi)^{d+1}) = ((L',\varphi')^{d+1}).\]
Thus, if the assertion of Step~1 holds for $K'$, then one has
\[
\widehat{\vol}_{\chi}(L,\varphi) \leqslant  \widehat{\vol}_{\chi}(L',\varphi') \leqslant ((L',\varphi')^{d+1}) = ((L,\varphi)^{d+1}).
\]
Therefore we may assume that $K$ is perfect.

By taking a tensor power of $L$ we may assume that $L$ is very ample and the canonical $K$-linear map
\begin{equation}S^n(H^0(X,L))\longrightarrow H^0(X,L^{\otimes n})\label{Equ: surjective map}\end{equation}
is surjective for any integer $n\geqslant 1$. Moreover, by Remark \ref{Rem: difference is constant}, the difference 
\[\widehat{\vol}_{\chi}(L,\varphi)-((L,\varphi)^{d+1})\]
does not depend on the choice of the metric family $\varphi$. Therefore, we may assume that $\varphi$ identifies with the quotient metric family induced by the norm family $\xi_1=(\|\ndot\|_{\varphi_\omega})_{\omega\in\Omega}$. By \cite[Proposition 2.2.22 (2)]{CMArakelovAdelic}, for any positive integer $n$, the metric $n\varphi$ identifies with the quotient metric family induced by the norm family $\xi_n=(\|\ndot\|_{n\varphi_\omega})_{\omega\in\Omega}$. Moreover, by changing metrics we may also assume that the minimal slope of $(H^0(X,L),\xi_1)$ is non-negative. Since the $K$-linear map \eqref{Equ: surjective map} is surjective, by \cite[Proposition 6.3.25]{CMArakelovAdelic}, we obtain that the minimal slope of $(H^0(X,L^{\otimes n}),\xi_n)$ is non-negative for any positive integer $n$. By \cite[Theorem 4.1.26]{CMArakelovAdelic}, there exists a Hermitian norm family $\xi_n'=(\|\ndot\|_{n,\omega}')$ of $H^0(X,L^{\otimes n})$ such that $\|\ndot\|_{n,\omega}=\|\ndot\|_{n\varphi_\omega}$ when $|\ndot|_\omega$ is non-Archimedean and 
\begin{equation}\label{Equ: encadrement norm hermitien}\|\ndot\|_{n,\omega}'\leqslant\|\ndot\|_{n\varphi_\omega}\leqslant (2r_n)^{1/2}\|\ndot\|_{n,\omega}'\end{equation}
when $|\ndot|_\omega$ is Archimedean, where $r_n$ denotes the dimension of $H^0(X,L^{\otimes n})$. Note that 
\[\Big|\widehat{\deg}(H^0(X,L^{\otimes n}),\xi_n)-\widehat{\deg}(H^0(X,L^{\otimes n}),\xi_n')\Big|\leqslant \frac 12\nu(\Omega_\infty)r_n\ln(2r_n),\]
so that 
\begin{equation}\label{Equ: vol chi computed by xin'}\operatorname{\widehat{vol}}_{\chi}(L,\varphi)=\lim_{n\rightarrow+\infty}\frac{\widehat{\deg}(H^0(X,L^{\otimes n}),\xi_n')}{n^{d+1}/(d+1)!}.\end{equation}

For any positive integer $n$, let $\|\ndot\|_n$ be the ultrametric norm on $H^0(X,L^{\otimes n})$ corresponding to the Harder-Narasimhan $\mathbb R$-filtration of 
$(H^0(X,L^{\otimes n}),\xi_n')$,
where we consider the trivial absolute value $|\ndot|_0$ on $K$. Let $\widetilde\varphi_n$ be the continuous metric on $L$ (where we still consider the trivial absolute value on $K$) such that $n\widetilde\varphi_n$ identifies with the quotient metric on $L^{\otimes n}$ induced by $\|\ndot\|_n$. By \cite[Proposition 2.2.22 (2)]{CMArakelovAdelic}, one has $\|\ndot\|_{n\widetilde\varphi_n}=\|\ndot\|_n$ on $H^0(X,L^{\otimes n})$ and hence
\begin{equation}\label{Equ: equality of degree}\widehat{\deg}(H^0(X,L^{\otimes n}),\|\ndot\|_{n\widetilde\varphi_n})=\widehat{\deg}(H^0(X,L^{\otimes n}),\|\ndot\|_{n})=\widehat{\deg}(H^0(X,L^{\otimes n}),\xi_n').\end{equation}

By Proposition \ref{Pro: comparison of intersection numbers} and the second inequality of \eqref{Equ: encadrement norm hermitien} we obtain that 
\begin{equation}\label{Equ: lower bound of intersection number}\begin{split}&\quad\;((nL,n\varphi)^{d+1})+\frac{1}{2}\nu(\Omega_\infty)(d+1)n^d(L^d)\ln(2r_n)\\
&\geqslant((nL,n\widetilde \varphi_n)^{d+1})-\nu(\Omega_\infty)n^d(L^d)\Big((d+1)\ln(r_n)+\ln(n^d(L^d))\Big),
\end{split}\end{equation}
where  we consider $X$ as an arithmetic variety over the adelic curve $S$ (resp. as an arithmetic variety over the adelic curve consisting of a single copy of the trivial absolute value on $K$) in the computation of the arithmetic intersection number on the left-hand side (resp. right-hand). Moreover, by Theorem \ref{Thm: trivial valuation case}, the following equality holds:
\begin{equation}\label{Equ: trivial valuation equalithy}\widehat{\vol}_{\chi}(L,\widetilde\varphi_n)=((L,\widetilde\varphi_n)^{d+1}).\end{equation}
By  \cite[Corollary 5.2]{MR3299845} (see also the proof of Theorem 7.3 of \emph{loc. cit.}), there exists a positive constant $C$ such that, for any positive integer $n$, one has
\[\widehat{\deg}(H^0(X,L^{\otimes n}),\|\ndot\|_{n\widetilde{\varphi}_n})\leqslant \frac{\widehat{\vol}_\chi(nL,n\widetilde\varphi_n)}{(d+1)!}+ Cn^d.\]
The constant $C$ can be taken in the form an invariant of the graded linear series $\bigoplus_{m\in\mathbb N}H^0(X,L^{\otimes m})$ multiplied by \[\sup_{m\in\mathbb N,\,m\geqslant 1}\frac{\widehat{\mu}_{\max}(H^0(X,L^{\otimes m}),\xi_m')}{m}.\]
By \eqref{Equ: equality of degree}, \eqref{Equ: lower bound of intersection number} and \eqref{Equ: trivial valuation equalithy}, we deduce that 
\[\begin{split}\widehat{\deg}(H^0(&X,L^{\otimes n}),\xi_n')\leqslant \frac{n^{d+1}}{(d+1)!}((L,\varphi)^{d+1})+Cn^d \\
&+\frac 12\nu(\Omega_\infty)(d+1)n^d(L^d)\ln(2r_n^3)+\nu(\Omega_\infty)n^d(L^d)\ln(n^d(L^d)).
\end{split}\]
 Dividing the two sides of the inequality by $n^{d+1}/(d+1)!$ and then taking the limit when $n\rightarrow+\infty$, by \eqref{Equ: vol chi computed by xin'} we obtain  
\[\widehat{\vol}_{\chi}(L,\varphi)\leqslant ((L,\varphi)^{d+1}).\]

{\noindent\bf Step 2:} {\it the converse inequality $\widehat{\vol}_\chi(L,\varphi)\geqslant ((L,\varphi)^{d+1})$.}

By replacing $L$ by a tensor power, we may assume that $L$ is very ample. Moreover, by 
the normalization of Noether (cf. \cite[Proposition 1.7.4]{CMIntersection}), we may also assume that there is a finite $K$-morphism $\pi:X\rightarrow\mathbb P_K^d$ such that $L\cong \pi^*(\mathcal O_{\mathbb P_K^d}(1))$. By Remark \ref{Rem: difference is constant}, we may further assume that there exists an element $\psi=(\psi_\omega)_{\omega\in\Omega}$ of $\mathscr M(\mathcal O_{\mathbb P_K^d}(1))$ such that $\varphi$ equals the pull-back of $\psi$ by $\pi$. 
Then, by Corollary~\ref{coro:chi:volume:generically:finite}, Proposition~\ref{prop:Hilbert:Samuel:Projective:Space} and \cite[Theorem 4.4.9]{CMIntersection},
one has
\begin{multline*}
\avol_\chi(L, \varphi) \geqslant \deg(\pi) \avol_\chi(\mathcal O_{\mathbb P_K^d}(1),\psi) 
=  \deg(\pi) ((\mathcal O_{\mathbb P_K^d}(1),\psi)^{d+1}) = ((L,\varphi)^{d+1}),
\end{multline*}
as required.
\end{proof}

\begin{coro}\label{Cor:Hilbert:Samuel:with:auxiliary:locally:free:module}
Let $X$ be a geometrically integral projective scheme over $\Spec K$, $d$ be the dimension of $X$,  $\overline L=(L,\varphi)$ be an adelic line bundle on $X$ and $\overline E=(E,U,\psi)$ be  a birational adelic torsion free $\OO_X$-module. 
Assume that $L$ is ample and the metrics in $\varphi$ are semi-positive. Moreover we suppose that either $(E, \psi)$ is an adelic invertible $\OO_X$-module or $X$ is normal.
Then one has
\[\lim_{n\rightarrow+\infty}\frac{\operatorname{\widehat{\deg}}(H^0(X,L^{\otimes n}\otimes E),(\|\ndot\|_{n\varphi_\omega+\psi_\omega})_{\omega\in\Omega})}{n^{d+1}/(d+1)!}=\rank(E)(\overline L{}^{d+1}).\]
\end{coro}

\begin{proof}
This is a consequence of Theorem~\ref{thm:lower:bound} together with Theorem~\ref{thm:upper:estimate:chi:vol:rank:one}.
\end{proof}

%!TEX root = ./Hilbert_Samuel_Adelic_Curves.tex

%% Positivity %%
\renewcommand{\thesubsection}{\arabic{chapter}.\arabic{section}.\arabic{subsection}}
\chapter{Positivity of adelic line bundles}

Throughout the chapter, we assume that the underlying field $K$ of the adelic curve $S=(K,(\Omega,\mathcal A,\nu),\phi)$ is perfect.

\section{Relative ampleness and nefness}

In this section, we consider relative ampleness and nefness.

\subsection{Convergence of minimal slopes}

\begin{lemm}\label{Lem: surjectivity tensor product} 
Let $k$ be a field, $X$ and $Y$ be projective $k$-schemes and $g:Y\rightarrow X$ be a projective $k$-morphism such that $g_*(\mathcal O_Y)=\mathcal O_X$. 
Let $L$ be an  ample line bundle on $Y$ and $M$ be an  ample line bundle on $X$. Then there exists $N\in\mathbb N_{\geqslant 1}$ such that, for any $(n,m)\in\mathbb N^2$ satisfying $\min\{n,m\}\geqslant N$, the $k$-linear map
\begin{multline*}
H^0(Y,L^{\otimes n})\otimes_KH^0(X,M^{\otimes m})=H^0(Y,L^{\otimes n})\otimes_KH^0(Y,g^*(M^{\otimes m})) \\
\longrightarrow H^0(Y,L^{\otimes n}\otimes g^*(M)^{\otimes m})
\end{multline*}
defined by multiplication of sections is surjective.
\end{lemm}
\begin{proof}
Consider the graphe
\[\Gamma_g:Y\longrightarrow Y\times_kX\]
of the morphism $g:Y\rightarrow X$. It is a closed immersion since $g$ is separated. Denote by $I$ the ideal sheaf of the image of $\Gamma_g$.  Let $p:Y\times_kX\rightarrow Y$ and $q:Y\times_kX\rightarrow X$ be the two projections, and $A=p^*(L)\otimes q^*(M)$. Since $M$ and $L$ are both ample, the line bundle $A$ on $Y\times_kX$ is ample. Moreover, one has $\Gamma_g^*(A)=L\otimes g^*(M)$. The short exact sequence 
\[\xymatrix{0\ar[r]&I\ar[r]&\mathcal O_{Y\times_KX}\ar[r]&\mathcal O_{Y\times_kX}/I\ar[r]&0}\]
induces, by tensor product with the invertible sheaf $p^*(L^{\otimes n})\otimes q^*(M^{\otimes m})$ and  then by taking cohomology groups on $Y\times_kX$,  an exact sequence of $K$-vector spaces
\begin{multline*}
H^0(Y,L^{\otimes n})\otimes_kH^0(Y,g^*(M)^{\otimes m})\longrightarrow H^0(Y,L^{\otimes n}\otimes g^*(M)^{\otimes m})\\
\longrightarrow H^1(Y\times_kX,I\otimes p^*(L^{\otimes n})\otimes q^*(M^{\otimes m})).
\end{multline*}
By \cite[Example 1.4.4]{MR2095471}, the line bundles $p^*(L)$ and $q^*(M)$ are nef.
By Fujita's vanishing theorem (cf. \cite[Theorem~5.1]{MR722501}), there exists $N\in\mathbb N_{\geqslant 1}$ such that, for any $(n,m)\in\mathbb N^2$ such that $\min\{n,m\}> N$, one has\[
\begin{split}&\quad\;H^1(Y\times_kX,I\otimes p^*(L^{\otimes n})\otimes q^*(M^{\otimes m}))\\
&=H^1(Y\times_kX,I\otimes A^{\otimes N}\otimes p^*(L^{\otimes (n-N)})\otimes q^*(M^{\otimes (m-N)}))=\boldsymbol{0}.\end{split}\] 
Therefore the assertion follows.
\end{proof}

\begin{lemm}\label{lemma:restriction:semi-positive}
Let $(k, |\ndot|)$ be a filed equipped with a complete absolute value. Let $X$ be a projective scheme over $k$, $L$ be a semi-ample line bundle on $X$ and $\varphi$ be a semi-positive metric of $L$.
Then, for any projective $K$-morphism $g:Y\rightarrow X$, $g^*(\varphi)$ is also semi-positive.
\end{lemm}
\begin{proof}
Replacing $L$ by a tensor power, we may assume that $L$ is generated by global sections, and that there exists a sequence of quotient metric families $(\varphi_n)_{n\in\mathbb N}$ such that 
\[\lim_{n\rightarrow+\infty}\frac 1nd(n\varphi,\varphi_n)=0.\]
Note that for each $n\in\mathbb N$, the pull-back $g^*(\varphi_n)$ is still a quotient metric, and one has 
\[d(ng^*(\varphi),g^*(\varphi_n))\leqslant d(n\varphi,\varphi_n).\]
Therefore we obtain that $g^*(\varphi)$ is semi-positive. 
\end{proof}

In the remaining of the section, we let $f:X\rightarrow\Spec K$ be a \emph{non-empty and reduced} projective scheme over $\Spec K$. Since the base field $K$ is supposed to be perfect, the $K$-scheme $X$ is geometrically reduced.

\begin{prop}\label{Pro: convergence mu min }
Let $\overline L=(L,\varphi)$ be an adelic line bundle on $X$ such that $L$ is ample. Then the sequence
\begin{equation}\label{Equ: asymptotic minimal slope}\frac{\widehat{\mu}_{\min}(f_*(\overline L^{\otimes n}))}{n},\quad n\in\mathbb N,\;n\geqslant 1\end{equation}
converges in $\mathbb R$.
\end{prop}
\begin{proof}
For any $n\in \mathbb N_{\geqslant 1}$ 
let $\overline E_n=(E_n,\xi_n)$ be the adelic vector bundle $f_*(\overline L^{\otimes n})$. Since $L$ is ample,  by Lemma~\ref{Lem: surjectivity tensor product} there exists $N\in\mathbb N_{\geqslant 1}$ such that, for any $(n,m)\in\mathbb N_{\geqslant N}^2$, the map 
\[E_n\otimes_K E_m\longrightarrow E_{n+m},\quad s\otimes t\longmapsto st\]
is surjective. Moreover, if we equip $E_n\otimes E_m$ with the $\varepsilon,\pi$-tensor product of the norm families $\xi_n$ and $\xi_m$, the above map has height $\leqslant 0$. By \cite[Proposition 4.3.31]{CMArakelovAdelic}, one has
\[\widehat{\mu}_{\min}(\overline E_{n+m})\geqslant\widehat{\mu}_{\min}(\overline E_n\otimes_{\varepsilon,\pi}\overline E_m).\]
Moreover, since the field $K$ is assumed to be perfect, by \cite[Corollary 5.6.2]{CMArakelovAdelic} (see also Remark \ref{Rem: positive characteristic minimal slope}), one has
\begin{multline*}\widehat{\mu}_{\min}(\overline E_{n}\otimes_{\varepsilon,\pi}\overline E_m)\geqslant\widehat{\mu}_{\min}(\overline E_n)+\widehat{\mu}_{\min}(\overline E_m) \\
-\frac{3}{2}\nu(\Omega_\infty)(\ln(\dim_{K}(E_n))+\ln(\dim_K(E_m))).\end{multline*}
Note that 
\[\ln(\dim_K(E_n))=O(\ln(n)),\]
and, by \cite[Propositions 6.4.4 and 6.2.7]{CMArakelovAdelic}, there exists a constant $C>0$ such that
\[\widehat{\mu}_{\min}(\overline E_n)\leqslant\widehat{\mu}_{\max}(\overline E_n)\leqslant Cn.\]
Therefore, by \cite[Corollary 3.6]{MR2768967}, we obtain the convergence of the sequence \eqref{Equ: asymptotic minimal slope}. 
\end{proof}

\subsection{Relative ampleness and asymptotic minimal slope}
\label{Sec: Relative ampleness}

\begin{defi} 
Let $\overline L=(L, \varphi)$ be an adelic line bundle on $X$. 
If $L$ is ample, we define the \emph{asymptotic minimal slope of $\overline L$}\index{slope!asymptotic minimal ---} as 
\[\widehat{\mu}^{\mathrm{asy}}_{\min}(\overline L):=\lim_{n\rightarrow+\infty}\frac{\widehat{\mu}_{\min}(f_*(\overline L^{\otimes n}))}{n}.\]
By definition, for any $m\in\mathbb N$ such that $m\geqslant 1$, one has
\begin{equation}\label{Equ: homogeneous}\widehat{\mu}_{\min}^{\mathrm{asy}}(\overline L{}^{\otimes m})=m\operatorname{\widehat{\mu}_{\min}^{\mathrm{asy}}}(\overline L).\end{equation}
\end{defi}

\begin{prop}\label{Pro: suradditivite}
Let $\overline L=(L,\varphi)$ and $\overline M=(M,\psi)$ be adelic line bundles on $X$ such that $L$ and $M$ are ample. Then one has
\begin{equation}\label{Equ: suradditivite of mu min asy}\widehat{\mu}_{\min}^{\mathrm{asy}}(\overline L\otimes\overline M)\geqslant\widehat{\mu}_{\min}^{\mathrm{asy}}(\overline L)+\widehat{\mu}_{\min}^{\mathrm{asy}}(\overline M).\end{equation}
\end{prop}
\begin{proof}
By Lemma \ref{Lem: surjectivity tensor product}, for sufficiently large natural number $n$, the $K$-linear map
\[H^0(X,L^{\otimes n})\otimes_KH^0(X,M^{\otimes n})\longrightarrow H^0(X,(L\otimes M)^{\otimes n}),\quad s\otimes t\longmapsto st\]
is surjective. Moreover, for any $\omega\in\Omega$, the following inequality holds:
\[\forall\,(s,t)\in H^0(X_\omega,L_\omega^{\otimes n})\times H^0(X_\omega, M_\omega^{\otimes n}),\quad \|st\|_{n(\varphi_\omega+\psi_\omega)}\leqslant\|s\|_{n\varphi_\omega}\cdot\|t\|_{n\psi_\omega}.\] Therefore, if we equip $H^0(X,L^{\otimes n})\otimes_KH^0(X,M^{\otimes n})$ with the $\varepsilon,\pi$-tensor product norm family, then the above $K$-linear map has height $\leqslant 0$.
Hence, by \cite[Proposition 4.3.31 and Corollary 5.6.2]{CMArakelovAdelic} (see also Remark \ref{Rem: positive characteristic minimal slope}), we obtain 
\[\begin{split}\widehat{\mu}_{\min}(f_*(&\overline L^{\otimes n}\otimes\overline M^{\otimes n}))\geqslant\widehat{\mu}_{\min}(f_*(\overline L^{\otimes n}))+\widehat{\mu}_{\min}(f_*(\overline M^{\otimes n}))\\
&-\frac 32\nu(\Omega_\infty)\Big(\ln(\dim_K(H^0(X,L^{\otimes n})))+\ln(\dim_K(H^0(X,M^{\otimes n})))\Big).
\end{split}\]
We divide the two sides of the inequality by $n$ and then take the limit when $n\rightarrow+\infty$, using 
\[\lim_{n\rightarrow+\infty}\frac 1n\ln(\dim_K(H^0(X,L^{\otimes n})))=\lim_{n\rightarrow+\infty}\frac 1n\ln(\dim_K(H^0(X,M^{\otimes n}))) =0.\]
we obtain the inequality \eqref{Equ: suradditivite of mu min asy}. 
\end{proof}

\begin{prop}\label{Pro: comparison mu min asy}
Let $L$ be an ample line bundle on $X$ and $\varphi_1$ and $\varphi_2$ be metric families on $L$ such that $(L,\varphi_1)$ and $(L,\varphi_2)$ are both adelic line bundles. Then the following inequality holds:
\begin{equation}\label{Equ: comparison mu min varphi 12}\Big|\operatorname{\widehat{\mu}_{\min}^{\mathrm{asy}}}(L,\varphi_1)-\operatorname{\widehat{\mu}_{\min}^{\mathrm{asy}}}(L,\varphi_2)\Big|\leqslant d(\varphi_1,\varphi_2).\end{equation} 
\end{prop}
\begin{proof}
For any $n\in\mathbb N$, the identity maps 
\[f_*(L^{\otimes n},n\varphi_1)\longrightarrow f_*(L^{\otimes n},n\varphi_2)\]
and
\[f_*(L^{\otimes n},n\varphi_2)\longrightarrow f_*(L^{\otimes n},n\varphi_1)\]
have heights $\leqslant d(n\varphi_1,n\varphi_2)=n\,d(\varphi_1,\varphi_2)$.
By \cite[Proposition 4.3.31]{CMArakelovAdelic}, we obtain that 
\[\Big|\operatorname{\widehat{\mu}_{\min}}(f_*(L^{\otimes n},n\varphi_1))-\operatorname{\widehat{\mu}_{\min}}(f_*(L^{\otimes n},n\varphi_2))\Big|\leqslant n\,d(\varphi_1,\varphi_2).\]
Dividing the two sides of the inequality by $n$ and then taking the limit when $n\rightarrow+\infty$, we obtain \eqref{Equ: comparison mu min varphi 12}. 
\end{proof}

\subsection{Asymptotic slope and intersection number}

In this subsection, we assume that $X$ is integral. Let $\overline L$ be an adelic line bundle on $X$ such that $L$ is ample.  Note that 
\[\dim_K(H^0(X,L^{\otimes n}))=\frac{(L^d)}{d!}n^d+o(n^d),\quad n\rightarrow+\infty.\]
Therefore, one has
\begin{equation}\label{Equ: asymptique sloep}\lim_{n\rightarrow+\infty}\frac{\widehat{\mu}(f_*(\overline L^{\otimes n}))}{n}=\frac{\widehat{\vol}_{\chi}(\overline L)}{(d+1)(L^d)}.\end{equation}
We denote by $\widehat{\mu}^{\mathrm{asy}}(\overline L)$ the value ${\displaystyle \frac{\widehat{\vol}_{\chi}(\overline L)}{(d+1)(L^d)}}$ and call it the \emph{asymptotic slope}\index{slope!asymptotic ---} of $\overline L$. By Theorem \ref{thm:lower:bound}, if $\overline L = (L, \varphi)$ is an adelic line bundle on $X$ such that $L$ is ample and $\varphi$ is semi-positive, that is,
$\overline L$ is relatively ample in the sense of Definition~\ref{def:relatively:ample} below,  then $\widehat{\mathrm{vol}}_\chi(\overline L)=(\overline L^{d+1})_S$ and hence
\begin{equation} \label{eqn:thm:Hilbert:Samuel:01}
\widehat{\mu}^{\mathrm{asy}}(\overline L) = \frac{(\overline L^{d+1})_S}{(d+1)(L^d)}.
\end{equation}

\begin{rema}\label{rema:mu:asy:geq:mu:asy:min}
Let $\overline L$ be an adelic line bundle on $X$ such that $L$ is ample. By definition the following inequality holds:
\begin{equation}\label{Equ: asymptotic slope and minimal slope}\widehat{\mu}^{\mathrm{asy}}(\overline L)\geqslant\widehat{\mu}^{\mathrm{asy}}_{\min}(\overline L).\end{equation}
\end{rema}

\subsection{Lower bound of intersection number for relatively ample adelic line bundles}

\begin{defi}\label{def:relatively:ample}
Let $(L,\varphi)$ be an adelic line bundle on $X$. We say $(L, \varphi)$ is \emph{relatively ample}\index{ample!relatively ---} if $L$ is ample and $\varphi$ is semi-positive. By  \cite[Proposition 2.3.5]{CMArakelovAdelic}, if $\overline L$ and $\overline M$ are relatively ample adelic line bundle, then the tensor product $\overline L\otimes\overline M$ is relatively ample.
\end{defi}

\begin{theo}\label{Thm: lower bound intersection number} 
Let $\overline L_i=(L_i,\varphi_i)$ be a family of relatively ample adelic line bundles on $X$, where $i\in\{0,\ldots,d\}$. For any $i\in\{0,\ldots,d\}$, let 
\[\delta_i=(L_0\cdots L_{i-1}L_{i+1}\cdots L_d).\]
Then the following inequality holds:
\begin{equation}\label{Equ: lower bound intersection prod}(\overline L_0\cdots\overline L_d)_S\geqslant\sum_{i=0}^d\delta_i\operatorname{\widehat{\mu}_{\min}^{\mathrm{asy}}}(\overline L_i).\end{equation}
\end{theo}
\begin{proof}
Without loss of generality, we may assume that $L_0,\ldots,L_d$ are very ample. For any $n\in\mathbb N_{\geqslant 1}$ 
and any $i\in\{0,\ldots,d\}$, we denote by $E_{i,n}$ the $K$-vector space $H^0(X,L_i^{\otimes n})$, and set $r_{i,n}=\dim_K(E_{i,n})-1$. We denote by $\xi_{n\varphi_i}$ the norm family $(\|\ndot\|_{n\varphi_{i,\omega}})_{\omega\in\Omega}$ on $E_{i,n}$, and 
let $\xi_{i,n}$ be a Hermitian norm family on $E_{i,n}$ such that $(E_{i,n},\xi_{i,n})$ forms an adelic vector bundle and that
\[d_\omega(\xi_{i,n},\xi_{n\varphi_i})\leqslant\frac 12\indic_{\Omega_\infty}(\omega)\ln(r_{i,n}+2).\]
The existence of such a Hermitian norm family is ensured by \cite[Theorem~4.1.26]{CMArakelovAdelic}.
Let $\varphi_i^{(n)}$ be the metric family on $L_i$ such that $n\varphi_i^{(n)}$ identifies with the quotient metric family induced by the closed embedding 
$X\rightarrow\mathbb P(E_{i,n})$
and the norm family $\xi_{i,n}$. Since 
\[\lim_{n\rightarrow+\infty}\frac{1}{n}\ln(r_{i,n}+2)=0\]
and the metric families $\varphi_i$ are semi-positive, by \cite[Proposition 3.3.12]{CMIntersection}, we obtain that 
\[\lim_{n\rightarrow+\infty}d(\varphi_i^{(n)},\varphi_i)=\lim_{n\rightarrow+\infty}\int_{\Omega}d_\omega(\varphi_i^{(n)},\varphi_i)\,\nu(\mathrm{d}\omega)=0.\]

For any $n\in\mathbb N_{\geqslant 1}$,
let $R_n$ be the one-dimensional vector space of 
\begin{equation}\label{Equ: resultant space}S^{n^d\delta_0}(E_{0,n}^\vee)\otimes_k\cdots\otimes_kS^{n^d\delta_d}(E_{d,n}^\vee)\end{equation}
spanned by any resultant of the closed embeddings $X\rightarrow\mathbb P(E_{i,n})$. We equip each $S^{n^d\delta_0}(E_{i,n}^\vee)$ with the orthogonal symmetric power norm family of $\xi_{i,n}^\vee$, and the tensor product space \eqref{Equ: resultant space} with the orthogonal tensor product norm family.   By \cite[Remark 4.2.14]{CMIntersection} and \cite[Corollary 1.4.3 and Lemma 4.3.8]{MR1260106}, we obtain that 
\begin{equation}\label{Equ: lower bound inter prod}\begin{split}((L_{0},\varphi_0^{(n)})\cdots&(L_d,\varphi_d^{(n)}))_S\geqslant-\frac{1}{n^{d+1}}\bigg(\operatorname{\widehat{\deg}}(\overline R_n)+\nu(\Omega_\infty)\sum_{i=0}^d\ln\binom{r_{i,n}+n^d\delta_i}{n^d\delta_i}\bigg)\\
&\geqslant-\frac{1}{n^{d+1}}\big(\operatorname{\widehat{\deg}}(\overline R_n)+\nu(\Omega_\infty)\sum_{i=0}^d n^d\delta_i\ln(r_{i,n}+1)\big),
\end{split}\end{equation}
where the second inequality comes from 
\[\forall\,(a,b)\in\mathbb N_{\geqslant 1}^2,\quad \binom{a+b}{b}\leqslant (a+1)^b.\]
Note that 
\begin{equation}\label{Equ: upper bound degree Rn}\widehat{\deg}(\overline R_n)\leqslant\widehat{\mu}_{\max}(S^{n^d\delta_0}(E_{0,n}^\vee,\xi_{0,n}^\vee)\otimes\cdots\otimes S^{n^d\delta_d}(E_{d,n}^\vee,\xi_{d,n}^\vee)).\end{equation}
In the case where $K$ is of characteristic $0$, by Remark \ref{Rem: operator norm symmetri cpower} and \cite[Proposition 4.3.31]{CMArakelovAdelic}, we obtain 
\begin{equation}\label{Equ: upper bound mu max Sym}\begin{split}&\quad\;\widehat{\mu}_{\max}\big(S^{n^d\delta_0}(E_{0,n}^\vee,\xi_{0,n}^\vee)\otimes\cdots\otimes S^{n^d\delta_d}(E_{d,n}^\vee,\xi_{d,n}^\vee)\big)\\
&\leqslant\widehat{\mu}_{\max}\big((E_{0,n}^\vee,\xi_{0,n}^\vee)^{\otimes n^d\delta_0}\otimes\cdots\otimes(E_{d,n}^\vee,\xi_{d,n}^\vee)^{\otimes n^d\delta_d}\big)+\nu(\Omega_\infty)\sum_{i=0}^dn^d\delta_i\ln(n^d\delta_i).
\end{split}\end{equation}
By \cite[Corollaries 4.3.27 and 5.6.2]{CMArakelovAdelic}, we have  
\begin{equation}\label{Equ: majoration mu max}\begin{split}
&\quad\;\widehat{\mu}_{\max}\big((E_{0,n}^\vee,\xi_{0,n}^\vee)^{\otimes n^d\delta_0}\otimes\cdots\otimes(E_{d,n}^\vee,\xi_{d,n}^\vee)^{\otimes n^d\delta_d}\big)\\
&\leqslant\sum_{i=0}^d n^d\delta_i\Big(\widehat{\mu}_{\max}(E_{i,n}^\vee,\xi_{i,n}^\vee)+\frac 12\nu(\Omega_\infty)\ln(r_{i,n}+1)\Big)\\
&=\sum_{i=0}^dn^d\delta_i\Big(-\widehat{\mu}_{\min}(E_{i,n},\xi_{i,n})+\frac 12\nu(\Omega_\infty)\ln(r_{i,n}+1)\Big)
\end{split}\end{equation}
Combining \eqref{Equ: lower bound inter prod}, \eqref{Equ: upper bound degree Rn}, \eqref{Equ: upper bound mu max Sym} and \eqref{Equ: majoration mu max}, we obtain 
\begin{equation}\label{Equ: lower bound of int numb}\begin{split}((L_{0},\varphi_0^{(n)})\cdots&(L_d,\varphi_d^{(n)}))_S\geqslant \sum_{i=0}^d\delta_i\frac{\widehat{\mu}_{\min}(E_{i,n},\xi_{i,n})}{n}\\
&-\frac 32\nu(\Omega_\infty)\sum_{i=0}^d\frac{\delta_i}{n}
\ln(r_{i,n}+1)-\nu(\Omega_\infty)\sum_{i=0}^d\frac{\delta_i}{n}\ln(n^d\delta_i).\end{split}\end{equation}
In the case where $K$ is of positive characteristic, by 
Corollary \ref{Cor: tensorial semistability} and Theorem \ref{Thm: maximal slope of symm power}, we obtain 
\[\widehat{\mu}_{\max}\big(S^{n^d\delta_0}(E_{0,n}^\vee,\xi_{0,n}^\vee)\otimes\cdots\otimes S^{n^d\delta_d}(E_{d,n}^\vee,\xi_{d,n}^\vee)\big)\leqslant\sum_{i=0}^dn^d\delta_i\operatorname{\widehat{\mu}_{\max}}(E_{i,n}^\vee,\xi_{i,n}^\vee).\]
Hence the inequality \eqref{Equ: lower bound of int numb} still holds in this case.
Since $r_{i,n}=O(n^d)$, taking the limit when $n$ goes to the infinity, we obtain the inequality \eqref{Equ: lower bound intersection prod}. 
\end{proof}

\subsection{Relative nefness and continuous extension of $\widehat{\mu}_{\min}^{\mathrm{asy}}$}

\begin{prop}\label{Pro: convergence and independence}
Let $\overline L$ and $\overline A$ be adelic line bundle on $X$. Assume that $L$ is  nef and $A$ is ample. Then the sequence \begin{equation}\label{Equ: sequence defiing mu min asy}\frac 1n\operatorname{\widehat{\mu}_{\min}^{\mathrm{asy}}
}(\overline L^{\otimes n}\otimes\overline A),\quad n\in\mathbb N_{\geqslant 1}\end{equation}
converges
in $\mathbb R\cup\{-\infty\}$, 
and the limit does not depend on the choice of $\overline A$. In particular, in the case where $L$ is ample, the following equality holds:
\begin{equation}\label{Equ: limit is overline L}\lim_{n\rightarrow+\infty}\frac{1}{n}\operatorname{\widehat{\mu}_{\min}}(\overline L^{\otimes n}\otimes\overline A)=\operatorname{\widehat{\mu}_{\min}}(\overline L).\end{equation}
\end{prop}
\begin{proof}
Let $p$ be a positive integer. By Proposition \ref{Pro: suradditivite}, for any $\ell\in\mathbb N_{\geqslant 1}$ and any $r\in \{1,\ldots, p\}$, one has 
\[\begin{split}&\quad\;\widehat{\mu}_{\min}^{\mathrm{asy}}(\overline L^{\otimes p}\otimes\overline A)=\frac{1}{\ell+1}\operatorname{\widehat{\mu}_{\min}^{\mathrm{asy}}}(\overline L^{\otimes(\ell+1)p}\otimes\overline A^{\otimes(\ell+1)})\\
&\geqslant\frac{1}{\ell+1}\Big(\operatorname{\widehat{\mu}_{\min}^{\mathrm{asy}}}(\overline L^{\otimes(\ell p+r)}\otimes\overline A)+\operatorname{\widehat{\mu}_{\min}^{\mathrm{asy}}}(\overline L^{\otimes(p-r)}\otimes \overline A)+(\ell-1)\operatorname{\widehat{\mu}_{\min}^{\mathrm{asy}}}(\overline A)\Big).
\end{split}\] 
Taking the limit superior when $\ell p+r\rightarrow+\infty$, we obtain 
\[\operatorname{\widehat{\mu}_{\min}^{\mathrm{asy}}}(\overline L^{\otimes p}\otimes\overline A)\geqslant p\limsup_{n\rightarrow+\infty}\frac 1n\operatorname{\widehat{\mu}_{\min}^{\mathrm{asy}}}(\overline L^{\otimes n}\otimes\overline A)+\operatorname{\widehat{\mu}_{\min}^{\mathrm{asy}}}(\overline A),\]
which leads to 
\[\liminf_{p\rightarrow+\infty}\frac 1p\operatorname{\widehat{\mu}_{\min}^{\mathrm{asy}}}(\overline L^{\otimes p}\otimes\overline A)\geqslant\limsup_{n\rightarrow+\infty}\frac 1n\operatorname{\widehat{\mu}_{\min}^{\mathrm{asy}}}(\overline L^{\otimes n}\otimes\overline A).\]
Therefore the sequence \eqref{Equ: sequence defiing mu min asy} converges in $[-\infty,+\infty]$.
Moreover, still by Proposition~\ref{Pro: suradditivite}, for any $p\in\mathbb N_{\geqslant 1}$, one has
\[\operatorname{\widehat{\mu}_{\min}^{\mathrm{asy}}}(\overline L\otimes\overline A)=\frac 1p\operatorname{\widehat{\mu}_{\min}^{\mathrm{asy}}}(\overline L^{\otimes p}\otimes\overline A^{\otimes p})\geqslant \frac 1p\operatorname{\widehat{\mu}_{\min}^{\mathrm{asy}}}(\overline L^{\otimes p}\otimes\overline A)+\frac{p-1}{p}\operatorname{\widehat{\mu}_{\min}^{\mathrm{asy}}}(\overline A),\]
which shows that 
\[\lim_{p\rightarrow+\infty}\frac 1p\operatorname{\widehat{\mu}_{\min}^{\mathrm{asy}}}(\overline L^{\otimes p}\otimes\overline A)\leqslant\operatorname{\widehat{\mu}_{\min}^{\mathrm{asy}}}(\overline L\otimes\overline A)-\operatorname{\widehat{\mu}_{\min}^{\mathrm{asy}}}(\overline A)<+\infty.\]

To prove the second assertion, we first show that the limit of the sequence does not depend on the choice of the metric family on $\overline A$. For this purpose, we consider two metric families $\varphi_1$ et $\varphi_2$ on $A$ such that both $(A,\varphi_1)$ and $(A,\varphi_2)$ are adelic line bundles on $X$. By Proposition \ref{Pro: comparison mu min asy}, for any $n\in\mathbb N$ one has
\[\mathopen{\Big |}\operatorname{\widehat{\mu}_{\min}^{\mathrm{asy}}}(\overline L^{\otimes n}\otimes(A,\varphi_1))-\operatorname{\widehat{\mu}_{\min}^{\mathrm{asy}}}(\overline L^{\otimes n}\otimes(A,\varphi_1))\mathclose{\Big |}\leqslant d(\varphi_1,\varphi_2),\]
so that 
\begin{equation}\label{Equ: do not depend on metric}\lim_{n\rightarrow+\infty}\frac 1n \operatorname{\widehat{\mu}_{\min}^{\mathrm{asy}}}(\overline L^{\otimes n}\otimes(A,\varphi_1))=\lim_{n\rightarrow+\infty}\frac 1n\operatorname{\widehat{\mu}_{\min}^{\mathrm{asy}}}(\overline L^{\otimes n}\otimes(A,\varphi_2)).\end{equation}

We then show that, for any $p\in\mathbb N_{\geqslant 2}$, the following inequality holds:
\begin{equation}\label{Equ: limit equality}\lim_{n\rightarrow+\infty}\frac{1}{n}\operatorname{\widehat{\mu}_{\min}^{\mathrm{asy}}}(\overline L^{\otimes n}\otimes\overline A)=\lim_{n\rightarrow+\infty}\frac{1}{n}\operatorname{\widehat{\mu}_{\min}^{\mathrm{asy}}}(\overline L^{\otimes n}\otimes\overline A^{\otimes p}).\end{equation}
In fact, by \eqref{Equ: homogeneous}, for any $n\in\mathbb N_{\geqslant 1}$ one has
\[\frac{1}{n}\operatorname{\widehat{\mu}_{\min}^{\mathrm{asy}}}(\overline L^{\otimes n}\otimes\overline A)=\frac{1}{np}\operatorname{\widehat{\mu}_{\min}^{\mathrm{asy}}}(\overline L^{\otimes np}\otimes\overline A^{\otimes p}).\]
Taking the limit when $n\rightarrow+\infty$, we obtain the equality \eqref{Equ: limit equality}. 

Note that if $\overline B$ is another adelic line bundle such that $B$ is ample, then the following inequality holds:
\begin{equation}\label{Equ: inequality mult by B}
\lim_{n\rightarrow+\infty}\frac{1}{n}\operatorname{\widehat{\mu}_{\min}^{\mathrm{asy}}}(\overline L^{\otimes n}\otimes\overline A)\leqslant \lim_{n\rightarrow+\infty}\frac{1}{n}\operatorname{\widehat{\mu}_{\min}^{\mathrm{asy}}}(\overline L^{\otimes n}\otimes\overline A\otimes\overline B).
\end{equation}
In fact, by Proposition \ref{Pro: suradditivite}, for any $n\in\mathbb N_{\geqslant 1}$, one has
\[\frac{1}{n}\operatorname{\widehat{\mu}_{\min}^{\mathrm{asy}}}(\overline L^{\otimes n}\otimes\overline A\otimes\overline B)\geqslant \frac{1}{n}\operatorname{\widehat{\mu}_{\min}^{\mathrm{asy}}}(\overline L^{\otimes n}\otimes\overline A)+\frac 1n\operatorname{\widehat{\mu}_{\min}^{\mathrm{asy}}}(\overline B).\]
Taking the limit when $n\rightarrow+\infty$, we obtain \eqref{Equ: inequality mult by B}.

Finally, we show that, if $\overline B$ is an arbitrary adelic line bundle such that $B$ is ample, then the equality 
\begin{equation}\label{Equ: does not depend on line bundle}
\lim_{n\rightarrow+\infty}\frac{1}{n}\operatorname{\widehat{\mu}_{\min}^{\mathrm{asy}}}(\overline L^{\otimes n}\otimes\overline A)=\lim_{n\rightarrow+\infty}\frac{1}{n}\operatorname{\widehat{\mu}_{\min}^{\mathrm{asy}}}(\overline L^{\otimes n}\otimes\overline B)
\end{equation}
holds. In fact, there exists $p\in\mathbb N_{\geqslant 1}$ such that $N=B^{\otimes p}\otimes A^\vee$ is ample. We equip it with an arbitrary  metric family such that $\overline N$ forms an adelic line bundle.  By \eqref{Equ: inequality mult by B} we obtain 
\[\lim_{n\rightarrow+\infty}\frac{1}{n}\operatorname{\widehat{\mu}_{\min}^{\mathrm{asy}}}(\overline L^{\otimes n}\otimes\overline A)\leqslant\lim_{n\rightarrow+\infty}\frac{1}{n}\operatorname{\widehat{\mu}_{\min}^{\mathrm{asy}}}(\overline L^{\otimes n}\otimes\overline A\otimes\overline N). \]
Since $A\otimes N$ is isomorphic to $ B^{\otimes p}$, by \eqref{Equ: do not depend on metric} and \eqref{Equ: limit equality} we obtain 
\[\lim_{n\rightarrow+\infty}\frac{1}{n}\operatorname{\widehat{\mu}_{\min}^{\mathrm{asy}}}(\overline L^{\otimes n}\otimes\overline A\otimes\overline N)=\lim_{n\rightarrow+\infty}\frac{1}{n}\operatorname{\widehat{\mu}_{\min}^{\mathrm{asy}}}(\overline L^{\otimes n}\otimes\overline B).\]
Therefore, we deduce 
\[\lim_{n\rightarrow+\infty}\frac{1}{n}\operatorname{\widehat{\mu}_{\min}^{\mathrm{asy}}}(\overline L^{\otimes n}\otimes\overline A)\leqslant\lim_{n\rightarrow+\infty}\frac{1}{n}\operatorname{\widehat{\mu}_{\min}^{\mathrm{asy}}}(\overline L^{\otimes n}\otimes\overline B).\]
Interchanging the roles of $\overline A$ and $\overline B$ we obtain the converse inequality.

To obtain the equality \eqref{Equ: limit is overline L}, it suffices to apply the equality \eqref{Equ: does not depend on line bundle}
 in the particular case where $\overline A=\overline L$. The proposition is thus proved.
\end{proof}

\begin{defi}\label{def:relatively:nef}
Let $\overline L$ be an adelic line bundle on $X$ such that $L$ is nef, we define 
\begin{equation*}\widehat{\mu}_{\min}^{\mathrm{asy}}(\overline L):=\lim_{n\rightarrow+\infty}\frac{1}{n}\operatorname{\widehat{\mu}_{\min}^{\mathrm{asy}}}(\overline L{}^{\otimes n}\otimes\overline A),\]
where $\overline A$ is an arbitrary adelic line bundle such that $A$ is ample. The element $\widehat{\mu}_{\min}^{\mathrm{asy}}(\overline L)$ of $\mathbb R\cup\{-\infty\}$ is called \emph{asymptotic minimal slope}\index{slope!asymptotic minimal ---} of $\overline L$.
\end{defi}

\begin{rema}
It is an interesting question to ask when the asymptotic minimal slope is a real number. As we will show in Theorem \ref{Thm: mu min asy and pull-back}, the asymptotic minimal slope does not decrease if we replace the adelic line bundle by its pullback by a projective morphism. In particular, if $L$ is the pullback of an ample line bundle by a projective morphism, then $\widehat{\mu}_{\min}^{\mathrm{asy}}(\overline L)\in\mathbb R$. 
\end{rema}

\begin{prop}\label{Pro: super additivity nef}
Let $\overline L$ and $\overline M$ be adelic line bundles on $X$ such that $L$ and $M$ are nef. One has
\begin{equation}\label{Equ: mu min super additive}\operatorname{\widehat{\mu}_{\min}^{\mathrm{asy}}}(\overline L\otimes \overline M) \geqslant \operatorname{\widehat{\mu}_{\min}^{\mathrm{asy}}}(\overline L)+\operatorname{\widehat{\mu}_{\min}^{\mathrm{asy}}}(\overline M).\end{equation}
Moreover, one has 
\begin{equation}\label{Equ: equality limit if mu min not infinite}\lim_{n\rightarrow+\infty}\frac 1n\operatorname{\widehat{\mu}_{\min}^{\mathrm{asy}}}(\overline L^{\otimes n}\otimes\overline M)=\operatorname{\widehat{\mu}_{\min}^{\mathrm{asy}}}(\overline L)\end{equation}
provided that $\operatorname{\widehat{\mu}_{\min}^{\mathrm{asy}}}(\overline M)>-\infty$.
\end{prop}
\begin{proof}
Let $\overline A$ be an adelic line bundle on $X$ such that $A$ is ample. For any $n\in\mathbb N_{\geqslant 1}$, by Proposition \ref{Pro: suradditivite} one has 
\[\frac{1}{n}\operatorname{\widehat{\mu}_{\min}^{\mathrm{asy}}}(\overline L^{\otimes n}\otimes\overline M^{\otimes n}\otimes\overline A^{\otimes 2})\geqslant
\frac 1n\operatorname{\widehat{\mu}_{\min}^{\mathrm{asy}}}(\overline L^{\otimes n}\otimes\overline A)+\operatorname{\widehat{\mu}_{\min}^{\mathrm{asy}}}(\overline M^{\otimes n}\otimes\overline A). \]
Taking the limit when $n\rightarrow+\infty$, we obtain the inequality \eqref{Equ: mu min super additive}. 

By \eqref{Equ: mu min super additive}, we obtain that, for any positive integer $n$, the inequality 
\[\frac{1}{n}\operatorname{\widehat{\mu}_{\min}^{\mathrm{asy}}}(\overline L^{\otimes n}\otimes \overline M)\geqslant \operatorname{\widehat{\mu}_{\min}^{\mathrm{asy}}}(\overline L)+\frac 1n \operatorname{\widehat{\mu}_{\min}^{\mathrm{asy}}}(\overline M).\]
Since $\operatorname{\widehat{\mu}_{\min}^{\mathrm{asy}}}(\overline M)\in\mathbb R$, taking the limit inferior when $n\rightarrow+\infty$, we obtain
\[\liminf_{n\rightarrow+\infty}\frac 1n\operatorname{\widehat{\mu}_{\min}^{\mathrm{asy}}}(\overline L^{\otimes n}\otimes\overline M)\geqslant \operatorname{\widehat{\mu}_{\min}^{\mathrm{asy}}}(\overline L).\]
Pick an adelic line bundle $\overline A$ on $X$ such that $A$ is ample. Since $A\otimes M$ is ample, one has 
\begin{equation}\label{Equ: limit tensor M A}\lim_{n\rightarrow+\infty }\frac 1n\operatorname{\widehat{\mu}_{\min}^{\mathrm{asy}}}(\overline L^{\otimes n}\otimes\overline M\otimes\overline A)=\operatorname{\widehat{\mu}_{\min}^{\mathrm{asy}}}(\overline L).\end{equation}
Moreover, by  \eqref{Equ: mu min super additive} one has 
\[\frac 1n\operatorname{\widehat{\mu}_{\min}^{\mathrm{asy}}}(\overline L^{\otimes n}\otimes\overline M\otimes A)\geqslant \frac 1n\operatorname{\widehat{\mu}_{\min}^{\mathrm{asy}}}(\overline L^{\otimes n}\otimes\overline M)+\frac 1n\operatorname{\widehat{\mu}_{\min}^{\mathrm{asy}}}(\overline A).\]
Taking the limit superior, by \eqref{Equ: limit tensor M A}  we obtain 
\[\limsup_{n\rightarrow+\infty}\frac 1n\operatorname{\widehat{\mu}_{\min}^{\mathrm{asy}}}(\overline L^{\otimes n}\otimes\overline M)\leqslant\operatorname{\widehat{\mu}_{\min}^{\mathrm{asy}}}(\overline L).\] 
Hence the equality \eqref{Equ: equality limit if mu min not infinite} holds.
\end{proof}

\begin{defi}
Let $\overline L=(L,\varphi)$ be an adelic line bundle on $X$. We say that $\overline L$ is \emph{relatively nef}\index{nef!relatively ---} if there exists a relatively ample adelic line bundle $\overline A$ on $X$ and a positive integer $N$ such that,
for any $n\in\mathbb N_{\geqslant N}$, the tensor product $\overline L{}^{\otimes n}\otimes\overline A$ is relatively ample.
\end{defi}

\begin{prop}\label{Pro: semi-ample is nef}
Let $\overline L=(L,\varphi)$ be an adelic line bundle on $X$ such that $L$ is semi-ample and $\varphi$ is semi-positive. Then, for any adelic line bundle $\overline A=(A,\psi)$ on $X$ which is relatively ample and any $n\in\mathbb N$, the tensor product $\overline L^{\otimes n}\otimes\overline A$ is relatively ample. In particular, $\overline L$ is relatively nef. 
\end{prop}
\begin{proof}
Since $L$ is semi-ample, we obtain that, for any $n\in\mathbb N$, $L^{\otimes n}\otimes A$ is ample.   Moreover, by \cite[Proposition 2.3.5]{CMArakelovAdelic}, $n\varphi+\psi$ is semi-positive. Hence $\overline L^{\otimes n}\otimes \overline A$ is relatively ample.
\end{proof}

\begin{prop}
\label{Pro: tensor product ref nef} Let $\overline L$ and $\overline M$ be adelic line bundles on $X$ which are relatively nef. Then the tensor product $\overline L\otimes\overline M$ is also relatively nef.
\end{prop}
\begin{proof}
Let $\overline A$ and $\overline B$ be relatively ample adelic line bundles on $X$, and $N$ be a positive integer such that $\overline L^{\otimes n}\otimes\overline A$ and $\overline M^{\otimes n}\otimes\overline B$ are relatively ample for any integer $n\geqslant N$. We then obtain that $(\overline L\otimes\overline M)^{\otimes n}\otimes(\overline A\otimes\overline B)$ is relatively ample. Therefore $\overline L\otimes\overline M$ is relatively nef.
\end{proof}

\begin{prop}\label{Pro: lower bound ample} 
Let $\overline L_0,\ldots,\overline L_d$ be a family of relatively nef adelic line bundles on $X$. For any $i\in\{0,\ldots,d\}$, let
\[\delta_i=(L_0\cdots L_{i-1}L_{i+1}\cdots L_d).\]
Assume that $\delta_i>0$ for those $i\in\{0,\ldots,d\}$ such that $\widehat{\mu}_{\min}^{\mathrm{asy}}(\overline L_i)=-\infty$.
Then the following inequality holds:
\begin{equation}\label{Equ: minormation intersection product}(\overline L_0\cdots\overline L_d)_S\geqslant\sum_{i=0}^d\delta_i\,\widehat{\mu}_{\min}^{\mathrm{asy}}(\overline L_i).\end{equation}
\end{prop}
\begin{proof}
If there is $i \in \{ 0, \ldots, d \}$ such that $\widehat{\mu}_{\min}^{\mathrm{asy}}(\overline L_i)=-\infty$, then the assertion is obvious, so that we may assume that
$\widehat{\mu}_{\min}^{\mathrm{asy}}(\overline L_i)> -\infty$ for all $i \in \{ 0, \ldots, d \}$.

Let $\overline A_i$ be a relatively ample adelic line bundle on $X$ such that $\overline L_i^{\otimes n}\otimes\overline A_i$ is
relatively ample for sufficiently large positive integer $n$.
For any $i\in\{0,\ldots,d\}$ and any positive integer $n$, let 
\begin{gather*}\overline{L}_{i,n}=\overline L_i^{\otimes n}\otimes\overline A_i,\\
\delta_{i,n}=(L_{0,n}\cdots L_{i-1,n}L_{i+1,n}\cdots L_{d,n}).\end{gather*}
By the multi-linearity of intersection product, we obtain that 
\[\lim_{n\rightarrow+\infty}\frac{\delta_{i,n}}{n^{d}}=\delta_i,\quad\lim_{n\rightarrow+\infty}\frac{(\overline L_{0,n}\cdots\overline L_{d,n})_S}{n^{d+1}}=(\overline L_0\cdots\overline L_d)_S.\]
Note that Theorem \ref{Thm: lower bound intersection number} leads to 
\[(\overline L_{0,n}\cdots\overline L_{d,n})_S\geqslant\sum_{i=0}^d\delta_{i,n}\operatorname{\widehat{\mu}_{\min}^{\mathrm{asy}}}(\overline L_{i,n})\]
for sufficiently large positive integer $n$.
Dividing the two sides by $n^{d+1}$ and then taking the limit when $n\rightarrow+\infty$, we obtain the inequality \eqref{Equ: minormation intersection product}.
\end{proof}

\subsection{Generalized Hodge index theorem}

\begin{theo}[Generalized Hodge index theorem]\label{coro:lower:bound:nef}
Let $(L, \varphi)$ be a relatively nef adelic invertible $\mathcal O_X$-module.
Then one has
\begin{equation}
\label{Equ: lower bound of vol}
\widehat{\vol}(L,\varphi)\geqslant ((L,\varphi)^{d+1}).
\end{equation}
\end{theo}

\begin{proof}
Let $(A, \psi)$ be a relatively ample adelic invertible $\mathcal O_X$-module and $n_0\in\mathbb N$ such that 
$(L, \varphi)^{\otimes n} \otimes (A, \psi)$ is relatively ample for $n \in\mathbb N_{\geqslant n_0}$. 
Then, by Theorem~\ref{thm:lower:bound},
\[
\avol((L, \varphi)^{\otimes n} \otimes (A, \psi)) \geqslant \avol_\chi((L, \varphi)^{\otimes n} \otimes (A, \psi)) \geqslant (((L, \varphi)^{\otimes n} \otimes (A, \psi))^{d+1})_S
\] 
for $n \geqslant n_0$, 
and hence by \cite[Theorem~6.4.14 and Theorem~6.4.24]{CMArakelovAdelic},
\begin{align*}
\avol(L, \varphi) & =  \lim_{n\to+\infty} \frac{1}{n^{d+1}}\avol((L, \varphi)^{\otimes n} \otimes (A, \psi))  \\
& \geqslant  \lim_{n\to+\infty} \frac{1}{n^{d+1}}  (((L, \varphi)^{\otimes n} \otimes (A, \psi))^{d+1})_S = ((L,\varphi)^{d+1}),
\end{align*}
as desired.
\end{proof}

\begin{coro}\label{Cor: bigness of line bundle}
Let $(L, \varphi)$ be a relatively nef adelic invertible $\mathcal O_X$-module.
If $((L,\varphi)^{d+1}) > 0$, then $L$ is big.
\end{coro}

\begin{proof}
By Corollary~\ref{coro:lower:bound:nef}, $\avol(L, \varphi) > 0$.
Let $(A, \psi)$ be a relatively ample adelic invertible $\mathcal O_X$-module.
By the continuity of $\avol$ (see \cite[Theorem 6.4.24]{CMArakelovAdelic}), there is a positive integer $n$ such that 
$\avol((L, \varphi)^{\otimes n} \otimes (A,\psi)^\vee) > 0$, so that, for some positive integer $m$,  $H^0(X, (L^{\otimes n} \otimes A^{\vee})^{\otimes m}) \not= \{ 0 \}$.
Therefore $L$ is big.
\end{proof}

\subsection{Pull-back by a projective morphism}

\begin{lemm}\label{lemma:restriction:relatively:nef}
If $\overline L=(L,\varphi)$ is a relatively nef adelic line bundle on $X$ and if  $g:Y\rightarrow X$ is projective morphism from a reduced $K$-scheme $Y$ to $X$, then the pull-back $g^*(\overline L)$ is a relatively nef adelic line bundle on $Y$.
\end{lemm}

\begin{proof}
Let $\overline{A}=(A,\psi)$ be a relatively ample line bundle on $X$ and $N$ be a positive integer such that $\overline{L}^{\otimes n} \otimes \overline{A}=(L^{\otimes n}\otimes A,n\varphi+\psi)$ is relatively ample for any $n\in\mathbb N_{\geqslant N}$.
Note that $L^{\otimes n}\otimes A$ is ample and hence $g^*(L)^{\otimes n}\otimes g^*(A)$ is semi-ample. Moreover, by Lemma~\ref{lemma:restriction:semi-positive}, $n\varphi+\psi$ is semi-positive. We choose an arbitrary relatively ample adelic line bundle $\overline B$ on $Y$. By Proposition \ref{Pro: semi-ample is nef}, we obtain that $g^*(\overline L)^{\otimes n}\otimes (g^*(\overline A)\otimes\overline B)$ is relatively ample for any $n\in\mathbb N_{\geqslant N}$.  
Thus the assertion follows.
\end{proof}

\begin{prop}\label{Pro: mu min of restriction}
Let $\overline L=(L,\varphi)$ be an adelic line bundle on $X$ such that $L$ is nef. For any non-empty and reduced closed subscheme $Y$ of $X$, the following inequality holds:
\begin{equation}\label{Equ: minoration mu muni asym}\widehat{\mu}_{\min}^{\mathrm{asy}}(\overline L|_Y)\geqslant\widehat{\mu}_{\min}^{\mathrm{asy}}(\overline L).\end{equation}
\end{prop}
\begin{proof}
We first consider the case where $L$ is ample. Clearly the restriction of $L$ to $Y$ is ample, and there exists $n_0\in\mathbb N$ such that the restriction map
\[\pi_n:H^0(X,L^{\otimes n})\longrightarrow H^0(Y,L|_Y^{\otimes n})\]
is surjective for any $n\in\mathbb N_{\geqslant n_0}$. Moreover, if we denote by $\varphi_\omega^{Y}$ the restriction of the metric $\varphi_\omega$ to $L_\omega|_{Y_\omega}$, then, for any $s\in H^0(X_\omega,L_{\omega}^{\otimes n})$, the inequality
\begin{equation*}\|s\|_{n\varphi_\omega}\geqslant\|\pi_{n,\omega}(s)\|_{n\varphi^Y_\omega}\end{equation*} holds, so that, by
\cite[Proposition 4.3.31]{CMArakelovAdelic}, we obtain
\[\operatorname{\widehat{\mu}_{\min}}(H^0(Y,L|_Y^{\otimes n}),(\|\ndot\|_{n\varphi^Y_\omega})_{\omega\in\Omega})\geqslant\operatorname{\widehat{\mu}_{\min}}(H^0(X,L^{\otimes n}),(\|\ndot\|_{n\varphi_\omega})_{\omega\in\Omega})\] 
for any $n\in\mathbb N_{\geqslant n_0}$.
Dividing the two sides of the inequality by $n$ and taking the limit when $n\rightarrow+\infty$, we obtain the inequality \eqref{Equ: minoration mu muni asym}.

In general, let $\overline A$ be an adelic line bundle on $X$ such that $A$ is ample. By the above argument, one has $\widehat{\mu}_{\min}^{\mathrm{asy}}(\overline A|_Y)\geqslant \widehat{\mu}_{\min}^{\mathrm{asy}}(\overline A)>-\infty$. Since $L$ is nef, $L|_Y$ is also nef (see \cite[Example 1.4.4]{MR2095471}) and therefore $\widehat{\mu}_{\min}^{\mathrm{asy}}(\overline L|_Y)$ is well defined. By \eqref{Equ: equality limit if mu min not infinite} and the above case, one has 
\[\widehat{\mu}_{\min}^{\mathrm{asy}}(\overline L|_Y)=\lim_{n\rightarrow+\infty}\frac 1n \widehat{\mu}_{\min}^{\mathrm{asy}}(\overline L|_Y^{\otimes n}\otimes \overline A|_Y)\geqslant\lim_{n\rightarrow+\infty}\frac 1n \widehat{\mu}_{\min}^{\mathrm{asy}}(\overline L^{\otimes n}\otimes\overline A)= \widehat{\mu}_{\min}^{\mathrm{asy}}(\overline L),\]
as required.
\end{proof}

\begin{prop}\label{Pro: minormation int prod}
Let $Y$ be a reduced  and non-empty closed subscheme of $X$ and $r$ be the dimension of $Y$. Let $\overline L_0,\ldots,\overline L_r$ be a family of relatively nef adelic line bundles on $X$. For any $i\in\{0,\ldots,r\}$, let 
\[\delta_i=(L_0|_Y\cdots L_{i-1}|_YL_{i+1}|_Y\cdots L_r|_Y).\]
Assume that, for any $i\in\{0,\ldots,r\}$, $\delta_i>0$ once $\operatorname{\widehat{\mu}_{\min}^{\mathrm{asy}}}(\overline L_i|_Y)=-\infty$. Then the following inequality holds:
\begin{equation}\label{Equ: restriction interesection product}(\overline L_0|_{Y}\cdots\overline L_r|_{Y})_S\geqslant\sum_{i=0}^r\delta_i\,\operatorname{\widehat{\mu}_{\min}^{\mathrm{asy}}}(\overline L_i|_Y).\end{equation} 
\end{prop}
\begin{proof}
This is a consequence of Proposition~\ref{Pro: lower bound ample} and Lemma~\ref{lemma:restriction:relatively:nef}.
\end{proof}

\begin{prop}\label{Pro: lower bound intersection number} 
Let $\overline L_0=(L_0,\varphi_0),\ldots,\overline L_d=(L_d,\varphi_d)$ be a family of adelic line bundles on $X$. For any $i\in\{0,\ldots,d\}$, let   
\[\delta_i=(L_0\cdots L_{i-1}L_{i+1}\cdots L_d).\]
Assume that $\overline L_1,\ldots,\overline L_d$ are relatively nef, $L_0$ admits a global section $s$ which is a regular meromorphic section, and, for any $i\in\{1,\ldots,d\}$, $\delta_i>0$ once $\operatorname{\widehat{\mu}_{\min}^{\mathrm{asy}}}(\overline L_i)=-\infty$. Then the following inequality holds:
\begin{equation}\label{Equ: lower bound with section}\begin{split}(\overline L_0\cdots\overline L_d)_S\geqslant \sum_{i=1}^d&\delta_i\operatorname{\widehat{\mu}_{\min}^{\mathrm{asy}}}(\overline L_i)\\
&-\int_{\Omega}\int_{X_\omega^{\mathrm{an}}}\ln\|s\|_{\varphi_{0,\omega}}c_1(L_{1,\omega},\varphi_{1,\omega})\cdots c_1(L_{d,\omega},\varphi_{d,\omega})\,\nu(\mathrm{d}\omega).\end{split}\end{equation}
\end{prop}
\begin{proof}
If there exists $i\in\{1,\ldots,d\}$ such that $\operatorname{\widehat{\mu}_{\min}^{\mathrm{asy}}}(\overline L_i)=-\infty$, then the inequality \eqref{Equ: lower bound with section} is trivial. Therefore, we can assume that $\operatorname{\widehat{\mu}_{\min}^{\mathrm{asy}}}(\overline L_i)\in\mathbb R$ for any $i\in\{1,\ldots,n\}$.
Let $\operatorname{div}(s)=a_1Z_1+\cdots+a_nZ_n$ be the decomposition of $\operatorname{div}(s)$ as linear combination of prime divisors, where $a_1,\ldots,a_n$ are non-negative integers since $s$ is a global section. By proposition \ref{Pro: mu min of restriction}, for any $i\in\{1,\ldots,d\}$ and and any $j\in\{1,\ldots,n\}$, one has 
\begin{equation}\label{Equ: mu min asy Zj}\operatorname{\widehat{\mu}_{\min}^{\mathrm{asy}}}(\overline L_i|_{Z_j})\geqslant \operatorname{\widehat{\mu}_{\min}^{\mathrm{asy}}}(\overline L_i).\end{equation}
By \cite[Proposition 4.4.4]{CMIntersection}, one has 
\[\begin{split}(\overline L_0\cdots\overline L_d)_S&=\sum_{j=1}^n a_j(\overline L_1|_{Z_j}\cdots \overline L_d|_{Z_j})_S\\
&\quad-\int_{\Omega}\int_{X_\omega^{\mathrm{an}}}\ln |s|_{\varphi_{0,\omega}}(x)\,c_1(L_{1,\omega},\varphi_{1,\omega})\cdots c_1(L_{d,\omega},\varphi_{d,\omega})(\mathrm{d}x)\,\nu(\mathrm{d}\omega).
\end{split}\]
By Proposition \ref{Pro: minormation int prod}, one has 
\[\sum_{j=1}^na_j(\overline L_1|_{Z_j}\cdots\overline L_d|_{Z_j})_S\geqslant\sum_{j=1}^na_j\sum_{i=1}^d\delta_{i,j}\operatorname{\widehat{\mu}_{\min}^{\mathrm{asy}}}(\overline L_i|_{Z_j})\geqslant\sum_{j=1}^na_j\sum_{i=1}^d\delta_{i,j}\operatorname{\widehat{\mu}_{\min}^{\mathrm{asy}}}(\overline L_i),\]
where 
\[\delta_{i,j}:=(L_1|_{Z_j}\cdots L_{i-1}|_{Z_j}L_{i+1}|_{Z_j}\cdots L_d|_{Z_j}),\] 
and the second inequality comes from \eqref{Equ: mu min asy Zj}. Note that, for any $i\in\{1,\ldots,d\}$, one has 
\[\sum_{j=1}^na_j\delta_{i,j}=\delta_i.\]
Hence we obtain the desired inequality.
\end{proof}

\subsection{Asymptotic minimal slope of a quotient adelic line bundle}
\begin{prop}\label{Pro: minormation tautological}
Let $(E,\xi)$ be an adelic vector bundle on $S$, $L$ be a quotient line bundle of $f^*(E)$ and $\varphi$ be the quotient metric family induced by $\xi$. Then the adelic line bundle $(L,\varphi)$ is relatively nef. Moreover, the following inequality holds:
\begin{equation}\label{Equ: lower bound mu asy min}
\operatorname{\widehat{\mu}_{\min}^{\mathrm{asy}}}(\overline L)\geqslant\widehat{\mu}_{\min}(\overline E)-\frac 32\nu(\Omega_\infty)\ln(\dim_K(E))
\end{equation}
\end{prop}
\begin{proof}
By \cite[Propositions 6.1.8 and 6.1.2]{MR2095472}, $f^*(E)$ is a nef vector bundle on $X$ and hence $L$ is a nef line bundle. Moreover, since quotient metrics are semi-positive (see \cite[Remark 2.3.1]{CMArakelovAdelic}), the adelic line bundle $\overline L$ is relatively nef. 

In the following, we prove the inequality \eqref{Equ: lower bound mu asy min}. Let $p$ be an integer and $\overline A$ be a relatively ample adelic line bundle on $X$. Then $\overline L^{\otimes p}\otimes\overline A$ is relatively ample. Let $Y=\mathbb P(f^*(E)^{\otimes p})$ and $g:Y\rightarrow X$ be the structural morphism. The quotient homomorphism $f^*(E)\rightarrow L$ induces by taking the tensor product a surjective homomorphism
$f^*(E)^{\otimes p}\rightarrow L^{\otimes p}$,
which corresponds to a section $s:X\rightarrow Y$ such that $s^*(\mathcal O_Y(1))\cong L^{\otimes p}$. Hence 
\[s^*(\mathcal O_{Y}(1)\otimes g^*(A))\cong L^{\otimes p}\otimes A.\] By Proposition \ref{Pro: mu min of restriction}, one has 
\begin{equation}\label{Equ: minoration mu min asy}\widehat{\mu}_{\min}^{\mathrm{asy}}(\overline L^{\otimes p}\otimes\overline A)\geqslant\widehat{\mu}_{\min}^{\mathrm{asy}}(\overline{\mathcal O_{Y}(1)}\otimes g^*(\overline A)),\end{equation}
where we consider Fubini-Study metric fiber by fiber on $\mathcal O_Y(1)$. Note that, for any integer $n\in\mathbb N_{\geqslant 1}$, by the adjunction formula one has
\[H^0(Y,\mathcal O_Y(n)\otimes g^*(A)^{\otimes n})=H^0(X,S^n(f^*(E)^{\otimes p})\otimes A^{\otimes n})=S^n(E^{\otimes p})\otimes H^0(X,A^{\otimes n}).\]
Moreover, the projection map
\[ E^{\otimes np}\otimes H^0(X, A^{\otimes n})\longrightarrow S^n(E^{\otimes p})\otimes H^0(X,A^{\otimes n}),\]
where we consider the $\varepsilon,\pi$-tensor product norm family on the left hand side of the arrow, and the adelic vector bundle structure of $(fg)_*(\overline{\mathcal O_Y(n)}\otimes g^*(\overline A)^{\otimes n})$ on the right hand side. By \cite[Corollary 5.6.2]{CMArakelovAdelic} (see also Remark \ref{Rem: positive characteristic minimal slope}), we obtain 
\[\begin{split}\frac{\widehat{\mu}_{\min}\big(g_*(\overline{\mathcal O_{Y}(n)}\otimes g^*(\overline A^{\otimes n})\big)}{n}&\geqslant \frac{1}{n}\bigg(np\operatorname{\widehat{\mu}_{\min}}(\overline E)+\operatorname{\widehat{\mu}_{\min}}(f_*(\overline A^{\otimes n}))\\
&\quad-\frac 32\nu(\Omega_\infty)\ln\big(\dim_K(E)^{np}\cdot\dim_K(H^0(X,A^{\otimes n}))\big)\bigg).
\end{split}\]
Taking the limit when $n\rightarrow+\infty$, we obtain
\[\widehat{\mu}_{\min}^{\mathrm{asy}}(\overline{\mathcal O}_Y(1)\otimes g^*(\overline A))\geqslant p\operatorname{\widehat{\mu}_{\min}}(\overline E)+\widehat{\mu}_{\min}^{\mathrm{asy}}(\overline A)-\frac 32\nu(\Omega_\infty)p\ln(\dim_K(E)).\]
Combining this inequality with \eqref{Equ: minoration mu min asy}, we obtain 
\[\frac 1p\operatorname{\widehat{\mu}_{\min}^{\mathrm{asy}}}(\overline L^{\otimes p}\otimes\overline A)\geqslant\widehat{\mu}_{\min}(\overline E)+\frac 1p\operatorname{\widehat{\mu}_{\min}^{\mathrm{asy}}}(\overline A)-\frac 32\nu(\Omega_\infty)\ln(\dim_K(E)). \]
Thus, due to Definition~\ref{def:relatively:nef}, we obtain
\[\operatorname{\widehat{\mu}_{\min}^{\mathrm{asy}}}(\overline L)\geqslant \limsup_{p\to\infty} \frac 1p\operatorname{\widehat{\mu}_{\min}^{\mathrm{asy}}}(\overline L^{\otimes p}\otimes\overline A) \geqslant
\widehat{\mu}_{\min}(\overline E)-\frac 32\nu(\Omega_\infty)\ln(\dim_K(E)),\]
as required. 
\end{proof}

\subsection{Asymptotic minimal slope of a pull-back}
\begin{prop}\label{Pro: pull back higher mu min asy}
Let $g:Y\rightarrow X$ be a projective morphism of $K$-schemes, which is surjective and such that $g_*(\mathcal O_Y)=\mathcal O_X$. Let $\overline L$ be an adelic line bundle on $X$ such that $L$ is nef. Then the following inequality holds:
\[\widehat{\mu}_{\min}^{\mathrm{asy}}(g^*(\overline L))\geqslant \widehat{\mu}_{\min}^{\mathrm{asy}}(\overline L).\]
\end{prop}
\begin{proof} By \cite[Example 1.4.4]{MR2095471}, the line bundle $g^*(L)$ is nef, and hence $\operatorname{\widehat{\mu}_{\min}^{\mathrm{asy}}}(g^*(\overline L))$ is well defined. We first consider the case where $L$ is ample.
Let $p$ be a positive integer and $\overline A$ be an adelic line bundle on $Y$ such that $A$ is ample.  By Lemma \ref{Lem: surjectivity tensor product}, for sufficiently positive integer $n$, the $K$-linear map 
\[H^0(Y,A^{\otimes n})\otimes H^0(X,L^{\otimes pn})\longrightarrow H^0(Y,A^{\otimes n}\otimes g^*(L)^{\otimes pn})\] 
is surjective. Moreover, if we equip the left hand side of the arrow with the $\varepsilon,\pi$-tensor product norm family of those of $(fg)_*(\overline A^{\otimes n})$ and $f_*(\overline L^{\otimes pn})$, then the $K$-linear map has height $\leqslant 0$. Therefore, by \cite[Corollary 5.6.2]{CMArakelovAdelic} we obtain   
\[\begin{split}\widehat{\mu}_{\min}((fg)_*(\overline A^{\otimes n}\otimes& g^*(\overline L)^{\otimes pn}))\geqslant \widehat{\mu}_{\min}((fg)_*({\overline A}^{\otimes n}))+\widehat{\mu}_{\min}(f_*(\overline L^{\otimes pn}))\\
&-\frac 32\nu(\Omega_\infty)\ln(\dim_K(H^0(Y,A^{\otimes n}))\dim_K(H^0(X,L^{\otimes pn}))). \end{split}\]
Dividing the two sides of the inequality by $pn$ and then taking the limit when $n\rightarrow+\infty$, we obtain 
\[\frac 1p\operatorname{\widehat{\mu}_{\min}^{\mathrm{asy}}}(\overline A\otimes g^*(\overline L)^{\otimes p})\geqslant\frac{\widehat{\mu}_{\min}^{\mathrm{asy}}(\overline A)}{p}+\widehat{\mu}_{\min}^{\mathrm{asy}}(\overline L),\]
which leads to 
\[\widehat{\mu}^{\mathrm{asy}}_{\min}(g^*(\overline L))\geqslant \limsup_{p\rightarrow+\infty}\frac 1p\operatorname{\widehat{\mu}_{\min}^{\mathrm{asy}}}(\overline A\otimes g^*(\overline L)^{\otimes p})\geqslant\widehat{\mu}_{\min}^{\mathrm{asy}}(\overline L).\]

We now consider the general case. Let $\overline B$ be an adelic line bundle on $X$ such that $B$ is ample. By the above argument we obtain that $\operatorname{\widehat{\mu}_{\min}^{\mathrm{asy}}}(g^*(\overline B))\geqslant\operatorname{\widehat{\mu}_{\min}^{\mathrm{asy}}}(\overline B)>-\infty$ and, for any positive integer $n$,
\[\frac{1}{n}\operatorname{\widehat{\mu}_{\min}^{\mathrm{asy}}}(g^*(\overline L^{\otimes n})\otimes g^*(\overline B))\geqslant \frac 1n \operatorname{\widehat{\mu}_{\min}^{\mathrm{asy}}}(\overline L^{\otimes n}\otimes\overline B).\]
Taking the limit when $n\rightarrow+\infty$, by \eqref{Equ: equality limit if mu min not infinite} we obtain $\widehat{\mu}_{\min}^{\mathrm{asy}}(g^*(\overline L))\geqslant \widehat{\mu}_{\min}^{\mathrm{asy}}(\overline L)$. 
\end{proof}

\begin{theo}\label{Thm: mu min asy and pull-back}
Let $g:Y\rightarrow X$ be a projective morphism of $K$-schemes. We assume that $Y$ is non-empty and reduced. For any adelic line bundle $\overline L$ on $X$ such that $L$ is nef, one has $\operatorname{\widehat{\mu}_{\min}^{\mathrm{asy}}}(g^*(\overline L))\geqslant\operatorname{\widehat{\mu}_{\min}^{\mathrm{asy}}}(\overline L) $.
\end{theo}
\begin{proof}
The projective morphism $g$ can be written as the composition of a closed immersion from $Y$ into a projective bundle on $X$ and the projection from the projective bundle to $X$. Hence the inequality follows from Propositions \ref{Pro: pull back higher mu min asy} and  \ref{Pro: mu min of restriction}.
\end{proof}

\subsection{Comparison with the normalized height}

The following height estimate can be deduced from Theorem \ref{Pro: lower bound ample}. Here we provide an alternative proof in the particular case where $X$ is  integral by using the arithmetic Hilbert-Samuel formula.
\begin{prop}\label{Pro: inequality height and mu min asy }
Let $\overline L$ be a relatively nef adelic line bundle on $X$ such that $(L^d)>0$. Then the following inequality holds
\begin{equation}\label{Equ: height and mu min asy}\frac{(\overline L^{d+1})_S}{(d+1)(L^d)}\geqslant\widehat{\mu}_{\min}^{\mathrm{asy}}(\overline L).\end{equation}
\end{prop}
\begin{proof} We assume that $X$ is integral.
In the case where $\overline{L}$ is relatively ample, it is a consequence of Theorem~\ref{thm:lower:bound} and Remark~\ref{rema:mu:asy:geq:mu:asy:min}.

We now consider the general case where $\overline L$ is only relatively nef. Let $\overline A$ be a relatively ample adelic line bundle and $N$ be a positive integer such that $L^{\otimes n}\otimes\overline A$ is relatively ample for any $n\in\mathbb N_{\geqslant N}$. For any $n\in\mathbb N_{\geqslant N}$, the adelic line bundle $\overline L_n=\overline L^{\otimes n}\otimes\overline A$ is relatively ample. Hence the particular case of the proposition proved above shows that 
\[\forall\,n\in\mathbb N_{\geqslant N},\quad \frac{(\overline L_n^{d+1})_S}{(d+1)(\overline L_n^d)}\geqslant\widehat{\mu}_{\min}^{\mathrm{asy}}(\overline L_n).\]
Moreover, by the multi-linearity of intersection product, one has
\[\lim_{n\rightarrow+\infty}\frac{(\overline L_n^{d+1})_S}{n^{d+1}}=(\overline L^{d+1})_S,\quad \lim_{n\rightarrow+\infty}\frac{(L_n^d)}{n^d}=(L^d).\]
Therefore, one obtains
\[\widehat{\mu}_{\min}^{\mathrm{asy}}(\overline L)=\lim_{n\rightarrow+\infty}\frac 1n\operatorname{\widehat{\mu}_{\min}^{\mathrm{asy}}}(\overline L_n)\leqslant \lim_{n\rightarrow+\infty}\frac{(\overline L_n^{d+1})_S}{n(d+1)(\overline L_n^d)}=\frac{(\overline L^{d+1})}{(d+1)(L^d)}.\]
\end{proof}

\begin{coro}\label{Cor: mu min bounded by height}
Let $\overline L$ be a relatively nef adelic line bundle on $X$. For any non-empty and reduced closed subscheme $Y$ of $X$, the following inequality holds:
\begin{equation}\label{Equ: bound of mu min by height}\widehat{\mu}_{\min}^{\mathrm{asy}}(\overline L)\leqslant \frac{h_{\overline L}(Y)}{(\dim_L(Y)+1)\deg_L(Y)}. \end{equation}
In particular, for any closed point $x$, one has 
\[\widehat{\mu}_{\min}^{\mathrm{asy}}(\overline L)\leqslant h_{\overline L}(x).\]
\end{coro}
\begin{proof}
By Lemma \ref{lemma:restriction:relatively:nef}, the restriction of $\overline L$ to $Y$ is relatively nef. By Proposition \ref{Pro: inequality height and mu min asy } one has 
\[\widehat{\mu}_{\min}^{\mathrm{asy}}(\overline L|_Y)\leqslant \frac{h_{\overline L}(Y)}{(\dim_L(Y)+1)\deg_L(Y)}.\] 
By Proposition \ref{Pro: mu min of restriction}, we obtain \eqref{Equ: bound of mu min by height}. 
\end{proof}

\section{Geometrically big and pseudoeffective adelic line bundles}

As in the previous section, we fix a proper adelic curve $S=(K,(\Omega,\mathcal A,\nu),\phi)$ such that, either $K$ is countable, or $(\Omega,\mathcal A)$ is discrete. We assume that $K$ is perfect. 

\subsection{Convergence of maximal slopes}
 
\begin{prop}\label{Pro: super additivity of mu max}
Let $X$ be an integral projective scheme over $\Spec K$, and $\overline L=(L,\varphi)$ and $\overline M=(M,\psi)$ be adelic line bundles on $X$ such that $H^0(X,L)$ and $H^0(X,M)$ are non-zero. Then the following inequality holds:
\[\widehat{\mu}_{\max}(f_*(\overline L\otimes\overline M))\geqslant\widehat{\mu}_{\max}(f_*(\overline L))+\widehat{\mu}_{\max}(f_*(\overline M))-\frac 32\nu(\Omega_\infty)(\ln(h^0(L)\cdot h^0(M))),\] 
where $h^0(L)=\dim_K(H^0(X,L))$ and $h^0(M)=\dim_K(H^0(X,M))$.
\end{prop}
\begin{proof}
By \cite[Theorem 4.3.58]{CMArakelovAdelic}, there exist non-zero vector subspaces $E$ and $F$ of $H^0(X,L)$ and $H^0(X,M)$, respectively, such that
\begin{equation*}
\widehat{\mu}_{\min}(\overline E)=\widehat{\mu}_{\max}(f_*(\overline L)),\quad \widehat{\mu}_{\min}(\overline F)=\widehat{\mu}_{\max}(f_*(\overline M)),
\end{equation*}
where we consider restricted norm families on $E$ and $F$. Since $X$ is integral, the map 
\[E\otimes_K F\longrightarrow H^0(X,L\otimes M),\quad s\otimes t\longmapsto st\]
is non-zero. Moreover, for any $\omega\in\Omega$, one has
\[\forall\,(s,t)\in E_\omega\times F_\omega,\quad
\|st\|_{\varphi_\omega+\psi_\omega}\leqslant \|s\|_{\varphi_\omega}\cdot\|t\|_{\psi_\omega}.\]
Therefore, the height of the above $K$-linear map is $\leqslant 0$ if we consider the $\varepsilon,\pi$-tensor product norm family on $E\otimes_KF$. By \cite[Theorem 4.3.31 and Corollary 5.6.2]{CMArakelovAdelic}, we obtain 
\[\begin{split}&\quad\;\widehat{\mu}_{\max}(f_*(\overline L\otimes\overline M))\\
&\geqslant\widehat{\mu}_{\min}(\overline E\otimes_{\varepsilon,\pi}\overline F)\geqslant\widehat{\mu}_{\min}(\overline E)+\widehat{\mu}_{\min}(\overline F)-\frac 32\nu(\Omega_\infty)\ln(\dim_K(E)\cdot\dim_K(F))\\
&=\widehat{\mu}_{\max}(f_*(\overline L))+\widehat{\mu}_{\max}(f_*(\overline M))-\frac 32\nu(\Omega_\infty)\ln(\dim_K(E)\cdot\dim_K(F))\\
&\geqslant\widehat{\mu}_{\max}(f_*(\overline L))+\widehat{\mu}_{\max}(f_*(\overline M))-\frac 32\nu(\Omega_\infty)(\ln(h^0(L)\cdot h^0(M))),
\end{split}\] 
as required.   
\end{proof}

\begin{coro}\label{Cor: convergence maximal slopes}
Let $\overline L$ be an adelic line bundle on $X$ such that $H^0(X,L^{\otimes n})$ is non-zero for sufficiently large natural number $n$. The sequence 
\[\frac 1n\operatorname{\widehat{\mu}_{\max}}(f_*(\overline L^{\otimes n})),\quad n\in\mathbb N_{\geqslant 1}\]
converges in $\mathbb R$.
\end{coro}
\begin{proof}
The convergence of the sequence follows from Proposition \ref{Pro: super additivity of mu max}, using the same argument as in the proof of Proposition \ref{Pro: convergence mu min }.
\end{proof}

\subsection{Asymptotic maximal slope}
\label{Sec: asymptotic maximal slope}
In this subsection, we let $f:X\rightarrow\Spec K$ be an integral projective $K$-scheme.

\begin{defi}
Let $\overline L$ be an adelic line bundle on $X$ such that $L$ is big. We define
\[\widehat{\mu}_{\max}^{\mathrm{asy}}(\overline L):=\lim_{n\rightarrow+\infty}\frac{\widehat{\mu}_{\max}(f_*(\overline L{}^{\otimes n}))}{n}.\]
By definition, for any $p\in\mathbb N_{\geqslant 1}$, the following equality holds:
\[\widehat{\mu}_{\max}^{\mathrm{asy}}(\overline L^{\otimes p})=p\operatorname{\widehat{\mu}_{\max}^{\mathrm{asy}}}(\overline L).\]
\end{defi}

\begin{prop}\label{Pro: superadditivty}
Let $\overline L$ and $\overline M$ be adelic line bundles on $X$ such that $L$ and $M$ are both big. One has 
\begin{equation}\label{Equ: super additive mu max asy}\widehat{\mu}_{\max}^{\mathrm{asy}}(\overline L\otimes\overline M)\geqslant\widehat{\mu}_{\max}^{\mathrm{asy}}(\overline L)+\widehat{\mu}_{\max}^{\mathrm{asy}}(\overline M).\end{equation} 
\end{prop}
\begin{proof}
For any $n\in\mathbb N_{\geqslant 1}$, let $a_n=\dim_K(H^0(X,L^{\otimes n}))$ and $b_n=\dim_K(X,M^{\otimes n})$. One has 
\[\ln(a_n)=O(\ln (n)),\quad \ln(b_n)=O(\ln(n)),\quad n\rightarrow+\infty.\]
By Proposition \ref{Pro: super additivity of mu max}, for sufficiently large $n$, one has 
\[\frac{\widehat{\mu}_{\max}(f_*((\overline L\otimes\overline M)^{\otimes n}))}{n}\geqslant\frac{\widehat{\mu}_{\max}(f_*(\overline L^{\otimes n}))}{n}+\frac{\widehat{\mu}_{\max}(f_*(\overline M^{\otimes n}))}{n}-\frac{3}{2}\nu(\Omega_\infty)\frac{\ln(a_nb_n)}{n}.\]
Taking the limit when $n\rightarrow+\infty$, we obtain the inequality \eqref{Equ: super additive mu max asy}. 
\end{proof}

\begin{prop}\label{Pro: sequence convergence mu max asy}
Let $\overline L$ and $\overline A$ be adelic line bundle on $X$. We assume that $L$ is pseudo-effective and $A$ is big. Then the sequence
\[\frac 1n\operatorname{\widehat{\mu}_{\max}^{\mathrm{asy}}}(\overline L^{\otimes n}\otimes \overline A),\quad n\in\mathbb N_{\geqslant 1}\]
converges in $\mathbb R\cup\{-\infty\}$. Moreover, its limit does not depend on the choice of $\overline A$. In particular, in the case where $L$ is big, the following equality  holds:
\begin{equation}
\lim_{n\rightarrow+\infty}\frac 1n\operatorname{\widehat{\mu}_{\max}^{\mathrm{asy}}}(\overline L^{\otimes n}\otimes \overline A)=\operatorname{\widehat{\mu}_{\max}^{\mathrm{asy}}}(\overline L).
\end{equation}
\end{prop}
\begin{proof}
The proof relies on the super-additivity of the function $\widehat{\mu}_{\max}^{\mathrm{asy}}(\ndot)$ (see Proposition \ref{Pro: superadditivty}) and follows the same strategy as that of Proposition \ref{Pro: convergence and independence}. We omit the details.
\end{proof}

\begin{defi}
Let $\overline L$ be an adelic line bundle on $X$ such that $L$ is pseudo-effective. We define $\widehat{\mu}_{\max}^{\mathrm{asy}}(\overline L)$ as the limite
\[\lim_{n\rightarrow+\infty}\frac{1}{n}\operatorname{\widehat{\mu}_{\max}^{\mathrm{asy}}}(\overline L^{\otimes n}\otimes\overline A),\]
where $\overline A$ is an arbitrary adelic line bundle on $X$ such that $A$ is big. The element $\widehat{\mu}_{\max}^{\mathrm{asy}}(\overline L)$ of $\RR \cup \{ -\infty \}$ is called the \emph{asymptotic maximal slope}\index{slope!asymptotic maximal ---} of $\overline L$.
\end{defi}

\begin{prop}\label{Pro: suradiditivity mu pi max}
Let $\overline L$ and $\overline M$ be adelic line bundles on $X$ such that $L$ and $M$ are pseudo-effective. Then the following inequality holds:
\[\widehat{\mu}_{\max}^{\mathrm{asy}}(\overline L\otimes\overline M)\geqslant\widehat{\mu}_{\max}^{\mathrm{asy}}(\overline L)+\widehat{\mu}_{\max}^{\mathrm{asy}}(\overline M).\]
\end{prop}
\begin{proof}
Let $\overline A$ be an adelic line bundle on $X$ such that $A$ is big. For any $n\in\mathbb N$, 
\[(L\otimes M)^{\otimes n}\otimes A^{\otimes 2}= (L^{\otimes n}\otimes A)\otimes(M^{\otimes n}\otimes A)\]
is big. Moreover, by Proposition \ref{Pro: superadditivty}, one has 
\[\frac{1}{n}\widehat{\mu}_{\max}^{\mathrm{asy}}\big((\overline L\otimes M)^{\otimes n}\otimes\overline A^{\otimes 2}\big)\geqslant\frac{1}{n}\widehat{\mu}_{\max}^{\mathrm{asy}}(\overline L^{\otimes n}\otimes \overline A)+\frac{1}{n}\widehat{\mu}_{\max}^{\mathrm{asy}}(\overline M^{\otimes n}\otimes \overline A).\]
Taking the limit when $n\rightarrow+\infty$, we obtain \[\widehat{\mu}_{\max}^{\mathrm{asy}}(\overline L\otimes\overline M)\geqslant\widehat{\mu}_{\max}^{\mathrm{asy}}(\overline L)+\widehat{\mu}_{\max}^{\mathrm{asy}}(\overline M).\]
\end{proof}

\subsection{Pullback by a surjective projective morphism}

Let $X$ and $Y$ be integral projective $K$-schemes and $g:Y\rightarrow X$ be a surjective projective morphism.

\begin{lemm}\label{Lem: pullback of pseudo effective line bundle}
Let $L$ be an invertible $\mathcal O_X$-module. If $L$ is pseudo-effective, then the pullback $g^*(L)$ is also pseudo-effective.
\end{lemm}
\begin{proof}
Let $A$ be a big invertible $\mathcal O_X$-module and $B$ be a big invertible $\mathcal O_Y$-module. For any positive integer $p$, the invertible $\mathcal O_X$-module $L^{\otimes p}\otimes A$ is big and hence $g^*(L^{\otimes p}\otimes A)$ is pseudo-effective since it has a tensor power which is effective. Similarly, $g^*(A)$ is also pseudo-effective. Thus we obtain that $g^*(A)\otimes B$ and $g^*(L)^{\otimes p}\otimes g^*(A)\otimes B$ are big. In particular, $g^*(L)$ is pseudo-effective.
\end{proof}

\begin{prop}\label{Pro: pullback by surjective morphism}
Let $\overline L$ be an adelic line bundle on $X$ such that $L$ is pseudo-effective. Then  the following inequality holds:
\[\widehat{\mu}_{\max}^{\mathrm{asy}}(g^*(\overline L))\geqslant\widehat{\mu}_{\max}^{\mathrm{asy}}(\overline L).\]
\end{prop}
\begin{proof}
We have seen in Lemma \ref{Lem: pullback of pseudo effective line bundle} that the invertible $\mathcal O_X$-module $L$ is pseudo-effective, so that $\widehat{\mu}_{\max}^{\mathrm{asy}}(\overline L)$ is well defined. We choose an  adelic line bundle $\overline A$ on $X$  such that $A$ is big. 

We first assume that $L$ is big. Let $n$ and $p$ be positive integers. We consider the $K$-linear map
\begin{multline*}H^0(X,g_*(A^{\otimes n}))\otimes H^0(X,L^{\otimes np})\longrightarrow H^0(X,g_*(A^{\otimes n})\otimes L^{\otimes np}) \\ =H^0(Y,A^{\otimes n}\otimes g^*(L^{\otimes np}))\end{multline*}
induced by multiplication of sections. Let $E$ be the destabilizing vector subspace of $(fg)_*(\overline A^{\otimes n})$ and let $F$ be the destabilizing vector subspace of $f_*(\overline L^{\otimes np})$. By \cite[Proposition 4.3.31]{CMArakelovAdelic}, one has 
\[\widehat{\mu}_{\min}(\overline E\otimes_{\varepsilon,\pi}\overline F)\leqslant\widehat{\mu}_{\max}((fg)_*(\overline A^{\otimes n}\otimes g^*(\overline L^{\otimes np}))).\]
by \cite[Corollary 5.6.2 and Remark 4.3.48]{CMArakelovAdelic} (see also Remark \ref{Rem: positive characteristic minimal slope}), one deduces
\[\begin{split}&\quad\;\widehat{\mu}_{\max}((fg)_*(\overline A^{\otimes n}\otimes g^*(\overline L^{\otimes np})))\\
&\geqslant\widehat{\mu}_{\min}(\overline E)+\widehat{\mu}_{\min}(\overline F)-\frac 32\nu(\Omega_\infty)(\ln(\dim_K(E))+\ln(\dim_K(F)))\\
&\geqslant\widehat{\mu}_{\max}((fg)_*(\overline A^{\otimes n}))+\widehat{\mu}_{\max}(f_*(\overline L^{\otimes np}))-2\nu(\Omega_\infty)\ln\dim_K(H^0(Y,A^{\otimes n}))\\
&\qquad -2\nu(\Omega_\infty)\ln\dim_K(H^0(X,L^{\otimes{np}})).
\end{split}\]
If we divide the two sides by $np$, taking the limit when $n\rightarrow+\infty$, we obtain 
\[\frac 1p\widehat{\mu}_{\max}^{\mathrm{asy}}(\overline A\otimes g^*(\overline L^{\otimes p}))\geqslant\frac 1p\widehat{\mu}_{\max}^{\mathrm{asy}}(\overline A)+\operatorname{\widehat{\mu}}_{\max}(\overline L).\]
Taking the limit when $p
\rightarrow+\infty$, we obtain $\widehat{\mu}_{\max}^{\mathrm{asy}}(g^*(\overline L))\geqslant\widehat{\mu}_{\max}^{\mathrm{asy}}(\overline L)$, as required.

We then consider the general case where $L$ is only assumed to be pseudo-effective. Let $\overline B$ be an adelic line bundle on $X$ such that $B$ is big. Note that, for any positive integer $p$, $L^{\otimes p}\otimes B$ is big. Hence, by the particular case of the proposition shown above, one has 
\[\widehat{\mu}_{\max}^{\mathrm{asy}}(g^*(\overline L)^{\otimes p}\otimes g^*(\overline B))\geqslant \operatorname{\widehat{\mu}_{\max}^{\mathrm{asy}}}(\overline L^{\otimes p}\otimes B).\]
Therefore, by Proposition \ref{Pro: suradiditivity mu pi max}, we obtain
\[\frac 1p\widehat{\mu}_{\max}^{\mathrm{asy}}(g^*(\overline L)^{\otimes p}\otimes g^*(\overline B)\otimes\overline A)\geqslant\frac 1p\widehat{\mu}_{\max}^{\mathrm{asy}}(\overline L^{\otimes p}\otimes\overline B)+\frac 1p\widehat{\mu}_{\max}^{\mathrm{asy}}(\overline A).\]
Taking the limit when $p\rightarrow+\infty$, we obtain $\widehat{\mu}_{\max}^{\mathrm{asy}}(g^*(\overline L))\geqslant\widehat{\mu}_{\max}^{\mathrm{asy}}(\overline L)$.
\end{proof}

\begin{rema}
Let $\overline L$ be an adelic line bundle on $X$. Assume that $L$ is the pull-back of a big line bundle by a surjective projective morphism. Then Proposition \ref{Pro: pullback by surjective morphism} shows that $\widehat{\mu}_{\max}^{\mathrm{asy}}(\overline L)\in\mathbb R$. 
\end{rema}

\subsection{Relative Fujita approximation}

Let $f:X\rightarrow\Spec K$ be a projective $K$-scheme, $K(X)$ be the field of rational functions on $X$, and $\mathscr M_X$ be the sheaf of meromorphic function on $X$.

\begin{defi}
 Let $L$ be a big line bundle on $X$. Note that $\mathscr M_X\otimes_{\mathcal O_X}L$ is isomorphic to the trivial invertible $\mathscr M_X$-module. In particular, if $s$ and $t$ are two global sections of $L$ such that $s\neq 0$, then there existe a unique rational function $\lambda\in K(X)$ such that $t=\lambda s$. We denote by $t/s$ this rational function.   If $E$ is a $K$-vector subspace of $H^0(X,L)$. We denote by $K(E)$ the sub-extension of $K(X)/K$ generated by elements of the form
$t/s$, where $t$ and $s$ are non-zero sections in $K(E)$. We say that $E$ is \emph{birational}\index{birational} if $K(E)=K(X)$. Moreover $L$ is said to be \emph{birational} if $K(H^0(X,L))=K(X)$.
\end{defi}

\begin{rema}\label{Rem: birational linear series}
Let $L$ and $M$ be line bundle on $X$, $E$ be a vector subspace of $H^0(X,L)$, $s$ be a non-zero global section of $M$ and \[F=\{ts\,|\,t\in E\}\subseteq H^0(X,L\otimes M).\] Then by definition one has $K(F)=K(E)$. In particular, if $E$ is birational,  so is $F$;  if $L$ is birational, so is $L\otimes M$.
\end{rema}

\begin{prop}
Let $L$ be a big line bundle on $X$. For sufficiently positive integer $p$, the line bundle $L^{\otimes p}$ is birational.
\end{prop}
\begin{proof}
Since $L$ is big, there exists a positive integer $q$, an ample line bundle $A$ and an effective line bundle $M$ on $X$ such that $L^{\otimes q}\cong A\otimes M$. By replacing $q$ by a multiple, we may assume that the graded $K$-algebra
\[\bigoplus_{n\in\mathbb N}H^0(X,A^{\otimes n})\]
is generated by $H^0(X,A)$ and that $L^{\otimes(q+1)}$ is effective. For any $a\in\mathbb N_{\geqslant 1}$, one has  
\[X=\operatorname{Proj}\Big(\bigoplus_{n\in\mathbb N}H^0(X,A^{\otimes an})\Big),\]
which implies that $A^{\otimes a}$ is birational and hence $L^{\otimes aq}$ is birational. Moreover, since $L^{\otimes (q+1)}$ is effective, for any $b\in\mathbb N_{\geqslant 1}$, the line bundle $L^{\otimes b(q+1)}$ is also effective. Therefore, for any $(a,b)\in\mathbb N_{\geqslant 1}^2$, the line bundle $L^{\otimes aq+b(q+1)}$ is birational. Since $q$ and $q+1$ are coprime, we obtain that $L^{\otimes p}$ is birational for sufficiently large $p\in\mathbb N$.
\end{proof}

\begin{defi}
Let $\overline L=(L,\varphi)$ be an adelic line bundle on $X$. If $s\in H^0(X,L)$ is a non-zero global section such that $\|s\|_{\varphi_\omega}\leqslant 1$ for any $\omega\in\Omega$, we say that the global section $s$ is \emph{effective}\index{effective}. We say that $\overline L$ is \emph{effective}\index{effective} if it admits at least an effective global section.
\end{defi}

\begin{lemm}
Let $\overline L$ be an adelic line bundle such that $L$ is big. For any $t<\widehat{\mu}_{\max}^{\mathrm{asy}}(\overline L)$ and any $N\in\mathbb N_{\geqslant 1}$, there exists an integer $p\geqslant N$ and a vector subspace $E$ of $H^0(X,L^{\otimes p})$ such that $K(E)=K(X)$ and 
$\widehat{\mu}_{\min}(\overline E)> pt$.
\end{lemm}
\begin{proof}
By replacing $\overline L$ by one of its tensor powers, we may assume without loss of generality that $L$ is birational. For any $n\in\mathbb N$, let $r_n=\dim_K(H^0(X,L^{\otimes n}))$. Since $t<\operatorname{\widehat{\mu}_{\max}^{\mathrm{asy}}}(\overline L)$, for sufficiently large $n\in\mathbb N$, one has 
\[\widehat{\mu}_{\max}(f_*(\overline L^{\otimes n}))> (n+1)t-\widehat{\mu}_{\min}(f_*(\overline L))+\frac 32\nu(\Omega_\infty)\ln\big(r_n\cdot r_1\big).\]
Let $F$ be a vector subspace of $H^0(X,L^{\otimes n})$ such that 
\[\widehat{\mu}_{\min}(\overline F)=\widehat{\mu}_{\max}(f_*(\overline L^{\otimes n})).\]
The existence of $F$ is ensured by \cite[Theorem 4.3.58]{CMArakelovAdelic}. Let $E$ be the image of $F\otimes_KH^0(X,L)$ by the $K$-linear map
\[H^0(X,L^{\otimes n})\otimes H^0(X,L)\longrightarrow H^0(X,L^{\otimes n+1}),\quad s\otimes t\longmapsto st.\]
Since $L$ is birational and $F$ is non-zero, we obtain that $E$ is birational.
By \cite[Corollary 5.6.2]{CMArakelovAdelic}, one has 
\[\begin{split}&\quad\;\widehat{\mu}_{\min}(\overline E)\geqslant\widehat{\mu}_{\min}(\overline F\otimes_{\varepsilon,\pi}f_*(\overline L))\\
&\geqslant\widehat{\mu}_{\min}(\overline F)+\widehat{\mu}_{\min}(f_*(\overline L))-\frac 32\nu(\Omega_\infty)\ln(\dim_K(F)\cdot r_1)\\
&=\widehat{\mu}_{\max}(f_*(\overline L^{\otimes n}))+\widehat{\mu}_{\min}(f_*(\overline L))-\frac 32\nu(\Omega_\infty)\ln(\dim_K(F)\cdot r_1)>(n+1)t.
\end{split}\]
\end{proof}

\begin{theo}[Relative Fujita approximation]\label{Thm: relative fujita}
Let $\overline L$ be an adelic line bundle on $X$ such that $L$ is big. For any real number $t<\widehat{\mu}_{\max}^{\mathrm{asy}}(\overline L)$, there exist
a positive integer $p$, a birational projective $K$-morphism $g:X'\rightarrow X$, a relatively nef adelic line bundle $\overline A$ and an effective adelic line bundle $\overline M$ on $X'$ such that $A$ is big, $g^*(\overline L{}^{\otimes p})$ is isomorphic to $\overline A\otimes\overline M$ and  
$\widehat{\mu}_{\min}^{\mathrm{asy}}(\overline A)\geqslant pt$.
\end{theo}
\begin{proof}
We pick a positive integer $p$ and a birational vector subspace $V$ of $H^0(X,L^{\otimes p})$ such that 
\[\widehat{\mu}_{\min}(\overline V)=\widehat{\mu}_{\max}(f_*(\overline L^{\otimes p}))\geqslant pt+\frac 32\nu(\Omega_\infty)\ln(\dim_K(H^0(X,L^{\otimes p}))).
\] 

Let $g:X'\rightarrow X$ be the blow-up of $L$ along the base locus of $V$, namely
\[X'=\operatorname{Proj}\bigg(\operatorname{Im}\bigg(\bigoplus_{n\in\mathbb N}S^n(f^*(V))\longrightarrow\bigoplus_{n\in\mathbb N}L^{\otimes np}\bigg)\bigg).\]
Denote by $E$ the exceptional divisor and by $s_E$ the global section of $\mathcal O_X(E)$ which trivializes $\mathcal O_X(E)$ outside of the exceptional divisor. One has \[\mathcal O_{X'}(1)\cong g^*(L^{\otimes p})\otimes\mathcal O_X(-E). \]
Moreover, the canonical surjective homomorphism 
\begin{equation}\label{Equ: pull back}g^*(f^*(V))\longrightarrow\mathcal O_{X'}(1)\end{equation}
induces a $K$-morphism $i:X'\rightarrow \mathbb P(V)$ such that $i^*(\mathcal O_V(1))=\mathcal O_{X'}(1)$, where $\mathcal O_V(1)$ denotes the universal invertible sheaf on $\mathbb P(V)$. Since $V$ is birational, the line bundle $\mathcal O_{X'}(1)$ is big.

We equip $V$ with the induced norm family of $(\|\ndot\|_{p\varphi_\omega})_{\omega\in\Omega}$ and $\mathcal O_{X'}(1)$ with the quotient metric family $\varphi'=(\varphi_\omega')_{\omega\in\Omega}$ induced by $(\|\ndot\|_{p\varphi_\omega})_{\omega\in\Omega}$ and the surjective homomorphism \eqref{Equ: pull back}. We identify $\mathcal O_X(E)$ with $g^*(L^{\otimes p})\otimes\mathcal O_{X'}(1)^\vee$ and equip it with the tensor product metric family. Then the section $s_E$ is effective.   Moreover, by Proposition \ref{Pro: minormation tautological}, the adelic line bundle $\overline{\mathcal O_{X'}(1)}$ is relatively nef, and the following inequality holds
\[\widehat{\mu}_{\min}^{\mathrm{asy}}(\overline{\mathcal O_{X'}(1)})\geqslant\operatorname{\widehat{\mu}_{\min}}(\overline V)-\frac 32\nu(\Omega_\infty)\ln(\dim_K(V))\geqslant pt,\]
as required.

\end{proof}

\begin{rema}\label{Rem: strong relative Fujita}
 Let $\overline L$ be an adelic line bundle on $X$ such that $L$ is big. Let $\overline B$ be a relatively ample adelic line bundle. There exists a positive integer $N$ such that $L^{\otimes m}\otimes B^\vee$ is  big for any $m\in\mathbb N_{\geqslant N}$. Let $t$ be a real number such that $t<\widehat{\mu}_{\max}^{\mathrm{asy}}(\overline L)$. There exists $m\in\mathbb N_{\geqslant N}$ such that 
\[mt-\operatorname{\widehat{\mu}_{\min}^{\mathrm{asy}}}(\overline B)<(m-N)\operatorname{\widehat{\mu}_{\max}^{\mathrm{asy}}}(\overline L)+\widehat{\mu}_{\max}^{\mathrm{asy}}(\overline L^{\otimes N}\otimes\overline B^\vee)\leqslant\widehat{\mu}_{\max}^{\mathrm{asy}}(\overline L^{\otimes m}\otimes\overline B^\vee),\]
where the second inequality comes from 
If we apply Theorem \ref{Thm: relative fujita} to $\overline L^{\otimes m}\otimes\overline B^\vee$, we obtain the existence of a positive integer $p$, a birational projective $K$-morphism $g:X'\rightarrow X$, a relatively nef adelic line bundle $\overline A$ and an effective adelic line bundle $\overline M$ on $X'$ such that $A$ is big, $g^*(\overline L^{\otimes mp}\otimes\overline B^{\vee\otimes p})$ is isomorphic to $\overline A\otimes\overline M$ and  
\begin{equation}\label{Equ: lower bound mu min A}\widehat{\mu}_{\min}^{\mathrm{asy}}(\overline A)\geqslant p(mt-\operatorname{\widehat{\mu}_{\min}^{\mathrm{asy}}}(\overline B)).\end{equation}
Let $\overline N=\overline A\otimes g^*(\overline B)^{\otimes p}$. This is a relatively ample line bundle, and one has 
\[\overline N\otimes\overline M\cong \overline A\otimes\overline M\otimes g^*(\overline B)^{\otimes p}\cong g^*(\overline L^{\otimes mp}).\]
Moreover, one has 
\[\widehat{\mu}_{\min}^{\mathrm{asy}}(\overline N)\geqslant\widehat{\mu}_{\min}^{\mathrm{asy}}(\overline A)+p\operatorname{\widehat{\mu}_{\min}^{\mathrm{asy}}}(g^*(\overline B))\geqslant \widehat{\mu}_{\min}^{\mathrm{asy}}(\overline A)+p\operatorname{\widehat{\mu}_{\min}^{\mathrm{asy}}}(\overline B),\]
where the first inequality comes from Proposition \ref{Pro: super additivity nef}, and the second comes from Theorem \ref{Thm: mu min asy and pull-back}. By \eqref{Equ: lower bound mu min A}, we obtain
\[\widehat{\mu}_{\min}^{\mathrm{asy}}(\overline N)\geqslant pmt.\] 
Therefore, in Theorem \ref{Thm: relative fujita}, the adelic line bundle $\overline A$ can be taken to be relatively ample.
\end{rema}

\subsection{Lower bound of intersection product}
\begin{theo} \label{Thm: lower bound intersection number bigness}
Let $X$ be an integral projective $K$-scheme, and $\overline L_0,\ldots,\overline L_d$ be adelic line bundles on $X$. For any $i\in\{0,\ldots,d\}$, let   
\[\delta_i=(L_0\cdots L_{i-1}L_{i+1}\cdots L_d).\] Suppose that
\begin{enumerate}[label=\rm(\arabic*)]
\item  $\overline L_1,\ldots,\overline L_d$ are relatively nef and $L_0$ is pseudo-effective. 
\item if $\delta_0=0$, then $\operatorname{\widehat{\mu}_{\max}^{\mathrm{asy}}}(\overline L_0)>-\infty$,
\item for any $i\in\{1,\ldots,d\}$, if $\delta_i=0$, then $\operatorname{\widehat{\mu}_{\min}^{\mathrm{asy}}}(\overline L_i)>-\infty$.
\end{enumerate}
Then the following inequality holds:
\begin{equation}\label{Equ: minoration dintersection pseudoeffective}(\overline L_0\cdots\overline L_d)_S\geqslant\delta_0\operatorname{\widehat{\mu}_{\max}^{\mathrm{asy}}}(\overline L_0)+\sum_{i=1}^d\delta_i\operatorname{\widehat{\mu}_{\min}^{\mathrm{asy}}}(\overline L_i).\end{equation}
\end{theo}
\begin{proof} If the set 
\[\{\operatorname{\widehat{\mu}_{\max}^{\mathrm{asy}}}(\overline L_0), \operatorname{\widehat{\mu}_{\min}^{\mathrm{asy}}}(\overline L_1),\ldots,\operatorname{\widehat{\mu}_{\min}^{\mathrm{asy}}}(\overline L_d)\}\]
contains $-\infty$, then the inequality \eqref{Equ: minoration dintersection pseudoeffective} is trivial. So we may assume without loss of generality that 
\[\{\operatorname{\widehat{\mu}_{\max}^{\mathrm{asy}}}(\overline L_0), \operatorname{\widehat{\mu}_{\min}^{\mathrm{asy}}}(\overline L_1),\ldots,\operatorname{\widehat{\mu}_{\min}^{\mathrm{asy}}}(\overline L_d)\}\subseteq\mathbb R.\]

Let $\overline M$ be an adelic line bundle on $X$ such that $M$ is big. For any $n\in\mathbb N_{\geqslant 1}$, let 
\[\overline L_{0,n}=\overline L_0^{\otimes n}\otimes\overline M.\]
For any $i\in\{1,\ldots,n\}$, let
\begin{gather*}
\delta_{i}'=(ML_1\cdots L_{i-1}L_{i+1}\cdots L_d)\\
\delta_{i,n}=(L_{0,n}L_1\cdots L_{i-1}L_{i+1}\cdots L_d)=n\delta_i+\delta_{i,n}'.\end{gather*}
By Theorem \ref{Thm: relative fujita} (see also Remark \ref{Rem: strong relative Fujita}), for any real number $t<\operatorname{\widehat{\mu}}_{\max}^{\mathrm{asy}}(\overline L_{0,n})$, there exists a positive integer $p$, a birational projective morphism $g:X'\rightarrow X$, a relatively ample adelic line bundle $\overline A$ and an effective adelic line bundle $E$ on $X'$ such that
\[g^*(\overline L_{0,n}^{\otimes p})=\overline A\otimes\overline E,\quad \widehat{\mu}_{\min}^{\mathrm{asy}}(\overline A)\geqslant pt.\]
 By Theorem \ref{Thm: mu min asy and pull-back}
 for any $i\in\{1,\ldots,d\}$, one has
\[\widehat{\mu}_{\min}^{\mathrm{asy}}(g^*(\overline L_i))\geqslant \widehat{\mu}_{\min}^{\mathrm{asy}}(\overline L_i).\]
Therefore, by Proposition \ref{Pro: lower bound intersection number} and Proposition \ref{Pro: lower bound ample}, we obtain 
\begin{gather*}(\overline E \cdot g^*(\overline L_1)\cdots g^*(\overline L_d))_S\geqslant\sum_{i=1}^d(E \cdot L_1\cdots L_{i-1} \cdot L_i\cdots L_{d})\operatorname{\widehat{\mu}_{\min}^{\mathrm{asy}}}(\overline L_i),\\
(\overline A \cdot g^*(\overline L_1)\cdots g^*(\overline L_d))_S\geqslant \delta_0\operatorname{\widehat{\mu}_{\min}^{\mathrm{asy}}}(\overline A)+\sum_{i=1}^d(A \cdot L_1\cdots L_{i-1} \cdot L_i\cdots L_d)\operatorname{\widehat{\mu}_{\min}^{\mathrm{asy}}}(\overline L_i).
\end{gather*}
Taking the sum, we obtain 
\[(\overline L_{0,n}^{\otimes p} \cdot \overline L_1\cdots\overline L_d)_S\geqslant\delta_0\operatorname{\widehat{\mu}}_{\min}^{\mathrm{asy}}(\overline A)+\sum_{i=1}^dp\delta_{i,n}\operatorname{\widehat{\mu}_{\min}^{\mathrm{asy}}}(\overline L_i)\geqslant\delta_0pt+\sum_{i=1}^dp\delta_{i,n}\operatorname{\widehat{\mu}_{\min}^{\mathrm{asy}}}(\overline L_i).
\]
Since $t$ is arbitrary, we deduce 
\[(\overline L_{0,n}^{\otimes p} \cdot \overline L_1\cdots\overline L_d)_S\geqslant\delta_0\operatorname{\widehat{\mu}_{\max}^{\mathrm{asy}}}(\overline L_{0,n})+\sum_{i=1}^d\delta_{i,n}\operatorname{\widehat{\mu}_{\min}^{\mathrm{asy}}}(\overline L_i).\]
Dividing the two sides by $n$ and then taking the limit when $n\rightarrow+\infty$, we obtain
\[(\overline L_0\cdots\overline L_d)_S\geqslant\delta_0\operatorname{\widehat{\mu}_{\max}^{\mathrm{asy}}}(\overline L_{0})+\sum_{i=1}^d\delta_{i}\operatorname{\widehat{\mu}_{\min}^{\mathrm{asy}}}(\overline L_i).\]
 
\end{proof}

\subsection{Convergence of the first minimum}

In this section, we let $f:X\rightarrow\Spec K$ be an integral projective scheme over $\Spec K$.

\begin{defi}
Let $\overline E=(E,(\|\ndot\|_\omega)_{\omega\in\Omega})$ be an adelic vector bundle on $S$. For any non-zero element $s$ in $E$, let 
\[\widehat{\deg}(s):=-\int_{\Omega}\ln\|s\|_\omega\,\nu(\mathrm{d}\omega).\]
If $E$ is non zero, we define \[\lambda_{\max}(\overline E):=\sup_{s\in E\setminus\{0\}}\widehat{\deg}(s).\] 
Clearly one has 
\begin{equation}\label{Equ: maximal minimum small than maximal slope}\lambda_{\max}(\overline E)\leqslant\widehat{\mu}_{\max}(\overline E).\end{equation}
\end{defi}

\begin{prop}\label{Pro: lambda max super additive}
Let $\overline L=(L,\varphi)$ and $\overline M=(M,\psi)$ be adelic line bundles on $X$ such that both $H^0(X,L)$ and $H^0(X,M)$ are non-zero. Then the following inequality holds:
\[\lambda_{\max}(f_*(\overline L\otimes\overline M))\geqslant\lambda_{\max}(f_*(\overline L))+\lambda_{\max}(f_*(\overline M)).\] 
\end{prop}
\begin{proof}
Let $s$ and $t$ be respectively non-zero elements of $H^0(X,L)$ and $H^0(X,M)$. For any $\omega\in\Omega$, one has 
\[\|st\|_{\varphi_\omega+\psi_\omega}\leqslant\|s\|_{\varphi_\omega}\cdot\|t\|_{\psi_\omega},\]
which leads to 
\[\lambda_{\max}(f_*(\overline L\otimes\overline M))\geqslant \widehat{\deg}(st)\geqslant\widehat{\deg}(s)+\widehat{\deg}(t).\]
Taking the supremum with respect to $s$ and $t$, we obtain the required inequality.
\end{proof}

Let $\overline L$ be an adelic line bundle on $X$ such that $L$ is big. Similarly to Corollary \ref{Cor: convergence maximal slopes}, the sequence 
\[\frac 1n\lambda_{\max}(f_*(\overline L^{\otimes n})),\quad n\in\mathbb N_{\geqslant 1}\]
converges to a real number, which we denote by $\lambda_{\max}^{\mathrm{asy}}(\overline L)$ and called the \emph{asymptotic first minimum}\index{first minimum!asymptotic ---} of $\overline L$. By definition, for any $p\in\mathbb N_{\geqslant 1}$ one has 
\[\lambda_{\max}^{\mathrm{asy}}(\overline L^{\otimes p})=p\lambda_{\max}^{\mathrm{asy}}(\overline L).\]
Proposition \ref{Pro: lambda max super additive} also implies that, if $\overline L$ and $\overline M$ are adelic line bundles on $X$ such that both $L$ and $M$ are big, one has 
\begin{equation}\label{Equ: super additivity lambda}\lambda_{\max}^{\mathrm{asy}}(\overline L\otimes\overline M)\geqslant\lambda_{\max}^{\mathrm{asy}}(\overline L)+\lambda_{\max}^{\mathrm{asy}}(\overline M).\end{equation}
Similarly to Proposition \ref{Pro: sequence convergence mu max asy}, this inequality allows to extend continuously the function $\lambda_{\max}^{\mathrm{asy}}(\ndot)$ to the cone of adelic line bundles $\overline L$ such that $L$ is pseudo-effective: if $\overline L$ is an adelic line bundle on $X$ such that $L$ is pseudo-effective, then, for any adelic line bundle $\overline A$ on $X$, the sequence 
\begin{equation}\label{Equ: convergence of lambda max asy}\frac 1n\lambda_{\max}^{\mathrm{asy}}(\overline L^{\otimes n}\otimes\overline A),\quad n\in\mathbb N_{\geqslant 1}\end{equation}
converges in $\mathbb R\cup\{-\infty\}$ and its limit does not depend on the choice of $\overline A$. For the proof of this statement one can following the strategy of the proof of Proposition \ref{Pro: convergence and independence} in using the inequality \ref{Equ: super additivity lambda} and the fact that, if $A$ is a big line bundle and $B$ is a line bundle on $X$, then there exists a positive integer $p$ such that $B^\vee\otimes A^{\otimes p}$ is big.
We denote the limit of the sequence \eqref{Equ: convergence of lambda max asy} by $\lambda_{\max}^{\mathrm{asy}}(\overline L)$. By \eqref{Equ: maximal minimum small than maximal slope} we obtain that 
\begin{equation}
\lambda_{\max}^{\mathrm{asy}}(\overline L)\leqslant\widehat{\mu}_{\max}^{\mathrm{asy}}(\overline L)
\end{equation}
for any adelic line bundle $\overline L$ such that $L$ is pseudo-effective.

\subsection{Height inequalities}

\begin{prop}\label{Pro: comparison height and mu max}
Let $f:X\rightarrow\Spec K$ be an integral projective scheme over $\Spec K$ and $\overline L$ be an adelic line bundle on $X$ which is relatively nef and such that $(L^d)>0$. Then the following inequality holds:
\begin{equation}
\widehat{\mu}^{\mathrm{asy}}(\overline L) =  \frac{(\overline L^{d+1})_S}{(d+1)(L^d)}\leqslant\widehat{\mu}_{\max}^{\mathrm{asy}}(\overline L).
\end{equation} 
\end{prop}
\begin{proof}
We first consider the case where $L$ is relatively ample.  As in the proof of Proposition \ref{Pro: inequality height and mu min asy }, one has 
\[\frac{(\overline L^{d+1})_S}{(d+1)(L^d)}=\lim_{n\rightarrow+\infty}\frac{\widehat{\mu}(f_*(\overline L^{\otimes n}))}{n}\leqslant\lim_{n\rightarrow+\infty}\frac{\widehat{\mu}_{\max}(f_*(\overline L^{\otimes n}))}{n}=\widehat{\mu}_{\max}^{\mathrm{asy}}(\overline L).\]

We now consider the general case. Let $\overline A$ be a relatively ample adelic line bundle on $X$ such that $\overline L^{\otimes n}\otimes\overline A$ is relatively ample for sufficiently large positive integer $n$. For any $n\in\mathbb N_{\geqslant 1}$, let $\overline L_n=\overline L^{\otimes n}\otimes\overline A$. The particular case of the proposition proved above shows that 
\[\frac{(\overline L_n^{d+1})_S}{(d+1)(L_n^d)}\leqslant\widehat{\mu}_{\max}^{\mathrm{asy}}(\overline L_n)\]
if $n$ is sufficiently large.
Taking the limit when $n\rightarrow+\infty$, by the relations 
\[\lim_{n\rightarrow+\infty}\frac{(\overline L_n^{d+1})_S}{n^{d+1}}=(\overline L^{d+1})_S,\quad \lim_{n\rightarrow+\infty}\frac{(L_n^d)}{n^d}=(L^d)\]
and Proposition \ref{Pro: sequence convergence mu max asy} we obtain the desired result.
\end{proof}

\begin{rema}
Combining Propositions \ref{Pro: comparison height and mu max} and \ref{Pro: inequality height and mu min asy }, we obtain that, if $\overline L$ is relatively nef and if $(L^d)>0$, then the inequality $\widehat{\mu}_{\min}^{\mathrm{asy}}(\overline L)\leqslant \widehat{\mu}_{\max}^{\mathrm{asy}}(\overline L)$. This inequality also holds for relatively nef adelic line bundle $\overline L$ with $(L^d)=0$. It suffices to choose an auxiliary relatively ample adelic line bundle $\overline M$ and deduce the inequality from
\[\frac 1n\widehat{\mu}_{\min}^{\mathrm{asy}}(\overline L^{\otimes n}\otimes\overline M)\leqslant\frac 1n\widehat{\mu}_{\max}^{\mathrm{asy}}(\overline L^{\otimes n}\otimes\overline M)\]
by taking the limit when $n\rightarrow+\infty$.
\end{rema}

\begin{theo}\label{Thm: minoration mu min asy L} Let $X$ be a non-empty and reduced projective $K$-scheme and $\Theta_X$ be the set of all integral closed subschemes of $X$.
Let $\overline L=(L,\varphi)$ be a relatively ample adelic line bundle on $X$. Then the following equality holds:
\begin{equation}\label{Equ: equality mu min asy and abs infimum}\widehat{\mu}_{\min}^{\mathrm{asy}}(\overline L)= \inf_{Y\in \Theta_X}\widehat{\mu}_{\max}^{\mathrm{asy}}(\overline L|_Y) = \inf_{Y\in \Theta_X}\frac{h_{\overline L}(Y)}{(\dim(Y)+1)\deg_L(Y)}.\end{equation}
\end{theo}
\begin{proof}  

For any $Y\in\Theta_X$ and any $n\in\mathbb N$, let $V_{Y,n}(L)$ be the image of the restriction map
\[H^0(X,L^{\otimes n})\longrightarrow H^0(Y,L|_Y^{\otimes n}).\]
We equip $V_{Y,n}(L)$ with the quotient norm family $\xi^Y_{n}$ induced by $\xi_{n\varphi}=(\|\ndot\|_{n\varphi_\omega})_{\omega\in\Omega}$ to obtain adelic vector bundle on $S$.

\ifsmf\begin{enonce}{Claim}\fi
\ifams\begin{clai}\fi
\label{Claim: quotient asy max slope}
For any $Y\in\Theta_X$, the following equality holds
\[\lim_{n\rightarrow+\infty}\frac{\widehat{\mu}_{\max}(V_{Y,n}(L),\xi_n^Y)}{n}=\widehat{\mu}_{\max}^{\mathrm{asy}}(\overline L|_Y).\]
\ifsmf\end{enonce}\fi
\ifams\end{clai}\fi
\begin{proof}
Since $L$ is ample, there exists $N\in\mathbb N_{\geqslant 1}$ such that, for any $n\in\mathbb N_{\geqslant N}$, the following conditions hold:
\begin{enumerate}[label=\rm(\alph*)]
\item $L^{\otimes n}$ is generated by global sections, and the canonical $K$-morphism \[j_n:X\longrightarrow\mathbb P(H^0(X,L^{\otimes n}))\] is a closed embedding,
\item  
one has $V_{Y,n}(L)=H^0(Y,L|_Y^{\otimes n})$. 
\end{enumerate}
For any $n\in\mathbb N_{\geqslant N}$, we denote by $\psi_n$ the metric family on $L$ such that $n\psi_n$ coincides with the quotient metric family induced by $\xi_{n\varphi}$. Similarly, we denote by $\psi^Y_n$ the metric family on $L|_Y$ such that $n\psi^Y_n$ coincides with the quotient metric family induced by $\xi^Y_n$. By definition, $\psi^Y_n$ identifies with the restriction of $\psi_n$ to $L|_Y$. Moreover, we denote by $\varphi^Y$ the restriction of the metric family $\varphi$ to $L|_Y$. By \cite[Proposition 2.2.22]{CMArakelovAdelic}, for any $\omega\in\Omega$, one has $\psi^Y_{n,\omega}\geqslant\varphi_\omega^Y$, and
\[\xi_n^Y=(\|\ndot\|_{n\psi^{Y}_{n,\omega}})_{\omega\in\Omega}. \]
Since $\varphi$ is semi-positive, for any $\omega\in\Omega$, one has 
\[\lim_{n\rightarrow+\infty}d_\omega(\psi_{n}^Y,\varphi^Y)=0.\]
Denote by $\xi_{n\varphi^Y}$ the norm family $(\|\ndot\|_{n\varphi^Y_\omega})_{\omega\in\Omega}$. Since 
\[d_\omega(\xi_n^Y,\xi_{n\varphi^Y})\leqslant d_\omega(n\psi_{n}^Y,n\varphi^Y)=nd_\omega(\psi_{n}^Y,\varphi^Y),\]
we obtain
\begin{equation}
\label{Equ: distance quotient and sup norms}\lim_{n\rightarrow+\infty}\frac 1nd_\omega(\xi^Y_n,\xi_{n\varphi^Y})=0.\end{equation}
By \cite[Proposition 2.2.22 (5)]{CMArakelovAdelic}, the function 
\[(\omega\in\Omega)\longmapsto d_\omega(\xi_n^Y,\xi_{n\varphi^Y})\] 
is dominated. Therefore, by Lebesgue's dominated convergence theorem we obtain 
\begin{equation}\label{Equ: extension property}\lim_{n\rightarrow+\infty}\frac 1n d(\xi_n^Y,\xi_{n\varphi^Y})=\lim_{n\rightarrow+\infty}\frac 1n \int_{\omega\in\Omega}d_\omega(\xi_n^Y,\xi_{n\varphi^Y})\,\nu(\mathrm{d}\omega).\end{equation}
Finally, by \cite[Proposition 4.3.31]{CMArakelovAdelic}, one has 
\[\bigg|\frac{\widehat{\mu}_{\max}(V_{Y,n}(L),\xi_{n}^Y)}{n}-\frac{\widehat{\mu}_{\max}(V_{Y,n}(L),\xi_{n\varphi^Y})}{n}\bigg|\leqslant\frac 1n d(\xi_n^Y,\xi_{n\varphi^Y}).\]
Passing to limit when $n\rightarrow+\infty$, we obtain the desired equality.
\end{proof}

Combining Claim \ref{Claim: quotient asy max slope} with \cite[Theorem 7.2.4]{CMArakelovAdelic}, we obtain  
\[\widehat{\mu}_{\min}^{\mathrm{asy}}(\overline L)\geqslant\inf_{Y\in\Theta_X}\widehat{\mu}_{\max}^{\mathrm{asy}}(\overline L|_Y).\]
By Proposition \ref{Pro: comparison height and mu max}, for any $Y\in\Theta_X$, one has 
\[\frac{h_{\overline L}(Y)}{(\dim(Y)+1)\deg_L(Y)}\leqslant\widehat{\mu}_{\max}^{\mathrm{asy}}(\overline L|_Y).\]
Finally, by Corollary \ref{Cor: mu min bounded by height}, for any $Y\in\Theta_X$, one has 
\[\frac{h_{\overline L}(Y)}{(\dim(Y)+1)\deg_L(Y)}\geqslant\widehat{\mu}_{\min}^{\mathrm{asy}}(\overline L).\]
Thus \eqref{Equ: equality mu min asy and abs infimum} is proved.

\end{proof}

\begin{coro}\label{Cor: monoration de hauteur normalisee} Let $X$ be an integral projective $K$-scheme and $\Theta_X$ be the set of all integral closed subschemes of $X$.
Let $\overline L=(L,\varphi)$ be a relatively ample adelic line bundle on $X$. Then the following inequality holds:
\begin{equation}\label{Equ: equality mu min asy and lambda max}\frac{(\overline L^{(d+1)})_S}{(d+1)(L^d)}\geqslant\frac{1}{d+1}\lambda_{\max}^{\mathrm{asy}}(\overline L)+\frac{d}{d+1}\inf_{Y\in\Theta_X\setminus\{X\}}\widehat{\mu}_{\max}^{\mathrm{asy}}(\overline L|_Y).\end{equation}
In particular, if 
\[\lambda_{\max}^{\mathrm{asy}}(\overline L)\geqslant\frac{(\overline L^{(d+1)})}{(d+1)(L^d)},\]
then the inequality 
\begin{equation}\frac{(\overline L^{(d+1)})_S}{(d+1)(L^d)}\geqslant\inf_{Y\in\Theta_X\setminus\{X\}}\widehat{\mu}_{\max}^{\mathrm{asy}}(\overline L|_Y)\end{equation}
holds. 
\end{coro}
\begin{proof}
The case where $d=0$ is trivial. In the following, we suppose that $d\geqslant 1$. By replacing $\overline L$ by a tensor power, we may assume that 
\[V_\sbullet(L)=\bigoplus_{n\in\mathbb N}H^0(X,L^{\otimes n})\]
is generated as $K$-algebra by $V_1(L)$.   For any $n\in\mathbb N$, we let $h^0(L^{\otimes n})$ be the dimension of $H^0(X,L^{\otimes n})$ over $K$. Let $s$ be a non-zero global section of $L$ and $I_\sbullet$ be the homogeneous ideal of  $V_\sbullet(L)$ generated by $s$. Then one can find a sequence 
\[I_{\sbullet}=I_{0,\sbullet}\subsetneq I_{1,\sbullet}\subsetneq\ldots\subsetneq I_{r,\sbullet}=V_\sbullet(L)\]
of homogeneous ideals of $R_\sbullet$ and non-zero homogeneous prime ideals $P_{i,\sbullet}$, $i\in\{1,\ldots,r\}$, of $V_\sbullet(L)$ such that 
\[\forall\,i\in\{1,\ldots,r\},\quad P_{i,\sbullet}\cdot I_{i,\sbullet}\subset I_{i-1,\sbullet}.\]
We denote by $Y_i$ the integral closed subscheme of $X$ defined by the homogeneous ideal $P_{i,\sbullet}$.

Consider the following sequence
\[
\begin{array}{ccccccc}
{V}_{0}(L)  & \overset{\cdot s}{\longrightarrow} & {I}_{0, 1} & \hookrightarrow \cdots \hookrightarrow &
{I}_{i, 1} & \hookrightarrow \cdots \hookrightarrow & {I}_{r, 1} = {V}_1(L) \\
 & \vdots &  \vdots & \vdots &  \vdots & \vdots & \vdots \\
& \overset{\cdot s}{\longrightarrow} & {I}_{0, j} & \hookrightarrow \cdots \hookrightarrow &
{I}_{i, j} & \hookrightarrow \cdots \hookrightarrow & {I}_{r, j} = {V}_j(\overline L) \\
& \overset{\cdot s}{\longrightarrow} & {I}_{0, j+1} & \hookrightarrow \cdots \hookrightarrow &
{I}_{i, j+1} & \hookrightarrow \cdots \hookrightarrow & {I}_{r, j+1} = {V}_{j+1}(L) \\
 & \vdots &  \vdots & \vdots &  \vdots & \vdots & \vdots \\
& \overset{\cdot s}{\longrightarrow} & {I}_{0, n} & \hookrightarrow \cdots \hookrightarrow &
{I}_{i, n} & \hookrightarrow \cdots \hookrightarrow & {I}_{r, n} = {V}_n(L)
\end{array}
\]
By \cite[Proposition 4.3.13]{CMArakelovAdelic}, one has 
\begin{equation}\label{Equ: minoration deg fstar Ln}\widehat{\deg}(f_*(\overline L^{\otimes n}))\geqslant\sum_{j=1}^n\sum_{i=1}^r\widehat{\deg}(\overline{I_{i,j}/I_{i-1,j}})+\widehat{\deg}(s)\sum_{k=0}^{n-1}h^0(L^{\otimes k}).\end{equation}
By \cite[Proposition 7.1.6]{CMArakelovAdelic} and \eqref{Equ: extension property}, one has 
\begin{equation}\label{Equ: minormation liminf Ii over I i-1}\liminf_{m\rightarrow+\infty}\frac{\widehat{\mu}_{\min}(\overline{I_{i,m}/I_{i-1,m}})}{m}\geqslant\widehat{\mu}_{\min}^{\mathrm{asy}}(\overline L|_{Y_i}).\end{equation}
Moreover, by the asymptotic Riemann-Roch formula, one has 
\[h^0(L^{\otimes k})=\frac{(L^d)}{d!}k^d+O(k^{d-1}),\]
which leads to 
\[\lim_{n\rightarrow+\infty}\frac{1}{nh^0(L^{\otimes n})}\sum_{j=0}^{n-1}h^0(L^{\otimes j})=\frac {1}{d+1}.\]
For any integers $n$ and $m$ such that $1\leqslant m\leqslant n$, we deduce from \eqref{Equ: minoration deg fstar Ln} that 
\[\begin{split}\widehat{\deg}(f_*(\overline L^{\otimes n}&))\geqslant\sum_{j=1}^m\sum_{i=1}^r\widehat{\deg}(\overline{I_{i,j}/I_{i-1,j}})\\
&+\min_{i\in\{1,\ldots,r\}}\inf_{\ell\in\mathbb N_{\geqslant m}}\frac{\widehat{\mu}_{\min}(\overline{I_{i,\ell}/I_{i-1,\ell}})}{\ell}\sum_{j=m+1}^nj(h^0(L^{\otimes j})-h^0(L^{\otimes(j-1)}))\\
&+\widehat{\deg}(s)\sum_{k=0}^{n-1}h^0(L^{\otimes k}).
\end{split}\]
Dividing the two sides by $nh^0(L^{\otimes n})$ and taking the limit when $n\rightarrow+\infty$, we obtain 
\[\frac{(\overline L^{(d+1)})_S}{(d+1)(L^d)}\geqslant\frac{d}{d+1}\min_{i\in\{1,\ldots,r\}}\inf_{\ell\in\mathbb N_{\geqslant m}}\frac{\widehat{\mu}_{\min}(\overline{I_{i,\ell}/I_{i-1,\ell}})}{\ell}+\frac{1}{d+1}\widehat{\deg}(s).\]
Since $m$ is arbitrary, taking the limit when $m\rightarrow+\infty$, by \eqref{Equ: minormation liminf Ii over I i-1} we obtain 
\[\frac{(\overline L^{(d+1)})_S}{(d+1)(L^d)}\geqslant\frac{d}{d+1}\min_{i\in\{1,\ldots,r\}}\widehat{\mu}_{\min}^{\mathrm{asy}}(\overline L|_{Y_i})+\frac 1{d+1}\widehat{\deg}(s).\]
By Theorem \ref{Thm: minoration mu min asy L}, for any $i\in\{1,\ldots,r\}$, one has 
\[\widehat{\mu}_{\min}^{\mathrm{asy}}(\overline L|_{Y_i})=\inf_{Z\in\Theta_{Y_i}}\widehat{\mu}_{\max}^{\mathrm{asy}}(\overline L|_Z)\geqslant\inf_{Z\in\Theta_X\setminus\{X\}}\widehat{\mu}_{\max}^{\mathrm{asy}}(\overline L|_Z).\] 
Since $s$ is arbitrary, we obtain
\[\frac{(\overline L^{(d+1)})_S}{(d+1)(L^d)}\geqslant\frac{1}{d+1}\widehat{\lambda}_{\max}(f_*(\overline L))+\frac{d}{d+1}\inf_{\begin{subarray}{c}Y\in\Theta_X\\
Y\neq X
\end{subarray}}\widehat{\mu}_{\max}^{\mathrm{asy}}(\overline L|_Y).\]
Finally, replacing $\overline L$ by $\overline L^{\otimes p}$ for $p\in\mathbb N_{\geqslant 1}$, we obtain 
\[\frac{(\overline L^{(d+1)})_S}{(d+1)(L^d)}\geqslant\frac{1}{p(d+1)}\widehat{\lambda}_{\max}(f_*(\overline L^{\otimes p}))+\frac{d}{d+1}\inf_{\begin{subarray}{c}Y\in\Theta_X\\
Y\neq X
\end{subarray}}\widehat{\mu}_{\max}^{\mathrm{asy}}(\overline L|_Y).\]
Taking the limite when $p\rightarrow+\infty$, we obtain the inequality \eqref{Equ: equality mu min asy and lambda max}.
\end{proof}

\begin{prop}
Let $X$ be an integral projective scheme over $\Spec K$ and $\overline L$ be an adelic line bundle on $X$ such that $L$ is big. Then the following inequality holds:
\[\widehat{\mu}_{\max}^{\mathrm{asy}}(\overline L)\leqslant\sup_{\begin{subarray}{c}Y\in\Theta_X\\
Y\neq X\end{subarray}}\inf_{x\in (X\setminus Y)^{(0)}}\frac{h_{\overline L}(x)}{[K(x):K]},\]
where $(X\setminus Y)^{(0)}$ denotes the set of closed points of $X\setminus Y$. 
\end{prop}
\begin{proof}
See \cite[Proposition 6.4.4]{CMArakelovAdelic}.
\end{proof}

\subsection{Minkowskian adelic line bundles}

\begin{defi}
Let $X$ be a reduced projective $K$-scheme and $\overline L$ be an adelic line bundle on $X$. We say that $\overline L$ is \emph{Minkowskian}\index{Minkowskian} 
if the inequality below holds:
\[{\lambda}_{\max}^{\mathrm{asy}}(\overline L)\geqslant\widehat{\mu}^{\mathrm{asy}}(\overline L)=\frac{h_{\overline L}(X)}{(\dim(X)+1)\deg_L(X)}.\]  
Moreover, $\overline L$ is said to be \emph{strongly Minkowskian}\index{Minkowskian!strongly ---} if for any integral closed sub-scheme $Y$ of $X$,  one has $\overline L|_Y$ is Minkowskian.  
Note that the strongly Minkowskian condition is satisfied in the following cases:
\begin{enumerate}[label=\rm(\arabic*)]
\item $S$ is the adelic curve associated with a number field, and the metrics of $\overline L$ over non-Archimedean places are almost everywhere induced  by a common integral model defined over the ring of algebraic integers in the number field;
\item $S$ is the adelic curve associated with a regular projective curve over a field, and the metrics of $\overline L$ are induced by an integral model of $L$ over the base curve;
\item $S$ is the adelic curve of a single copy of the trivial absolute value. 
\end{enumerate} 
The case (1) comes from the classic Minkowski theory of Euclidean lattices. The case (2) is a consequence of Riemann-Roch theorem on curves. The case (3) follows from \cite[Remark 4.3.63]{CMArakelovAdelic}.
\end{defi}

\begin{coro}\label{Cor: mu min equals absolute minimum Minkowskian}
Let $X$ be an integral projective $K$-scheme and $\overline L$ be a relatively ample adelic line bundle on $X$. Assume that $\overline L$ is strongly Minkowskian. Then the following inequality holds:
\begin{equation}\label{Equ: bounded by absolute minimum}\frac{(\overline L^{(d+1)})_S}{(d+1)(L^d)}\geqslant\inf_{x\in X^{(0)}}h_{\overline L}(x),\end{equation}
where $X^{(0)}$ denotes the set of closed points of $X$.
Moreover, one has 
\begin{equation}\label{Equ: mu min asy L}\widehat{\mu}_{\min}^{\mathrm{asy}}(\overline L)=\inf_{x\in X^{(0)}}\frac{h_{\overline L}(x)}{[K(x):K]}.\end{equation}
\end{coro}
\begin{proof}
We reason by induction on the dimension $d$ of $X$. The case where $d=0$ is trivial. Assume that $d\geqslant1$ and that the result is true for any integral projective $K$-scheme of dimension $<d$. By Corollary \ref{Cor: monoration de hauteur normalisee} one has 
\[\frac{(\overline L^{(d+1)})_S}{(d+1)(L^d)}\geqslant\inf_{Y\in\Theta_X\setminus\{X\}}\widehat{\mu}_{\max}^{\mathrm{asy}}(\overline L|_Y)\geqslant\inf_{Y\in\Theta_X\setminus\{X\}}\frac{h_{\overline L}(Y)}{(\dim(Y)+1)\deg_L(Y)},\]
where the second inequality comes from Proposition \ref{Pro: comparison height and mu max}. For any $Y\in\Theta_X$ such that $Y\neq X$, one has $\dim(Y)\leqslant 1$. Hence the induction hypothesis leads to 
\begin{equation}\label{Equ: lower bound height}\frac{h_{\overline L}(Y)}{(\dim(Y)+1)\deg_L(Y)}\geqslant\inf_{x\in Y^{(0)}}\frac{h_{\overline L}(x)}{[K(x):K]}\geqslant \inf_{x\in X^{(0)}}\frac{h_{\overline L}(x)}{[K(x):K]}.\end{equation}
The inequality \eqref{Equ: bounded by absolute minimum} is thus proved. 

By Corollary \ref{Cor: mu min bounded by height}, the inequality
\[\widehat{\mu}_{\min}^{\mathrm{asy}}(\overline L)\leqslant\inf_{x\in X^{(0)}}\frac{h_{\overline L}(x)}{[K(x):K]}\]
holds. Conversely, by Theorem \ref{Thm: minoration mu min asy L} and the inequality \eqref{Equ: lower bound height}, one has
\[\begin{split}\widehat{\mu}_{\min}^{\mathrm{asy}}(\overline L) &= \inf_{Y\in \Theta_X}\frac{h_{\overline L}(Y)}{(\dim(Y)+1)\deg_L(Y)}\\&\geqslant \inf_{Y\in\Theta_X}\inf_{x\in Y^{(0)}}\frac{h_{\overline L}(x)}{[K(x):K]}=\inf_{x\in X^{(0)}}\frac{h_{\overline L}(x)}{[K(x):K]}.
\end{split}\]
\end{proof}

\begin{lemm}\label{lemm:finite:morphism:Minkowskian}
Let $\pi : X \to Y$ be a generically finite and surjective morphism of $d$-dimensional projective integral schemes over $K$.
Let $\overline M$ be a relatively ample adelic line bundle on $Y$.
Then we have the following:
\begin{enumerate}[label=\rm(\arabic*)] 
\item $\widehat{\mu}^{\mathrm{asy}}(\pi^* (\overline M)) = \widehat{\mu}^{\mathrm{asy}}(\overline M)$.
\item  $\lambda^{\mathrm{asy}}_{\max}(\pi^* (\overline M))  \geqslant \lambda^{\mathrm{asy}}_{\max}(\overline M)$.
\item If $\overline M$ is Minkowskian, then $\pi^*( \overline M)$ is also Minkowskian.
\end{enumerate}
\end{lemm}

\begin{proof}
(1) By the Hilbert-Samual formula,
\[
\widehat{\mu}^{\mathrm{asy}}(\overline M) = \frac{(\overline M^{d+1})_S}{(M^d)(\dim Y + 1)}\quad\text{and}\quad
\widehat{\mu}^{\mathrm{asy}}(\pi^*(\overline M)) = \frac{(\pi^* (\overline M)^{d+1})_S}{(\pi^*(M)^d)(\dim X + 1)},
\]
and hence the assertion follows because
\[
(\pi^* (\overline M)^{d+1})_S =( \deg \pi) (\overline M^{d+1})_S\quad\text{and}\quad (\pi^*(M)^d) =( \deg\pi )(M^d).
\]

\medskip
(2) is obvious because $\adeg(s) = \adeg (\pi^*(s))$ for $s \in H^0(Y, M) \setminus \{ 0 \}$. Moreover, (3) is a consequence (1) and (2).
\end{proof}

\begin{prop}
Let $\pi : X \to Y$ be a finite morphism of projective integral schemes over $K$.
Let $\overline M$ be an adelic line bundle on $Y$ such that $M$ is ample and $\overline M$ is semi-positive.
If $\overline M$ is strongly Minkowskian, then $\pi^* (\overline M)$ is also strongly Minkowskian.
\end{prop}

\begin{proof}
Let $Z$ be a subvariety of $X$. Then $\rest{\pi}{Z} : Z \to \pi(Z)$ is a finite and surjective morphism, and hence,
by Lemma~\ref{lemm:finite:morphism:Minkowskian},
$\rest{\pi^* (\overline M)}{Z} = (\rest{\pi}{Z})^*(\rest{\overline M}{\pi(Z)})$ is Minkowskian, as required.
\end{proof}

\begin{rema}
Let $L$ be a very ample line bundle on $X$. Then there exist a finite and surjective morphism $\pi : X \to \PP^d_K$
such that $\pi^*(\OO_{\PP^d_K}(1)) \simeq L$. 
Let $\{ s_0, \ldots, s_d \}$ be a basis of $H^0(\PP^d_K, \OO_{\PP^d_K}(1))$. For $\omega \in \Omega$,
let $\|\ndot\|_{\omega}$ be the norm on $H^0(\PP^d_{K_{\omega}}, \OO_{\PP^d_{K_{\omega}}}(1)) = H^0(\PP^d_K, \OO_{\PP^d_K}(1)) \otimes K_{\omega}$ given by
\[
\| a_0 s_1 + \cdots + a_d s_d \|_{\omega} = \begin{cases}
\max \{ |a_0|_{\omega}, \ldots, |a_d|_{\omega} \} & \text{if $\omega \in \Omega_{\mathrm{fin}}$}, \\[1ex]
\sqrt{|a_0|^2_{\omega} + \cdots + |a_d|^2_{\omega}} & \text{if $\omega \in \Omega_{\infty}$},
\end{cases}
\]
where $a_0, \ldots, a_d \in K_{\omega}$.
Let $\psi_{\omega}$ be the Fubini-Study metric of $\OO_{\PP^d_{K_{\omega}}}(1)$ induced by $\|\ndot\|_{\omega}$.
Then it is not difficult to see that $(\OO_{\PP^d_K}(1), \psi = ( \psi_{\omega} )_{\omega \in \Omega})$  is
semipositive and Minkowskian,
so that, by Lemma~\ref{lemm:finite:morphism:Minkowskian},
$(L, \pi^*(\psi))$ is Minkowskian.
\end{rema}

\subsection{Successive minima}

Let $X$ be a reduce projective scheme over $\Spec K$ and $\overline L$ be a relatively ample adelic line bundle on $X$. For any $i\in\{1,\ldots,d+1\}$, let 
\[e_i(\overline L)=\sup_{\begin{subarray}{c}\text{$Y\subseteq X$ closed} \\
\operatorname{codim}(Y)\geqslant i\end{subarray}}\;\inf_{\begin{subarray}{c}\text{$Z\in\Theta_X$}\\
\text{$Z\not\subseteq Y$}
\end{subarray}}\operatorname{\widehat{\mu}}_{\max}^{\mathrm{asy}}(\overline L|_{Z})
.\]
By definition, the following inequalities hold:
\[e_1(\overline L)\geqslant\ldots\geqslant e_{d+1}(\overline L).\]
Moreover, by Theorem \ref{Thm: minoration mu min asy L}, one has 
\[e_{d+1}(\overline L)=\operatorname{\widehat{\mu}}_{\min}^{\mathrm{asy}}(\overline L).\]

\begin{prop}\label{Pro: e1 is equal to mu max asy}
Assume that the scheme $X$ is integral. For any relatively ample adelic line bundle $\overline L$ on $X$, the equality $e_1(\overline L)=\widehat{\mu}_{\max}^{\mathrm{asy}}(\overline L)$ holds.
\end{prop}
\begin{proof}
If $Y$ is a closed subscheme of codimension $1$ of $X$, then $X\not\subseteq Y$. Therefore, the inequality $e_1(\overline L)\leqslant\widehat{\mu}_{\max}^{\mathrm{asy}}(\overline L)$ holds. In the following, we show the converse inequality. Let $t$ be a real number such that $t>e_1(\overline L)$. By definition, there exists a family $(Z_i)_{i\in I}$ of integral closed subschemes of $Y$ such that $\widehat{\mu}_{\max}^{\mathrm{asy}}(\overline L|_{Z_i})\leqslant t$ for any $i\in I$ and that the generic points of $Z_i$ form a Zariski dense family in $X$. 

Let $m$ be a positive integer and $E_m$ be a vector subspace of $H^0(X,L^{\otimes m})$ such that 
\begin{equation}\label{Equ: mu min Em }\widehat{\mu}_{\min}(\overline E_m)=\widehat{\mu}_{\max}(f_*(\overline L^{\otimes m})).\end{equation} For any positive integer $n$, let $F_{m,n}$ be the image of $E_m^{\otimes n}$ by the multiplication map
\[H^0(X,L^{\otimes m})^{\otimes n}\longrightarrow H^0(X,L^{\otimes mn}).\]
By \cite[Proposition 4.3.31 and Corollary 5.6.2]{CMArakelovAdelic} (see also Remark \ref{Rem: positive characteristic minimal slope}), one has
\begin{equation}\label{Equ: lower bound Fmn}\widehat{\mu}_{\min}(\overline F_{m,n})\geqslant n\Big(\widehat{\mu}_{\min}(\overline E_m)-\frac 32\nu(\Omega_\infty)\ln(\dim_K(E_m))\Big).\end{equation}
Moreover, there exists $i\in I$ such that the generic point of $Z_i$ does not belong to the base locus of $E_m$ (namely the closed subscheme of $X$ defined by the ideal sheaf $\operatorname{Im}(E_m\otimes L^{\vee\otimes m}\rightarrow\mathcal O_X)$). Therefore the image of $F_{m,n}$ by the restriction map 
\[H^0(X,L^{\otimes mn})\longrightarrow H^0(Z_i,L^{\otimes mn}|_{Z_i})\]
is non-zero. By \cite[Proposition 4.3.31]{CMArakelovAdelic}, one has
\[\widehat{\mu}_{\min}(\overline F_{m,n})\leqslant\widehat{\mu}_{\max}((f|_{Z_i})_*(\overline L|_{Z_i}^{\otimes mn})).\]
Combining this inequality with \eqref{Equ: mu min Em } and \eqref{Equ: lower bound Fmn}, we obtain  
\[\frac{1}{m}\widehat{\mu}_{\max}(f_*(\overline L^{\otimes m}))\leqslant\frac{1}{mn}\widehat{\mu}_{\max}((f|_{Z_i})_*(\overline L|_{Z_i}^{\otimes mn}))+\frac{3}{2m}\nu(\Omega_\infty)\ln(\dim_K(E_m)).\]
Taking the limit when $n\rightarrow+\infty$, we obtain
\[\frac 1m\widehat{\mu}_{\max}(f_*(\overline L^{\otimes m}))\leqslant t+\frac{3}{2m}\nu(\Omega_\infty)\ln(\dim_K(E_m)).\]
Taking the limit when $m\rightarrow+\infty$, we obtain 
$\widehat{\mu}^{\mathrm{asy}}_{\max}(\overline L)\leqslant t$. Since $t>e_1(\overline L)$
is arbitrary, we get $\widehat{\mu}^{\mathrm{asy}}_{\max}(\overline L)\leqslant e_1(\overline L)$, as required.
\end{proof}

\begin{rema}\label{Rem: successive minima}
Let $\overline L$ be a relatively ample adelic line bundle on $X$. For any $t\in\mathbb R$ and any positive integer $n$, we let $V_n^t(\overline L)$ be the vector subspace of $H^0(X,L^{\otimes n})$ generated by non-zero vector subspaces of minimal slope $\geqslant nt$ and $r_n(t)$ be the dimension of the base locus of $V_n^t(\overline L)$. For $t\in\mathbb R$, let 
\[r(t)=\liminf_{n\rightarrow+\infty}r_n(t).\]
By using the method used in the proof of Proposition \ref{Pro: e1 is equal to mu max asy}, we can show that, for any $i\in\{1,\ldots,d+1\}$ 
\[\sup\{t\in\mathbb R\,|\,r(t)\leqslant i\}\leqslant e_i(\overline L).\]
It is a natural question to ask  if the equality holds. Moreover, we expect that the following inequality is true:
\begin{equation}\label{Equ: Zhang inequality}(d+1)\widehat{\mu}^{\mathrm{asy}}(\overline L)\geqslant\sum_{i=1}^{d+1}e_i(\overline L).\end{equation}

For any $i\in\{1,\ldots, d+1\}$, one has
\[e_i(\overline L)=\sup_{\begin{subarray}{c}\text{$Y\subseteq X$ closed} \\
\operatorname{codim}(Y)\geqslant i\end{subarray}}\;\inf_{\begin{subarray}{c}\text{$Z\in\Theta_X$}\\
\text{$Z\not\subseteq Y$}
\end{subarray}}\operatorname{\widehat{\mu}}_{\max}^{\mathrm{asy}}(\overline L|_Z)\leqslant\sup_{\begin{subarray}{c}\text{$Y\subseteq X$ closed} \\
\operatorname{codim}(Y)\geqslant i\end{subarray}}\;\inf_{x\in (X\setminus Y)^{(0)}}h_{\overline L}(x),
\]
where $(X\setminus Y)^{(0)}$ denotes the set of closed points of $X$ outside of $Y$.
 In the case where $S$ is the adelic curve consisting of places of a number field, by \cite[Theorem 1.5]{MR4271920}, for any integral closed subscheme $Z$ of $X$, one has 
\[\widehat{\mu}_{\max}^{\mathrm{asy}}(\overline L|_Z)=\sup_{\begin{subarray}{c}W\in\Theta_Z\\
W\neq Z
\end{subarray}}\inf_{x\in (Z\setminus W)^{(0)}}h_{\overline L}(x).\]
If $Z$ is not contained in $Y$, then 
\[\widehat{\mu}_{\max}^{\mathrm{asy}}(\overline L|_Z)\geqslant\inf_{x\in (Z\setminus Y)^{(0)}}h_{\overline L}(x)
\geqslant\inf_{x\in(X\setminus Y)^{(0)}}h_{\overline L}(x).\]
Therefore, in this case $e_i(\overline L)$ identifies with the $i$-th minimum of the height function $h_{\overline L}$ in the sense of Zhang. In particular, the inequality follows from \cite[Theorem 5.2]{MR1254133}.
\end{rema}

\section{Positivity conditions for adelic line bundles}

\subsection{Ampleness and nefness}

In this subsection, we let $X$ be a non-empty and reduced projective scheme over $\Spec K$, and let $d$ be the dimension of $X$.

\begin{defi}
We say that an adelic line bundle $\overline L$ on $X$ is \emph{ample}\index{ample} if it is relatively ample and if there exists $\varepsilon>0$ such that the inequality 
\[h_{\overline L}(Y)\geqslant\varepsilon\deg_L(Y)(\dim(Y)+1)\]
holds for any integral closed subscheme $Y$ of $X$.
\end{defi}

\begin{prop}\label{Pro: equivlance condition of ampleness}
Let $\overline L$ be an adelic line bundle which is relatively ample. Then the following statements are equivalent:
\begin{enumerate}[label=\rm(\arabic*)]
\item $\overline L$ is ample,
\item $\widehat{\mu}_{\min}^{\mathrm{asy}}(\overline L)>0$,
\item there exists $\varepsilon>0$ such that, for any integral closed subscheme $Y$ of $X$, one has $\widehat{\mu}_{\max}^{\mathrm{asy}}(\overline L|_Y)>\varepsilon$.

\end{enumerate}
\end{prop}
\begin{proof}
This is a consequence of Theorem \ref{Thm: minoration mu min asy L}.
\end{proof}

\begin{prop}\label{Pro: positivity intersection number}
If $\overline L_0,\ldots,\overline L_d$ are ample adelic line bundles on $X$, then the inequality 
\[(\overline L_0\cdots\overline L_d)_S> 0\]
holds.
\end{prop}
\begin{proof}
This is a consequence of Theorem \ref{Thm: lower bound intersection number} and Proposition \ref{Pro: equivlance condition of ampleness}. 
\end{proof}

\begin{prop}\label{Pro: strongly minkowskian and ampleness}
Let $\overline L$ be an adelic line bundle which is relatively ample and strongly Minkowskian. Then the following conditions are equivalent:
\begin{enumerate}[label=\rm(\arabic*)]
\item $\overline L$ is ample,
\item there exists $\varepsilon>0$ such that, for any closed point $x$ of $X$, one has $h_{\overline L}(x)>\varepsilon$.
\end{enumerate}
\end{prop}
\begin{proof}
This is a consequence of Corollary \ref{Cor: mu min equals absolute minimum Minkowskian}.
\end{proof}

\begin{defi}
We say that an adelic line bundle $\overline L$ on $X$ is \emph{nef}\index{nef} if there exists an ample adelic line bundle $\overline A$ and a positive integer $N$ such that $\overline L^{\otimes n}\otimes\overline A$ is ample for any $n\in\mathbb N_{\geqslant N}$.

\end{defi}

\begin{prop}\label{Pro: nef with positivity of mu min}
Let $\overline L$ be an adelic line bundle on $X$. The following conditions are equivalent:
\begin{enumerate}[label=\rm(\arabic*)]
\item $\overline L$ is nef,
\item $\overline L$ is relatively nef and $\widehat{\mu}_{\min}^{\mathrm{asy}}(\overline L)\geqslant 0$.
\end{enumerate}
\end{prop}
\begin{proof}
Assume that $\overline L$ is nef. By definition, it is relatively nef. Let $\overline A$ be an ample adelic line bundle and $N$ be a positive integer such that $\overline L^{\otimes n}\otimes\overline A$ is ample for any $n\in\mathbb N_{\geqslant N}$. Then one has $\widehat{\mu}_{\min}^{\mathrm{asy}}(\overline L^{\otimes n}\otimes\overline A)>0$, which leads to 
\[\widehat{\mu}_{\min}^{\mathrm{asy}}(\overline L)=\lim_{n\rightarrow+\infty}\frac 1n\widehat{\mu}_{\min}^{\mathrm{asy}}(\overline L^{\otimes n}\otimes\overline A)\geqslant 0.\]

Conversely, we assume that $\overline L$ is relatively nef and $\widehat{\mu}_{\min}^{\mathrm{asy}}(\overline L)\geqslant 0$. Since $\overline L$ is relatively nef, there exists a relatively ample line bundle $\overline A$ and a positive integer $N$ such that $\overline L^{\otimes n}\otimes\overline A$ is relatively ample for any $n\in\mathbb N_{\geqslant N}$. By dilating the metrics of $\overline A$, we may assume that $\widehat{\mu}_{\min}^{\mathrm{asy}}(\overline L)>0$. Then, by Proposition \ref{Pro: super additivity nef}
we obtain that 
\[\forall\,n\in\mathbb N_{\geqslant N},\quad\widehat{\mu}_{\min}^{\mathrm{asy}}(\overline L^{\otimes n}\otimes\overline A)\geqslant n\,\widehat{\mu}_{\min}^{\mathrm{asy}}(\overline L)+\widehat{\mu}_{\min}^{\mathrm{asy}}(\overline A)\geqslant\widehat{\mu}_{\min}^{\mathrm{asy}}(\overline A)>0. \]
Therefore $\overline L^{\otimes n}\otimes\overline A$ is ample.
\end{proof}

\begin{prop}\label{Pro: property of nef line bundles}
\begin{enumerate}[label=\rm(\arabic*)]
\item If $\overline L_0,\ldots,\overline L_d$ are nef adelic line bundles on $X$, then the inequality
$(\overline L_0\cdots\overline L_d)_S\geqslant 0$ holds.
\item If $\overline L$ is a nef adelic line bundle on $X$ and if $g:Y\rightarrow X$ is a projective $K$-morphism, then the pullback $g^*(\overline L) $ is nef.
\item If $\overline L$ is a nef
adelic line bundle on $X$, for any integral closed subscheme $Y$ of $X$, one has $h_{\overline L}(Y)\geqslant 0$.
\item If $\overline L$ is a relatively ample adelic line bundle on $X$ such that $h_{\overline L}(Y)\geqslant 0$ for any integral closed subscheme $Y$ of $X$, then $\overline L$ is nef.
\item If $\overline L$ is a relatively ample adelic line bundle on $X$ such that $\widehat{\mu}_{\max}^{\mathrm{asy}}(\overline L|_Y)\geqslant 0$ for any integral closed subscheme $Y$ of $X$, then $\overline L$ is nef.
\end{enumerate}
\end{prop}
\begin{proof}
The first statement is a consequence of Proposition \ref{Pro: lower bound ample} and Proposition \ref{Pro: nef with positivity of mu min}. The second statement follows from Lemma \ref{lemma:restriction:relatively:nef}, Theorem \ref{Thm: mu min asy and pull-back} and Proposition \ref{Pro: nef with positivity of mu min}. The third statement is a consequence of the first and the second ones. The last two statements are consequences of Theorem \ref{Thm: minoration mu min asy L} and Proposition \ref{Pro: nef with positivity of mu min}.
\end{proof}

\subsection{Bigness and pseudo-effectivity}

In this subsection, we let $X$ be an integral projective $K$-scheme $f:X\rightarrow\Spec K$ and let $d$ be its dimension.

\begin{defi}
Let $\overline L$ be an adelic line bundle on $X$. We define the \emph{arithmetic volume}\index{volume!arithmetic ---} of $\overline L$ as
\[\widehat{\vol}(\overline L):=\limsup_{n\rightarrow+\infty}\frac{\widehat{\deg}_+(f_*(\overline L^{\otimes n}))}{n^{d+1}/(d+1)!}.\]
If $\widehat{\vol}(\overline L)>0$, we say that $\overline L$ is \emph{big}\index{big}. It has been shown in \cite[Proposition 6.4.18]{CMArakelovAdelic} that $\overline L$ is big if and only if $L$ is big and $\widehat{\mu}_{\max}^{\mathrm{asy}}(\overline L)>0$.
\end{defi}

\begin{prop}
An ample adelic line bundle is big.
\end{prop}
\begin{proof}
Let $\overline L$ be an ample adelic line bundle on $X$. Then 
one has $\widehat{\mu}_{\min}^{\mathrm{asy}}(\overline L)>0$, namely for sufficiently large positive integer $n$ one has $\widehat{\mu}_{\min}(f_*(\overline L^{\otimes n}))>0$. By \cite[Proposition 4.3.13]{CMArakelovAdelic}, for such $n$ one has 
\[\widehat{\deg}(f_*(\overline L^{\otimes n}))=\widehat{\deg}_+(f_*(\overline L^{\otimes n})),\]
which leads to, by Theorem~\ref{thm:lower:bound},
\[\widehat{\vol}(\overline L)=(\overline L^{d+1})_S>0,\]
where the inequality comes from  Proposition \ref{Pro: positivity intersection number}. Hence $\overline L$ is big.  
\end{proof}

\begin{rema}
We expect that a variant of the method in the proof of Theorem \ref{Thm: relative fujita} leads to an arithmetic version of Fujita's approximation theorem for big adelic line bundles, which generalizes the results of \cite{MR2722508,MR2575090}.
\end{rema}

\begin{prop}
Let $\overline L_0,\ldots,\overline L_d$ be adelic line bundles on $X$. Assume that $\overline L_0$ is big and $\overline L_1,\ldots,\overline L_d$ are ample, then 
\[(\overline L_0\cdots\overline L_d)_S>0.\]
\end{prop}
\begin{proof}
This is a consequence of Theorem \ref{Thm: lower bound intersection number bigness}.
\end{proof}

\begin{defi}
Let $\overline L$ be an adelic line bundle on $X$. We say $\overline L$  is \emph{pseudo-effective}\index{pseudo-effective} if there exist a big adelic line bundle $\overline M$ on $X$ and a positive integer $n_0$ such that $\overline L^{\otimes n}\otimes\overline M$ is big for any $n\in\mathbb N_{\geqslant n_0}$.
\end{defi}

\begin{prop}\label{Pro: criterion pseudo-effectivity}
Let $\overline L$ be an adelic line bundle on $X$. The following assertions are equivalent:
\begin{enumerate}[label=\rm(\arabic*)]
\item $\overline L$ is pseudo-effective,
\item $L$ is pseudo-effective and $\widehat{\mu}_{\max}^{\mathrm{asy}}(\overline L)\geqslant 0$.
\end{enumerate}
\end{prop}
\begin{proof}
Assume that $\overline L$ is pseudo-effective. Let $\overline M$ be a big adelic line bundle and $n_0$ be a positive integer such that $\overline L^{\otimes n}\otimes\overline M$ is big for any integer $n\geqslant n_0$. In particular, $L^{\otimes n}\otimes M$ is big for any integer $n\geqslant n_0$. Hence $L$ is pseudo-effective. Moreover, for $n\geqslant n_0$, one has $\widehat{\mu}_{\max}^{\mathrm{asy}}(\overline L^{\otimes n}\otimes\overline M)>0$, which implies that
\[\widehat{\mu}_{\max}^{\mathrm{asy}}(\overline L)=\lim_{n\rightarrow+\infty}\frac{1}{n}\widehat{\mu}_{\max}^{\mathrm{asy}}(\overline L^{\otimes n}\otimes\overline M)\geqslant 0.\]

Conversely, assume that $L$ is pseudo-effective and $\widehat{\mu}_{\max}^{\mathrm{asy}}(\overline L)\geqslant 0$. Let $\overline M$ be a big adelic line bundle on $X$. Since $L$ is pseudo-effective, for any positive integer $n$, $L^{\otimes n}\otimes M$ is big. Moreover, by Proposition \ref{Pro: suradiditivity mu pi max} one has
\[\widehat{\mu}_{\max}^{\mathrm{asy}}(\overline L^{\otimes n}\otimes\overline M)\geqslant n\operatorname{\widehat{\mu}_{\max}^{\mathrm{asy}}}(\overline L)\otimes\operatorname{\widehat{\mu}_{\max}^{\mathrm{asy}}}(\overline M)>0.\]
Hence $\overline L^{\otimes n}\otimes\overline M$ is big for any $n\in\mathbb N$, which shows that $\overline L$ is pseudo-effective.
\end{proof}

\begin{prop}\label{Pro: nef pseudoeffective}
\begin{enumerate}[label=\rm(\arabic*)]
\item Let $\overline L_0,\ldots,\overline L_d$ be adelic line bundles on $X$. Assume that $\overline L_0$ is pseudo-effective and that $\overline L_1,\ldots,\overline L_d$ are nef, then the inequality $(\overline L_0\cdots\overline L_d)_S\geqslant 0$ holds.
\item If $\overline L$ is a pseudo-effective adelic line bundle on $X$ and if $g:Y\rightarrow X$ is a surjective and projective morphism, then the pullback $g^*(\overline L)$ is also pseudo-effective.
\item If $\overline L$ is nef, then it is pseudo-effective.
\end{enumerate}
\end{prop}
\begin{proof}
The first statement is a consequence of Theorem \ref{Thm: lower bound intersection number bigness}; the second one is a consequence of Proposition \ref{Pro: pullback by surjective morphism}.

(3) Since $\overline L$ is nef, we obtain that $L$ is nef, and hence is pseudo-effective. Let $\overline A$ be an ample adelic line bundle. For any positive integer $p$, by Proposition \ref{Pro: suradiditivity mu pi max} one has 
\[\frac 1p\widehat{\mu}_{\max}^{\mathrm{asy}}(\overline L^{\otimes p}\otimes\overline A)\geqslant\widehat{\mu}_{\max}^{\mathrm{asy}}(\overline L)+\frac 1p\widehat{\mu}_{\max}^{\mathrm{asy}}(\overline A).\]
Taking the limit when $p\rightarrow+\infty$, we obtain $\widehat{\mu}_{\max}^{\mathrm{asy}}(\overline L)\geqslant 0$. By Proposition \ref{Pro: criterion pseudo-effectivity}, $\overline L$ is pseudo-effective.
\end{proof}

%!TEX root = ./Hilbert_Samuel_Adelic_Curves.tex

\appendix
\renewcommand{\thesection}{\Alph{section}}
\renewcommand{\thesubsection}{\Alph{section}.\arabic{subsection}}
\renewcommand{\thetheo}{\Alph{section}.\arabic{theo}}
\renewcommand{\thechapter}{}
\renewcommand{\theequation}{\Alph{section}.\arabic{equation}}

\chapter*{Appendix}
\ifams\setcounter{theo}{0}\fi
\setcounter{section}{0}
\setcounter{subsection}{0}

\section{Tensorial semi-stability}

 We  recall some constructions and  facts of multi-linear algebra and classical invariant theory. Then we prove a lifting theorem for invariants in a symmetric power of a tensor product  under the action of the product of special linear groups. 
In subsections \ref{Subsec: symmetric power}--\ref{Subsec: several modules},  we fix a commutative ring $k$ with unit. 

\subsection{Symmetric power}
\label{Subsec: symmetric power}

Let $V$ be a free $k$-module of finite rank and $\delta$ be a natural number. We denote by $V^{\otimes\delta}$ the $\delta$-th tensor power of the $k$-module $V$. Note that the symmetric group $\mathfrak S_\delta$ acts $k$-linearly on $V^{\otimes\delta}$ by permuting the tensor factors. The quotient $k$-module of $V^{\otimes\delta}$ by this action of $\mathfrak S_\delta$ is denoted by $S^\delta(V)$. The class of $x_1\otimes\cdots\otimes x_\delta$ in $S^\delta(V)$ is denoted by $x_1\cdots x_\delta$. If $\boldsymbol{e}=(e_i)_{i=1}^d$ is a basis of $V$ over $k$, then 
\[\boldsymbol{e}^{\boldsymbol{a}}:=\prod_{i\in \{1,\ldots,d\}}e_i^{a_i},\quad \boldsymbol{a}=(a_i)_{i=1}^d\in\mathbb N^{d},\;|\boldsymbol{a}|:=a_1+\cdots+a_d=\delta\]
form a basis of $S^\delta(V)$ over $k$. In particular, $S^\delta(V)$ is a free $k$-module. More generally, if $V$ is decomposed into a direct sum 
\[V=V^{(1)}\oplus\cdots\oplus V^{(r)}\]
of free sub-$k$-modules,
then the $k$-linear map
\begin{equation}\label{Equ: symmetric power of a direct sum}\bigoplus_{\begin{subarray}{c}\boldsymbol{b}=(b_1,\ldots,b_r)\in\mathbb N^r\\
|\boldsymbol{b}|=b_1+\cdots+b_r=\delta
\end{subarray}}S^{b_1}(V^{(1)})\otimes\cdots\otimes S^{b_r}(V^{(r)})\longrightarrow S^{\delta}(V),\end{equation}
which sends
\[x_1\otimes\cdots\otimes x_r\in S^{b_1}(V^{(1)})\otimes\cdots\otimes S^{b_r}(V^{(r)})\]
to $x_1\cdots x_r\in S^\delta(V)$, is an isomorphism. 

We call $S^\delta(V)$ the \emph{$\delta$-th symmetric power}\index{symmetric power} of the free $k$-module $V$. Note that $S^\delta$ defines a functor from the category of free $k$-modules of finite rank to itself. Moreover, it preserves the extension of scalars, namely, for any commutative $k$-algebra $A$, one has \[S^\delta(V\otimes_kA)\cong S^\delta(V)\otimes_kA.\]

The graded $k$-algebra structure of the tensor algebra
\[T(V):=\bigoplus_{\delta\in\mathbb N}V^{\otimes\delta}\]
induces by passing to quotient a graded $k$-algebra structure on 
\[\operatorname{Sym}(V):=\bigoplus_{\delta\in\mathbb N} S^\delta(V).\]
This $k$-algebra is commutative and finitely generated, and it is isomorphic to the polynomial ring $k[X_1, \ldots, X_r]$ over $k$, where $r$ is the rank of $V$. 

\subsection{Exterior power}
Let $V$ be a free $k$-module. We denote by \[\Lambda(V)=\bigoplus_{\delta\in\mathbb N}\Lambda^\delta(V)\] the quotient graded $k$-algebra of $T(V)$ by the two-sided ideal generated by elements of the form $x\otimes x$, $x\in V$. If $x_1,\ldots,x_\delta$ are elements of $V$, the image of $x_1\otimes\cdots\otimes x_\delta$ in $\Lambda^\delta(V)$ is denoted by $x_1\wedge\cdots\wedge x_\delta$. Note that, if $(e_i)_{i=1}^d$ is a basis of $V$ over $k$, then 
\[e_{i_1}\wedge\cdots\wedge e_{i_\delta},\quad 1\leqslant i_1<\ldots<i_\delta\leqslant d\]
form a basis of $\Lambda^\delta(V)$ over $k$. The $k$-module $\Lambda^\delta(V)$ is called the \emph{$\delta$-th exterior power}\index{exterior power} of $V$. Note that $\Lambda^\delta$ also defines a functor from the category of free $k$-modules of finite rank to itself, and preserves extensions of scalars.

If $\delta$ is a natural number, we denote by $\iota_\delta:\Lambda^{\delta}(V)\rightarrow V^{\otimes \delta}$ the anti-symmetrization map which sends $x_1\wedge\cdots\wedge x_\delta\in\Lambda^{\delta}(V)$ to 
\[\sum_{\sigma\in\mathfrak S_\delta}\operatorname{sgn}(\sigma)x_{\sigma(1)}\otimes\cdots\otimes x_{\sigma(\delta)}.\]
This is an injective $\operatorname{GL}(V)$-linear map. It is however \emph{not} a section of the projection $E^{\otimes \delta}\rightarrow\Lambda^\delta(V)$. 

\subsection{Schur functor}

We denote by $\mathbb N_{\geqslant 1}$ the set of positive integers and by $\mathbb N^{\oplus\infty}$ the set of sequences $\lambda=(\lambda_i)_{i\in\mathbb N_{\geqslant 1}}$ of natural numbers indexed by $\mathbb N_{\geqslant 1}$ such that $\lambda_i=0$ except finitely many $i$.  
 For any $\lambda=(\lambda_i)_{i\in\mathbb N_{\geqslant 1}}\in\mathbb N^{\oplus\infty}$, we denote by $|\lambda|$ the sum 
\[\sum_{i\in\mathbb N_{\geqslant 1}}\lambda_i,\]
called the \emph{weight}\index{weight} of $\lambda$.
If $n$ is a natural number and $\lambda=(\lambda_1,\ldots,\lambda_n)\in\mathbb N^n$,  by abuse of notation we denote by $\lambda$ the sequence
\[(\lambda_1,\ldots,\lambda_n,0,\ldots,0,\ldots)\in\mathbb N^{\oplus\infty}.\]

We call \emph{partition}\index{partition} any sequence $\lambda=(\lambda_i)_{i\in\mathbb N_{\geqslant 1}}\in\mathbb N^{\oplus\infty}$ such that \[\lambda_1\geqslant\lambda_2\geqslant\ldots.\] The value $\sup\{i\in\mathbb N_{\geqslant 1}\,:\,\lambda_i>0\}$ (with the convention $\sup\varnothing=0$) is called the \emph{length}\index{length} of the partition $\lambda$. For any $\delta\in\mathbb N$, we denote by $\mathscr P_{\delta}$ the set of partitions of weight $\delta$.

If $\lambda=(\lambda_i)_{i\in\mathbb N_{\geqslant 1}}$ is a sequence in $\mathbb N^{\oplus\infty}$, we denote by $\widetilde{\lambda}=(\widetilde{\lambda}_n)_{n\in\mathbb N_{\geqslant 1}}$ the sequence defined as 
\[\widetilde{\lambda}_n=\sum_{i\in\mathbb N_{\geqslant 1},\,\lambda_i\geqslant n}1.\]
We call $\widetilde{\lambda}$ the \emph{transpose}\index{transpose} of $\lambda$. Clearly one has 
\[\widetilde\lambda_1\geqslant\widetilde\lambda_2\geqslant\ldots\]
and hence $\widetilde{\lambda}$ is a partition.
Moreover, the following equalities hold: 
\[\sum_{n\in\mathbb N_{\geqslant 1}}\widetilde{\lambda}_n=\sum_{\begin{subarray}{c}(i,n)\in\mathbb N_{\geqslant 1}^2\\
\lambda_i> n\end{subarray}}1=\sum_{i\in\mathbb N_{\geqslant 1}}\sum_{\begin{subarray}{c}n\in\mathbb N_{\geqslant 1}\\n\leqslant\lambda_i\end{subarray}}1=\sum_{i\in\mathbb N_{\geqslant 1}}\lambda_i.\]
Note that the double transpose $\widetilde{\widetilde\lambda}$ is equal to the sequence $\lambda$ sorted in the decreasing order. The following graph illustrates the transpose of a partition.

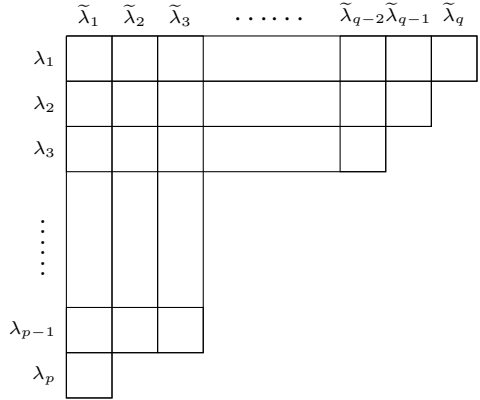
\begin{center}
\begin{tikzpicture}[scale=0.6]
\draw (0,8) -- (9,8) -- (9,7) -- (0,7) -- (0,8) ;
\draw (0,7) -- (8,7) -- (8,6) -- (0,6) -- (0,7) ;
\draw (0,6) -- (7,6) -- (7,5) -- (0,5) -- (0,6) ;
\draw (0,2) -- (3,2) -- (3,1) -- (0,1) -- (0,2) ;
\draw (0,1) -- (1,1) -- (1,0) -- (0,0) -- (0,1) ;

\draw (0,8) -- (1,8) -- (1,0) -- (0,0) -- (0,8) ;
\draw (1,8) -- (2,8) -- (2,1) -- (1,1) -- (1,8) ;
\draw (2,8) -- (3,8) -- (3,1) -- (2,1) -- (2,8) ;
\draw (6,8) -- (7,8) -- (7,5) -- (6,5) -- (6,8) ;
\draw (7,8) -- (8,8) -- (8,6) -- (7,6) -- (7,8) ;
\draw (8,8) -- (9,8) -- (9,7) -- (8,7) -- (8,8) ;

\filldraw  (0,0.5) node[left] {\tiny $\lambda_p$}; 
\filldraw  (0,1.5) node[left] {\tiny $\lambda_{p-1}$}; 
\filldraw  (0,5.5) node[left] {\tiny $\lambda_{3}$}; 
\filldraw  (0,6.5) node[left] {\tiny $\lambda_{2}$}; 
\filldraw  (0,7.5) node[left] {\tiny $\lambda_{1}$}; 
\filldraw  (-0.2,3.85) node[left] {$\vdots$}; 
\filldraw  (-0.2,3.15) node[left] {$\vdots$}; 

\filldraw  (0.5,8) node[above] {\tiny $\widetilde{\lambda}_{1}$}; 
\filldraw  (1.5,8) node[above] {\tiny $\widetilde{\lambda}_{2}$}; 
\filldraw  (2.5,8) node[above] {\tiny $\widetilde{\lambda}_{3}$}; 
\filldraw  (6.5,8) node[above] {\tiny $\widetilde{\lambda}_{q-2}$}; 
\filldraw  (7.5,8) node[above] {\tiny $\widetilde{\lambda}_{q-1}$}; 
\filldraw  (8.5,8) node[above] {\tiny $\widetilde{\lambda}_{q}$}; 
\filldraw  (4.1,8) node[above] {$\cdots$}; 
\filldraw  (4.9,8) node[above] {$\cdots$}; 
\end{tikzpicture}
\end{center}

Let $V$ be a free $k$-module of finite rank. For any  $\lambda=(\lambda_1,\ldots,\lambda_p,0,\ldots,0,\ldots)$ in $\mathbb N^{\oplus\infty}$, we let 
\begin{gather*}
V^{\otimes\lambda}:=V^{\otimes\lambda_1}\otimes\cdots\otimes V^{\otimes\lambda_p}\\
\Lambda^{\lambda}(V):=\Lambda^{\lambda_1}(V)\otimes\cdots\otimes\Lambda^{\lambda_p}(V),\\ S^{\lambda}(V):=S^{\lambda_1}(V)\otimes\cdots\otimes S^{\lambda_p}(V).
\end{gather*}
By the isomorphism \eqref{Equ: symmetric power of a direct sum}, we can identify $S^\lambda(V)$ with a direct summand of $S^\delta(V^{\oplus p})$. Moreover, if $\lambda$ and $\mu$ are two elements of $\mathbb N^{\oplus\infty}$, one has a commutative diagram of canonical $\operatorname{SL}(V)$-linear maps 
\[\xymatrix{V^{\otimes\lambda}\otimes V^{\otimes\mu}\ar[r]^-{\simeq}\ar[d]&V^{\otimes(\lambda+\mu)}\ar[d]\\ S^{\lambda}(V)\otimes S^{\mu}(V)\ar[r]& S^{\lambda+\mu}(V)}\]

If $\lambda=(\lambda_1,\ldots,\lambda_p,0,\ldots, 0,\ldots)$ is a partition  if its transpose is of the form \[\widetilde\lambda=(\widetilde\lambda_1,\ldots,\widetilde{\lambda}_q,0,\ldots,0,\ldots)\] with $\widetilde{\lambda}_q>0$, we denote by $L^\lambda(V)$ the image of the following composed map 
\[\Lambda^\lambda(V)=\Lambda^{\lambda_1}(V)\otimes\cdots\otimes\Lambda^{\lambda_p}(V)\stackrel{\alpha_\lambda}{\longrightarrow} V^{\otimes|\lambda|}\stackrel{\beta_\lambda}{\longrightarrow}S^{\widetilde{\lambda}_1}(V)\otimes\cdots\otimes S^{\widetilde{\lambda}_q}(V)=S^{\widetilde{\lambda}}(V),\]
where $\alpha_{\lambda}$ is induced by the anti-symmetrization maps $\Lambda^{\lambda_i}(V)\rightarrow V^{\otimes\lambda_i}$, and $\beta_\lambda$ sends $x_1\otimes\cdots\otimes x_{|\lambda|}$ to
\begin{multline*}(x_1x_{\lambda_1+1}x_{\lambda_1+\lambda_2+1}\cdots x_{\lambda_1+\cdots+\lambda_{\widetilde{\lambda}_1-1}+1}) \\
\otimes(x_2x_{\lambda_1+2}x_{\lambda_1+\lambda_2+2}\cdots x_{\lambda_1+\cdots+\lambda_{\widetilde{\lambda}_2-1}+2})\\
\otimes\cdots\otimes (x_{\lambda_1}x_{\lambda_1+\lambda_1}x_{\lambda_1+\lambda_2+\lambda_1}\cdots x_{\lambda_1+\cdots+\lambda_{\widetilde{\lambda}_q-1}+\lambda_1}) .
\end{multline*}
The following are two fundamental examples for partitions of $\delta\in\mathbb N$ 
\[L^{(\delta)}(V)=\Lambda^\delta(V),\quad L^{(1,\ldots,1)}(V)=S^\delta(V).\]
It can be shown that $L^\lambda(V)$ is a free $k$-module of finite rank, and $L^\lambda$ defines a functor from the category of free $k$-modules of finite rank to itself (see \cite[\S II.2]{MR658729}). It is called the \emph{Schur functor}\index{Schur functor} with respect to $\lambda$.

\subsection{First fundamental theorem of classical invariant theory}
In this subsection, we assume that $k$ is a field (of any characteristic). We recall the first fundamental theorem of classical invariant theory in a form proved by De Concini and Procesi. We refer to \cite[Theorem 3.3]{MR422314} for proof, see also \cite[Theorem 2.1]{MR2004511}.

\begin{theo}\label{Thm: FFCIT}
Let $V$ be a finite-dimensional vector space over $k$, $r$ be the dimension of $V$ over $k$, and $p$ be a positive integer. Let $V_1,\ldots, V_p$ be $p$ identical copies of $V$. We consider the canonical action of the special linear group $\operatorname{SL}(V)$ on the symmetric algebra $\operatorname{Sym}(V_1\oplus\cdots\oplus V_p)$. Then the invariant sub-$k$-algebra $\operatorname{Sym}(V_1\oplus\cdots\oplus V_p)^{\mathrm{SL}(V)} $ is generated by one-dimensional $k$-vector subspaces of the form
\[\operatorname{Im}\Big(\det(V)=\Lambda^r(V)\longrightarrow V_{i_1}\otimes\cdots\otimes V_{i_r}\Big)\]
in identifying $V_{i_1}\otimes\cdots\otimes V_{i_r}$ with a direct factor of $S^{r}(V_1\oplus\cdots\oplus V_p)$ via the isomorphism \eqref{Equ: symmetric power of a direct sum}, where $i_1,\ldots,i_r$ are distinct elements of $\{1,\ldots,p\}$, and in the above formula we consider the anti-symmetrization map. 
\end{theo}

\begin{rema}\label{Rem: invariant in symmetric power}
Let $\lambda$ be a partition, $p$ be the length of $\lambda$, and $\delta$ be the weight of $\lambda$. We identify $S^\lambda(V)$ with a vector subspace of $S^\delta(V^{\oplus p})$. Theorem \ref{Thm: FFCIT} shows that one can lift the invariant vectors of $S^\delta(V^{\oplus p})$ to tensor powers. More precisely, if $S^\lambda(V)^{\operatorname{SL}(V)}$ is not zero, then $\delta$ should be divisible by $r$, and $S^\lambda(V)^{\operatorname{SL}(V)}$ identifies with the image of
\[\bigoplus_{\begin{subarray}{c}(\mu_1,\ldots,\mu_{\delta/r})\in\mathscr D_{r}^{\delta/r}\\
\mu_1+\cdots+\mu_{\delta/r}=\lambda
\end{subarray}}\det(V)^{\otimes(\delta/r)}\longrightarrow S^\lambda(V).\]  
where $\mathscr D_{r}$ denotes the set of sequences in $\mathbb N^{\oplus\infty}$ of weight $r$ and with coordinates in $\{0,1\}$, and for any $(\mu_1,\ldots,\mu_{\delta/r})\in\mathscr D_r^{\delta/r}$ such that $\mu_1+\cdots+\mu_{\delta/r}=\lambda$, we consider the composed map
\[\det(V)^{\otimes(\delta/r)}\longrightarrow V^{\otimes\mu_1}\otimes\cdots\otimes V^{\otimes\mu_{\delta/r}}\longrightarrow V^{\otimes\lambda}\longrightarrow S^\lambda(V),\]
where the first arrow is induced by anti-symmetrization maps.
\end{rema}

\subsection{Cauchy decomposition}

In this subsection, we consider two free $k$-modules of finite rank $V$ and $W$. The symmetric algebra $\operatorname{Sym}(V\otimes W)$ is naturally equipped with a structure of graded $\mathrm{GL}(V)\times\mathrm{GL}(W)$-module. In the case where $k$ contains $\mathbb Q$, then it is known that $\operatorname{Sym}(V\otimes W)$ is isomorphic as $\mathrm{GL}(V)\times\mathrm{GL}(W)$-module to a direct sum
\[\bigoplus_{\lambda}L^\lambda(V)\otimes L^\lambda(W),\]
where $\lambda$ runs over the set of all partitions. In general, such decomposition is not always possible. 

We equip the set $\mathbb N^{\oplus\infty}$ with the lexicographic order. This is a total order. For any $\delta\in\mathbb N$, the $\mathrm{GL}(V)\times\mathrm{GL}(W)$-module $S^\delta(V\otimes W)$ admits a decreasing filtration indexed by $\mathscr P_\delta$
such that the sub-quotient indexed by $\lambda$ is isomorphic to $L^{\lambda}(V)\otimes  L^{\lambda}(W)$. In particular, $S^{\delta}(V\otimes W)$ admits a sub-$\mathrm{GL}(V)\times\mathrm{GL}(W)$-module which is isomorphic to \[L^{(\delta)}(V)\otimes L^{(\delta)}(W)=\Lambda^\delta(V)\otimes\Lambda^\delta(W)\] This result is called the \emph{Cauchy decomposition formula}\index{Cauchy decomposition formula} for symmetric power. We refer the readers to \cite[Theorem III.1.4]{MR658729} for more details. 

\subsection{Case of several modules}\label{Subsec: several modules}

In this subsection, we apply Cauchy decomposition formula to several $k$-modules. We first illustrate the case of three modules. Let $V_1$, $V_2$ and $V_3$ be three free $k$-modules of finite rank and $\delta$ be a natural number. By Cauchy decomposition formula, the symmetric power $S^\delta(V_1\otimes V_2\otimes V_3)$ admits a decreasing filtration of sub-$\mathrm{GL}(V_1)\times\mathrm{GL}(V_2)\times\mathrm{GL}(V_3)$-modules indexed by 
\[\mathscr P_\delta=\{\lambda\in\mathbb N^{\oplus\infty}\,:\,|\lambda|=\delta\}\]
such that the subquotient indexed by $\lambda$ of the filtration is isomorphic to 
\[L^\lambda(V_1)\otimes L^{\lambda}(V_2\otimes V_3).\] 
Let $\widetilde{\lambda}=(\widetilde{\lambda}_1,\ldots,\widetilde{\lambda}_q,0,\ldots,0,\ldots)$ be the transpose of $\lambda$. By definition $L^\lambda(V_2\otimes V_3)$ is a sub-$\operatorname{GL}(V_2)\times\operatorname{GL}(V_3)$-module of 
\begin{equation}\label{Equ: subquotient lambda tilde}S^{\widetilde{\lambda}_1}(V_2\otimes V_3)\otimes\cdots\otimes S^{\widetilde{\lambda}_q}(V_2\otimes V_3).\end{equation}
We now apply Cauchy decomposition formula to each of the tensor powers $S^{\widetilde{\lambda}_i}(V_2\otimes V_3)$. By passing to the tensor product of filtrations, we obtain a decreasing filtration of \eqref{Equ: subquotient lambda tilde} indexed by 
$\mathscr P_{\widetilde{\lambda}_1}\times\cdots\times\mathscr P_{\widetilde{\lambda}_q}$ (equipped with the lexicographic order), such that the sub-quotient indexed by 
\[(\mu_1,\ldots,\mu_q)\in\mathscr P_{\widetilde{\lambda}_1}\times\cdots\times\mathscr P_{\widetilde{\lambda}_q}\]
is isomorphic to 
\begin{equation}\label{Equ: subquotient of tensor product}L^{\mu_1}(V_2)\otimes\cdots\otimes L^{\mu_q}(V_2)\otimes L^{\mu_1}(V_3)\otimes\cdots\otimes L^{\mu_q}(V_3).\end{equation}
By combining all non-zero coordinates of $\widetilde{\mu}_1,\ldots,\widetilde{\mu}_q$ into a single partition, we obtain a partition $\eta$ and can identify \eqref{Equ: subquotient of tensor product} with a sub-$\operatorname{GL}(V_2)\times\operatorname{GL}(V_3)$-module of $S^\eta(V_2)\otimes S^\eta(V_3)$.
This filtration induces by restriction a decreasing filtration on $L^\lambda(V_2\otimes V_3)$. The sub-quotient of the latter indexed by $(\mu_1,\ldots,\mu_q)$ identifies with a sub-$\operatorname{GL}(V_2)\times\operatorname{GL}(V_3)$-module of \eqref{Equ: subquotient of tensor product}. 
By induction we obtain the following result.

\begin{prop}\label{Pro: filtration Theta delta}
Let $d\in\mathbb N_{\geqslant 2}$ and $(V_i)_{i=1}^d$ be a family of free $k$-modules of finite rank and $\delta$ be a natural number. Let $V$ be the tensor product $V_1\otimes\cdots\otimes V_d$ and $G=\operatorname{SL}(V_1)\times\cdots\times\operatorname{SL}(V_d)$. There exist a finite totally ordered set $\Theta_{\delta,d}$, a map 
\[h=(h_1,\ldots,h_d):\Theta_{\delta,d}\longrightarrow\mathscr P_\delta^d,\]
and a decreasing $\Theta_{\delta,d}$-filtration of $S^\delta(V)$
such that the subquotient indexed by $a\in\Theta_{\delta,d}$ is isomorphic to a sub-$G$-module of \[S^{h_1(a)}(V_1)\otimes \cdots\otimes S^{h_d(a)}(V_d).\]   
\end{prop}
\begin{proof}
We reason by induction on $d$. The case where $d=2$ comes from Cauchy decomposition formula. Assume that $d\geqslant 3$ and that the proposition has been proved for $d-1$ free $k$-modules of finite rank. We apply the induction hypothesis to $V_1,\ldots,V_{d-2}$ and $V_{d-1}\otimes V_d$ to obtain a finite totally ordered set $\Theta_{\delta,d-1}$ and a map 
\[h=(h_1,\ldots,h_{d-1}):\Theta_{\delta,d-1}\longrightarrow\mathscr P_\delta^{d-1}\]
together with a decreasing $\Theta_{\delta,d-1}$-filtration of $S^\delta(V)$
such that the subquotient indexed by $a\in\Theta_{\delta,d-1}$ is isomorphic to a sub-$G$-module of \[S^{h_1(a)}(V_1)\otimes \cdots\otimes S^{h_{d-2}(a)}(V_{d-2})\otimes S^{h_{d-1}(a)}(V_{d-1}\otimes V_d).\] 
Assume that $h_{d-1}(a)$ is of the form $(\lambda_1,\ldots,\lambda_p,0,\ldots,0,\ldots)$. We apply Cauchy decomposition formula to $S^{\lambda_i}(V_{d-1}\otimes V_d)$ to obtain a decreasing filtration of 
\[S^{h_{d-1}(a)}(V_{d-1}\otimes V_d)\]
indexed by $\mathscr P_{\lambda_1}\times\cdots\times\mathscr P_{\lambda_p}$ such that the sub-quotient indexed by 
\[(\mu_1,\ldots,\mu_p)\in\mathscr P_{\lambda_1}\times\cdots\times\mathscr P_{\lambda_p}\]
is isomorphic to $L^{\mu_1}(V_{d-1})\otimes\cdots\otimes L^{\mu_p}(V_{d-1})\otimes L^{\mu_1}(V_d)\otimes\cdots\otimes L^{\mu_p}(V_d)$, which identifies a sub-$\operatorname{GL}(V_{d-1})\times \operatorname{GL}(V_{d})$-module of some $S^{\eta}(V_{d-1})\otimes S^{\eta}(V_d)$, where $\eta$ is a partition of weight $\delta$. In this way we obtain a refinement of the $\Theta_{\delta,d-1}$-filtration of $S^\delta(V)$ which satisfies the required property. The proposition is thus proved. 
\end{proof}

\begin{theo}\label{Thm: invariant element in sym prod}
We keep the notation and the assumptions of Proposition \ref{Pro: filtration Theta delta}, and assume in addition that $k$ is a field. For any $i\in\{1,\ldots,d\}$, let $r_i$ be the dimension of $V_i$ over $k$. If the space $S^\delta(V)^G$ of $G$-invariant vectors in $S^\delta(V)$ is non-zero, then $\delta$ is divisible by $\operatorname{lcm}(r_1,\ldots,r_d)$. Moreover, $S^\delta(V)^G$ identifies with the image of the following $k$-linear map
\begin{equation}\label{Equ: invariant vectors description}\bigoplus_{(\sigma_1,\ldots,\sigma_d)\in\mathfrak S_{\delta}^d}\det(V_1)^{\otimes(\delta/r_1)}\otimes\cdots\otimes\det(V_d)^{\otimes(\delta/r_d)}\longrightarrow S^\delta(V),\end{equation}
where for each $(\sigma_1,\ldots,\sigma_d)\in\mathfrak S_\delta^d$, we consider the composed map
{\small \[\xymatrix{\det(V_1)^{\otimes(\delta/r_1)}\otimes\cdots\otimes\det(V_d)^{\otimes(\delta/r_d)}\ar[r]& V_1^{\otimes\delta}\otimes\cdots\otimes V_d^{\otimes\delta}\ar[d]^-{\sigma_1\otimes\cdots\otimes\sigma_d}\\
&V_1^{\otimes\delta}\otimes\cdots\otimes V_d^{\otimes\delta}\cong V^{\otimes\delta}\ar[r]&S^\delta(V)}\]}%
where the first arrow is induced by the anti-symmetrization map.
\end{theo}
\begin{proof}
By Proposition \ref{Pro: filtration Theta delta}, there exists a finite totally ordered  set $\Theta_{\delta,d}$, a map 
\[h=(h_1,\ldots,h_d):\Theta_{\delta,d}\longrightarrow\mathscr P_\delta^d,\]
and a decreasing $\Theta_{\delta,d}$-filtration $\mathcal F$ of $S^\delta(V)$
such that the subquotient indexed by $a\in\Theta_{\delta,d}$ is isomorphic to a sub-$G$-module of \[S^{h_1(a)}(V_1)\otimes \cdots\otimes S^{h_d(a)}(V_d).\] 
Let $s$ be a non-zero element of $S^\delta(V)^G$, $a$ be the greatest element of $\Theta_{\delta,d}$ such that $s\in\mathcal F^a(S^\delta(V))$. Let $\operatorname{sq}^a(S^\delta(V))$ be the subquotient of the filtration $\mathcal F$ at $a$. By definition, the canonical image of $s$ in $\operatorname{sq}^a(S^\delta(V))$ is a non-zero element of $\operatorname{sq}^a(S^\delta(V))^G$, which is contained in \[\big(S^{h_1(a)}(V_1)\otimes \cdots\otimes S^{h_d(a)}(V_d)\big)^G=S^{h_1(a)}(V_1)^{\operatorname{SL}(V_1)}\otimes\cdots\otimes S^{h_1(a)}(V_d)^{\operatorname{SL}(V_d)}.\]
By Remark \ref{Rem: invariant in symmetric power}, there exists an element $s'$ in the image of \eqref{Equ: invariant vectors description} such that $s-s'$ belongs to 
\[\bigcup_{b\in\Theta_{\delta,d},\;b>a}\mathcal F^b(S^{\delta}(V)).\]
Iterating this procedure we obtain that $s$ actually belongs to the image of \eqref{Equ: invariant vectors description}. The assertion is thus proved.
\end{proof}

\section{Symmetric power norm}

Throughout the section, we let $k$ be a field and $V$ be a finite-dimensional vector space over $k$. We assume that the field $k$ is equipped with an absolute value $|\ndot|$ such that $k$ is complete with respect to the topology induced by $|\ndot|$. We also assume that the vector space $V$ is equipped with a norm $\|\ndot\|$, which is either ultrametric (when $|\ndot|$ is non-Archimedean) or induced by an inner product  (when $|\ndot|$ is Archimedean). We denote by $\|\ndot\|_*$ the dual norm of $\|\ndot\|$ on the dual vector space $V^\vee$. Recall that the norm $\|\ndot\|_*$ is defined as 
\[\forall\, f\in V^\vee,\quad \|f\|_*=\sup_{x\in V}\frac{|f(x)|}{\|x\|}.\]
It is also ultrametric or induced by an inner product.

\subsection{Orthogonal basis}

Let $\alpha$ be an element of $\mathopen{]}0,1\mathclose{]}$. We say that a basis $(e_i)_{i=1}^d$ of $V$ is \emph{$\alpha$-orthogonal}\index{orthogonal!$\alpha$- ---} if the following inequality holds:
\[\forall\,(\lambda_1,\ldots,\lambda_d)\in k^n,\quad\|\lambda_1e_1+\cdots+\lambda_de_d\|\geqslant \alpha\max_{i\in\{1,\ldots, d\}}|\lambda_i|\cdot\|e_i\|.\]
A $1$-orthogonal basis is also called an \emph{orthogonal basis}\index{orthogonal!--- basis}. It is not hard to check that, in the case where $|\ndot|$ is Archimedean (and $\|\ndot\|$ is induced by an inner product $\emptyinnprod$), the orthogonality is equivalent to the usual definition (cf. \cite[Proposition~1.2.3]{CMArakelovAdelic}): a basis $(e_i)_{i=1}^d$ is $1$-orthogonal if and only if 
\[\forall\,(i,j)\in\{1,\ldots,d\}^2,\quad \text{if $i\neq j$ then $\langle e_i,e_j\rangle=0$.}\]   

Assume that $|\ndot|$ is non-Archimedean. Let $\boldsymbol{e}=(e_i)_{i=1}^d$ be a basis of $V$ and let $\|\ndot\|_{\boldsymbol{e}}$ be the norm of $V$ defined as 
\[\forall\,(\lambda_1,\ldots,\lambda_d)\in k^d,\quad \|\lambda_1e_1+\cdots+\lambda_de_d\|_{\boldsymbol{e}}=\max_{i\in\{1,\ldots,d\}}|\lambda_i|\cdot\|e_i\|.\]
Note that $\|\ndot\|_{\boldsymbol{e}}$ is an ultrametric norm of $V$, and one has $\|\ndot\|\leqslant\|\ndot\|_{\boldsymbol{e}}$ (since the norm $\|\ndot\|$ is ultrametric). Moreover, $\boldsymbol{e}$ is an orthogonal of $(V,\|\ndot\|_{\boldsymbol{e}})$. For any $\alpha\in\mathopen{]}0,1\mathclose{]}$ the basis $(e_i)_{i=1}^d$ is $\alpha$-orthogonal with respect to $\|\ndot\|$ if and only if 
\[d(\|\ndot\|,\|\ndot\|_{\boldsymbol{e}}):=\sup_{x\in V\setminus\{0\}}\Big|\ln\|x\|-\ln\|x\|_{\boldsymbol{e}}\Big|\leqslant |\ln(\alpha)|.\]
By the ultrametric Gram-Schimdt procedure (see for example \cite[Proposition 1.2.30]{CMArakelovAdelic}), for any $\alpha\in\mathrm{]}0,1\mathrm{[}$, the ultrametrically normed vector space $(V,\|\ndot\|)$ admits an $\alpha$-orthogonal basis. Therefore, there exists a a sequence of ultrametric norms $(\|\ndot\|_{n\in\mathbb N})$ such that $(V,\|\ndot\|_n)$ admits an orthogonal basis for any $n$, and that 
\[\lim_{n\rightarrow+\infty}d(\|\ndot\|,\|\ndot\|_n)=0.\]

\subsection{Direct sum}

Let $(V_i,\|\ndot\|_i)$, $i\in\{1,\ldots,\delta\}$ be a family of finite-di\-men\-sion\-al normed vector space over $k$.  We assume that, for any $i\in\{1,\ldots,\delta\}$, the norm $\|\ndot\|_i$ is either ultrametric or induced by an inner product. In the case where $|\ndot|$ is non-Archimedean, we equip $V_1\oplus\cdots\oplus V_\delta$ with the \emph{ultrametric direct sum norm}\index{direct sum norm!ultrametric ---}, defined as
\[\forall\, (x_1,\ldots,x_\delta)\in V_1\oplus\cdots\oplus V_{\delta},\quad \|(x_1,\ldots,x_\delta)\|=\max_{i\in\{1,\ldots,\delta\}}\|x_i\|_i.\]
In the case where $|\ndot|$ is Archimedean, we equip $V_1\oplus\cdots\oplus V_\delta$ with the \emph{orthogonal direct sum norm}\index{direct sum norm!orthogonal ---}, such that 
\[\|(x_1,\ldots,x_\delta)\|^2=\sum_{i=1}^\delta\|x_i\|_i^2.\]

\subsection{Symmetric power norm}
Let $\delta$ be a natural number and let 
\[\pi:V^{\otimes\delta}\longrightarrow S^\delta(V)\] 
be the surjective $k$-linear map which sends $x_1\otimes\cdots\otimes x_\delta\in V^{\otimes\delta}$ to $x_1\cdots x_\delta$. We equip $V^{\otimes\delta}$ with the $\varepsilon$-tensor product norm or the orthogonal tensor product norm according to whether $|\ndot|$ is non-Archimedean or Archimedean, respectively. We then equip $S^\delta(V)$ with the quotient norm. 

\begin{prop}\label{Pro: quotient metric on Sdelta}
Assume that the absolute value $|\ndot|$ is non-Archimedean. Let $\alpha\in\mathopen{]}0,1\mathclose{]}$ and let $\boldsymbol{e}=(e_i)_{i=1}^d$ be an $\alpha$-orthogonal basis of $V$. Then the elements 
\[\boldsymbol{e}^{\boldsymbol{a}},\quad \boldsymbol{a}\in\mathbb N^d,\;|\boldsymbol{a}|=\delta\]
form an $\alpha^\delta$-orthogonal basis of $S^\delta(V)$. Moreover, for any $\boldsymbol{a}=(a_1,\ldots,a_d)\in\mathbb N^d$ such that $|\boldsymbol{a}|=\delta$, one has 
\begin{equation}\label{Equ: encadrement norme}\alpha^{\delta}\prod_{i=1}^d\|e_i\|^{a_i}\leqslant\|\boldsymbol{e}^{\boldsymbol{a}}\|\leqslant\prod_{i=1}^d\|e_i\|^{a_i}.\end{equation}
\end{prop}
\begin{proof}
Denote by $f:\{1,\dots,d\}^\delta\rightarrow\mathbb N^d$ the map which sends $(b_1,\ldots,b_\delta)$ to the vector 
\[\Big(\operatorname{card}\big(\{j\in\{1,\ldots,\delta\}\,|\,b_j=i\}\big)\Big)_{i=1}^d.\] 
Let $\pi:V^{\otimes\delta}\rightarrow S^\delta(V)$ be the projection map. For any \[b=(b_1,\ldots,b_\delta)\in\{1,\ldots,d\}^\delta,\] denote by $e_{b}$ the split tensor  $e_{b_1}\otimes\cdots\otimes e_{b_\delta}\in V^{\otimes\delta}$. 

For $\boldsymbol{a}=(a_1,\ldots,a_n)\in\mathbb N^n$ such that $|\boldsymbol{a}|=\delta$, one has
\[\|\boldsymbol{e}^{\boldsymbol{a}}\|=\inf\bigg\{\bigg\|\sum_{b\in f^{-1}(\{\boldsymbol{a}\})}\lambda_be_b\bigg\|\;:\;\sum_{b\in f^{-1}(\{\boldsymbol{a}\})}\lambda_b=1,\;\bigg\}.\]
Hence (see \cite[Remark~1.1.56]{CMArakelovAdelic})
\[\|\boldsymbol{e}^{\boldsymbol{a}}\|\leqslant\|e_1\|^{a_1}\cdots\|e_n\|^{a_n}.\]
Since $(e_i)_{i=1}^d$ is an $\alpha$-orthogonal basis, $(e_b)_{b\in\{1,\ldots,d\}^\delta}$ is an $\alpha^\delta$-orthogonal basis of $V^{\otimes\delta}$ (see \cite[Proposition 1.2.19]{CMArakelovAdelic}). For any $(\lambda_b)_{b\in f^{-1}(\{\boldsymbol{a}\})}\in k^{f^{-1}(\{\boldsymbol{a}\})}$ such that 
\[\sum_{b\in f^{-1}(\{\boldsymbol{a}\})}\lambda_b=1,\]
one has 
\[\|e_1\|^{i_1}\cdots\|e_d\|^{i_d}\leqslant \|e_1\|^{i_1}\cdots\|e_d\|^{i_d}\max_{b\in f^{-1}(\{\boldsymbol{a}\})}|\lambda_b|\leqslant \alpha^{-\delta}\bigg\|\sum_{b\in f^{-1}(\{\boldsymbol{a}\})}\lambda_be_b\bigg\|,\]
which leads to $\|\boldsymbol{e}^{\boldsymbol{a}}\|\geqslant\alpha^{-\delta}\|e_1\|^{a_1}\cdots\|e_d\|^{a_d}$.

For any
\[s=\sum_{b\in\{1,\ldots,n\}^\delta}\mu_b e_{b}\in E^{\otimes\delta},\]
one has
\[\pi(s)=\sum_{\boldsymbol{a}\in\mathbb N^n,\,|\boldsymbol{a}|=\delta}\bigg(\sum_{b\in f^{-1}(\{\boldsymbol{a}\})}\mu_b\bigg)\boldsymbol{e}^{\boldsymbol{a}}.\]
Moreover,
\[\begin{split}\|s\|&\geqslant\alpha^\delta\max_{\begin{subarray}{c}\boldsymbol{a}=(a_1,\ldots,a_d)\in\mathbb N^d\\|\boldsymbol{a}|=\delta\end{subarray}}\|e_{1}\|^{a_1}\cdots\|e_d\|^{a_d}\max_{b\in f^{-1}(\{\boldsymbol{a}\})}|\mu_b|\\
&\geqslant\alpha^\delta\max_{\begin{subarray}{c}\boldsymbol{a}=(a_1,\ldots,a_d)\in\mathbb N^d\\|\boldsymbol{a}|=\delta\end{subarray}}\|e_{1}\|^{a_1}\cdots\|e_d\|^{a_d}\bigg|\sum_{b\in f^{-1}(\{\boldsymbol{a}\})}\mu_b\bigg|\\
&\geqslant \alpha^\delta\max_{\begin{subarray}{c}\boldsymbol{a}=(a_1,\ldots,a_d)\in\mathbb N^d\\|\boldsymbol{a}|=\delta\end{subarray}}\|\boldsymbol{e}^{\boldsymbol{a}}\|\cdot\bigg|\sum_{b\in f^{-1}(\{\boldsymbol{a}\})}\mu_b\bigg|.
\end{split}\]
Therefore, we obtain that $(e^{\boldsymbol{a}})_{\boldsymbol{a}\in\mathbb N^d,\,|\boldsymbol{a}|=\delta}$ forms an $\alpha^\delta$-orthogonal basis of $S^\delta(V)$, as required.
\end{proof}

\begin{rema}\label{Rem: filtration on symmetric product}
Consider the case where $|\ndot|$ is the trivial absolute value. In this case, the ultrametric norm $\|\ndot\|$ corresponds to a decreasing $\mathbb R$-filtration $\mathcal F$ on $V$ such that 
\[\forall\,t\in\mathbb R,\quad \mathcal F^{t}(V)=\{x\in V\,:\,\|x\|\leqslant\mathrm{e}^{-t}\}.\]
We can also express this $\mathbb R$-filtration as an increasing sequence
\[0=V_0\subsetneq V_1\subsetneq\ldots\subsetneq V_r=V\]
together with a decreasing sequence
\[\mu_1>\ldots>\mu_r,\]
with $\mathcal F^t(V)=V_i$ when $t\in \mathopen{]}\mu_{i+1},\mu_i\mathclose{]}\cap\mathbb R$, where by convention $\mu_0=+\infty$ and $\mu_{r+1}=-\infty$. For any $t\in\mathbb R$, the subquotient
\[\operatorname{sq}^t(V):=\mathcal F^{t}(V)\bigg/\bigcup_{\varepsilon>0}\mathcal F^{t+\varepsilon}(V)\]
is either the zero vector space when $t\not\in\{\mu_1,\ldots,\mu_r\}$, or is equal to $V_{i}/V_{i-1}$ when $t=\mu_i$. 

By \cite[Proposition 1.2.30]{CMArakelovAdelic}, there exists an orthogonal basis $\boldsymbol{e}$ such that $\boldsymbol{e}\cap \mathcal F^t(V)$ forms a basis of $\mathcal F^t(V)$ for any $t\in\mathbb R$. By Proposition \ref{Pro: quotient metric on Sdelta}, we obtain that the elements 
\[\boldsymbol{e}^{\boldsymbol{a}},\quad \boldsymbol{a}\in\mathbb N^d,\;|\boldsymbol{a}|=\delta\]
form an orthogonal basis of $S^\delta(V)$. Moreover, for any $\boldsymbol{a}=(a_1,\ldots,a_n)\in\mathbb N^n$ such that $|\boldsymbol{a}|=\delta$, one has 
\[\|\boldsymbol{e}^{\boldsymbol{a}}\|=\prod_{i=1}^d\|e_i\|^{a_i}.\]
Therefore, if we equip $S^\delta(V)$ with the $\mathbb R$-filtration induced by the symmetric product norm of $\|\ndot\|$, for any $t\in\mathbb R$ one has a natural isomorphism
\[\operatorname{sq}^t(S^\delta(V))\cong\bigoplus_{\begin{subarray}{c}\boldsymbol{b}=(b_1,\ldots,b_r)\in\mathbb N^r\\
|\boldsymbol{b}|=b_1+\cdots+b_r=\delta\\
b_1\mu_1+\cdots+b_r\mu_r=t
\end{subarray}}S^{b_1}(V_1/V_0)\otimes\cdots\otimes S^{b_r}(V_r/V_{r-1}).\]
\end{rema}

\subsection{Subquotient metric on symmetric power} \label{Subsec: subquotient on symmetric power}

In this subsection, we assume that $|\ndot|$ is non-Archimedean and $\|\ndot\|$ is  ultrametric. We let $|\ndot|_0$ be the trivial absolute value on $k$ and $\|\ndot\|_0$ be an ultrametric norm on $V$ (with respect to the trivial absolute value), which corresponds to an $\mathbb R$-filtration $\mathcal F$ on $V$,  or an increasing sequence
\[0=V_0\subsetneq V_1\subsetneq\ldots\subsetneq V_r=V\]
of vector subspaces of $V$ together with a decreasing sequence
\[\mu_1>\ldots>\mu_r\]
of real numbers, as explained in Remark \ref{Rem: filtration on symmetric product}. 

Let $\delta$ be a natural integer. We equip $S^\delta(V)$ with the $\mathbb R$-filtration corresponding to the symmetric product norm of $\|\ndot\|_0$. As we have seen in Remark \ref{Rem: filtration on symmetric product}, for any $t\in\mathbb R$, the subquotient $\operatorname{sq}^t(S^\delta(V))$ is isomorphic to 
\begin{equation}\label{Equ: tensor product of Sa }\bigoplus_{\begin{subarray}{c}\boldsymbol{b}=(b_1,\ldots,b_r)\in\mathbb N^r\\
|\boldsymbol{b}|=b_1+\cdots+b_r=\delta\\
b_1\mu_1+\cdots+b_r\mu_r=t
\end{subarray}}S^{b_1}(V_1/V_0)\otimes\cdots\otimes S^{b_r}(V_r/V_{r-1}).\end{equation}
Note that the norm $\|\ndot\|$ induces by passing to subquotient a norm on each $V_i/V_{i-1}$, which leads to a symmetric power norm on $S^b(V_i/V_{i-1})$ for any $b\in\mathbb N$. For $(b_1,\ldots,b_r)\in\mathbb N^r$, we equip $S^{b_1}(V_1/V_0)\otimes\cdots\otimes S^{b_r}(V_r/V_{r-1})$ with the tensor product of symmetric power norms ($\varepsilon$-tensor product when $|\ndot|$ is non-Archimedean and orthogonal tensor product when $|\ndot|$ is Archimedean), and the vector space \eqref{Equ: tensor product of Sa } with the direct sum norm (ultrametric direct sum if $|\ndot|$ is non-Archimedean and orthogonal direct sum if $|\ndot|$ is Archimedean).

Here we are interested in the comparison between the subquotient norm on $\operatorname{sq}^t(S^\delta(V))$ induced by the symmetric tensor power norm and the direct sum of tensor product norm on \eqref{Equ: tensor product of Sa } described above.

\begin{prop}\label{Pro:subquotient of filtration}
Assume that the absolute value $|\ndot|$ is non-Archimedean. Then, for any $t\in\mathbb R$ the isomorphism 
\begin{equation}\label{Equ: sq t iso}\operatorname{sq}^t(S^\delta(V))\cong \bigoplus_{\begin{subarray}{c}\boldsymbol{b}=(b_1,\ldots,b_r)\in\mathbb N^r\\
|\boldsymbol{b}|=b_1+\cdots+b_r=\delta\\
b_1\mu_1+\cdots+b_r\mu_r=t
\end{subarray}}S^{b_1}(V_1/V_0)\otimes\cdots\otimes S^{b_r}(V_r/V_{r-1})\end{equation}
is an isometry.
\end{prop}
\begin{proof}
Let $\alpha$ be an element of $\mathopen{]}0,1\mathclose{[}$. By \cite[Proposition 1.2.30]{CMArakelovAdelic}, for any $i\in\{1,\ldots,r\}$ there exists 
\[\boldsymbol{e}^{(i)}=(e^{(i)}_1,\ldots,e^{(i)}_{d_i})\in (V_i\setminus V_{i-1})^{d_i}\]
such that 
\begin{enumerate}[label=\rm(\alph*)]
\item the images of $e^{(i)}_1,\ldots,e^{(i)}_{d_i}$ in $V_i/V_{i-1}$ form a basis of the latter, where $d_i=\dim_k(V_i/V_{i-1})$,
\item $(\boldsymbol{e}^{(1)},\ldots,\boldsymbol{e}^{(r)})$ forms an $\alpha$-orthogonal basis of $V$.
\end{enumerate} 
By Proposition \ref{Pro: quotient metric on Sdelta}, the elements
\[(\boldsymbol{e}^{(1)})^{\boldsymbol{a}^{(1)}}\cdots (\boldsymbol{e}^{(r)})^{\boldsymbol{a}^{(r)}},\quad (\boldsymbol{a}^{(1)},\ldots,\boldsymbol{a}^{(r)})\in\mathbb N^{d_1}\times\cdots\times\mathbb N^{d_r},\;\big|\boldsymbol{a}^{(1)}\big|+\cdots\big|\boldsymbol{a}^{(r)}\big|=\delta\]
form an $\alpha^\delta$-orthogonal basis of $S^\delta(V)$. We let $\boldsymbol{x}^{(i)}=(x_1^{(i)},\ldots,x_{d_i}^{(i)})\in (V_i/V_{i-1})^{d_i}$, where $x_j^{(i)}$ denotes the class of $e_j^{(i)}$ in $V_i/V_{i-1}$. Since $(\boldsymbol{e}^{(1)},\ldots,\boldsymbol{e}^{(r)})$ forms an $\alpha$-orthogonal basis of $V$, we obtain that \begin{equation}\label{Equ: estimation norm xji}\alpha\|e_j^{(i)}\|\leqslant\|x_j^{(i)}\|\leqslant\|e_j^{(i)}\|.\end{equation} For any $\boldsymbol{b}=(b_1,\ldots,b_r)\in\mathbb N^r$ such that $b_1+\cdots+b_r=\delta$, the vectors
{\footnotesize \[(\boldsymbol{x}^{(1)})^{\boldsymbol{a}^{(1)}}\otimes\cdots\otimes(\boldsymbol{x}^{(r)})^{\boldsymbol{a}^{(r)}},\quad (\boldsymbol{a}^{(1)},\ldots,\boldsymbol{a}^{(r)})\in\mathbb N^{d_1}\times\cdots\times\mathbb N^{d_r},\;\forall\,i\in\{1,\ldots,r\},\,\big|\boldsymbol{a}^{i}\big|=b_i,\]}%
form an $\alpha^{\delta}$-orthogonal basis of \[S^{b_1}(V_1/V_0)\otimes\cdots\otimes S^{b_r}(V_{r}/V_{r-1}).\]
Moreover, by \eqref{Equ: encadrement norme} and \eqref{Equ: estimation norm xji} we obtain that 
\[\left|\ln\frac{\big\|(\boldsymbol{e}^{(1)})^{\boldsymbol{a}^{(1)}}\cdots (\boldsymbol{e}^{(r)})^{\boldsymbol{a}^{(r)}}\big\|}{\|(\boldsymbol{x}^{(1)})^{\boldsymbol{a}^{(1)}}\otimes\cdots\otimes(\boldsymbol{x}^{(r)})^{\boldsymbol{a}^{(r)}}\|}\right|\leqslant 2\delta|\ln(\alpha)|.\]
Therefore, under the isomorphism \eqref{Equ: sq t iso}, the distance between the norms on the left hand side and the right hand side is bounded from above by $3\delta|\ln(\alpha)|$. Since $\alpha\in\mathopen{]}0,1\mathclose{[}$ is arbitrary, we obtain that \eqref{Equ: sq t iso} is actually an isometry.
\end{proof}

\subsection{Symmetric tensor}
 Let $\delta$ be a positive integer. We denote by \begin{equation*}\operatorname{sym}:V^{\otimes\delta}\longrightarrow V^{\otimes\delta}\end{equation*}
sending $x_1\otimes\cdots\otimes x_\delta$ to 
\[\sum_{\sigma\in\mathfrak S_{\delta}}x_{\sigma(1)}\otimes\cdots\otimes x_{\sigma(\delta)}.\]

\begin{prop}\label{Pro: norm antisym p}Assume that the absolute value $|\ndot|$ is non-Archimedean. 
Let $\delta\in\mathbb N_{\geqslant 1}$. We equip $V^{\otimes\delta}$ with the $\varepsilon$-tensor power norm of $\|\ndot\|$. Then the $K$-linear map $\operatorname{sym}:V^{\otimes\delta}\longrightarrow V^{\otimes\delta}$ has operator norm $\leqslant 1$. 
\end{prop}
\begin{proof}
Let $T$ be an element of $V^{\otimes\delta}$. If we consider $T$ as a $\delta$-multilinear form on $V^\vee$, then the $\varepsilon$-tensor power norm of $T$ is given by 
\[\|T\|_{\varepsilon}=\sup_{(\alpha_1,\ldots,\alpha_\delta)\in (V^\vee\setminus\{0\})^\delta}\frac{|T(\alpha_1,\ldots,\alpha_\delta)|}{\|\alpha_1\|_*\cdots\|\alpha_\delta\|_*}.\]
Note that the element $\operatorname{sym}(T)$, viewed as a $\delta$-multilinear form on $V^\vee$, is given by 
\[\operatorname{sym}(T)(\alpha_1,\ldots,\alpha_\delta)=\sum_{\sigma\in\mathfrak S_\delta}T(\alpha_{\sigma(1)},\ldots,\alpha_{\sigma(\delta)}).\]
Since the absolute value $|\ndot|$ is non-Archimedean, we obtain 
\[|\operatorname{sym}(T)(\alpha_1,\ldots,\alpha_\delta)|\leqslant\max_{\sigma\in\mathfrak S_\delta}|T(\alpha_{\sigma(1)},\ldots,\alpha_{\sigma(\delta)})|\leqslant\|T\|_\varepsilon\cdot\|\alpha_1\|_{*}\cdots\|\alpha_\delta\|_*,\]
which shows $\|\operatorname{sym}(T)\|_{\varepsilon}\leqslant\|T\|_{\varepsilon}$. 
\end{proof}

\begin{prop}
\label{Pro: norm antisym 0}Assume that the absolute value $|\ndot|$ is Archimedean. 
Let $\delta\in\mathbb N_{\geqslant 1}$. We equip $V^{\otimes\delta}$ with the orthogonal tensor power norm
of $\|\ndot\|$. Then the $K$-linear map $\operatorname{sym}:V^{\otimes\delta}\rightarrow V^{\otimes\delta}$
has operator norm $\leqslant \delta!$. 
\end{prop}
\begin{proof}
Let $(e_j)_{j=1}^d$ be an orthonormal basis of $V$. Recall that an orthonormal basis of $(V^{\otimes\delta},\|\ndot\|)$ is given by 
\[e_{j_1}\otimes\cdots\otimes e_{j_\delta},\quad (j_1,\ldots,j_\delta)\in\{1,\ldots,d\}^\delta.\] Let 
\[T=\sum_{\lambda=(\lambda_1,\ldots,\lambda_\delta)\in\{1,\ldots,d\}^\delta}a_\lambda e_{\lambda_1}\otimes\cdots\otimes e_{\lambda_\delta}\in V^{\otimes\delta}.\]
One has
\[\|T\|^2=\sum_{\lambda\in\{1,\ldots,d\}^\delta}|a_\lambda|^2.\]
Let $\mathscr P_{\delta,d}$ be the set of vectors $(a_1,\ldots,a_d)\in\mathbb N^d$ such that $a_1+\cdots+a_d=\delta$. For each $\boldsymbol{a}=(a_1,\ldots,a_d)\in S_\delta$, let $I_{\boldsymbol{a}}$ be the set of $(\lambda_1,\ldots,\lambda_\delta)\in\{1,\ldots,d\}^\delta$ such that 
\[\forall\,j\in\{1,\ldots,d\},\quad a_j=\sum_{\begin{subarray}{c}
i\in\{1,\ldots,\delta\}\\
\lambda_i=j
\end{subarray}}1.\]
Then the following equality holds
\[\operatorname{sym}(T)=\sum_{\boldsymbol{a}=(a_1,\ldots,a_d)\in \mathscr P_{\delta,d}}\bigg(\sum_{\lambda\in I_{\boldsymbol{a}}}a_\lambda\bigg)a_1!\cdots a_d!\sum_{\lambda=(\lambda_1,\ldots,\lambda_\delta)\in I_{\boldsymbol{a}}}e_{\lambda_1}\otimes\cdots\otimes e_{\lambda_\delta},\]
which leads to 
\[\begin{split}\|\operatorname{sym}(T)\|^2&=\sum_{\boldsymbol{a}=(a_1,\ldots,a_d)\in \mathscr P_{\delta,d}}(a_1!\cdots a_d!)^2\bigg|\sum_{\lambda\in I_{\boldsymbol{a}}}a_\lambda\bigg|^2\frac{\delta!}{a_1!\cdots a_d!}\\
&\leqslant\sum_{\boldsymbol{a}=(a_1,\ldots,a_d)\in \mathscr P_{\delta,d}}(a_1!\cdots a_d!)^2\bigg(\sum_{\lambda\in I_{\boldsymbol{a}}}|a_\lambda|^2\bigg)\bigg(\frac{\delta!}{a_1!\cdots a_d!}\bigg)^2\\
&=(\delta!\|T\|)^2.
\end{split}\]
\end{proof}

\begin{rema}\label{Rem: operator norm symmetri cpower}
Note that the $k$-linear map $\operatorname{sym}:V^{\otimes\delta}\longrightarrow V^{\otimes\delta}$ factors though the symmetric power $S^\delta(V)$. Moreover, in the case where $k$ is  of characteristic $0$, the unique $k$-linear map $\operatorname{sym}':S^{\delta}(V)\rightarrow V^{\otimes \delta}$  such that the composition
\[\xymatrix{V^{\otimes\delta}\ar@{->>}[r]&S^\delta(V)\ar[r]^-{\operatorname{sym}'}&V^{\otimes\delta}}\]
identifies with $\operatorname{sym}:V^{\otimes\delta}\rightarrow V^{\otimes\delta}$ is injective.  The above propositions show that, if we equip $V^{\otimes\delta}$ with the $\varepsilon$-tensor power norm (resp. orthogonal tensor power norm) of $\|\ndot\|$ in the case where $|\ndot|$ is non-Archimedean (resp. Archimedean) and equip $S^\delta(V)$ with the quotient norm, then the operator norm of $\operatorname{sym}'$ is bounded from above by $1$ (resp. $\delta!$).  
\end{rema}

\subsection{Determinant norm}\label{subsection:determinant:norm}
Recall that we have fixed a finite-dimension normed vector space $(V,\|\ndot\|)$ over $k$. Let $r$ be the dimension of $V$ over $k$. We denote by $\det(V)$ the exterior power $\Lambda^r(V)$. This is a one-dimensional vector space over $k$. We equip it with the \emph{determinant norm}\index{determinant norm} $\|\ndot\|_{\det}$, which is defined as 
\[\forall\,\eta\in\det(V),\quad \|\eta\|_{\det}=\inf_{\begin{subarray}{c}(x_i)_{i=1}^r\in V^r\\
\eta=x_1\wedge\cdots\wedge x_r\end{subarray}}\|x_1\|\cdots\|x_r\|.\]

\begin{prop}
Assume that the absolue value $|\ndot|$ is non-Archimedean and the norm $\|\ndot\|$ is ultrametric. Then the anti-symmetrization map
$\det(V)\rightarrow V^{\otimes r}$ is an isometry from $(\det(V),\|\ndot\|_{\det})$ to its image, where we consider the $\varepsilon$-tensor power norm on $V^{\otimes r}$. 
\end{prop} 
\begin{proof}
By the ultrametric Gram-Schmidt procedure, one can approximate the norm $\|\ndot\|$ by a sequence of norms for which $V$ admits an orthogonal basis. Therefore, we may assume without loss of generality that $(V,\|\ndot\|)$ has an orthogonal basis $(e_i)_{i=1}^r$. By \cite[Proposition 1.2.25]{CMArakelovAdelic}, one has 
\[\|e_1\wedge\cdots\wedge e_r\|=\|e_1\|\cdots\|e_r\|.\]
Moreover, the anti-symmetrization of $e_1\wedge\cdots\wedge e_r$ is given by 
\begin{equation}\label{Equ: anti-sym of det}\sum_{\sigma\in\mathfrak S_r}\operatorname{sgn}(\sigma)e_{\sigma(1)}\otimes\cdots\otimes e_{\sigma(r)}.\end{equation}
Since 
\[\{e_{\sigma(1)}\otimes\cdots\otimes e_{\sigma(r)}\,:\,\sigma\in\mathfrak S_r\}\]
is a subset of the orthogonal basis 
\[e_{j_1}\otimes\cdots\otimes e_{j_r},\quad (j_1,\ldots,j_r)\in \{1,\ldots, r\}^r,\]
we obtain that the norm of \eqref{Equ: anti-sym of det} is equal to $\|e_1\|\cdots\|e_r\|$. The proposition is thus proved. 
\end{proof}

\begin{rema}
In the case where $|\ndot|$ is Archimedean and $\|\ndot\|$ is induced by an inner product, the above result is no longer true. The orthogonal tensor power norm of the anti-symmetrization of an element $\eta$ of $\det(V)$ is equal to $\sqrt{r!}\,\|\eta\|_{\det}$.
\end{rema}

\section{Maximal slopes of symmetric power}

In this section, we prove that, on an adelic curve of perfect underlying field and without non-Archimedean places, the tensor product of semi-stable Hermitian adelic vector bundles remains semi-stable. This allows to justify that the argument used in the proof of \cite[Proposition 5.3.1]{CMArakelovAdelic} is still valid in the positive characteristic case.

We fix a proper adelic curve $S=(K,(\Omega,\mathcal A,\nu),\phi)$ with a perfect underlying field $K$. We assume in addition that, either the $\sigma$-algebra $\mathcal A$ is discrete, or the field $K$ is countable. 

\subsection{Tensorial semi-stability}

In this subsection, we assume that $\Omega_\infty$ is empty, namely the absolute value $|\ndot|_{\omega}$ is non-Archimedean for any $\omega\in\Omega$. We let $(\overline E_i)_{i=1}^d$ be a family of Hermitian adelic vector bundles on $S$ and $\overline E$ be the orthogonal tensor product
\[\overline E_1\otimes\cdots\otimes\overline E_d.\]
For any $i\in\{1,\ldots,d\}$, we let $r_i$ be the dimension $E_i$ over $K$.
The purpose is to prove the following estimate.

\begin{theo}
Let $Q$ be a one-dimensional quotient vector space of $E$, equipped with the quotient norm family. Then the following inequality holds:
\begin{equation}\label{Equ: quotient line bundle}\widehat{\deg}(\overline Q)\geqslant\sum_{i=1}^d\widehat{\mu}_{\min}(\overline E_i).\end{equation} 
\end{theo}
\begin{proof}
Let $G$ be the product of special linear group schemes \[\mathbb{SL}(E_1)\times\cdots\times\mathbb{SL}(E_d).\] Note that the algebraic group $G$ acts on the scheme $\mathbb P(E)$ and the tautological line bundle $\mathcal O_{E}(1)$ is naturally equipped with a $G$-linear structure. In particular,  the group
\[G(K)=\operatorname{SL}(E_1)\times \cdots \times\operatorname{SL}(E_d)\]
acts naturally on the sectional $K$-algebra
\[\bigoplus_{n\in\mathbb N}H^0(\mathbb P(E),\mathcal O_E(n))=\bigoplus_{n\in\mathbb N}S^n(E)=\operatorname{Sym}(E).\]
Let $x$ be the rational point of $\mathbb P(E)$ which is represented by the one-dimensional quotient space $Q$.

{\bf Step 1:}  We suppose firstly that $x$ is semi-stable in the sense of geometric invariant theory with respect to the $G$-linear line bundle $\mathcal O_E(1)$. In other words, we assume that there exists a positive integer $\delta$ and a section in $S^\delta(E)=H^0(\mathbb P(E),\mathcal O_E(\delta))$ invariant by the action of $G(K)$, which does not vanish at $x$. By Theorem \ref{Thm: invariant element in sym prod},  we obtain that $\delta$ is divisible by $\operatorname{lcm}(r_1,\ldots,r_d)$ and there exists $(\sigma_1,\ldots,\sigma_d)\in\mathfrak S_\delta^d$ such that the following composed map is non-zero
{\small \[\xymatrix{\det(E_1)^{\otimes(\delta/r_1)}\otimes\cdots\otimes\det(E_d)^{\otimes(\delta/r_d)}\ar[r]& E_1^{\otimes\delta}\otimes\cdots\otimes E_d^{\otimes\delta}\ar[d]^-{\sigma_1\otimes\cdots\otimes\sigma_d}\\
&E_1^{\otimes\delta}\otimes\cdots\otimes E_d^{\otimes\delta}\cong E^{\otimes\delta}\ar[d]\\
&S^\delta(E)\ar[r]&F^{\otimes\delta}}\]}%
Therefore, we obtain
\[\widehat{\deg}(\overline Q)\geqslant\sum_{i=1}^d\frac{\delta}{r_i}\operatorname{\widehat{\deg}}(\overline E_i)=\delta\sum_{i=1}^d\operatorname{\widehat{\mu}}(\overline E_i)\geqslant\delta\sum_{i=1}^d\operatorname{\widehat{\mu}_{\min}}(\overline E_i).\]

{\bf Step 2:} In this step, we assume that $x$ is not semi-stable under the action of $G$ with respect to the $G$-linear line bundle $\mathcal O_E(1)$. Note that this condition is equivalent to the following: $x$ is not semi-stable under the action of 
\[\mathbb{GL}(E_1)\times\cdots\times\mathbb{GL}(E_d)\]
with respect to 
\[\mathcal O_E(r_1\cdots r_d)\otimes\pi^*(\det(E_1^\vee)^{\otimes b_1}\otimes\cdots\otimes\det(E_d^\vee)^{\otimes b_d}),\]
where $\pi:\mathbb P(E)\rightarrow\Spec K$ is the structural morphism and, for any $i\in\{1,\ldots,d\}$,
\[b_i:=\frac{r_1\cdots r_d}{r_i}.\]
Then the inequality \eqref{Equ: quotient line bundle} can be obtained following the same argument as in the proof of \cite[Theorem 5.6.1]{CMArakelovAdelic}.
\end{proof}

\begin{coro}\label{Cor: tensorial semistability}
Let $(\overline E_i)_{i=1}^d$ be a family of Hermitian adelic line bundles on $S$. For any vector subspace $F$ of $E_1\otimes\cdots\otimes E_d$, one has
\[\widehat{\mu}(\overline F)\leqslant\sum_{i=1}^d\widehat{\mu}_{\max}(\overline E_i).\]
In particular, if $\overline E_1,\ldots,\overline E_d$ are all semi-stable, then $\overline E_1\otimes\cdots\otimes \overline E_d$ is also semi-stable.
\end{coro}
\begin{proof}
We first treat the case where $F$ is of dimension $1$. We identify $F^\vee$ with a quotient vector space of $E_1^\vee\otimes\cdots\otimes E_d^\vee$. By Theorem \ref{Equ: quotient line bundle} we obtain 
\[\widehat{\mu}(\overline F)=\widehat{\deg}(\overline F)=-\widehat{\deg}(\overline F^\vee)\leqslant -\sum_{i=1}^d\widehat{\mu}_{\min}(\overline E_i^\vee)=\sum_{i=1}^d\widehat{\mu}_{\max}(\overline E_i),\]
where the last equality comes from \cite[Corollary 4.3.27]{CMArakelovAdelic}.

In the following, we consider the general case. Without loss of generality, we may assume that $F$ is the destabilizing vector subspace of $\overline E_1\otimes\cdots\otimes\overline E_d$. In particular, $\overline F$ is semi-stable. Let $s$ be the element of $F^\vee\otimes E_1\otimes \cdots\otimes E_d$ be the element which corresponds to the inclusion map $f:F\rightarrow E_1\otimes\cdots\otimes E_d$. Let $L$ be the one-dimension vector subspace of $F\rightarrow E_1\otimes\cdots\otimes E_d$ spanned by $s$.  By the one-dimensional case of the statement proved above, one has
\begin{align*}\widehat{\deg}(\overline L) &\leqslant\widehat{\mu}_{\max}(\overline F^\vee)+\sum_{i=1}^d\widehat{\mu}_{\max}(\overline E_d)=\widehat{\mu}(\overline F^\vee)+\sum_{i=1}^d\widehat{\mu}_{\max}(\overline E_d) \\
& =-\widehat{\mu}(\overline F)+\sum_{i=1}^d\widehat{\mu}_{\max}(\overline E_d)\end{align*}
since $\overline F$ is assumed to be semi-stable. Moreover, for any $\omega\in\Omega$, if we denote by $\|\ndot\|_\omega$ the $\varepsilon$-tensor product norm on $F_\omega^\vee\otimes E_{1,\omega}\otimes\cdots\otimes E_{d,\omega}$, then $\|s\|_\omega$ identifies with the operator norm of $f_\omega$, which is bounded from above by $1$. Therefore one has $\widehat{\deg}(\overline L)\geqslant 0$, which shows that 
\[\widehat{\mu}(\overline F)\leqslant\sum_{i=1}^d\widehat{\mu}(\overline E_d),\]
as required.
\end{proof}

\begin{rema}\label{Rem: positive characteristic minimal slope}
By passing to dual, we obtain from Corollary \ref{Cor: tensorial semistability} that the inequality \eqref{Equ: quotient line bundle} actually holds for quotient vector subspace of $\overline E_1\otimes\cdots\otimes\overline E_d$ of arbitrary rank. In other words, the following inequality holds:
\[\widehat{\mu}_{\min}(\overline E_1\otimes\cdots\otimes\overline E_d)\geqslant\sum_{i=1}^d\widehat{\mu}_{\min}(\overline E_i).\]
Therefore, the results of \cite[Chapter 5]{CMArakelovAdelic} still hold in the case where $K$ is a perfect field of positive characteristic. In particular, if $\overline F$ is a vector subspace of $E_1\otimes\cdots\otimes E_d$, equipped with the restriction of the orthogonal tensor product norm family, then the dual statement of \cite[Theorem 5.6.1]{CMArakelovAdelic} leads to the following inequality
\[\widehat{\mu}(\overline F)\leqslant \widehat{\mu}(F,\|\ndot\|_{0,F})+\frac 12\nu(\Omega_\infty)\sum_{j=1}^d\ln(\dim_K(E_j)), \]
where $\|\ndot\|_{0,F}$ denotes the restriction of the ultrametric norm (where we consider the trivial absolute value on $K$) on $E_1\otimes\cdots\otimes E_d$ by taking the $\varepsilon$-tensor product of norms associated with Harder-Narasimhan $\mathbb R$-filtrations of $\overline E_1,\ldots,\overline E_d$. 
\end{rema}

\subsection{Slope of a symmetric power}

\begin{prop}\label{Pro: mu S delta}
Assume that $\Omega_\infty$ is empty. Let $\overline{E}$ be a Hermitian adelic vector bundle on $S$ and $\delta$ be a positive number. The following equality holds
\[\widehat{\mu}(S^\delta(\overline E))=\delta\operatorname{\widehat{\mu}}(\overline E).\]
Moreover, if $\overline E$ is semi-stable, then $S^\delta(\overline E)$ is also semi-stable.
\end{prop}
\begin{proof}Let $r$ be the dimension of $E$ over $K$.
Without loss of generality, we may assume that $r\geqslant 2$. Let 
\[0=E_0\subsetneq E_1\subsetneq \ldots\subsetneq E_r=E\]
be a complete flag of vector subspaces of $E$. By Proposition \ref{Pro:subquotient of filtration}, one has 
\[\begin{split}\widehat{\deg}(S^\delta(\overline E))&=\sum_{\begin{subarray}{c}(a_1,\ldots,a_r)\in\mathbb N^r\\
a_1+\cdots+a_r=\delta\end{subarray}}\sum_{i=1}^ra_i\operatorname{\widehat{\deg}}(\overline{E_{i}/E_{i-1}})\\
&=\sum_{i=1}^r\widehat{\deg}(\overline{E_i/E_{i-1}})\sum_{a=0}^\delta a\binom{r+\delta-a-2}{r-2}.
\end{split}\]
Note that 
\[\begin{split}&\quad\;\sum_{a=0}^\delta a\binom{r+\delta-a-2}{r-2}= \delta\sum_{a=0}^\delta\binom{r+\delta-a-2}{r-2} \\ & {\kern20em-\sum_{a=0}^\delta(\delta-a)\binom{r+\delta-a-2}{\delta-a}}\\
&=\delta\binom{\delta+r-1}{r-1}-\sum_{a=0}^{\delta-1}(r-1)\binom{r+\delta-a-2}{\delta-a-1}\\
&=\delta\binom{\delta+r-1}{r-1}-(r-1)\binom{\delta+r-1}{r}=\Big(\delta-\frac{\delta(r-1)}{r}\Big)\binom{\delta+r-1}{r-1}.
\end{split}\]
Since 
\[\dim_K(S^\delta(E))=\binom{r+\delta-1}{\delta},\]
we obtain 
\[\widehat{\mu}(S^\delta(\overline E))=\frac{\delta}{r}\sum_{i=1}^r\widehat{\deg}(\overline E_i/E_{i-1})=\delta\operatorname{\widehat{\mu}}(\overline E).\]

In the case where $\overline E$ is semi-stable, by Corollary \ref{Cor: tensorial semistability} we obtain that $\overline E^{\otimes\delta}$ is also semi-stable. Moreover, its slope is also equal to $\delta\operatorname{\widehat{\mu}}(\overline E)$. Since any quotient vector space of $S^\delta(E)$ is also a quotient vector space of $E^{\otimes\delta}$, we obtain that, for any quotient vector space $Q$ of $E^{\otimes\delta}$, one has
\[\widehat{\mu}(\overline Q)\geqslant \widehat{\mu}(\overline E^{\otimes\delta})=\delta\widehat{\mu}(\overline E)=\widehat{\mu}(S^\delta(\overline E)).\]
Therefore $S^{\delta}(\overline E)$ is also semi-stable.
\end{proof}

\subsection{Symmetric power}

In this subsection, we fix a Hermitian adelic vector bundle $\overline E$ on $S$.

\begin{theo}\label{Thm: maximal slope of symm power}
For any positive integer $\delta$, the following inequality holds:
\begin{equation}
\widehat{\mu}_{\max}(S^\delta(\overline E))\leqslant\delta\operatorname{\widehat{\mu}_{\max}}(\overline E)+\nu(\Omega_\infty)\ln(\delta !).
\end{equation}
Moreover, in the case where $\Omega_\infty$ is empty, the norm \textup{(}where we consider the trivial absolute value on $K$\textup{)} on $S^\delta(E)$ associated with the Harder-Narasimhan $\mathbb R$-filtration of $S^\delta(\overline E)$ coincides with the $\varepsilon$-symmetric power of that associated with the Harder-Narasimhan $\mathbb R$-filtration of $\overline E$.
\end{theo}
\begin{proof}
We first treat the case where $\Omega_\infty=\varnothing$. Let  $\mathcal F$ be the Harder-Narasimhan $\mathbb R$-filtration of $\overline E$, which correspond to a sequence
\[0=E_0\subsetneq E_1\subsetneq \ldots\subsetneq E_r=E\]
of vector subspaces of $E$, together with a decreasing sequence
\[\mu_1>\ldots>\mu_r\]
of successive slopes. We equip $S^\delta(E)$ with the symmetric power of the $\mathbb R$-filtration $\mathcal F$. Note that the subquotient $\operatorname{sq}^t(S^\delta(E))$ of index $t$ is given by (see \S\ref{Subsec: subquotient on symmetric power})
\[\bigoplus_{\begin{subarray}{c}
\boldsymbol{b}=(b_1,\ldots,b_r)\in\mathbb N^r\\
|\boldsymbol{b}|=b_1+\cdots+b_r=\delta\\
b_1\mu_1+\cdots+b_r\mu_r=t
\end{subarray}}S^{b_1}(V_1/V_0)\otimes\cdots\otimes S^{b_r}(V_r/V_{r-1}).\]
By Corollary \ref{Cor: tensorial semistability} and  Proposition \ref{Pro: mu S delta}, each Hermitian adelic vector bundle 
\[S^{b_1}(\overline{E_1/E_0})\otimes\cdots\otimes S^{b_r}(\overline{E_r/E_{r-1}})\]
is semi-stable of slope
\[b_1\mu_1+\cdots+b_r\mu_r=t.\]
Therefore, the symmetric power of the $\mathbb R$-filtration $\mathcal F$ identifies with the Harder-Narasimhan $\mathbb R$-filtration of $S^\delta(\overline E)$ and  the maximal slope of $S^\delta(\overline E)$ is equal to $\delta\operatorname{\widehat{\mu}}_{\max}(\overline E)$.

In the case where $\Omega_\infty$ is not empty, the field $K$ is necessarily of characteristic $0$. Let $\operatorname{sym}':S^\delta(E)\rightarrow E^{\otimes\delta}$ be the $K$-linear map induced by the symmetrization map (see Remark \ref{Rem: operator norm symmetri cpower}). Since $K$ is of characteristic $0$, this map is injective and hence
\[\widehat{\mu}_{\max}(S^{\delta}(\overline E))\leqslant\widehat{\mu}_{\max}(\overline E^{\otimes\delta})+h(\operatorname{sym}')\leqslant\widehat{\mu}_{\max}(\overline E^{\otimes\delta})+\nu(\Omega_\infty)\ln(\delta!),\]
as required.
\end{proof}

%\lquery{95pt}{\textcolor{mred}{In the following, I keep the color `darkblue' for the parts related to this revision. In this sense, minor typos are not subject to the colored items.}}

%\include{AppendixA.tex}
%\include{AppendixB.tex}

\ifsmf
\backmatter
\fi

\bibliography{Hilbert-Samuel}
\bibliographystyle{smfplain_nobysame}

\printindex

%remerciement: Laurent Berger
\end{document}